%% file: Takato-Vallejo-final.tex
\documentclass{mathincs}

\usepackage{graphicx}
\usepackage{amsmath,amssymb}
\usepackage{color}
\usepackage{url}
\usepackage{bm}
\usepackage{tikz}
\usepackage{pict2e}
\usepackage{hyperref}
\usepackage{ketpic2e,ketlayer2e}

 \theoremstyle{definition}
 
 \theoremstyle{remark}
 
 \newtheorem{ex}{Example}

\begin{document}
%
%
%
\title[Using Oshima splines]
 {Using Oshima splines to produce accurate numerical results and
 high quality graphical output}
\author[S. Takato]{Setsuo Takato}

\address{%
T\=oh\=o University \\
2-2-1, Miyama\\
Funabashi\\
Japan}
\email{takato@phar.toho-u.ac.jp}

\author[J. A. Vallejo]{Jos\'e A. Vallejo}
\address{
Universidad Aut\'onoma de San Luis Potos\'i\\
Av. Salvador Nava s/n\\ 
78290 San Luis Potos\'i (SLP)\\ 
M\'exico
}
\email{jvallejo@fc.uaslp.mx }



\date{\today}

\begin{abstract}

We illustrate the use of Oshima splines in producing high-quality \LaTeX\ output in two
cases: first, the numerical computation of derivatives and integrals, and second, the
display of silhouettes and wireframe surfaces, using the macros package \ketcindy . Both cases are of particular interest for
college and university teachers wanting to create handouts to be used by students, or
drawing figures for a research paper. When dealing with numerical computations, \ketcindy\ can make a call to the CAS Maxima to check for accuracy; in the case of surface graphics,
it is particularly important to be able to detect intersections of projected curves, and we show how to do it in a seamlessly manner using Oshima splines in \ketcindy . A C compiler can be called in this case to speed up computations. 
\end{abstract}

\maketitle
\section{Introduction}

Cinderella is a dynamic geometry software (DGS) comprising two main components:
CindyScreen, an interactive screen where geometric elements can be constructed
like in any other DGS, and CindyScript, a scripting language which can manipulate
not only  geometric objects but more general constructions.
The following figure shows both components of Cinderella, CindyScreen (left), and the
CindyScript editor (right).

\begin{center}
\includegraphics[scale=0.26]{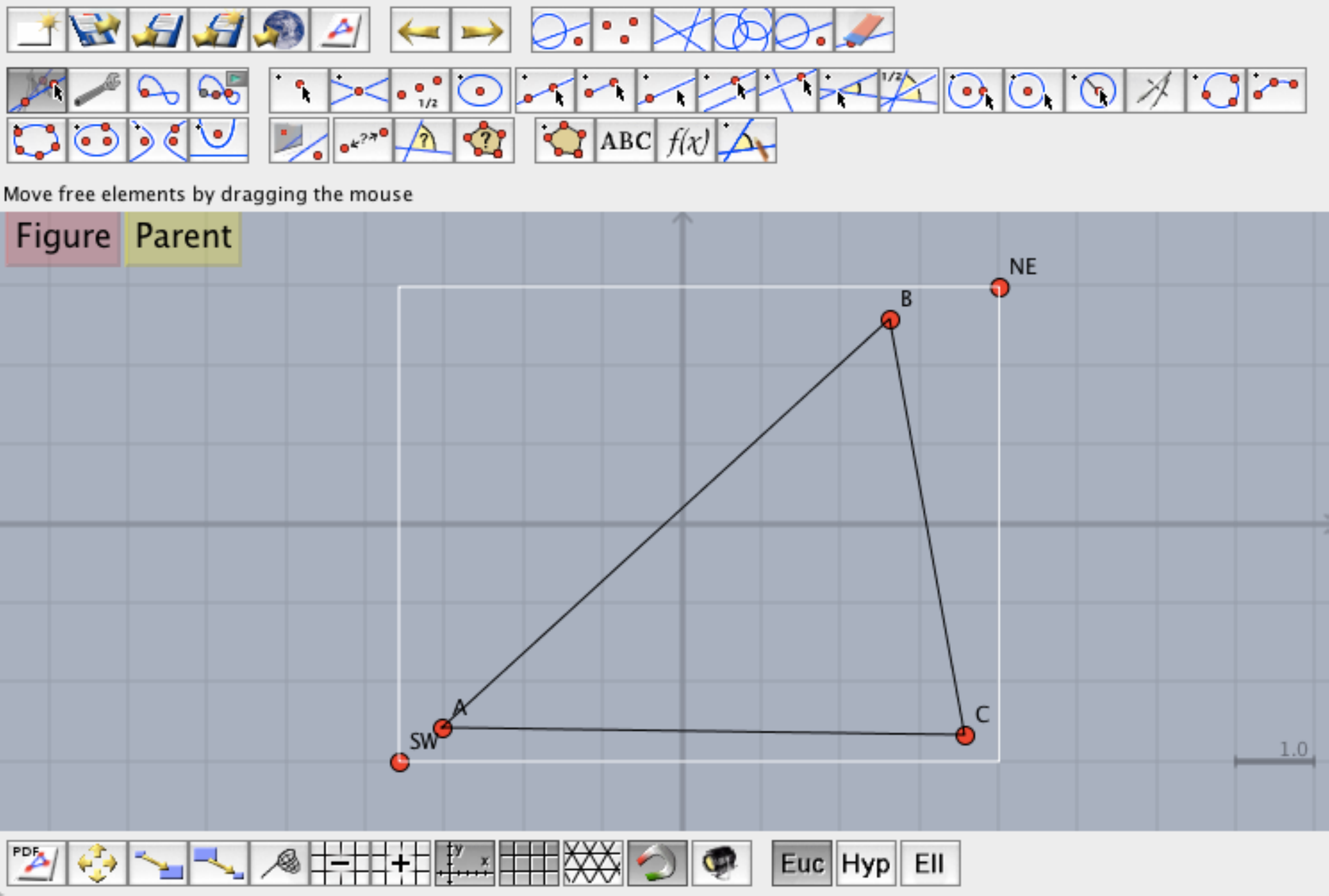}\hspace{10mm}
\includegraphics[scale=0.25]{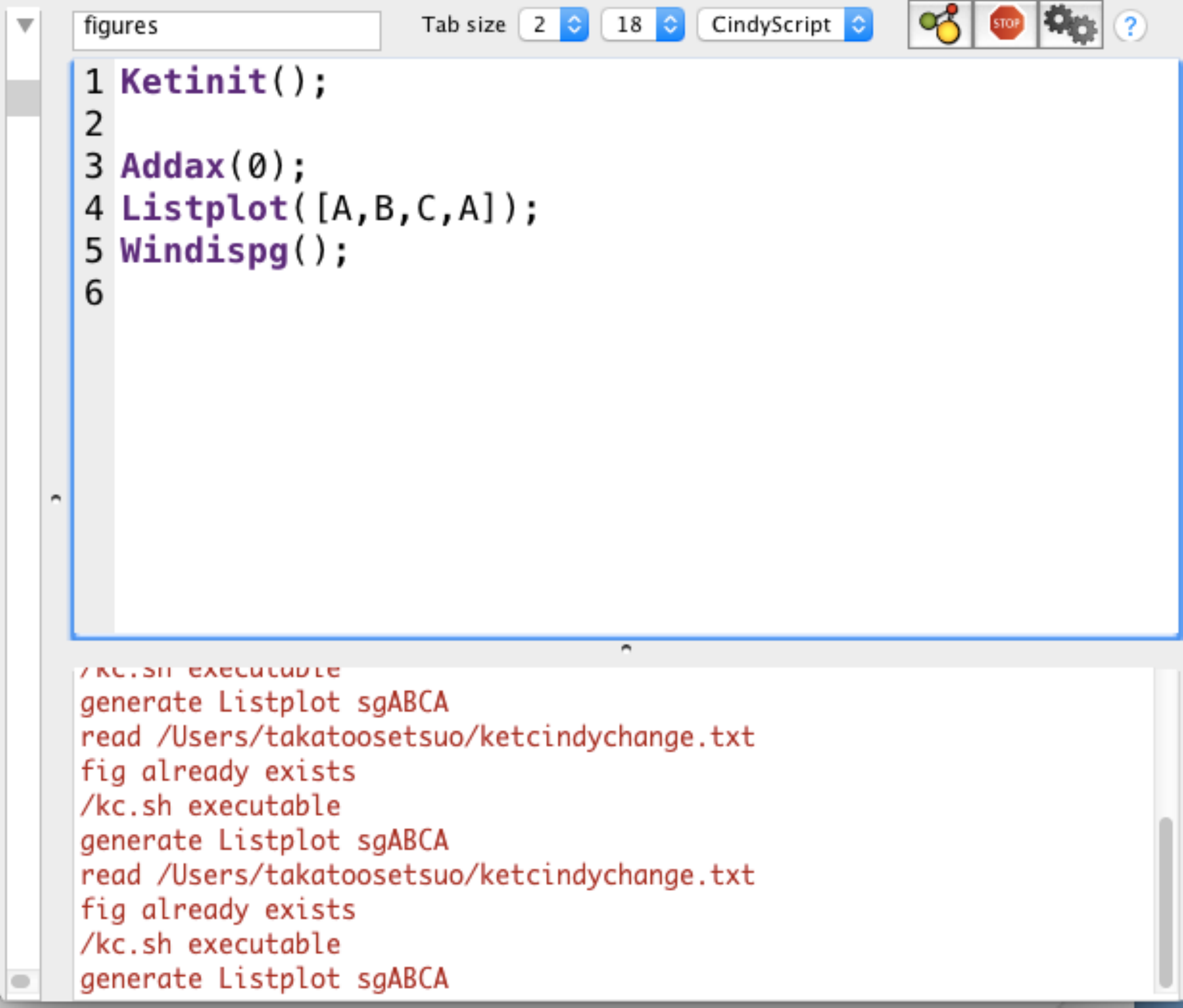}
\\
{\bf Fig.\thefigure\ }CindyScreen and CindyScript\addtocounter{figure}{1}
\end{center}

\ketcindy\ is a package of CindyScript macros designed to produce high-quality \LaTeX\ 
figures and animations (see \cite{T2016,TMVK}). It is well suited for particularly
complex graphics such as the ones appearing in dynamical systems theory \cite{VT2}, and
flexible enough to satisfy the demands of a wide class of topics ranging from linear
algebra to the calculus of variations, including the preprocessing of data for 3D printing \cite{KYKMNKT}--\cite{kor}. It can be accessed from
the Comprehensive \TeX\ Archive Network (CTAN) repository \cite{ctan}, requiring
\verb|Cinderella| \cite{cindy} and \verb|R| \cite{R} as dependences. In this text we will also require the optional Computer Algebra System (CAS) \verb|Maxima| \cite{maxima}.
The package comes with its own documentation explaining the installation process.

As an example of its use, we now show how to draw a freehand smooth curve. For this
task, \ketcindy\ uses B\'ezier curves and has several commands to create splines from them: 
\verb|Bezier|,  \verb|CRspline|, \verb|Ospline|, \verb|Bspline| and \verb|Mkbezierptcrv| are 
already implemented. Here we use \verb|CRspline|, the command to draw a Catmull-Rom spline.
The steps to generate a plotting data file suitable for being included in a \LaTeX\ 
document are as follows:

\begin{enumerate}
\item Open a template in the work folder of ketcindyfolder or any \ketcindy\ file. The rectangle in the screen shows the drawing range for the \verb|picture| environment in \LaTeX.
\item Change points SW, NE to fix the drawing area and add some other points A, B, C, D
that will serve as control points for the spline.

\begin{center}
\noindent
\begin{minipage}{0.485\textwidth}
\centering\noindent
\includegraphics[scale=0.1925]{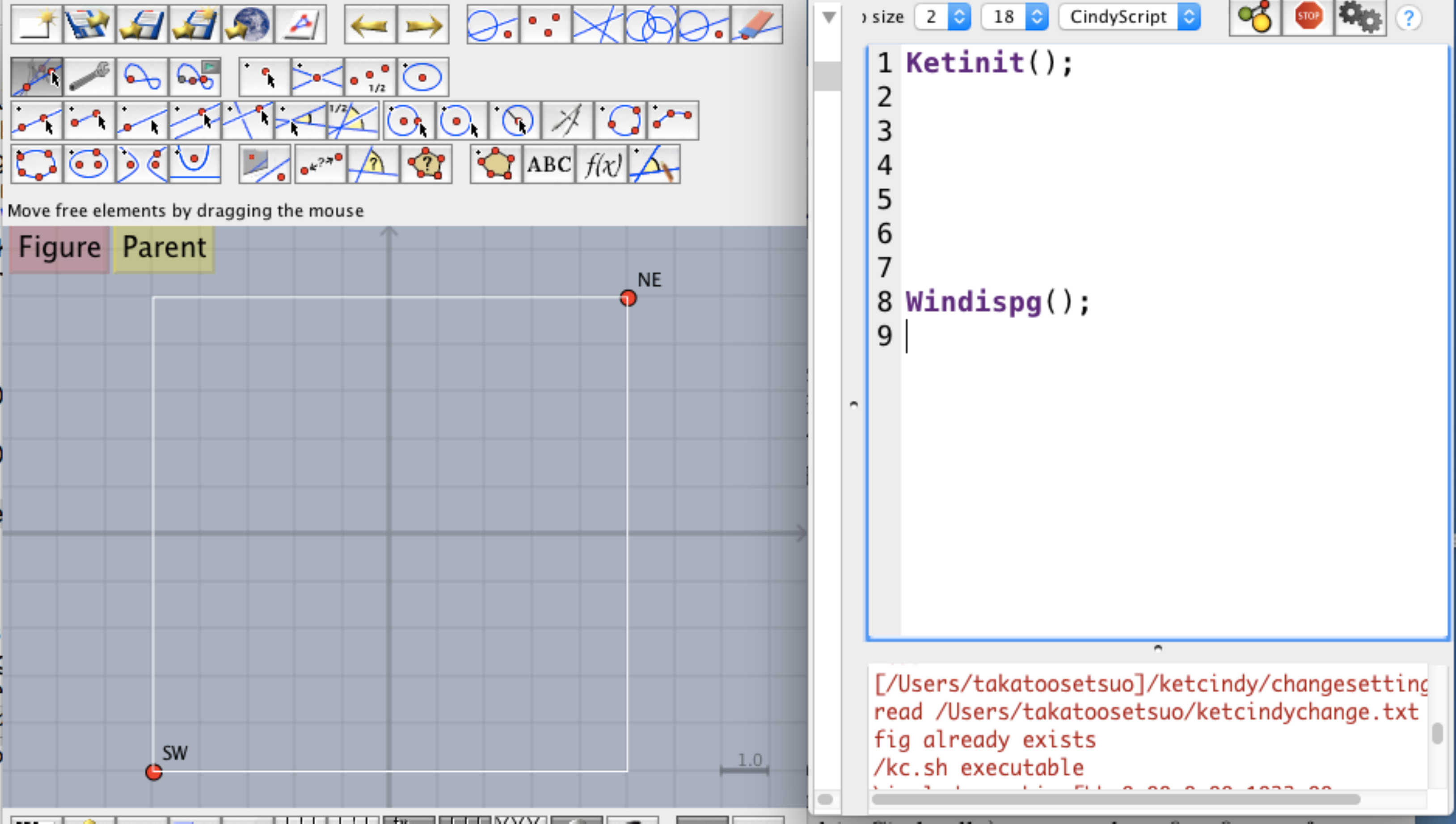}
\\
{\bf Fig.\thefigure\ }Initial screen\addtocounter{figure}{1}
\end{minipage}%
\begin{minipage}{0.485\textwidth}
\centering
\includegraphics[scale=0.1925]{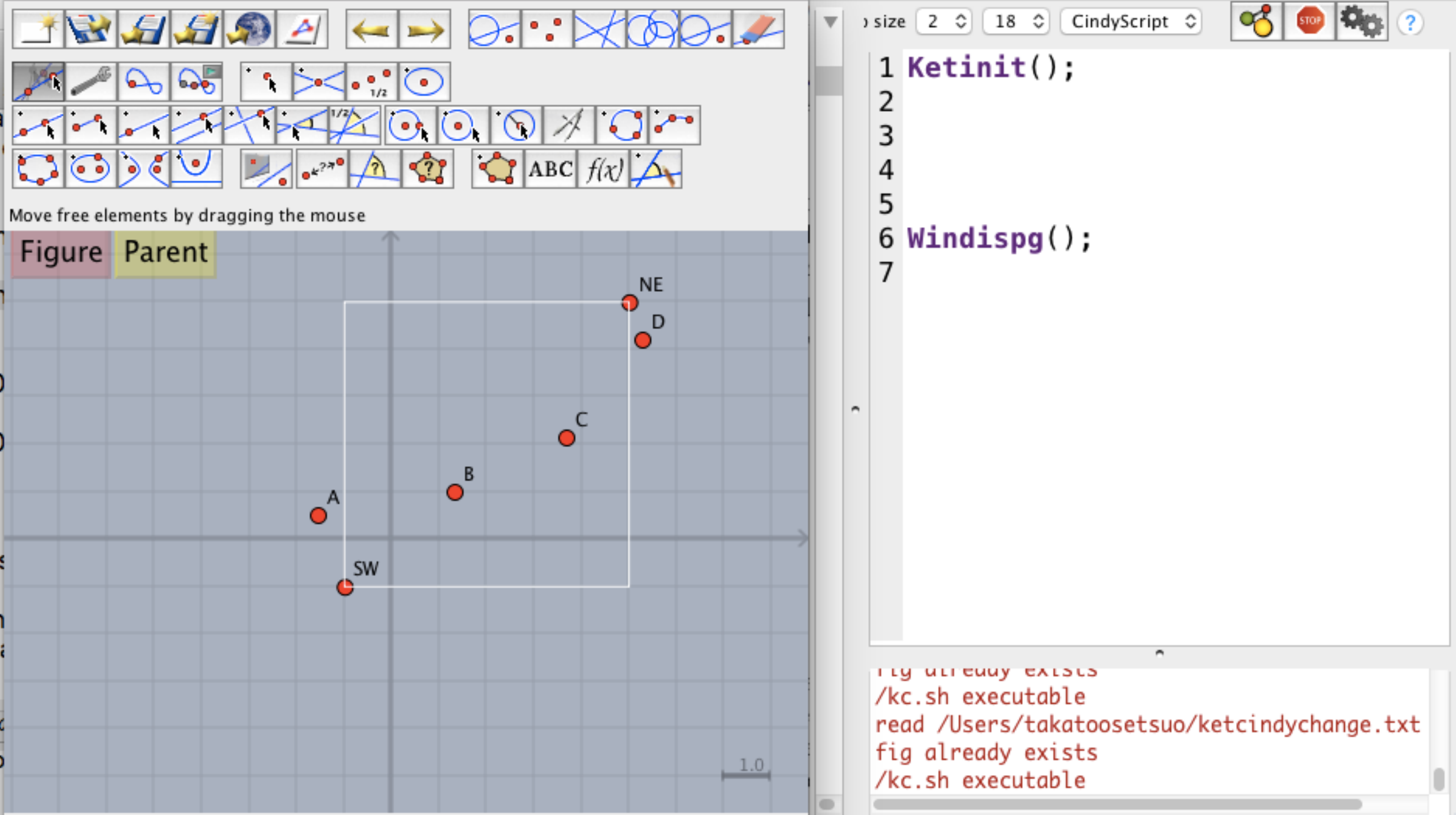}
\\
{\bf Fig.\thefigure \ }Adding points\addtocounter{figure}{1}
\end{minipage}
\end{center}

\item In the script editor, write the command \verb|CRspline("1",[A,B,C,D])|: the curve is displayed in the screen. One can change the shape of the curve by moving any of the control
points A,B,C,D in the CindyScreen.
\item Pressing the button \verb|Figure| in the screen will generate the plotting data for \LaTeX\ . The output is represented in Figure 5 below.

\begin{center}
\noindent
\begin{minipage}{0.49\textwidth}
\centering
\includegraphics[scale=0.2]{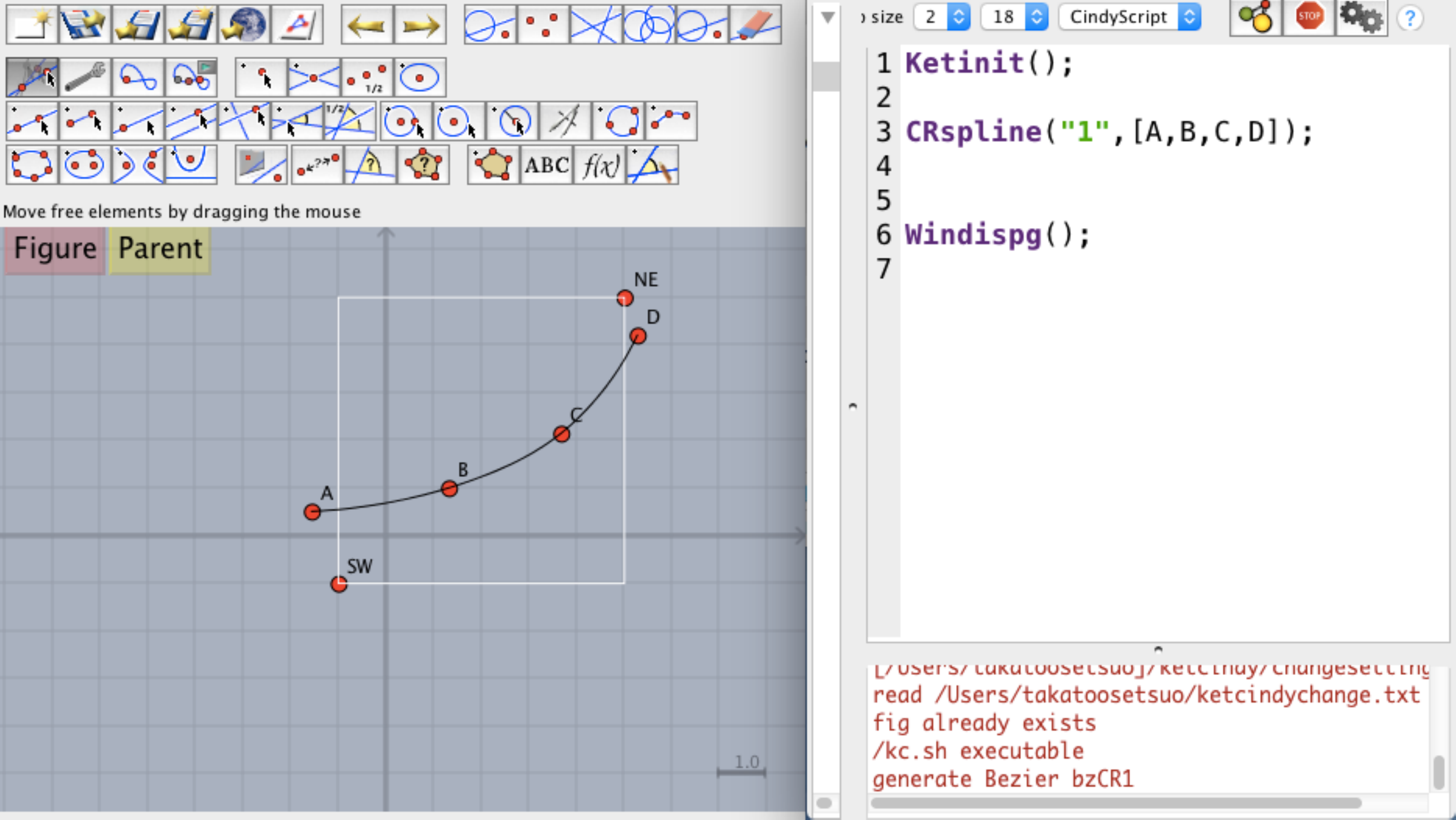}
\\
{\bf Fig.\thefigure\ }Drawing CRspline curve\addtocounter{figure}{1}
\end{minipage}%
\begin{minipage}{0.49\textwidth}
\centering\noindent
\includegraphics[scale=0.55]{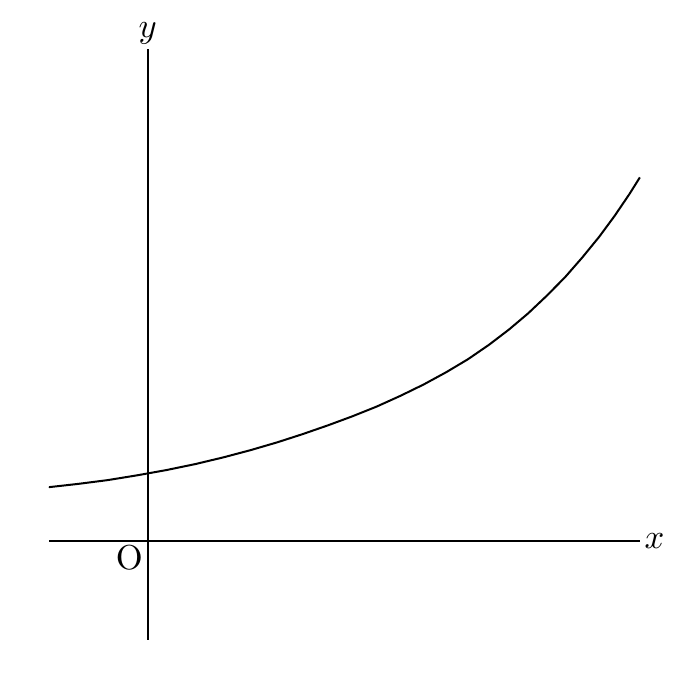}
\\
{\bf Fig.\thefigure \ }Plotting \LaTeX\ data \addtocounter{figure}{1}
\end{minipage}
\end{center}

\noindent\hspace*{2mm}%

\vspace{2mm}

\item Further embellishment of the figure is, of course, possible. Le us complete it
to illustrate the notion of differences quotient. It suffices to write the following commands in the CindyScript editor:

\begin{verbatim}
  CRspline("1",[A,B,C,D]);
  pt=[C.x,B.y];
  Listplot("1",[B,pt,C,B]);
  Bowdata("1",[B,pt],["do","E=\Delta x"]);
  Bowdata("2",[pt,C],[1.5,"do","E=\Delta y"]);
  Letter(B,''n'',''P'');
  Htickmark([B.x,"x",C.x,"\xi"]);
  Vtickmark([B.y,"f(x)",C.y,"f(\xi)"]);
\end{verbatim}

\noindent
Again, pressing \verb|Figure| will generate the \LaTeX\ code and the corresponding
\verb|pdf| file.

\begin{center}
\includegraphics[scale=0.3]{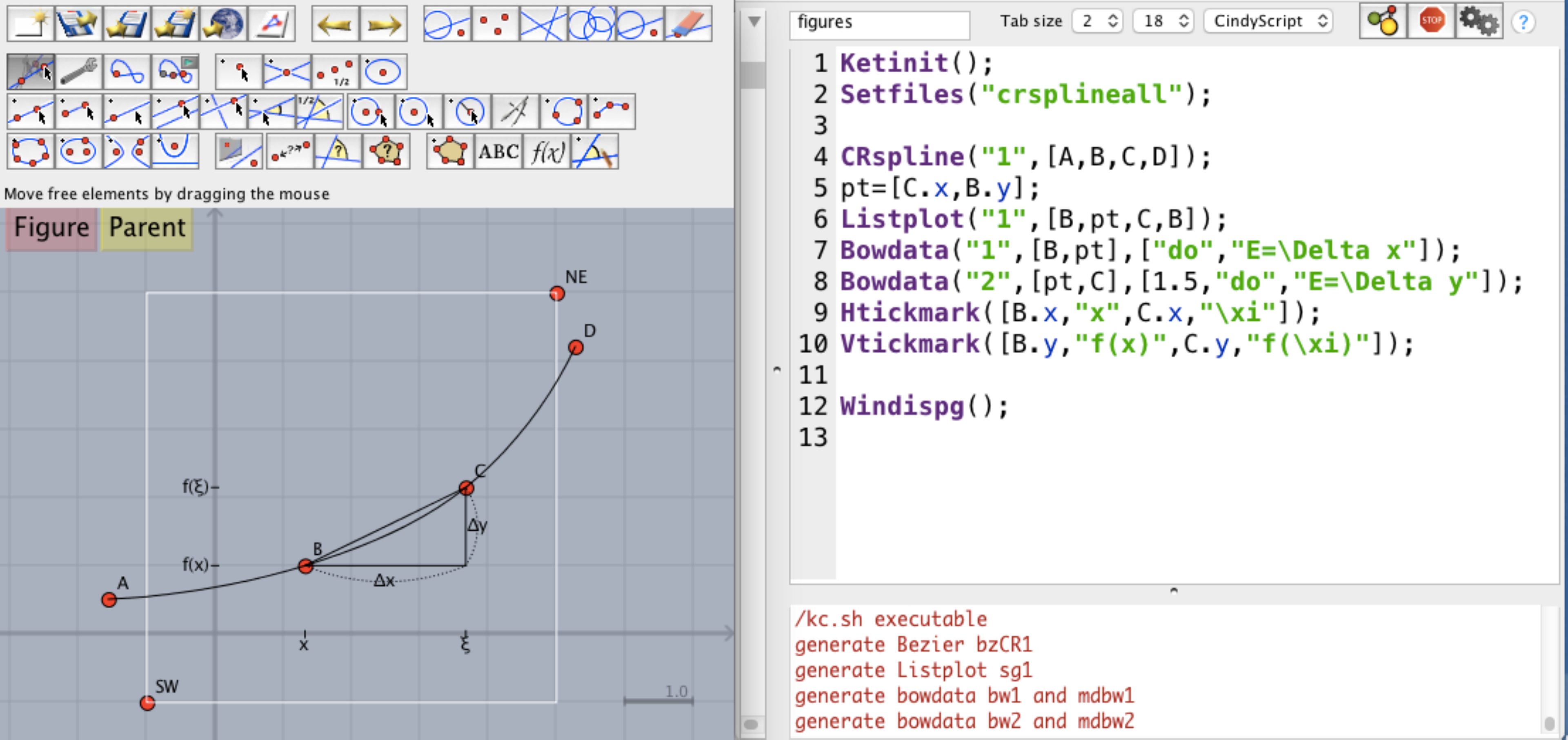} 
\\
{\bf Fig.\thefigure \ }Adding other components to the figure\addtocounter{figure}{1}
\end{center}
\end{enumerate}

As a side remark, sometimes more flexibility than that provided by \LaTeX\ is required
to insert figures in a document. This is particularly the case when creating handouts to
be distributed to students, quizzes, cheat sheets, and the like. To deal with these cases,
the \verb|layer| environment, which is defined in \verb|ketlayer.sty|, can be useful. If
\verb|graph1.tex| is the result of steps $1-5$ above, the code

\begin{verbatim}
    \begin{layer}{140}{50}
    \putnotes{115}{2}{$\frac{\varDelta y}{\varDelta x}\to f'(x)$}
    \putnotese{80}{5}{\inputf{fig/graph1.tex}}
    \end{layer}
\end{verbatim}

\noindent
when inserted in a \LaTeX\ document, will produce the following:

\vspace{3mm}

\begin{layer}{140}{50}
\putnotese{85}{3}{\scalebox{0.8}{\input{fig/crsplineall.tex}}}
\putnotes{115}{2}{$\frac{\varDelta y}{\varDelta x}\to f'(x)$}
\boxframese{7}{-1}{130}{53}{}
\end{layer}

\noindent
\hspace*{10mm}\begin{minipage}{85mm}
\vspace{1cm}
Let $\mathrm{P}(x,\ f(x))$ be a point on the graph of $y=f(x)$.\\
Let $(\xi,\ f(\xi))$ be a point close to P, and put
$$
\Delta x=\xi-x, \Delta y=f(\xi)-f(x)\,.
$$
Then $\dfrac{dy}{dx}=\displaystyle\lim_{\xi \to x}\dfrac{\Delta y}{\Delta x}$ is the derivative of the function.\\
Question : Calculate $\dfrac{dy}{dx}$ when $f(x)=x^2$.
\vspace{1cm}
\end{minipage}
\vspace{5mm}

\begin{center}
{\bf Fig.\thefigure \ }Example of a material for use in teaching\addtocounter{figure}{1}
\end{center}

Change the arguments of \verb|\putnote| to move each component.
If they are placed properly, set the second argument of \verb|layer| to 0, then grids will  disappear
without moving the position of all components.

\begin{verbatim}
    \begin{layer}{140}{0}
    \putnotese{80}{5}{\inputf{fig/graph1.tex}}
    \putnotes{115}{2}{$\frac{\varDelta y}{\varDelta x}\to f'(x)$}
    \end{layer}
\end{verbatim}

\vspace{3mm}

\begin{layer}{140}{0}
\putnotese{85}{3}{\scalebox{0.8}{\input{fig/crsplineall.tex}}}
\putnotes{115}{2}{$\frac{\varDelta y}{\varDelta x}\to f'(x)$}
\boxframese{7}{-1}{130}{53}{}
\end{layer}

\noindent
\hspace*{10mm}\begin{minipage}{85mm}
\vspace{1cm}
Let $\mathrm{P}(x,\ f(x))$ be a point on the graph of $y=f(x)$.\\
Let $(\xi,\ f(\xi))$ be a point close to P, and put
$$
\Delta x=\xi-x, \Delta y=f(\xi)-f(x)\,.
$$
Then $\dfrac{dy}{dx}=\displaystyle\lim_{\xi \to x}\dfrac{\Delta y}{\Delta x}$ is the derivative of the function.\\
Question : Calculate $\dfrac{dy}{dx}$ when $f(x)=x^2$.
\vspace{1cm}
\end{minipage}
\\
\begin{center}
{\bf Fig.\thefigure \ }Example of a material for use in teaching\addtocounter{figure}{1}
\end{center}

\verb|\includegraphics| of \verb|graphicx| cannot do the same thing by itself. Moreover, 
\verb|\input| might be more convienient when inputting the file cosist of graphical codes such as \verb|pict2e| and \verb|TikZ|.

\vspace{3mm}

The Catmull-Rom spline which we used in this example is well known and often used to draw 
freehand shapes or to make a nonlinear interpolation. However, as with any other class of 
splines, there are many cases in which the results it gives are not as good as one desires.
This is particularly true for conic sections, which can not be exactly generated using
B\'ezier curves: if we take four points A, B, C, D on an ellipse and use the Catmull-Rom 
spline to approximate it, the corresponding code is
\begin{verbatim}
   Paramplot("1","[3*cos(t),2*sin(t)]","t=[0,2*pi]",["da"]);
     // shows the ellipse with a dashed line.
   A=[3,0]; B=[0,2]; C=[-3,2]; D=[0,-2];
   CRspline("1",[A,B,C,D,A]);
\end{verbatim}

\noindent
and the resulting figure is the following:

\begin{center}
\scalebox{0.8}{\input{fig/crsplineclosed.tex}}

{\bf Fig.\thefigure \ }Catmull-Rom spline for an ellipse\addtocounter{figure}{1}
\end{center}

\vspace{2mm}

Even recognizing the impossibility of exactly reproducing conic sections, this result is
less than sub-optimal. Due to the shortcomings of Catmull-Rom splines, Oshima developed a 
new idea to determine the spline control points from the given points on the curve 
\cite{Oshima2016}, in such a way that the interpolation is nearly optimal in the case of
conics (compare Figure 9 with Figure 12 below). 
This opens the possibility of constructing numerical schemes for derivation and integration
based on this spline, as well as a promising technique to attack the problem
of determining the intersection points of surfaces and curves, a basic problem in Computer
Aided Geometric Design (or CAGD) that must be solved in order to construct wireframes.
Thus, our contribution in this paper is to describe these techniques, and illustrate both
its practical implementation and its application to some typical problems. In this regard,
we have used the set of macros \ketcindy\ because it allows us to exploit the very precise
numerical computation of wireframes resulting from Oshima techniques (when those computations
are done through an external C compiler), and as a result, we can produce high-quality
graphical output, as shown in the examples below. The resulting framework, integrating
the DGS Cinderella, the CAS Maxima and a C compiler with \ketcindy\ as the interface, has proven
to be very flexible and powerful, allowing the user to carry on the numerical tasks as well
as the generation of graphics in an intuitive and unified environment.

\section{Oshima spline curve}

Let points $\mathrm{P}_{j-1}\,, \mathrm{P}_{j}\,, \mathrm{P}_{j+1}\,, \mathrm{P}_{j+2}$ be on a certain curve which we want to approximate, and $\mathrm{Q}_{j},\ \mathrm{R}_{j}$ be the
control points corresponding to an interval $\mathrm{P}_{j}\mathrm{P}_{j+1}$. In the case of
the Catmull-Rom spline curve, these control points are defined solely by $\mathrm{P}_{j},\mathrm{P}_{j+1}$ using the conditions
$$
\overrightarrow{\mathrm{P_j Q_j}}=\frac{1}{6}\overrightarrow{\mathrm{P_{j-1} P_{j+1}}},\quad 
\overrightarrow{\mathrm{P_{j+1}R_j}}=\frac{1}{6}\overrightarrow{\mathrm{P_{j+2} P_{j}}}\,.
$$
Note that the coefficients appearing here are contants. Figure \thefigure, shows the shape of the resulting spline, where it is to be noticed that the curve bends rapidly. 

Oshima's definition is
$$
\overrightarrow{\mathrm{P_j Q_j}}=c\,\overrightarrow{\mathrm{P_{j-1} P_{j+1}}},\quad 
\overrightarrow{\mathrm{P_{j+1}R_j}}=c\,\overrightarrow{\mathrm{P_{j+2} P_{j}}}\,,
$$
where the value of the coefficient $c$ is determined from
$\mathrm{P}_{j-1}\,, \mathrm{P}_{j}\,, \mathrm{P}_{j+1}\,, \mathrm{P}_{j+2}$
as follows:
$$ 
c=\frac{4\overrightarrow{\mathrm{P_j P_{j+1}}}}
{3(\overrightarrow{\mathrm{P_{j-1}P_{j+1}}}+\overrightarrow{\mathrm{P_{j}P_{j+2}}})}
\times \frac{1}{1+\sqrt{\frac{1}{2}\left( 1+\cos\theta\right)}}\,,
$$
where $\theta$ is the angle between $\overrightarrow{\mathrm{P_{j-1}P_{j+1}}}$ and $\overrightarrow{\mathrm{P_{j}P_{j+2}}}$.
See Figure \addtocounter{figure}{1}\thefigure\addtocounter{figure}{-1}, and notice that 
this time the curve is  smoother, as well as the different placement of the control points. Figure \addtocounter{figure}{2}\thefigure\addtocounter{figure}{-2} represents the Oshima spline for the ellipse of the preceding section.

\vspace{-3mm}

\noindent\hspace*{10mm}%
\begin{minipage}{60mm}
\begin{center}
\scalebox{0.8}{\input{fig/crospline1}}

\vspace{-5mm}

{\bf Fig.\thefigure \ }Case of Catmull-Rom spline\addtocounter{figure}{1}
\end{center}
\end{minipage}%
\hspace{5mm}
\begin{minipage}{60mm}
\begin{center}
\scalebox{0.8}{\input{fig/crospline2}}

\vspace{-5mm}

{\bf Fig.\thefigure \ }Case of Oshima spline\addtocounter{figure}{1}
\end{center}
\end{minipage}

\vspace{2mm}

\begin{center}
\scalebox{0.8}{\input{fig/osplineclosed.tex}}

{\bf Fig.\thefigure \ }Oshima spline for an ellipse\addtocounter{figure}{1}

\end{center}

\vspace{2mm}

Though we can not say that Oshima spline gives a better interpolation in all cases,
in general it certainly does when it comes to drawing smooth freehand curves. As an
application of this fact, we will apply next the Oshima technique to the problem of
numerical differentiation and integration. \ketcindy\ provides two commands 
(\verb|Integrate|, \verb|Derivative|) for these tasks.

\subsection{Numerical integration}

\begin{layer}{140}{0}
\putnotes{102}{2}{\input{fig/int1.tex}}
\putnotes{102}{40}{{\bf Fig.\thefigure \ }Definite integral of a B\'ezier curve\addtocounter{figure}{1}}
\end{layer}

\noindent
\begin{minipage}{65mm}
Let $C$ be a B\'ezier curve determined by the control
points $\mathrm{P}_j$,\ $\mathrm{P}_{j+1}$,\ 
$\mathrm{P}_{j+2}$,\ $\mathrm{P}_{j+3}$, and let $(x_k,y_k)$ denote the coordinates of $\mathrm{P}_k$.

Since the parametric equations of $C$ are
\begin{align}\label{eq1}
\mathrm{P}=\mathrm{P}_{j}(1-t)^3+3\mathrm{P}_{j+1}(1-t)^2 t\nonumber\\
+3\mathrm{P}_{j+2}(1-t) t^2+\mathrm{P}_{j+3}t^3\\
(0\leqq t \leqq 1),\nonumber
\end{align}
the definite integral $\displaystyle\int_{x_j}^{x_{j+3}}y\,\mathrm{d}x$ becomes
$$
\int_{0}^{1}y\frac{\mathrm{d}x}{\mathrm{d}t}\mathrm{d}t\,.
$$
\end{minipage}

\vspace{2mm}

\ketcindy\ can call Maxima using \verb|Mxfun| to execute a single command and 
\verb|CalcbyM| to execute several commands \cite{TMVK}. Thus, we can evaluate
the preceding integral with the following CindyScript code:

\begin{verbatim}
    cmdL=[
     "P:[x1,y1]*(1-t)^3+3*[x2,y2]*(1-t)^2*t
              +3*[x3,y3]*(1-t)*t^2+[x4,y4]*t^3",[],
     "f:P[2]*diff(P[1],t)",[],
     "ans:integrate",["f","t",0,1],
     "ans",[]
    ];
    CalcbyM("ans",cmdL);
    println(ans);
\end{verbatim}

\noindent
The output is displayed in the console as 
\begin{verbatim}
    ((10*x4-6*x3-3*x2-x1)*y4+(6*x4-3*x2-3*x1)*y3
            +(3*x4+3*x3-6*x1)*y2+(x4+3*x3+6*x2-10*x1)*y1)/20.
\end{verbatim}

This procedure is implemented in \ketcindy 's command \verb|Integrate|. Here we 
illustrate its use through several examples.

\begin{ex}
To compute the area surrounded by ellipse $\dfrac{x^2}{3^2}+\dfrac{y^2}{2^2}=1$
we would execute the script

\begin{verbatim}
    Paramplot("1","[3*cos(t),2*sin(t)]","t=[0,pi]",["Num=25"]);
    Paramplot("2","[3*cos(t),2*sin(t)]","t=[pi,2*pi]",["Num=25"]);
    ans=Integrate("gp1",[-pi,pi])-Integrate("gp2",[-pi,pi]);
    println(Sprintf(ans/(6*pi),6));
\end{verbatim}

\noindent
The result is $0.999937$, which gives a good approximation to $\dfrac{S}{\pi a b}=1$.
\end{ex}

\begin{ex}
The definite integral of $y=x^2\sin x$ from $0$ to $\pi$ is computed by the script

\begin{verbatim}
    Plotdata("1","x^2*sin(x)","x=[-pi,pi]",["Num=50"]);
    ans=Integrate("gr1",[0,pi]);
    println(Sprintf(ans,6));
\end{verbatim}

\noindent
The result is $5.869063$. We can check this result with Maxima:

\begin{verbatim}
    Mxfun("1","integrate",["x^2*sin(x)","x",0,"%pi"]);
    Mxfun("2","float",[mx1]);
\end{verbatim}
\noindent
The output is \verb|%pi^2-4|, whose numerical value in double precision is $5.869604401089358$.
\end{ex}

Notice that \verb|Integrate| can be applied to a list of points. Actually, \verb|gr1| and 
\verb|gp1|, \verb|gp2| in the preceding examples are lists of points.

\begin{ex}\label{example3}
The definite integral of the implicit function
$8x^2-4\sqrt{2}xy+y^2-3x-6\sqrt{2}y+2=0$ in the region
$[-2,2]\times [-2,2.5]$ is performed by

\begin{verbatim}
    Implicitplot("1","8*x^2-4*sqrt(2)*x*y+y^2-3*x-6*sqrt(2)*y+2=0",
            "x=[-2,2]","y=[-2,2.5]");
    P=Ptstart("imp1");
    Q=Ptend("imp1");
    Letter([P,"n","P",Q,"n","Q"]);
    ans=Integrate("imp1",[Q_1,P_1]);
    println(Sprintf(ans,6));
\end{verbatim}

\noindent
The result is $1.698725$.

\begin{center}
\input{fig/intimp.tex}

{\bf Fig.\thefigure \ }Integration of an implicit function\addtocounter{figure}{1}
\end{center}
\end{ex}

\subsection{Numerical differentiation}

\noindent
We can also calculate the derivative of a function given by a list of points. To find the derivative of $C$ at the point $\mathrm{P}_j$ in Figure \addtocounter{figure}{-2}\thefigure\addtocounter{figure}{2}, we differentiate \eqref{eq1} and put $t=0$, then
$$
\frac{d\mathrm{P}}{dt}=-3\mathrm{P}_j+\mathrm{P}_{j+1}\,.
$$

\noindent
Using this procedure, we have improved the command \verb|Derivative| and implemented \verb|Tangentplot| to draw the tangent line at a point as an application.
The following is an example of the use of \verb|Tangentplot| with the implicit function
of Example \ref{example3}.

\begin{verbatim}
    dx=(Q_1-P_1)/10;
    forall(0..10,ii,
      v=P_1+ii*dx;
      Tangentplot(text(ii),"imp1","x="+format(v,6),["dr,0.2"]);
    );
\end{verbatim}

\vspace{1mm}

\begin{center}
\input{fig/derivimp.tex}

\vspace{1mm}

{\bf Fig.\thefigure \ }Tangent lines to an implicit function\addtocounter{figure}{1}
\end{center}

\section{Drawing 3D surfaces}
\subsection{Intersection of silhouette lines}
Teaching materials displaying 3D figures are often required in mathematics courses.
For such printed materials, figures presented as line drawings are better suited, because students can write their own remarks over them on the paper (Figure \thefigure). \ketcindy\ supports the line drawing of 3D figures as explained below, but let us first present a couple of examples:

\noindent\begin{minipage}{0.5\textwidth}
\begin{center}
\input{fig/mant.tex}

{\bf Fig.\thefigure}\ Line drawing of a surface\addtocounter{figure}{1}%
\end{center}
\end{minipage}
\begin{minipage}{0.5\textwidth}
\begin{center}
\input{fig/mant2.tex}

{\bf Fig.\thefigure}\ Segmenting curves\addtocounter{figure}{1}%
\end{center}
\end{minipage}

\vspace{2mm}

To produce these 3D figures, \ketcindy\ follows the steps:

\begin{enumerate}
\item The silhouette lines of the surface are determined. To this end, numerical 
data are obtained from an implicit function of the form
$$
J(u,v)=\frac{dX}{du}\frac{dY}{dv}-\frac{dX}{dv}\frac{dY}{du}=0\,,
$$
where $(X,Y)=\mathrm{Proj}(x,y,z)$ is the map onto the plane of projection.
\item The intersections of silhouette lines and projected curves (see Figure \addtocounter{figure}{-1}\thefigure\addtocounter{figure}{1}) are computed.
\item The curve is segmented by these intersections, and it is determined
whether each segment is hidden by the surface or not.
\end{enumerate}

\noindent
Of the above, the second item is of fundamental importance, but it represents a
difficult task in the case of contacting curves because, numerically, they are
polygonal lines (see Figure \thefigure\ for an example, where the small window
in red shows the contact region that will be further discussed below).
 
{\color{red}
\begin{layer}{140}{0}
\dashboxframec[16]{57}{32}{8}{7}{}
\end{layer}
}
\begin{center}
\includegraphics[scale=0.75]{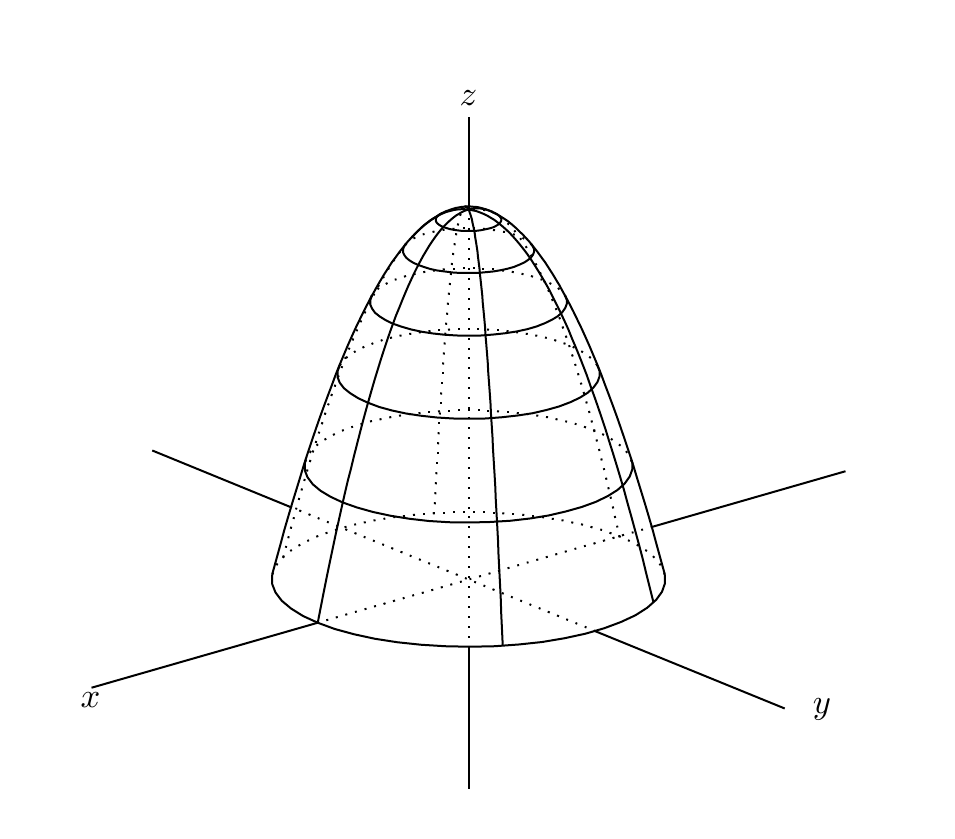}
\\
{\bf Fig.\thefigure}\ Case of contacting curves\addtocounter{figure}{1}
\end{center}

The following graphical examples illustrate this setting: the right figure shows 
a magnification of the contact region in the left one. The parametric equations of
the surface are $x=u\cos v$, $y=u\sin v$, and $z=4-u^2$.

\begin{center}
\includegraphics[scale=0.5]{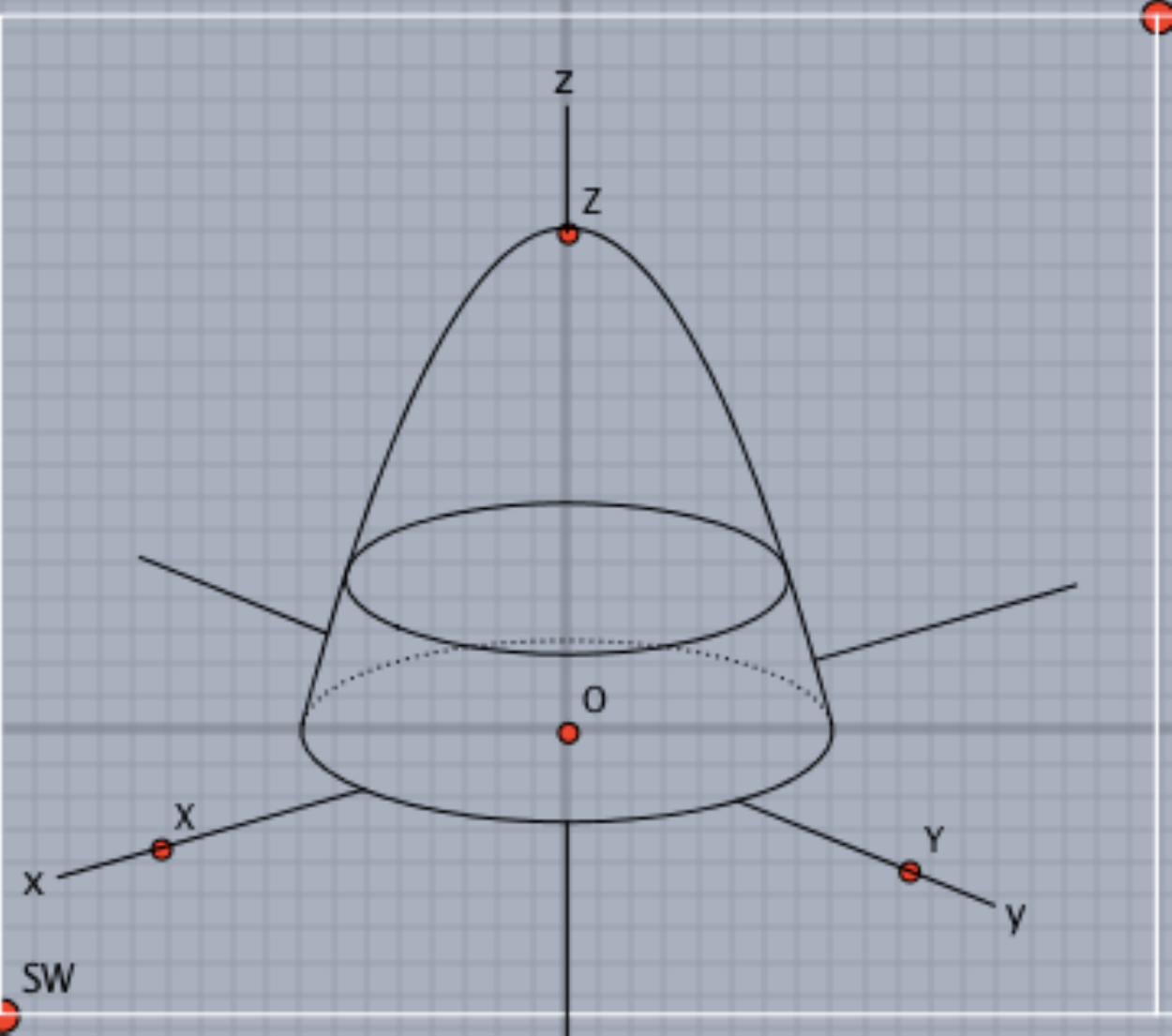}%
\hspace{5mm}%
\includegraphics[scale=0.5]{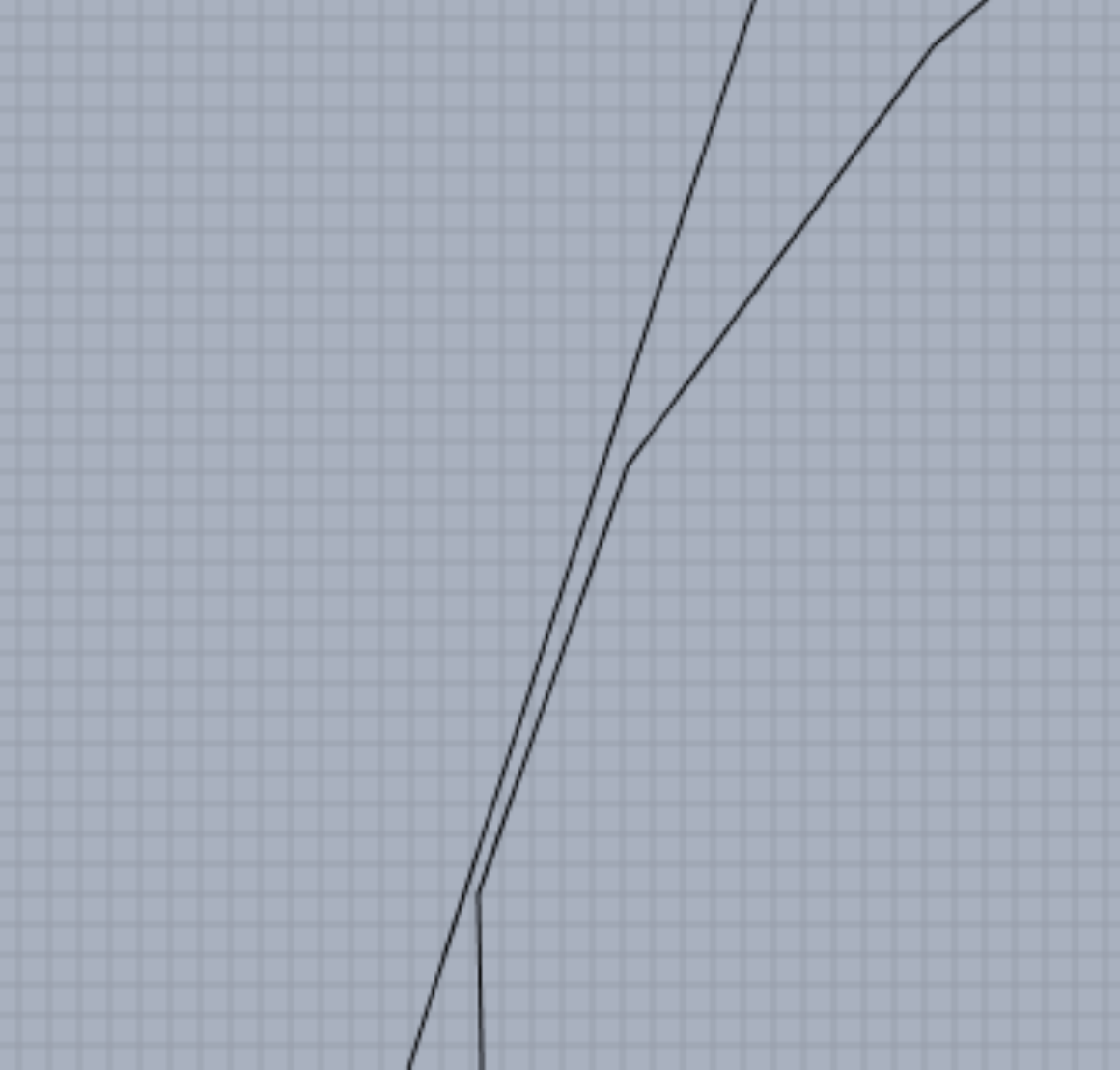} 
\\
{\bf Fig.\thefigure}\ Magnifying around the contact point\addtocounter{figure}{1}
\end{center}

To refine the calculation in item 2, we have adopted an interpolatory scheme using 
Oshima splines around the contact point. 

In the figure below, the left pane shows a further magnification. The right one
shows the Oshima splines of the silhouette line of the surface and the curve in red.

\begin{center}
\includegraphics[scale=0.42]{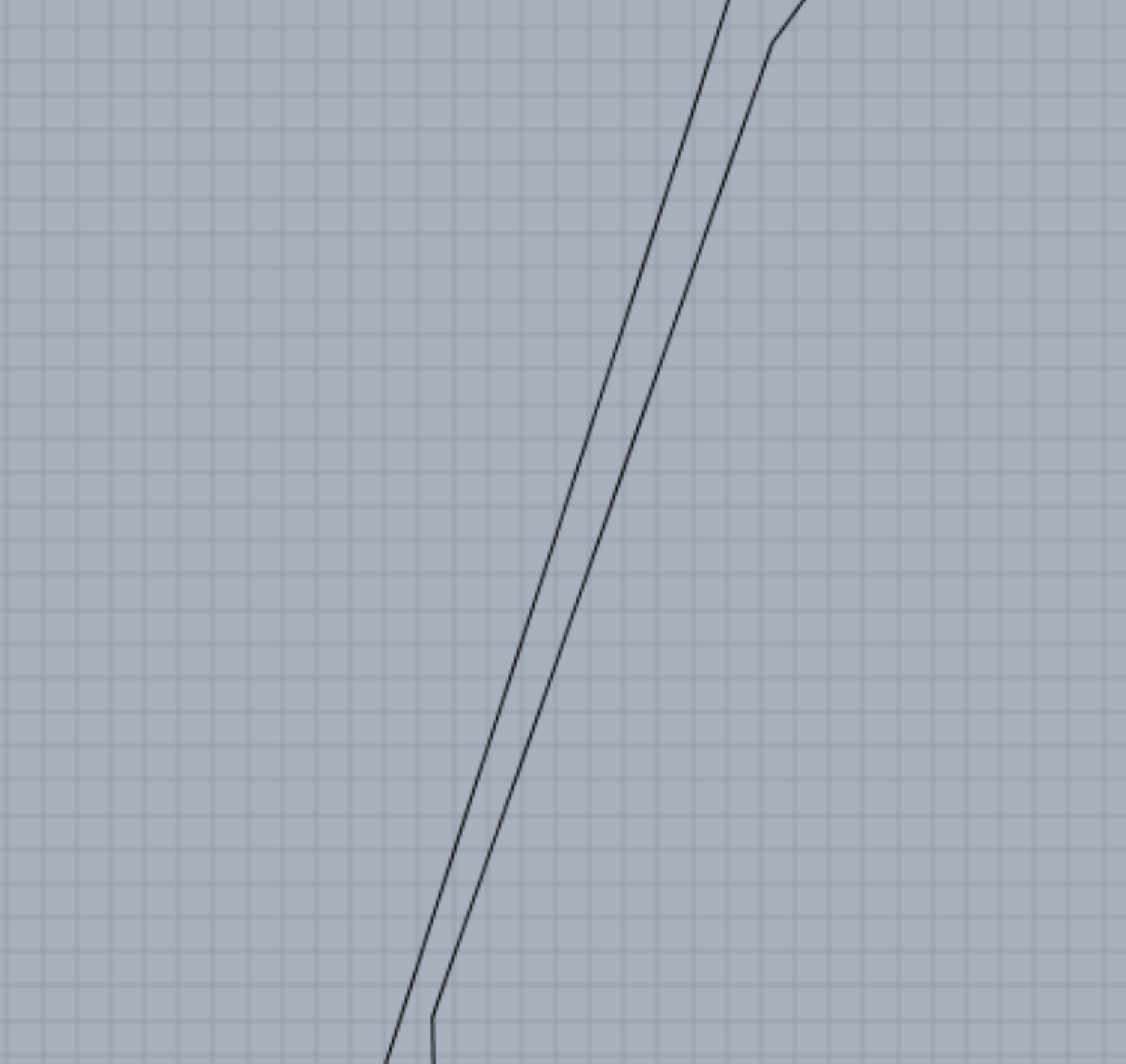}%
\hspace{5mm}%
\includegraphics[scale=0.425]{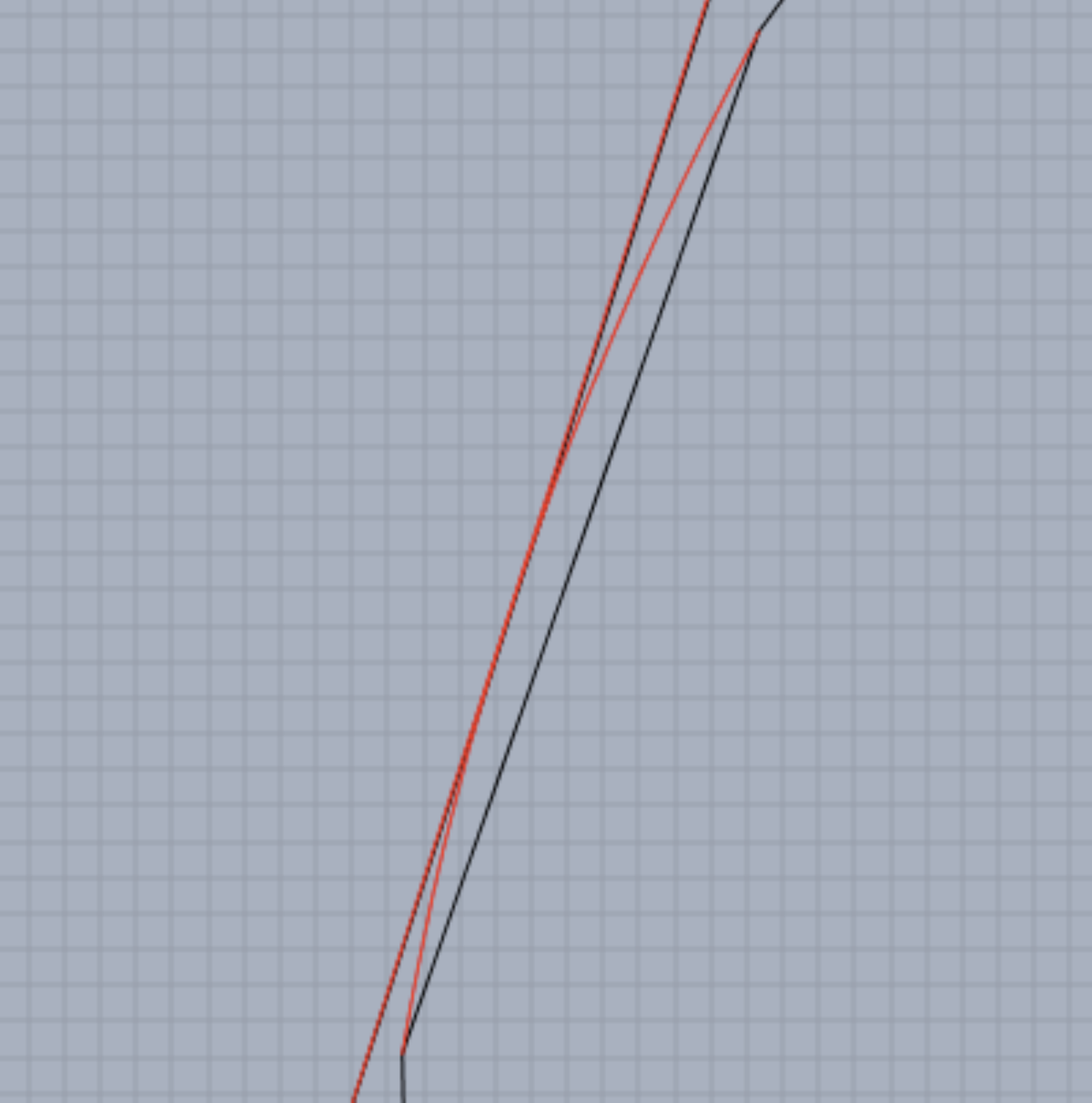} 
\\
{\bf Fig.\thefigure}\ Use of Oshima spline curve\addtocounter{figure}{1}
\end{center}

We can see that the intersection is still represented by a cluster of points, but a
much narrow one. A randomly chosen point in this intersection has coordinates
\begin{equation}\label{eq2}
\mathrm{P}=[-1.65827,1.20578]\,.
\end{equation}

The exact values for the coordinates of the intersecting point (in double
precision) computed by Maxima are:
\begin{equation*}
\mathrm{P}=[-1.656701299244927,1.210755779027779]\,,
\end{equation*}
confirming that \eqref{eq2} is a good approximation to the contact point.

Once the lines representing the surface and the curves on it have been determined, the next
step is to hide those portions that lie behind the surface from the perspective of the observer. However, it takes a long time to apply the algorithm for the elimination of those hidden parts if only CindyScript is used. To speed up computations, \ketcindy\ can call
\texttt{gcc}(Gnu C Compiler). The following is the CindyScript code used to generate 
Figure\ \addtocounter{figure}{-3}\thefigure\addtocounter{figure}{3}\:

\begin{verbatim}
    fd=[
      "z=4-(x^2+y^2)",
      "x=R*cos(T)","y=R*sin(T)",
      "R=[0,2]","T=[0,2*pi]","e"
    ];
    Startsurf();
    Sfbdparadata("1",fd);
    Crvsfparadata("1","ax3d","sfbd3d1",fd);
    Wireparadata("1","sfbd3d1",fd,5,2*pi/6*(0..5));
    ExeccmdC("1",[""],[]);
\end{verbatim}

\noindent
Its execution takes only a few seconds in a standard desktop computer.

\subsection{3D animations}

\ketcindy\ can produce both, a flip movie which displays slides step by step, and \LaTeX\ animations suitable to be included in \texttt{pdf} documents (or translated into graphics formats such as \texttt{APNG}), through the use of \verb|animate.sty|. The use of the 
\texttt{gcc} compiler makes it feasible to animate 3D surfaces. We show an example of an
animation illustrating the construction of the M\"{o}bius band. A snapshot of the
corresponding CindyScreen appears in Figure \thefigure, and the CindyScript code is
the following:

\begin{center}
\includegraphics[scale=0.45]{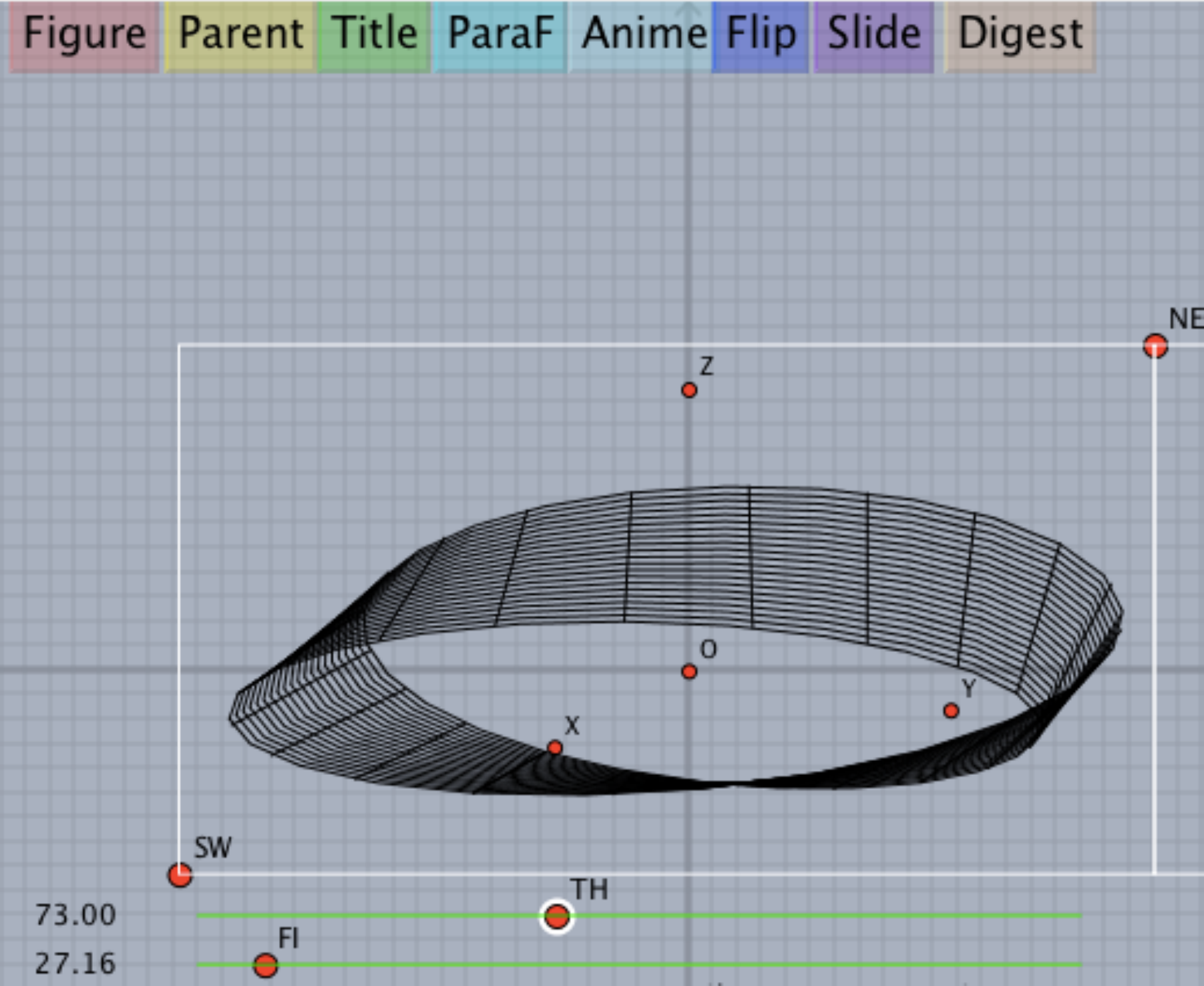} 
\\
{\bf Fig.\thefigure}\ Screen for 3D movies\addtocounter{figure}{1}
\end{center}

\begin{verbatim}
    Start3d([A,B,C,S,Sl,Sr]);
    mf(s):=(
      regional(tmp);
      Startsurf();
      fd=["p",
       "x=2*cos(t)*(2+r*cos(t/2))",
       "y=2*sin(t)*(2+r*cos(t/2))",
       "z=2*r*sin(t/2)","r=[-0.4,0.4]","t=[0,"+text(s)+"]","nsew"
      ];
      Sfbdparadata("1",fd);
      tmp=select((1..12)/12*2*pi,#<=s);
      Wireparadata("1","sfbd3d1",fd,5,tmp);
      ExeccmdC("",[""],["nodisp"]);
    );
    Setpara("mobius","mf(s)","s=[0,2*pi]",["Div=24"]);
\end{verbatim}

\noindent
Pressing the button \verb|Flip| generates the flip movie, while \verb|Anime| generates the animation. The separate slides of the flip movie are shown below.

\begin{center}
\fbox{\scalebox{0.4}{\input{fig/p001.tex}}}
\fbox{\scalebox{0.4}{\input{fig/p004.tex}}}
\fbox{\scalebox{0.4}{\input{fig/p007.tex}}}

\fbox{\scalebox{0.4}{\input{fig/p010.tex}}}
\fbox{\scalebox{0.4}{\input{fig/p013.tex}}}
\fbox{\scalebox{0.4}{\input{fig/p016.tex}}}

\fbox{\scalebox{0.4}{\input{fig/p019.tex}}}
\fbox{\scalebox{0.4}{\input{fig/p022.tex}}}
\fbox{\scalebox{0.4}{\input{fig/p025.tex}}}

\vspace{2mm}

{\bf Fig.\thefigure}\ Slides for the animated M\"{o}bius band\addtocounter{figure}{1}

\end{center}

\section{Conclusions and future work}

It is quite common that teachers, and researchers as well, need added functionalities
in their software of choice for producing high-quality graphical output or accurate
numerical computations. It is very important that these functionalities do not come
at the cost of a significant increase in the time required to learn their use.
CindyScript is a scripting language that allows to build them easily and, moreover,
one that can interact with other software components such as Maxima and the \texttt{gcc}
compiler. This opens a whole new world of possibilities. For example, consider how easy it
is to write the following code to interactively compute the area surrounded by a closed curve:

\begin{verbatim}
    Findarea(pd):=(
     regional(p0,p1,p2,p3,s);
     if(isstring(pdstr),pd=parse(pdstr),pd=pdstr);
     s=0;
     forall(1..(length(pd)-1),
       p1=pd_#;
       p2=pd_(#+1);
       if(#==1,p0=pd_(length(pd)-1),p0=pd_(#-1));
       if(#==length(pd)-1,p3=pd_2,p3=pd_(#+2));
       s=s+IntegrateO(p0,p1,p2,p3);
    );
\end{verbatim}

\noindent
It is then possible to put a slider on the screen, and the results obtained by sliding the
control point are shown below.\\

\begin{center}
\includegraphics[scale=0.4]{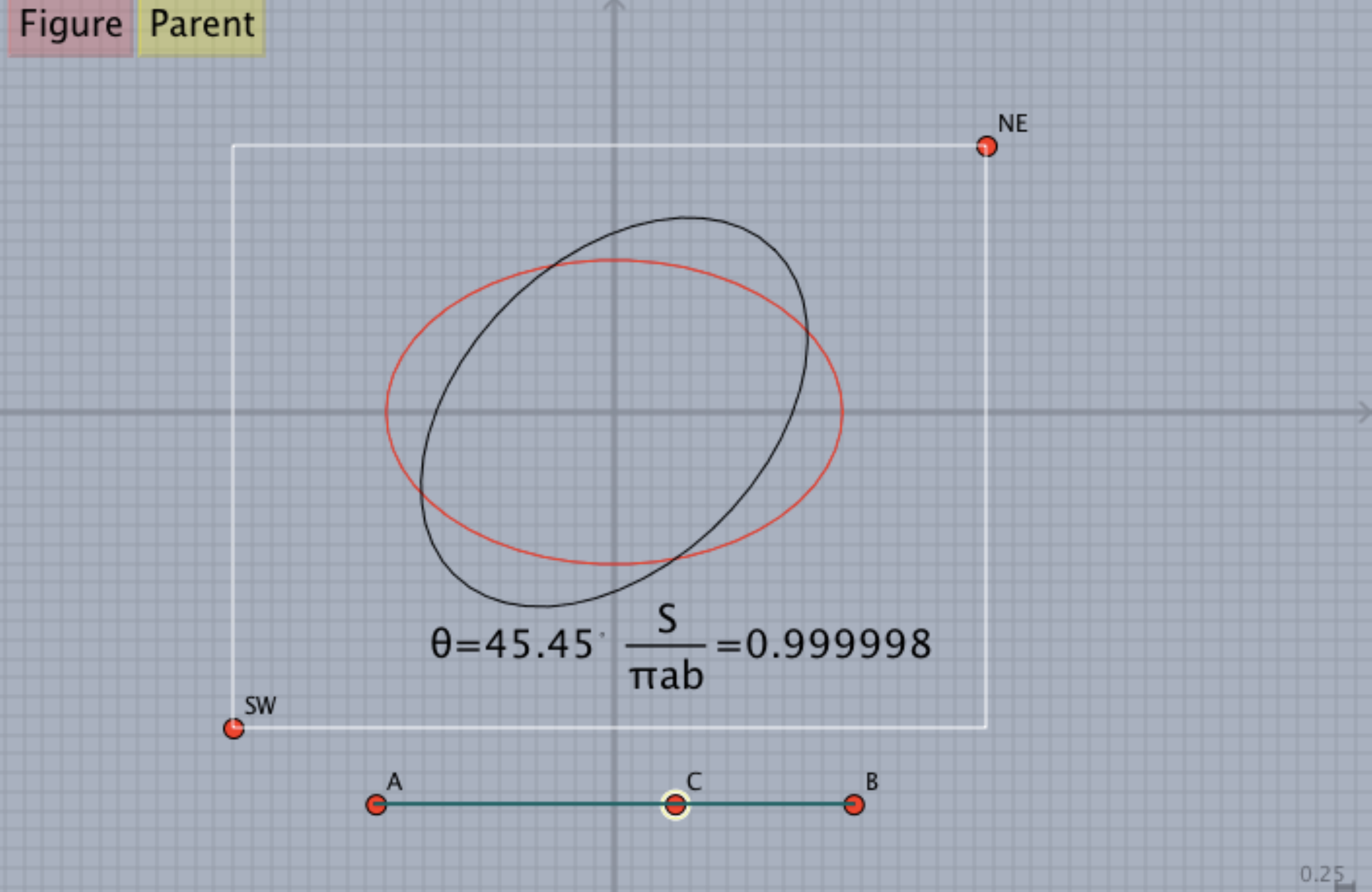}
\\
{\bf Fig.\thefigure}\ Adding a command with CindyScript \addtocounter{figure}{1}\\
\end{center}

Programming is an important factor when creating more appealing materials, and DGSs such as Cinderella make it easy to visualize the output and modify it as necessary.
Free CASs such as Maxima are also useful for symbolic computations, and it is remarkable
that \ketcindy\ can work as their user interface. Moreover, C compiler is very efficient
in speeding up the calculations required for drawing 3D figures with \ketcindy. 
We could say that the combined use of \ketcindy, Cinderella, Maxima and C 
is a powerful tool to develop programs. Finally, let us comment that \ketcindy\ also has built-in commands to generate files in \texttt{obj} format, suitable to print 3D models.
As a future work, we will develop a C library to speed up the process.


\subsection*{Acknowledgments}
This work was supported by JSPS KAKENHI Grant Number 16K0115,18K02948,18K02872.


\subsection*{Acknowledgment}
This work was supported by JSPS KAKENHI Grant Numbers 
16K01152, 18K02948, 18K02872.

\end{document}

%% file: fig/crsplineall.tex
{\unitlength=1cm%
\begin{picture}%
(6,6)(-1,-1)%
\linethickness{0.008in}
\polyline(-1.00000,0.55998)(-0.74163,0.59811)(-0.48660,0.64809)(-0.23725,0.70986)%
(0.00869,0.78405)(0.25349,0.87128)(0.49944,0.97217)(0.74880,1.08735)(1.00386,1.21745)%
(1.26644,1.36982)(1.53521,1.54791)(1.80767,1.74646)(2.08135,1.96022)(2.35375,2.18393)%
(2.62239,2.41235)(2.88479,2.64020)(3.13845,2.86225)(3.38090,3.07323)(3.60964,3.26790)%
(3.82103,3.44544)(4.01561,3.61085)(4.19720,3.76760)(4.36962,3.91914)(4.53668,4.06892)%
(4.70219,4.22041)(4.86998,4.37706)(5.00000,4.50064)%
\polyline(1.00386,1.21745)(3.60964,1.21745)(3.60964,3.26790)(1.00386,1.21745)%
\settowidth{\Width}{$\Delta x$}\setlength{\Width}{-0.5\Width}%
\settoheight{\Height}{$\Delta x$}\settodepth{\Depth}{$\Delta x$}\setlength{\Height}{-0.5\Height}\setlength{\Depth}{0.5\Depth}\addtolength{\Height}{\Depth}%
\put(2.3100000,0.9600000){\hspace*{\Width}\raisebox{\Height}{$\Delta x$}}%
\put(1.00386,1.21745){\circle*{0.040640}}\put(1.09587,1.18074){\circle*{0.040640}}
\put(1.18890,1.14672){\circle*{0.040640}}\put(1.28288,1.11542){\circle*{0.040640}}
\put(1.37774,1.08686){\circle*{0.040640}}\put(1.47338,1.06107){\circle*{0.040640}}
\put(1.56973,1.03807){\circle*{0.040640}}\put(1.66671,1.01790){\circle*{0.040640}}
\put(1.76426,1.00066){\circle*{0.040640}}\put(1.86226,0.98626){\circle*{0.040640}}
\put(1.96065,0.97471){\circle*{0.040640}}\put(2.05932,0.96602){\circle*{0.040640}}
\put(2.15821,0.96020){\circle*{0.040640}}\put(2.25722,0.95726){\circle*{0.040640}}
\put(2.35628,0.95726){\circle*{0.040640}}\put(2.45530,0.96020){\circle*{0.040640}}
\put(2.55418,0.96602){\circle*{0.040640}}\put(2.65286,0.97471){\circle*{0.040640}}
\put(2.75124,0.98626){\circle*{0.040640}}\put(2.84925,1.00066){\circle*{0.040640}}
\put(2.94679,1.01790){\circle*{0.040640}}\put(3.04378,1.03807){\circle*{0.040640}}
\put(3.14013,1.06107){\circle*{0.040640}}\put(3.23577,1.08686){\circle*{0.040640}}
\put(3.33062,1.11542){\circle*{0.040640}}\put(3.42460,1.14672){\circle*{0.040640}}
\put(3.51763,1.18074){\circle*{0.040640}}\put(3.60964,1.21745){\circle*{0.040640}}
\linethickness{0.008in}
\settowidth{\Width}{$\Delta y$}\setlength{\Width}{-0.5\Width}%
\settoheight{\Height}{$\Delta y$}\settodepth{\Depth}{$\Delta y$}\setlength{\Height}{-0.5\Height}\setlength{\Depth}{0.5\Depth}\addtolength{\Height}{\Depth}%
\put(3.9200000,2.2400000){\hspace*{\Width}\raisebox{\Height}{$\Delta y$}}%
\put(3.60964,1.21745){\circle*{0.040640}}\put(3.66167,1.30130){\circle*{0.040640}}
\put(3.70926,1.38775){\circle*{0.040640}}\put(3.75226,1.47657){\circle*{0.040640}}
\put(3.79049,1.56754){\circle*{0.040640}}\put(3.82374,1.66045){\circle*{0.040640}}
\put(3.85210,1.75497){\circle*{0.040640}}\put(3.87548,1.85084){\circle*{0.040640}}
\put(3.89365,1.94784){\circle*{0.040640}}\put(3.90663,2.04566){\circle*{0.040640}}
\put(3.91449,2.14403){\circle*{0.040640}}\put(3.91721,2.24267){\circle*{0.040640}}
\put(3.91449,2.34132){\circle*{0.040640}}\put(3.90663,2.43969){\circle*{0.040640}}
\put(3.89365,2.53751){\circle*{0.040640}}\put(3.87548,2.63450){\circle*{0.040640}}
\put(3.85210,2.73038){\circle*{0.040640}}\put(3.82374,2.82489){\circle*{0.040640}}
\put(3.79049,2.91781){\circle*{0.040640}}\put(3.75226,3.00878){\circle*{0.040640}}
\put(3.70926,3.09760){\circle*{0.040640}}\put(3.66167,3.18405){\circle*{0.040640}}
\put(3.60964,3.26790){\circle*{0.040640}}
\linethickness{0.008in}
\settowidth{\Width}{P}\setlength{\Width}{-0.5\Width}%
\settoheight{\Height}{P}\settodepth{\Depth}{P}\setlength{\Height}{\Depth}%
\put(1.0000000,1.2700000){\hspace*{\Width}\raisebox{\Height}{P}}%
\polyline(1.00386,0.05000)(1.00386,-0.05000)%
\settowidth{\Width}{$x$}\setlength{\Width}{-0.5\Width}%
\settoheight{\Height}{$x$}\settodepth{\Depth}{$x$}\setlength{\Height}{-\Height}%
\put(1.0000000,-0.1000000){\hspace*{\Width}\raisebox{\Height}{$x$}}%
\polyline(3.60964,0.05000)(3.60964,-0.05000)%
\settowidth{\Width}{$\xi$}\setlength{\Width}{-0.5\Width}%
\settoheight{\Height}{$\xi$}\settodepth{\Depth}{$\xi$}\setlength{\Height}{-\Height}%
\put(3.6100000,-0.1000000){\hspace*{\Width}\raisebox{\Height}{$\xi$}}%
\polyline(0.05000,1.21745)(-0.05000,1.21745)%
\settowidth{\Width}{$f(x)$}\setlength{\Width}{-1\Width}%
\settoheight{\Height}{$f(x)$}\settodepth{\Depth}{$f(x)$}\setlength{\Height}{-0.5\Height}\setlength{\Depth}{0.5\Depth}\addtolength{\Height}{\Depth}%
\put(-0.1000000,1.2200000){\hspace*{\Width}\raisebox{\Height}{$f(x)$}}%
\polyline(0.05000,3.26790)(-0.05000,3.26790)%
\settowidth{\Width}{$f(\xi)$}\setlength{\Width}{-1\Width}%
\settoheight{\Height}{$f(\xi)$}\settodepth{\Depth}{$f(\xi)$}\setlength{\Height}{-0.5\Height}\setlength{\Depth}{0.5\Depth}\addtolength{\Height}{\Depth}%
\put(-0.1000000,3.2700000){\hspace*{\Width}\raisebox{\Height}{$f(\xi)$}}%
\polyline(-1.00000,0.00000)(5.00000,0.00000)%
\polyline(0.00000,-1.00000)(0.00000,5.00000)%
\settowidth{\Width}{$x$}\setlength{\Width}{0\Width}%
\settoheight{\Height}{$x$}\settodepth{\Depth}{$x$}\setlength{\Height}{-0.5\Height}\setlength{\Depth}{0.5\Depth}\addtolength{\Height}{\Depth}%
\put(5.0500000,0.0000000){\hspace*{\Width}\raisebox{\Height}{$x$}}%
\settowidth{\Width}{$y$}\setlength{\Width}{-0.5\Width}%
\settoheight{\Height}{$y$}\settodepth{\Depth}{$y$}\setlength{\Height}{\Depth}%
\put(0.0000000,5.0500000){\hspace*{\Width}\raisebox{\Height}{$y$}}%
\settowidth{\Width}{O}\setlength{\Width}{-1\Width}%
\settoheight{\Height}{O}\settodepth{\Depth}{O}\setlength{\Height}{-\Height}%
\put(-0.0500000,-0.0500000){\hspace*{\Width}\raisebox{\Height}{O}}%
\end{picture}}%

%% file: fig/crsplineclosed.tex
{\unitlength=1cm%
\begin{picture}%
(7,5)(-3.5,-2.5)%
\linethickness{0.008in}
\polyline(3.00000,0.00000)(2.99069,0.09866)\polyline(2.98138,0.19731)(2.97634,0.25067)(2.96383,0.29441)%
\polyline(2.93657,0.38968)(2.90931,0.48495)\polyline(2.86800,0.57484)(2.82458,0.66391)%
\polyline(2.77859,0.75147)(2.72144,0.83242)\polyline(2.66430,0.91338)(2.62892,0.96351)(2.60290,0.99084)%
\polyline(2.53458,1.06262)(2.46625,1.13439)\polyline(2.39416,1.20207)(2.31701,1.26425)%
\polyline(2.23985,1.32643)(2.18691,1.36909)(2.16054,1.38560)\polyline(2.07655,1.43818)(1.99256,1.49076)%
\polyline(1.90834,1.54293)(1.81916,1.58613)\polyline(1.72997,1.62932)(1.64079,1.67252)%
\polyline(1.54919,1.71002)(1.45615,1.74412)\polyline(1.36311,1.77822)(1.27734,1.80965)(1.26985,1.81163)%
\polyline(1.17404,1.83692)(1.07822,1.86221)\polyline(0.98241,1.88750)(0.92705,1.90211)(0.88581,1.90917)%
\polyline(0.78814,1.92589)(0.69047,1.94261)\polyline(0.59279,1.95933)(0.56214,1.96457)(0.49439,1.97028)%
\polyline(0.39564,1.97860)(0.29690,1.98691)\polyline(0.19815,1.99523)(0.18837,1.99605)(0.09909,1.99605)%
\polyline(0.00000,1.99605)(-0.09909,1.99605)\polyline(-0.19815,1.99523)(-0.29690,1.98691)%
\polyline(-0.39564,1.97860)(-0.49439,1.97028)\polyline(-0.59279,1.95933)(-0.69047,1.94261)%
\polyline(-0.78814,1.92589)(-0.88581,1.90917)\polyline(-0.98241,1.88750)(-1.07822,1.86221)%
\polyline(-1.17404,1.83692)(-1.26985,1.81163)\polyline(-1.36311,1.77822)(-1.45615,1.74412)%
\polyline(-1.54919,1.71002)(-1.60748,1.68866)(-1.64079,1.67252)\polyline(-1.72997,1.62932)(-1.81916,1.58613)%
\polyline(-1.90834,1.54293)(-1.91227,1.54103)(-1.99256,1.49076)\polyline(-2.07655,1.43818)(-2.16054,1.38560)%
\polyline(-2.23985,1.32643)(-2.31701,1.26425)\polyline(-2.39416,1.20207)(-2.42705,1.17557)(-2.46625,1.13439)%
\polyline(-2.53458,1.06262)(-2.60290,0.99084)\polyline(-2.66430,0.91338)(-2.72144,0.83242)%
\polyline(-2.77859,0.75147)(-2.78933,0.73625)(-2.82458,0.66391)\polyline(-2.86800,0.57484)(-2.90575,0.49738)(-2.90931,0.48495)%
\polyline(-2.93657,0.38968)(-2.96383,0.29441)\polyline(-2.98138,0.19731)(-2.99069,0.09866)%
\polyline(-3.00000,0.00000)(-2.99069,-0.09866)\polyline(-2.98138,-0.19731)(-2.97634,-0.25067)(-2.96383,-0.29441)%
\polyline(-2.93657,-0.38968)(-2.90931,-0.48495)\polyline(-2.86800,-0.57484)(-2.82458,-0.66391)%
\polyline(-2.77859,-0.75147)(-2.72144,-0.83242)\polyline(-2.66430,-0.91338)(-2.62892,-0.96351)(-2.60290,-0.99084)%
\polyline(-2.53458,-1.06262)(-2.46625,-1.13439)\polyline(-2.39416,-1.20207)(-2.31701,-1.26425)%
\polyline(-2.23985,-1.32643)(-2.18691,-1.36909)(-2.16054,-1.38560)\polyline(-2.07655,-1.43818)(-1.99256,-1.49076)%
\polyline(-1.90834,-1.54293)(-1.81916,-1.58613)\polyline(-1.72997,-1.62932)(-1.64079,-1.67252)%
\polyline(-1.54919,-1.71002)(-1.45615,-1.74412)\polyline(-1.36311,-1.77822)(-1.27734,-1.80965)(-1.26985,-1.81163)%
\polyline(-1.17404,-1.83692)(-1.07822,-1.86221)\polyline(-0.98241,-1.88750)(-0.92705,-1.90211)(-0.88581,-1.90917)%
\polyline(-0.78814,-1.92589)(-0.69047,-1.94261)\polyline(-0.59279,-1.95933)(-0.56214,-1.96457)(-0.49439,-1.97028)%
\polyline(-0.39564,-1.97860)(-0.29690,-1.98691)\polyline(-0.19815,-1.99523)(-0.18837,-1.99605)(-0.09909,-1.99605)%
\polyline(-0.00000,-1.99605)(0.09909,-1.99605)\polyline(0.19815,-1.99523)(0.29690,-1.98691)%
\polyline(0.39564,-1.97860)(0.49439,-1.97028)\polyline(0.59279,-1.95933)(0.69047,-1.94261)%
\polyline(0.78814,-1.92589)(0.88581,-1.90917)\polyline(0.98241,-1.88750)(1.07822,-1.86221)%
\polyline(1.17404,-1.83692)(1.26985,-1.81163)\polyline(1.36311,-1.77822)(1.45615,-1.74412)%
\polyline(1.54919,-1.71002)(1.60748,-1.68866)(1.64079,-1.67252)\polyline(1.72997,-1.62932)(1.81916,-1.58613)%
\polyline(1.90834,-1.54293)(1.91227,-1.54103)(1.99256,-1.49076)\polyline(2.07655,-1.43818)(2.16054,-1.38560)%
\polyline(2.23985,-1.32643)(2.31701,-1.26425)\polyline(2.39416,-1.20207)(2.42705,-1.17557)(2.46625,-1.13439)%
\polyline(2.53458,-1.06262)(2.60290,-0.99084)\polyline(2.66430,-0.91338)(2.72144,-0.83242)%
\polyline(2.77859,-0.75147)(2.78933,-0.73625)(2.82458,-0.66391)\polyline(2.86800,-0.57484)(2.90575,-0.49738)(2.90931,-0.48495)%
\polyline(2.93657,-0.38968)(2.96383,-0.29441)\polyline(2.98138,-0.19731)(2.99069,-0.09866)%
\polyline(3.00000,0.00000)(2.94300,0.21800)(2.78400,0.46400)(2.54100,0.72600)(2.23200,0.99200)%
(1.87500,1.25000)(1.48800,1.48800)(1.08900,1.69400)(0.69600,1.85600)(0.32700,1.96200)%
(0.00000,2.00000)(-0.32700,1.96200)(-0.69600,1.85600)(-1.08900,1.69400)(-1.48800,1.48800)%
(-1.87500,1.25000)(-2.23200,0.99200)(-2.54100,0.72600)(-2.78400,0.46400)(-2.94300,0.21800)%
(-3.00000,0.00000)(-2.94300,-0.21800)(-2.78400,-0.46400)(-2.54100,-0.72600)(-2.23200,-0.99200)%
(-1.87500,-1.25000)(-1.48800,-1.48800)(-1.08900,-1.69400)(-0.69600,-1.85600)(-0.32700,-1.96200)%
(0.00000,-2.00000)(0.32700,-1.96200)(0.69600,-1.85600)(1.08900,-1.69400)(1.48800,-1.48800)%
(1.87500,-1.25000)(2.23200,-0.99200)(2.54100,-0.72600)(2.78400,-0.46400)(2.94300,-0.21800)%
(3.00000,0.00000)%
\linethickness{0.003in}
\put(3.00000,0.66667){\circle*{0.120000}}\put(1.00000,2.00000){\circle*{0.120000}}
\linethickness{0.008in}
\linethickness{0.003in}
\put(-1.00000,2.00000){\circle*{0.120000}}\put(-3.00000,0.66667){\circle*{0.120000}}
\linethickness{0.008in}
\linethickness{0.003in}
\put(-3.00000,-0.66667){\circle*{0.120000}}\put(-1.00000,-2.00000){\circle*{0.120000}}
\linethickness{0.008in}
\linethickness{0.003in}
\put(1.00000,-2.00000){\circle*{0.120000}}\put(3.00000,-0.66667){\circle*{0.120000}}
\linethickness{0.008in}
\linethickness{0.003in}
\put(3.00000,0.00000){\circle*{0.120000}}\put(0.00000,2.00000){\circle*{0.120000}}
\put(-3.00000,0.00000){\circle*{0.120000}}\put(0.00000,-2.00000){\circle*{0.120000}}
\linethickness{0.008in}
\settowidth{\Width}{A}\setlength{\Width}{0\Width}%
\settoheight{\Height}{A}\settodepth{\Depth}{A}\setlength{\Height}{\Depth}%
\put(3.0500000,0.0500000){\hspace*{\Width}\raisebox{\Height}{A}}%
\settowidth{\Width}{B}\setlength{\Width}{0\Width}%
\settoheight{\Height}{B}\settodepth{\Depth}{B}\setlength{\Height}{\Depth}%
\put(0.0500000,2.0500000){\hspace*{\Width}\raisebox{\Height}{B}}%
\settowidth{\Width}{C}\setlength{\Width}{-1\Width}%
\settoheight{\Height}{C}\settodepth{\Depth}{C}\setlength{\Height}{\Depth}%
\put(-3.0500000,0.0500000){\hspace*{\Width}\raisebox{\Height}{C}}%
\settowidth{\Width}{D}\setlength{\Width}{0\Width}%
\settoheight{\Height}{D}\settodepth{\Depth}{D}\setlength{\Height}{-\Height}%
\put(0.0500000,-2.0500000){\hspace*{\Width}\raisebox{\Height}{D}}%
\polyline(-3.50000,0.00000)(3.50000,0.00000)%
\polyline(0.00000,-2.50000)(0.00000,2.50000)%
\settowidth{\Width}{$x$}\setlength{\Width}{0\Width}%
\settoheight{\Height}{$x$}\settodepth{\Depth}{$x$}\setlength{\Height}{-0.5\Height}\setlength{\Depth}{0.5\Depth}\addtolength{\Height}{\Depth}%
\put(3.5500000,0.0000000){\hspace*{\Width}\raisebox{\Height}{$x$}}%
\settowidth{\Width}{$y$}\setlength{\Width}{-0.5\Width}%
\settoheight{\Height}{$y$}\settodepth{\Depth}{$y$}\setlength{\Height}{\Depth}%
\put(0.0000000,2.5500000){\hspace*{\Width}\raisebox{\Height}{$y$}}%
\settowidth{\Width}{O}\setlength{\Width}{-1\Width}%
\settoheight{\Height}{O}\settodepth{\Depth}{O}\setlength{\Height}{-\Height}%
\put(-0.0500000,-0.0500000){\hspace*{\Width}\raisebox{\Height}{O}}%
\end{picture}}%

%% file: fig/crospline1.tex
{\unitlength=1cm%
\begin{picture}%
(7,5)(-3.5,-2.5)%
\linethickness{0.008in}
\linethickness{0.003in}
\put(-3.27911,-1.83186){\circle*{0.120000}}\put(-1.81922,0.02107){\circle*{0.120000}}
\put(1.28942,1.15938){\circle*{0.120000}}\put(3.14406,0.41072){\circle*{0.120000}}
\put(1.75564,-1.08320){\circle*{0.120000}}

\linethickness{0.008in}
\polyline(-3.27911,-1.83186)(-3.19605,-1.77747)\polyline(-3.11298,-1.72309)(-3.02992,-1.66870)%
\polyline(-2.94685,-1.61432)(-2.86379,-1.55993)\polyline(-2.78072,-1.50554)(-2.69766,-1.45116)%
\polyline(-2.61460,-1.39677)(-2.53153,-1.34238)\polyline(-2.44847,-1.28800)(-2.36540,-1.23361)%
\polyline(-2.28234,-1.17923)(-2.19928,-1.12484)\polyline(-2.11621,-1.07045)(-2.03315,-1.01607)%
\polyline(-1.95008,-0.96168)(-1.86702,-0.90729)\polyline(-1.78395,-0.85291)(-1.70089,-0.79852)%
\polyline(-1.61783,-0.74414)(-1.53476,-0.68975)\polyline(-1.45170,-0.63536)(-1.36863,-0.58098)%
\polyline(-1.28557,-0.52659)(-1.20251,-0.47221)\polyline(-1.11944,-0.41782)(-1.03638,-0.36343)%
\polyline(-0.95331,-0.30905)(-0.87025,-0.25466)\polyline(-0.78718,-0.20027)(-0.70412,-0.14589)%
\polyline(-0.62106,-0.09150)(-0.53799,-0.03712)\polyline(-0.45493,0.01727)(-0.37186,0.07166)%
\polyline(-0.28880,0.12604)(-0.20574,0.18043)\polyline(-0.12267,0.23481)(-0.03961,0.28920)%
\polyline(0.04346,0.34359)(0.12652,0.39797)\polyline(0.20959,0.45236)(0.29265,0.50675)%
\polyline(0.37571,0.56113)(0.45878,0.61552)\polyline(0.54184,0.66990)(0.62491,0.72429)%
\polyline(0.70797,0.77868)(0.79103,0.83306)\polyline(0.87410,0.88745)(0.95716,0.94184)%
\polyline(1.04023,0.99622)(1.12329,1.05061)\polyline(1.20636,1.10499)(1.28942,1.15938)%
\polyline(-1.81922,0.02107)(-1.72190,0.02871)\polyline(-1.62458,0.03635)(-1.52726,0.04399)%
\polyline(-1.42994,0.05163)(-1.33262,0.05927)\polyline(-1.23530,0.06691)(-1.13799,0.07455)%
\polyline(-1.04067,0.08219)(-0.94335,0.08983)\polyline(-0.84603,0.09747)(-0.74871,0.10511)%
\polyline(-0.65139,0.11275)(-0.55407,0.12039)\polyline(-0.45675,0.12803)(-0.35943,0.13567)%
\polyline(-0.26211,0.14331)(-0.16479,0.15095)\polyline(-0.06747,0.15859)(0.02985,0.16623)%
\polyline(0.12716,0.17387)(0.22448,0.18151)\polyline(0.32180,0.18915)(0.41912,0.19679)%
\polyline(0.51644,0.20443)(0.61376,0.21207)\polyline(0.71108,0.21972)(0.80840,0.22736)%
\polyline(0.90572,0.23500)(1.00304,0.24264)\polyline(1.10036,0.25028)(1.19768,0.25792)%
\polyline(1.29499,0.26556)(1.39231,0.27320)\polyline(1.48963,0.28084)(1.58695,0.28848)%
\polyline(1.68427,0.29612)(1.78159,0.30376)\polyline(1.87891,0.31140)(1.97623,0.31904)%
\polyline(2.07355,0.32668)(2.17087,0.33432)\polyline(2.26819,0.34196)(2.36551,0.34960)%
\polyline(2.46283,0.35724)(2.56014,0.36488)\polyline(2.65746,0.37252)(2.75478,0.38016)%
\polyline(2.85210,0.38780)(2.94942,0.39544)\polyline(3.04674,0.40308)(3.14406,0.41072)%
\polyline(1.28942,1.15938)(1.30969,1.06188)\polyline(1.32996,0.96437)(1.35023,0.86687)%
\polyline(1.37050,0.76937)(1.39077,0.67186)\polyline(1.41104,0.57436)(1.43131,0.47686)%
\polyline(1.45158,0.37935)(1.47185,0.28185)\polyline(1.49212,0.18435)(1.51239,0.08684)%
\polyline(1.53267,-0.01066)(1.55294,-0.10817)\polyline(1.57321,-0.20567)(1.59348,-0.30317)%
\polyline(1.61375,-0.40068)(1.63402,-0.49818)\polyline(1.65429,-0.59568)(1.67456,-0.69319)%
\polyline(1.69483,-0.79069)(1.71510,-0.88819)\polyline(1.73537,-0.98570)(1.75564,-1.08320)%
\settowidth{\Width}{$\mathrm{P}_{j-1}$}\setlength{\Width}{-1\Width}%
\settoheight{\Height}{$\mathrm{P}_{j-1}$}\settodepth{\Depth}{$\mathrm{P}_{j-1}$}\setlength{\Height}{\Depth}%
\put(-3.3300000,-1.7800000){\hspace*{\Width}\raisebox{\Height}{$\mathrm{P}_{j-1}$}}%
\settowidth{\Width}{$\mathrm{P}_{j}$}\setlength{\Width}{-1\Width}%
\settoheight{\Height}{$\mathrm{P}_{j}$}\settodepth{\Depth}{$\mathrm{P}_{j}$}\setlength{\Height}{\Depth}%
\put(-1.8700000,0.0700000){\hspace*{\Width}\raisebox{\Height}{$\mathrm{P}_{j}$}}%
\settowidth{\Width}{$\mathrm{P}_{j+1}$}\setlength{\Width}{-0.5\Width}%
\settoheight{\Height}{$\mathrm{P}_{j+1}$}\settodepth{\Depth}{$\mathrm{P}_{j+1}$}\setlength{\Height}{\Depth}%
\put(1.2900000,1.2100000){\hspace*{\Width}\raisebox{\Height}{$\mathrm{P}_{j+1}$}}%
\settowidth{\Width}{$\mathrm{P}_{j+2}$}\setlength{\Width}{0\Width}%
\settoheight{\Height}{$\mathrm{P}_{j+2}$}\settodepth{\Depth}{$\mathrm{P}_{j+2}$}\setlength{\Height}{-0.5\Height}\setlength{\Depth}{0.5\Depth}\addtolength{\Height}{\Depth}%
\put(3.1900000,0.4100000){\hspace*{\Width}\raisebox{\Height}{$\mathrm{P}_{j+2}$}}%
\settowidth{\Width}{$\mathrm{P}_{j+3}$}\setlength{\Width}{-0.5\Width}%
\settoheight{\Height}{$\mathrm{P}_{j+3}$}\settodepth{\Depth}{$\mathrm{P}_{j+3}$}\setlength{\Height}{-\Height}%
\put(1.7600000,-1.1300000){\hspace*{\Width}\raisebox{\Height}{$\mathrm{P}_{j+3}$}}%
\polyline(-3.27911,-1.83186)(-3.18430,-1.58520)(-3.08267,-1.35795)(-2.97250,-1.14794)%
(-2.85210,-0.95301)(-2.71975,-0.77100)(-2.57375,-0.59973)(-2.41241,-0.43706)(-2.23401,-0.28081)%
(-2.03685,-0.12882)(-1.81922,0.02107)(-1.56949,0.17233)(-1.28295,0.32466)(-0.96831,0.47453)%
(-0.63428,0.61842)(-0.28957,0.75283)(0.05711,0.87422)(0.39706,0.97909)(0.72157,1.06392)%
(1.02192,1.12519)(1.28942,1.15938)(1.54027,1.16429)(1.79250,1.14234)(2.04014,1.09695)%
(2.27723,1.03155)(2.49781,0.94956)(2.69590,0.85442)(2.86554,0.74953)(3.00075,0.63834)%
(3.09558,0.52426)(3.14406,0.41072)(3.13428,0.27671)(3.06583,0.10701)(2.94995,-0.08628)%
(2.79789,-0.29105)(2.62090,-0.49521)(2.43023,-0.68666)(2.23712,-0.85328)(2.05282,-0.98298)%
(1.88857,-1.06365)(1.75564,-1.08320)%
\linethickness{0.003in}
\put(-1.05780,0.51961){\circle*{0.120000}}\put(0.46221,1.09444){\circle*{0.120000}}
\linethickness{0.008in}
\linethickness{0.003in}
\put(2.11664,1.22432){\circle*{0.120000}}\put(3.06636,0.78448){\circle*{0.120000}}
\linethickness{0.008in}
\settowidth{\Width}{$\mathrm{Q}_{j}$}\setlength{\Width}{-0.5\Width}%
\settoheight{\Height}{$\mathrm{Q}_{j}$}\settodepth{\Depth}{$\mathrm{Q}_{j}$}\setlength{\Height}{\Depth}%
\put(-1.0600000,0.5700000){\hspace*{\Width}\raisebox{\Height}{$\mathrm{Q}_{j}$}}%
\settowidth{\Width}{$\mathrm{R}_{j}$}\setlength{\Width}{-0.5\Width}%
\settoheight{\Height}{$\mathrm{R}_{j}$}\settodepth{\Depth}{$\mathrm{R}_{j}$}\setlength{\Height}{\Depth}%
\put(0.4600000,1.1400000){\hspace*{\Width}\raisebox{\Height}{$\mathrm{R}_{j}$}}%
\end{picture}}%

%% file: fig/crospline2.tex
{\unitlength=1cm%
\begin{picture}%
(7,5)(-3.5,-2.5)%
\linethickness{0.008in}
\linethickness{0.003in}
\put(-3.27911,-1.83186){\circle*{0.120000}}\put(-1.81922,0.02107){\circle*{0.120000}}
\put(1.28942,1.15938){\circle*{0.120000}}\put(3.14406,0.41072){\circle*{0.120000}}
\put(1.75564,-1.08320){\circle*{0.120000}}

\linethickness{0.008in}
\polyline(-3.27911,-1.83186)(-3.19605,-1.77747)\polyline(-3.11298,-1.72309)(-3.02992,-1.66870)%
\polyline(-2.94685,-1.61432)(-2.86379,-1.55993)\polyline(-2.78072,-1.50554)(-2.69766,-1.45116)%
\polyline(-2.61460,-1.39677)(-2.53153,-1.34238)\polyline(-2.44847,-1.28800)(-2.36540,-1.23361)%
\polyline(-2.28234,-1.17923)(-2.19928,-1.12484)\polyline(-2.11621,-1.07045)(-2.03315,-1.01607)%
\polyline(-1.95008,-0.96168)(-1.86702,-0.90729)\polyline(-1.78395,-0.85291)(-1.70089,-0.79852)%
\polyline(-1.61783,-0.74414)(-1.53476,-0.68975)\polyline(-1.45170,-0.63536)(-1.36863,-0.58098)%
\polyline(-1.28557,-0.52659)(-1.20251,-0.47221)\polyline(-1.11944,-0.41782)(-1.03638,-0.36343)%
\polyline(-0.95331,-0.30905)(-0.87025,-0.25466)\polyline(-0.78718,-0.20027)(-0.70412,-0.14589)%
\polyline(-0.62106,-0.09150)(-0.53799,-0.03712)\polyline(-0.45493,0.01727)(-0.37186,0.07166)%
\polyline(-0.28880,0.12604)(-0.20574,0.18043)\polyline(-0.12267,0.23481)(-0.03961,0.28920)%
\polyline(0.04346,0.34359)(0.12652,0.39797)\polyline(0.20959,0.45236)(0.29265,0.50675)%
\polyline(0.37571,0.56113)(0.45878,0.61552)\polyline(0.54184,0.66990)(0.62491,0.72429)%
\polyline(0.70797,0.77868)(0.79103,0.83306)\polyline(0.87410,0.88745)(0.95716,0.94184)%
\polyline(1.04023,0.99622)(1.12329,1.05061)\polyline(1.20636,1.10499)(1.28942,1.15938)%
\polyline(-1.81922,0.02107)(-1.72190,0.02871)\polyline(-1.62458,0.03635)(-1.52726,0.04399)%
\polyline(-1.42994,0.05163)(-1.33262,0.05927)\polyline(-1.23530,0.06691)(-1.13799,0.07455)%
\polyline(-1.04067,0.08219)(-0.94335,0.08983)\polyline(-0.84603,0.09747)(-0.74871,0.10511)%
\polyline(-0.65139,0.11275)(-0.55407,0.12039)\polyline(-0.45675,0.12803)(-0.35943,0.13567)%
\polyline(-0.26211,0.14331)(-0.16479,0.15095)\polyline(-0.06747,0.15859)(0.02985,0.16623)%
\polyline(0.12716,0.17387)(0.22448,0.18151)\polyline(0.32180,0.18915)(0.41912,0.19679)%
\polyline(0.51644,0.20443)(0.61376,0.21207)\polyline(0.71108,0.21972)(0.80840,0.22736)%
\polyline(0.90572,0.23500)(1.00304,0.24264)\polyline(1.10036,0.25028)(1.19768,0.25792)%
\polyline(1.29499,0.26556)(1.39231,0.27320)\polyline(1.48963,0.28084)(1.58695,0.28848)%
\polyline(1.68427,0.29612)(1.78159,0.30376)\polyline(1.87891,0.31140)(1.97623,0.31904)%
\polyline(2.07355,0.32668)(2.17087,0.33432)\polyline(2.26819,0.34196)(2.36551,0.34960)%
\polyline(2.46283,0.35724)(2.56014,0.36488)\polyline(2.65746,0.37252)(2.75478,0.38016)%
\polyline(2.85210,0.38780)(2.94942,0.39544)\polyline(3.04674,0.40308)(3.14406,0.41072)%
\polyline(1.28942,1.15938)(1.30969,1.06188)\polyline(1.32996,0.96437)(1.35023,0.86687)%
\polyline(1.37050,0.76937)(1.39077,0.67186)\polyline(1.41104,0.57436)(1.43131,0.47686)%
\polyline(1.45158,0.37935)(1.47185,0.28185)\polyline(1.49212,0.18435)(1.51239,0.08684)%
\polyline(1.53267,-0.01066)(1.55294,-0.10817)\polyline(1.57321,-0.20567)(1.59348,-0.30317)%
\polyline(1.61375,-0.40068)(1.63402,-0.49818)\polyline(1.65429,-0.59568)(1.67456,-0.69319)%
\polyline(1.69483,-0.79069)(1.71510,-0.88819)\polyline(1.73537,-0.98570)(1.75564,-1.08320)%
\settowidth{\Width}{$\mathrm{P}_{j-1}$}\setlength{\Width}{-1\Width}%
\settoheight{\Height}{$\mathrm{P}_{j-1}$}\settodepth{\Depth}{$\mathrm{P}_{j-1}$}\setlength{\Height}{\Depth}%
\put(-3.3300000,-1.7800000){\hspace*{\Width}\raisebox{\Height}{$\mathrm{P}_{j-1}$}}%
\settowidth{\Width}{$\mathrm{P}_{j}$}\setlength{\Width}{-1\Width}%
\settoheight{\Height}{$\mathrm{P}_{j}$}\settodepth{\Depth}{$\mathrm{P}_{j}$}\setlength{\Height}{\Depth}%
\put(-1.8700000,0.0700000){\hspace*{\Width}\raisebox{\Height}{$\mathrm{P}_{j}$}}%
\settowidth{\Width}{$\mathrm{P}_{j+1}$}\setlength{\Width}{-0.5\Width}%
\settoheight{\Height}{$\mathrm{P}_{j+1}$}\settodepth{\Depth}{$\mathrm{P}_{j+1}$}\setlength{\Height}{\Depth}%
\put(1.2900000,1.2100000){\hspace*{\Width}\raisebox{\Height}{$\mathrm{P}_{j+1}$}}%
\settowidth{\Width}{$\mathrm{P}_{j+2}$}\setlength{\Width}{0\Width}%
\settoheight{\Height}{$\mathrm{P}_{j+2}$}\settodepth{\Depth}{$\mathrm{P}_{j+2}$}\setlength{\Height}{-0.5\Height}\setlength{\Depth}{0.5\Depth}\addtolength{\Height}{\Depth}%
\put(3.1900000,0.4100000){\hspace*{\Width}\raisebox{\Height}{$\mathrm{P}_{j+2}$}}%
\settowidth{\Width}{$\mathrm{P}_{j+3}$}\setlength{\Width}{-0.5\Width}%
\settoheight{\Height}{$\mathrm{P}_{j+3}$}\settodepth{\Depth}{$\mathrm{P}_{j+3}$}\setlength{\Height}{-\Height}%
\put(1.7600000,-1.1300000){\hspace*{\Width}\raisebox{\Height}{$\mathrm{P}_{j+3}$}}%
\polyline(-3.27911,-1.83186)(-3.17835,-1.59544)(-3.06753,-1.37038)(-2.94667,-1.15669)%
(-2.81575,-0.95435)(-2.67479,-0.76338)(-2.52377,-0.58377)(-2.36271,-0.41551)(-2.19160,-0.25862)%
(-2.01043,-0.11310)(-1.81922,0.02107)(-1.52253,0.20679)(-1.22148,0.37811)(-0.91652,0.53444)%
(-0.60810,0.67519)(-0.29669,0.79976)(0.01727,0.90756)(0.33331,0.99802)(0.65099,1.07054)%
(0.96984,1.12452)(1.28942,1.15938)(1.59149,1.17095)(1.87237,1.15800)(2.13029,1.12251)%
(2.36349,1.06649)(2.57018,0.99192)(2.74860,0.90081)(2.89697,0.79513)(3.01352,0.67690)%
(3.09647,0.54810)(3.14406,0.41072)(3.14776,0.31118)(3.11979,0.20056)(3.06014,0.07886)%
(2.96881,-0.05391)(2.84581,-0.19777)(2.69113,-0.35270)(2.50477,-0.51870)(2.28674,-0.69579)%
(2.03703,-0.88396)(1.75564,-1.08320)%
\linethickness{0.003in}
\put(-0.83803,0.66350){\circle*{0.120000}}\put(0.22346,1.07570){\circle*{0.120000}}
\linethickness{0.008in}
\linethickness{0.003in}
\put(2.32964,1.24104){\circle*{0.120000}}\put(3.04635,0.88072){\circle*{0.120000}}
\linethickness{0.008in}
\settowidth{\Width}{$\mathrm{Q}_{j}$}\setlength{\Width}{-0.5\Width}%
\settoheight{\Height}{$\mathrm{Q}_{j}$}\settodepth{\Depth}{$\mathrm{Q}_{j}$}\setlength{\Height}{\Depth}%
\put(-0.8400000,0.7100000){\hspace*{\Width}\raisebox{\Height}{$\mathrm{Q}_{j}$}}%
\settowidth{\Width}{$\mathrm{R}_{j}$}\setlength{\Width}{-0.5\Width}%
\settoheight{\Height}{$\mathrm{R}_{j}$}\settodepth{\Depth}{$\mathrm{R}_{j}$}\setlength{\Height}{\Depth}%
\put(0.2200000,1.1300000){\hspace*{\Width}\raisebox{\Height}{$\mathrm{R}_{j}$}}%
\settowidth{\Width}{$\theta$}\setlength{\Width}{-0.5\Width}%
\settoheight{\Height}{$\theta$}\settodepth{\Depth}{$\theta$}\setlength{\Height}{-0.5\Height}\setlength{\Depth}{0.5\Depth}\addtolength{\Height}{\Depth}%
\put(0.4000000,0.3700000){\hspace*{\Width}\raisebox{\Height}{$\theta$}}%
\polyline(0.23508,0.18235)(0.23360,0.20596)(0.22183,0.26763)(0.20243,0.32735)(0.17569,0.38417)%
(0.15507,0.41667)%
\end{picture}}%

%% file: fig/osplineclosed.tex
{\unitlength=1cm%
\begin{picture}%
(7,5)(-3.5,-2.5)%
\linethickness{0.008in}
\polyline(3.00000,0.00000)(2.99069,0.09866)\polyline(2.98138,0.19731)(2.97634,0.25067)(2.96383,0.29441)%
\polyline(2.93657,0.38968)(2.90931,0.48495)\polyline(2.86800,0.57484)(2.82458,0.66391)%
\polyline(2.77859,0.75147)(2.72144,0.83242)\polyline(2.66430,0.91338)(2.62892,0.96351)(2.60290,0.99084)%
\polyline(2.53458,1.06262)(2.46625,1.13439)\polyline(2.39416,1.20207)(2.31701,1.26425)%
\polyline(2.23985,1.32643)(2.18691,1.36909)(2.16054,1.38560)\polyline(2.07655,1.43818)(1.99256,1.49076)%
\polyline(1.90834,1.54293)(1.81916,1.58613)\polyline(1.72997,1.62932)(1.64079,1.67252)%
\polyline(1.54919,1.71002)(1.45615,1.74412)\polyline(1.36311,1.77822)(1.27734,1.80965)(1.26985,1.81163)%
\polyline(1.17404,1.83692)(1.07822,1.86221)\polyline(0.98241,1.88750)(0.92705,1.90211)(0.88581,1.90917)%
\polyline(0.78814,1.92589)(0.69047,1.94261)\polyline(0.59279,1.95933)(0.56214,1.96457)(0.49439,1.97028)%
\polyline(0.39564,1.97860)(0.29690,1.98691)\polyline(0.19815,1.99523)(0.18837,1.99605)(0.09909,1.99605)%
\polyline(0.00000,1.99605)(-0.09909,1.99605)\polyline(-0.19815,1.99523)(-0.29690,1.98691)%
\polyline(-0.39564,1.97860)(-0.49439,1.97028)\polyline(-0.59279,1.95933)(-0.69047,1.94261)%
\polyline(-0.78814,1.92589)(-0.88581,1.90917)\polyline(-0.98241,1.88750)(-1.07822,1.86221)%
\polyline(-1.17404,1.83692)(-1.26985,1.81163)\polyline(-1.36311,1.77822)(-1.45615,1.74412)%
\polyline(-1.54919,1.71002)(-1.60748,1.68866)(-1.64079,1.67252)\polyline(-1.72997,1.62932)(-1.81916,1.58613)%
\polyline(-1.90834,1.54293)(-1.91227,1.54103)(-1.99256,1.49076)\polyline(-2.07655,1.43818)(-2.16054,1.38560)%
\polyline(-2.23985,1.32643)(-2.31701,1.26425)\polyline(-2.39416,1.20207)(-2.42705,1.17557)(-2.46625,1.13439)%
\polyline(-2.53458,1.06262)(-2.60290,0.99084)\polyline(-2.66430,0.91338)(-2.72144,0.83242)%
\polyline(-2.77859,0.75147)(-2.78933,0.73625)(-2.82458,0.66391)\polyline(-2.86800,0.57484)(-2.90575,0.49738)(-2.90931,0.48495)%
\polyline(-2.93657,0.38968)(-2.96383,0.29441)\polyline(-2.98138,0.19731)(-2.99069,0.09866)%
\polyline(-3.00000,0.00000)(-2.99069,-0.09866)\polyline(-2.98138,-0.19731)(-2.97634,-0.25067)(-2.96383,-0.29441)%
\polyline(-2.93657,-0.38968)(-2.90931,-0.48495)\polyline(-2.86800,-0.57484)(-2.82458,-0.66391)%
\polyline(-2.77859,-0.75147)(-2.72144,-0.83242)\polyline(-2.66430,-0.91338)(-2.62892,-0.96351)(-2.60290,-0.99084)%
\polyline(-2.53458,-1.06262)(-2.46625,-1.13439)\polyline(-2.39416,-1.20207)(-2.31701,-1.26425)%
\polyline(-2.23985,-1.32643)(-2.18691,-1.36909)(-2.16054,-1.38560)\polyline(-2.07655,-1.43818)(-1.99256,-1.49076)%
\polyline(-1.90834,-1.54293)(-1.81916,-1.58613)\polyline(-1.72997,-1.62932)(-1.64079,-1.67252)%
\polyline(-1.54919,-1.71002)(-1.45615,-1.74412)\polyline(-1.36311,-1.77822)(-1.27734,-1.80965)(-1.26985,-1.81163)%
\polyline(-1.17404,-1.83692)(-1.07822,-1.86221)\polyline(-0.98241,-1.88750)(-0.92705,-1.90211)(-0.88581,-1.90917)%
\polyline(-0.78814,-1.92589)(-0.69047,-1.94261)\polyline(-0.59279,-1.95933)(-0.56214,-1.96457)(-0.49439,-1.97028)%
\polyline(-0.39564,-1.97860)(-0.29690,-1.98691)\polyline(-0.19815,-1.99523)(-0.18837,-1.99605)(-0.09909,-1.99605)%
\polyline(-0.00000,-1.99605)(0.09909,-1.99605)\polyline(0.19815,-1.99523)(0.29690,-1.98691)%
\polyline(0.39564,-1.97860)(0.49439,-1.97028)\polyline(0.59279,-1.95933)(0.69047,-1.94261)%
\polyline(0.78814,-1.92589)(0.88581,-1.90917)\polyline(0.98241,-1.88750)(1.07822,-1.86221)%
\polyline(1.17404,-1.83692)(1.26985,-1.81163)\polyline(1.36311,-1.77822)(1.45615,-1.74412)%
\polyline(1.54919,-1.71002)(1.60748,-1.68866)(1.64079,-1.67252)\polyline(1.72997,-1.62932)(1.81916,-1.58613)%
\polyline(1.90834,-1.54293)(1.91227,-1.54103)(1.99256,-1.49076)\polyline(2.07655,-1.43818)(2.16054,-1.38560)%
\polyline(2.23985,-1.32643)(2.31701,-1.26425)\polyline(2.39416,-1.20207)(2.42705,-1.17557)(2.46625,-1.13439)%
\polyline(2.53458,-1.06262)(2.60290,-0.99084)\polyline(2.66430,-0.91338)(2.72144,-0.83242)%
\polyline(2.77859,-0.75147)(2.78933,-0.73625)(2.82458,-0.66391)\polyline(2.86800,-0.57484)(2.90575,-0.49738)(2.90931,-0.48495)%
\polyline(2.93657,-0.38968)(2.96383,-0.29441)\polyline(2.98138,-0.19731)(2.99069,-0.09866)%
\polyline(3.00000,0.00000)(2.96162,0.32972)(2.85021,0.64055)(2.67135,0.92876)(2.43062,1.19062)%
(2.13363,1.42241)(1.78594,1.62041)(1.39314,1.78090)(0.96083,1.90014)(0.49459,1.97441)%
(0.00000,2.00000)(-0.49459,1.97441)(-0.96083,1.90014)(-1.39314,1.78090)(-1.78594,1.62041)%
(-2.13363,1.42241)(-2.43062,1.19062)(-2.67135,0.92876)(-2.85021,0.64055)(-2.96162,0.32972)%
(-3.00000,0.00000)(-2.96162,-0.32972)(-2.85021,-0.64055)(-2.67135,-0.92876)(-2.43062,-1.19062)%
(-2.13363,-1.42241)(-1.78594,-1.62041)(-1.39314,-1.78090)(-0.96083,-1.90014)(-0.49459,-1.97441)%
(0.00000,-2.00000)(0.49459,-1.97441)(0.96083,-1.90014)(1.39314,-1.78090)(1.78594,-1.62041)%
(2.13363,-1.42241)(2.43062,-1.19062)(2.67135,-0.92876)(2.85021,-0.64055)(2.96162,-0.32972)%
(3.00000,0.00000)%
\linethickness{0.003in}
\put(3.00000,1.12644){\circle*{0.120000}}\put(1.68967,2.00000){\circle*{0.120000}}
\linethickness{0.008in}
\linethickness{0.003in}
\put(-1.68967,2.00000){\circle*{0.120000}}\put(-3.00000,1.12644){\circle*{0.120000}}
\linethickness{0.008in}
\linethickness{0.003in}
\put(-3.00000,-1.12644){\circle*{0.120000}}\put(-1.68967,-2.00000){\circle*{0.120000}}
\linethickness{0.008in}
\linethickness{0.003in}
\put(1.68967,-2.00000){\circle*{0.120000}}\put(3.00000,-1.12644){\circle*{0.120000}}
\linethickness{0.008in}
\linethickness{0.003in}
\put(3.00000,0.00000){\circle*{0.120000}}\put(0.00000,2.00000){\circle*{0.120000}}
\put(-3.00000,0.00000){\circle*{0.120000}}\put(0.00000,-2.00000){\circle*{0.120000}}
\linethickness{0.008in}
\settowidth{\Width}{A}\setlength{\Width}{0\Width}%
\settoheight{\Height}{A}\settodepth{\Depth}{A}\setlength{\Height}{\Depth}%
\put(3.0500000,0.0500000){\hspace*{\Width}\raisebox{\Height}{A}}%
\settowidth{\Width}{B}\setlength{\Width}{0\Width}%
\settoheight{\Height}{B}\settodepth{\Depth}{B}\setlength{\Height}{\Depth}%
\put(0.0500000,2.0500000){\hspace*{\Width}\raisebox{\Height}{B}}%
\settowidth{\Width}{C}\setlength{\Width}{-1\Width}%
\settoheight{\Height}{C}\settodepth{\Depth}{C}\setlength{\Height}{\Depth}%
\put(-3.0500000,0.0500000){\hspace*{\Width}\raisebox{\Height}{C}}%
\settowidth{\Width}{D}\setlength{\Width}{0\Width}%
\settoheight{\Height}{D}\settodepth{\Depth}{D}\setlength{\Height}{-\Height}%
\put(0.0500000,-2.0500000){\hspace*{\Width}\raisebox{\Height}{D}}%
\polyline(-3.50000,0.00000)(3.50000,0.00000)%
\polyline(0.00000,-2.50000)(0.00000,2.50000)%
\settowidth{\Width}{$x$}\setlength{\Width}{0\Width}%
\settoheight{\Height}{$x$}\settodepth{\Depth}{$x$}\setlength{\Height}{-0.5\Height}\setlength{\Depth}{0.5\Depth}\addtolength{\Height}{\Depth}%
\put(3.5500000,0.0000000){\hspace*{\Width}\raisebox{\Height}{$x$}}%
\settowidth{\Width}{$y$}\setlength{\Width}{-0.5\Width}%
\settoheight{\Height}{$y$}\settodepth{\Depth}{$y$}\setlength{\Height}{\Depth}%
\put(0.0000000,2.5500000){\hspace*{\Width}\raisebox{\Height}{$y$}}%
\settowidth{\Width}{O}\setlength{\Width}{-1\Width}%
\settoheight{\Height}{O}\settodepth{\Depth}{O}\setlength{\Height}{-\Height}%
\put(-0.0500000,-0.0500000){\hspace*{\Width}\raisebox{\Height}{O}}%
\end{picture}}%

%% file: fig/int1.tex
{\unitlength=1cm%
\begin{picture}%
(6.72,3.65)(-0.24,-0.2)%
\linethickness{0.008in}
\linethickness{0.012in}
\polyline(0.97979,1.38061)(1.04045,1.43760)(1.10477,1.49365)(1.17260,1.54872)(1.24381,1.60278)%
(1.31826,1.65579)(1.39581,1.70773)(1.47632,1.75855)(1.55965,1.80822)(1.64566,1.85672)%
(1.73422,1.90399)(1.82518,1.95001)(1.91841,1.99475)(2.01377,2.03818)(2.11111,2.08024)%
(2.21031,2.12092)(2.31121,2.16018)(2.41369,2.19798)(2.51759,2.23430)(2.62279,2.26909)%
(2.72915,2.30232)(2.83652,2.33395)(2.94476,2.36397)(3.05375,2.39232)(3.16333,2.41897)%
(3.27337,2.44390)(3.38374,2.46706)(3.49428,2.48842)(3.60487,2.50795)(3.71536,2.52562)%
(3.82561,2.54139)(3.93550,2.55522)(4.04487,2.56708)(4.15358,2.57694)(4.26151,2.58476)%
(4.36850,2.59051)(4.47443,2.59415)(4.57915,2.59566)(4.68252,2.59499)(4.78441,2.59211)%
(4.88467,2.58699)(4.98317,2.57959)(5.07977,2.56989)(5.17432,2.55783)(5.26670,2.54340)%
(5.35675,2.52656)(5.44435,2.50726)(5.52935,2.48549)(5.61161,2.46120)(5.69100,2.43436)%
(5.76737,2.40493)%
\polyline(0.97979,1.38061)(0.97979,0.00000)%
\polyline(5.76737,2.40493)(5.76737,0.00000)%
\polyline(5.73764,0.00000)(5.76737,0.02973)%
\polyline(5.38409,-0.00000)(5.76737,0.38328)%
\polyline(5.03053,0.00000)(5.76737,0.73684)%
\polyline(4.67698,0.00000)(5.76737,1.09039)%
\polyline(4.32343,0.00000)(5.76737,1.44394)%
\polyline(3.96987,0.00000)(5.76737,1.79750)%
\polyline(3.61632,0.00000)(5.76737,2.15105)%
\polyline(3.26277,0.00000)(5.69542,2.43265)%
\polyline(2.90921,0.00000)(5.42151,2.51229)%
\polyline(2.55566,0.00000)(5.12037,2.56471)%
\polyline(2.20211,0.00000)(4.79374,2.59163)%
\polyline(1.84855,0.00000)(4.44158,2.59302)%
\polyline(1.49500,0.00000)(4.06379,2.56879)%
\polyline(1.14145,0.00000)(3.65788,2.51643)%
\polyline(0.97979,0.19190)(3.21961,2.43172)%
\polyline(0.97979,0.54545)(2.73979,2.30545)%
\polyline(0.97979,0.89900)(2.19573,2.11495)%
\polyline(0.97979,1.25256)(1.49975,1.77252)%
\linethickness{0.003in}
\put(0.97979,1.38061){\circle*{0.120000}}\put(1.95957,2.33812){\circle*{0.120000}}
\put(4.52037,2.91709){\circle*{0.120000}}\put(5.76737,2.40493){\circle*{0.120000}}
\linethickness{0.008in}
\settowidth{\Width}{$\mathrm{P}_j$}\setlength{\Width}{-1\Width}%
\settoheight{\Height}{$\mathrm{P}_j$}\settodepth{\Depth}{$\mathrm{P}_j$}\setlength{\Height}{\Depth}%
\put(0.9300000,1.4300000){\hspace*{\Width}\raisebox{\Height}{$\mathrm{P}_j$}}%
\settowidth{\Width}{$\mathrm{P}_{j+1}$}\setlength{\Width}{-0.5\Width}%
\settoheight{\Height}{$\mathrm{P}_{j+1}$}\settodepth{\Depth}{$\mathrm{P}_{j+1}$}\setlength{\Height}{\Depth}%
\put(1.9600000,2.3900000){\hspace*{\Width}\raisebox{\Height}{$\mathrm{P}_{j+1}$}}%
\settowidth{\Width}{$\mathrm{P}_{j+2}$}\setlength{\Width}{-0.5\Width}%
\settoheight{\Height}{$\mathrm{P}_{j+2}$}\settodepth{\Depth}{$\mathrm{P}_{j+2}$}\setlength{\Height}{\Depth}%
\put(4.5200000,2.9700000){\hspace*{\Width}\raisebox{\Height}{$\mathrm{P}_{j+2}$}}%
\settowidth{\Width}{$\mathrm{P}_{j+3}$}\setlength{\Width}{0\Width}%
\settoheight{\Height}{$\mathrm{P}_{j+3}$}\settodepth{\Depth}{$\mathrm{P}_{j+3}$}\setlength{\Height}{\Depth}%
\put(5.8200000,2.4500000){\hspace*{\Width}\raisebox{\Height}{$\mathrm{P}_{j+3}$}}%
\settowidth{\Width}{$C$}\setlength{\Width}{-0.5\Width}%
\settoheight{\Height}{$C$}\settodepth{\Depth}{$C$}\setlength{\Height}{-0.5\Height}\setlength{\Depth}{0.5\Depth}\addtolength{\Height}{\Depth}%
\put(3.2400000,2.6300000){\hspace*{\Width}\raisebox{\Height}{$C$}}%
\polyline(-0.24000,0.00000)(6.48000,0.00000)%
\polyline(0.00000,-0.20000)(0.00000,3.45000)%
\settowidth{\Width}{$x$}\setlength{\Width}{0\Width}%
\settoheight{\Height}{$x$}\settodepth{\Depth}{$x$}\setlength{\Height}{-0.5\Height}\setlength{\Depth}{0.5\Depth}\addtolength{\Height}{\Depth}%
\put(6.5300000,0.0000000){\hspace*{\Width}\raisebox{\Height}{$x$}}%
\settowidth{\Width}{$y$}\setlength{\Width}{-0.5\Width}%
\settoheight{\Height}{$y$}\settodepth{\Depth}{$y$}\setlength{\Height}{\Depth}%
\put(0.0000000,3.5000000){\hspace*{\Width}\raisebox{\Height}{$y$}}%
\settowidth{\Width}{O}\setlength{\Width}{-1\Width}%
\settoheight{\Height}{O}\settodepth{\Depth}{O}\setlength{\Height}{-\Height}%
\put(-0.0500000,-0.0500000){\hspace*{\Width}\raisebox{\Height}{O}}%
\end{picture}}%

%% file: fig/intimp.tex
{\unitlength=1cm%
\begin{picture}%
(3.5,3.5)(-1,-0.5)%
\linethickness{0.008in}
\polyline(2.00000,1.53299)(1.98064,1.51000)(1.92000,1.43679)(1.90553,1.42000)(1.84000,1.34269)%
(1.82879,1.33000)(1.76000,1.25078)(1.75022,1.24000)(1.68000,1.16119)(1.66957,1.15000)%
(1.60000,1.07404)(1.58652,1.06000)(1.52000,0.98946)(1.50072,0.97000)(1.44000,0.90759)%
(1.41172,0.88000)(1.36000,0.82858)(1.31899,0.79000)(1.28000,0.75261)(1.22185,0.70000)%
(1.20000,0.67984)(1.12000,0.61048)(1.11940,0.61000)(1.04000,0.54505)(1.00703,0.52000)%
(0.96000,0.48352)(0.88529,0.43000)(0.88000,0.42613)(0.80000,0.37365)(0.74308,0.34000)%
(0.72000,0.32605)(0.64000,0.28412)(0.56424,0.25000)(0.56000,0.24805)(0.48000,0.21901)%
(0.40000,0.19733)(0.32000,0.18403)(0.24000,0.18030)(0.16000,0.18761)(0.08000,0.20770)%
(0.00000,0.24272)(-0.01164,0.25000)(-0.08000,0.29642)(-0.12739,0.34000)(-0.16000,0.37277)%
(-0.20515,0.43000)(-0.24000,0.47869)(-0.26448,0.52000)(-0.31273,0.61000)(-0.32000,0.62482)%
(-0.35171,0.70000)(-0.38577,0.79000)(-0.40000,0.83136)(-0.41474,0.88000)(-0.43961,0.97000)%
(-0.46188,1.06000)(-0.48000,1.14163)(-0.48167,1.15000)(-0.49822,1.24000)(-0.51301,1.33000)%
(-0.52618,1.42000)(-0.53787,1.51000)(-0.54817,1.60000)(-0.55720,1.69000)(-0.56000,1.72162)%
(-0.56476,1.78000)(-0.57113,1.87000)(-0.57654,1.96000)(-0.58105,2.05000)(-0.58473,2.14000)%
(-0.58761,2.23000)(-0.58975,2.32000)(-0.59119,2.41000)(-0.59196,2.50000)%
\polyline(-0.59196,2.50000)(-0.59196,0.00000)(2.00000,0.00000)(2.00000,1.53299)%
\polyline(1.97487,0.00000)(2.00000,0.02513)%
\polyline(1.62132,0.00000)(2.00000,0.37868)%
\polyline(1.26777,0.00000)(2.00000,0.73223)%
\polyline(0.91421,0.00000)(2.00000,1.08579)%
\polyline(0.56066,0.00000)(2.00000,1.43934)%
\polyline(0.20711,0.00000)(0.40609,0.19898)%
\polyline(-0.14645,0.00000)(0.06696,0.21341)%
\polyline(-0.50000,0.00000)(-0.14366,0.35634)%
\polyline(-0.59145,0.26210)(-0.28859,0.56497)%
\polyline(-0.59076,0.61635)(-0.39379,0.81332)%
\polyline(-0.59007,0.97059)(-0.46892,1.09174)%
\polyline(-0.58938,1.32484)(-0.52210,1.39211)%
\polyline(-0.58869,1.67908)(-0.55887,1.70889)%
\polyline(-0.58799,2.03333)(-0.58059,2.04073)%
\settowidth{\Width}{P}\setlength{\Width}{-0.5\Width}%
\settoheight{\Height}{P}\settodepth{\Depth}{P}\setlength{\Height}{\Depth}%
\put(2.0000000,1.5800000){\hspace*{\Width}\raisebox{\Height}{P}}%
\settowidth{\Width}{Q}\setlength{\Width}{-0.5\Width}%
\settoheight{\Height}{Q}\settodepth{\Depth}{Q}\setlength{\Height}{\Depth}%
\put(-0.5900000,2.5500000){\hspace*{\Width}\raisebox{\Height}{Q}}%
\polyline(-1.00000,0.00000)(2.50000,0.00000)%
\polyline(0.00000,-0.50000)(0.00000,3.00000)%
\settowidth{\Width}{$x$}\setlength{\Width}{0\Width}%
\settoheight{\Height}{$x$}\settodepth{\Depth}{$x$}\setlength{\Height}{-0.5\Height}\setlength{\Depth}{0.5\Depth}\addtolength{\Height}{\Depth}%
\put(2.5500000,0.0000000){\hspace*{\Width}\raisebox{\Height}{$x$}}%
\settowidth{\Width}{$y$}\setlength{\Width}{-0.5\Width}%
\settoheight{\Height}{$y$}\settodepth{\Depth}{$y$}\setlength{\Height}{\Depth}%
\put(0.0000000,3.0500000){\hspace*{\Width}\raisebox{\Height}{$y$}}%
\settowidth{\Width}{O}\setlength{\Width}{-1\Width}%
\settoheight{\Height}{O}\settodepth{\Depth}{O}\setlength{\Height}{-\Height}%
\put(-0.0500000,-0.0500000){\hspace*{\Width}\raisebox{\Height}{O}}%
\end{picture}}%

%% file: fig/derivimp.tex
{\unitlength=1cm%
\begin{picture}%
(3.48,3.8)(-1.08,-0.67)%
\linethickness{0.008in}
\linethickness{0.008in}
\polyline(2.00000,1.53299)(1.98064,1.51000)(1.92000,1.43679)(1.90553,1.42000)(1.84000,1.34269)%
(1.82879,1.33000)(1.76000,1.25078)(1.75022,1.24000)(1.68000,1.16119)(1.66957,1.15000)%
(1.60000,1.07404)(1.58652,1.06000)(1.52000,0.98946)(1.50072,0.97000)(1.44000,0.90759)%
(1.41172,0.88000)(1.36000,0.82858)(1.31899,0.79000)(1.28000,0.75261)(1.22185,0.70000)%
(1.20000,0.67984)(1.12000,0.61048)(1.11940,0.61000)(1.04000,0.54505)(1.00703,0.52000)%
(0.96000,0.48352)(0.88529,0.43000)(0.88000,0.42613)(0.80000,0.37365)(0.74308,0.34000)%
(0.72000,0.32605)(0.64000,0.28412)(0.56424,0.25000)(0.56000,0.24805)(0.48000,0.21901)%
(0.40000,0.19733)(0.32000,0.18403)(0.24000,0.18030)(0.16000,0.18761)(0.08000,0.20770)%
(0.00000,0.24272)(-0.01164,0.25000)(-0.08000,0.29642)(-0.12739,0.34000)(-0.16000,0.37277)%
(-0.20515,0.43000)(-0.24000,0.47869)(-0.26448,0.52000)(-0.31273,0.61000)(-0.32000,0.62482)%
(-0.35171,0.70000)(-0.38577,0.79000)(-0.40000,0.83136)(-0.41474,0.88000)(-0.43961,0.97000)%
(-0.46188,1.06000)(-0.48000,1.14163)(-0.48167,1.15000)(-0.49822,1.24000)(-0.51301,1.33000)%
(-0.52618,1.42000)(-0.53787,1.51000)(-0.54817,1.60000)(-0.55720,1.69000)(-0.56000,1.72162)%
(-0.56476,1.78000)(-0.57113,1.87000)(-0.57654,1.96000)(-0.58105,2.05000)(-0.58473,2.14000)%
(-0.58761,2.23000)(-0.58975,2.32000)(-0.59119,2.41000)(-0.59196,2.50000)%
\settowidth{\Width}{P}\setlength{\Width}{-0.5\Width}%
\settoheight{\Height}{P}\settodepth{\Depth}{P}\setlength{\Height}{\Depth}%
\put(2.0000000,1.5800000){\hspace*{\Width}\raisebox{\Height}{P}}%
\settowidth{\Width}{Q}\setlength{\Width}{-0.5\Width}%
\settoheight{\Height}{Q}\settodepth{\Depth}{Q}\setlength{\Height}{\Depth}%
\put(-0.5900000,2.5500000){\hspace*{\Width}\raisebox{\Height}{Q}}%
\linethickness{0.002in}
\polyline(2.40000,2.00819)(0.14563,-0.67000)%
\linethickness{0.002in}
\polyline(2.40000,1.96779)(0.04483,-0.67000)%
\linethickness{0.002in}
\polyline(2.40000,1.88650)(-0.10803,-0.67000)%
\linethickness{0.002in}
\polyline(2.40000,1.77315)(-0.28220,-0.67000)%
\linethickness{0.002in}
\polyline(2.40000,1.54982)(-0.59808,-0.67000)%
\linethickness{0.002in}
\polyline(2.40000,1.21157)(-1.08000,-0.62264)%
\linethickness{0.002in}
\polyline(2.40000,0.74464)(-1.08000,-0.20789)%
\linethickness{0.002in}
\polyline(2.40000,-0.08788)(-1.08000,0.34139)%
\linethickness{0.002in}
\polyline(1.17497,-0.67000)(-1.08000,1.06755)%
\linethickness{0.002in}
\polyline(0.21585,-0.67000)(-1.08000,2.45991)%
\linethickness{0.002in}
\polyline(-0.50677,-0.67000)(-0.61274,3.13000)%
\polyline(-1.08000,0.00000)(2.40000,0.00000)%
\polyline(0.00000,-0.67000)(0.00000,3.13000)%
\settowidth{\Width}{$x$}\setlength{\Width}{0\Width}%
\settoheight{\Height}{$x$}\settodepth{\Depth}{$x$}\setlength{\Height}{-0.5\Height}\setlength{\Depth}{0.5\Depth}\addtolength{\Height}{\Depth}%
\put(2.4500000,0.0000000){\hspace*{\Width}\raisebox{\Height}{$x$}}%
\settowidth{\Width}{$y$}\setlength{\Width}{-0.5\Width}%
\settoheight{\Height}{$y$}\settodepth{\Depth}{$y$}\setlength{\Height}{\Depth}%
\put(0.0000000,3.1800000){\hspace*{\Width}\raisebox{\Height}{$y$}}%
\settowidth{\Width}{O}\setlength{\Width}{-1\Width}%
\settoheight{\Height}{O}\settodepth{\Depth}{O}\setlength{\Height}{-\Height}%
\put(-0.0500000,-0.0500000){\hspace*{\Width}\raisebox{\Height}{O}}%
\end{picture}}%

%% file: fig/mant.tex
{\unitlength=1cm%
\begin{picture}%
(6.2,5.86)(-2.96,-1.54)%
\linethickness{0.008in}
\settowidth{\Width}{$x$}\setlength{\Width}{-0.5\Width}%
\settoheight{\Height}{$x$}\settodepth{\Depth}{$x$}\setlength{\Height}{-0.5\Height}\setlength{\Depth}{0.5\Depth}\addtolength{\Height}{\Depth}%
\put(-2.3800000,-0.7600000){\hspace*{\Width}\raisebox{\Height}{$x$}}%
\settowidth{\Width}{$y$}\setlength{\Width}{-0.5\Width}%
\settoheight{\Height}{$y$}\settodepth{\Depth}{$y$}\setlength{\Height}{-0.5\Height}\setlength{\Depth}{0.5\Depth}\addtolength{\Height}{\Depth}%
\put(2.2200000,-0.8200000){\hspace*{\Width}\raisebox{\Height}{$y$}}%
\settowidth{\Width}{$z$}\setlength{\Width}{-0.5\Width}%
\settoheight{\Height}{$z$}\settodepth{\Depth}{$z$}\setlength{\Height}{-0.5\Height}\setlength{\Depth}{0.5\Depth}\addtolength{\Height}{\Depth}%
\put(0.0000000,3.9500000){\hspace*{\Width}\raisebox{\Height}{$z$}}%
\polyline(-1.58905,1.76317)(-1.54846,1.68504)%
\polyline(-1.54846,1.68504)(-1.51761,1.62566)(-1.50925,1.60901)%
\polyline(-1.50925,1.60901)(-1.48666,1.56408)(-1.47226,1.53494)%
\polyline(-1.47226,1.53494)(-1.43736,1.46439)(-1.40094,1.39011)(-1.34718,1.28069)%
(-1.32366,1.23328)(-1.25215,1.09278)(-1.24463,1.07842)(-1.18598,0.97026)(-1.12623,0.86833)%
(-1.06147,0.76948)(-1.00211,0.69104)(-0.98721,0.67332)(-0.94490,0.62749)(-0.88868,0.57712)%
(-0.83261,0.53896)(-0.77776,0.51316)(-0.72283,0.49836)(-0.66734,0.49403)(-0.65601,0.49440)%
(-0.61207,0.49956)(-0.55618,0.51423)(-0.49907,0.53765)(-0.44167,0.56867)(-0.38285,0.60710)%
(-0.34881,0.63196)(-0.32264,0.65221)(-0.26107,0.70313)(-0.21911,0.74001)(-0.19747,0.75959)%
(-0.13195,0.82063)(-0.10257,0.84865)(-0.06384,0.88595)(0.00284,0.95064)(0.08190,1.02705)%
(0.09966,1.04405)(0.16135,1.10230)(0.18867,1.12761)(0.24710,1.18056)(0.27076,1.20147)%
(0.33522,1.25657)%
\polyline(0.33522,1.25657)(0.34196,1.26234)(0.41727,1.32292)(0.45520,1.35178)(0.48236,1.37172)%
(0.54264,1.41362)(0.59731,1.44854)(0.63601,1.47126)(0.64726,1.47750)(0.69239,1.50081)%
(0.71682,1.51207)(0.73194,1.51846)(0.76617,1.53095)(0.79324,1.53818)(0.80179,1.53970)%
\polyline(-1.58905,1.76317)(-1.51761,1.62566)(-1.48666,1.56408)(-1.43736,1.46439)%
(-1.40094,1.39011)(-1.34718,1.28069)(-1.32366,1.23328)(-1.25215,1.09278)(-1.24463,1.07842)%
(-1.18598,0.97026)(-1.12623,0.86833)(-1.06147,0.76948)(-1.00211,0.69104)(-0.98721,0.67332)%
(-0.94490,0.62749)(-0.88868,0.57712)(-0.83261,0.53896)(-0.77776,0.51316)(-0.72283,0.49836)%
(-0.66734,0.49403)(-0.65601,0.49440)(-0.61207,0.49956)(-0.55618,0.51423)(-0.49907,0.53765)%
(-0.44167,0.56867)(-0.38285,0.60710)(-0.34881,0.63196)(-0.32264,0.65221)(-0.26107,0.70313)%
(-0.21911,0.74001)(-0.19747,0.75959)(-0.13195,0.82063)(-0.10257,0.84865)(-0.06384,0.88595)%
(0.00284,0.95064)(0.08190,1.02705)(0.09966,1.04405)(0.16135,1.10230)(0.18867,1.12761)%
(0.24710,1.18056)(0.27076,1.20147)(0.34196,1.26234)(0.41727,1.32292)(0.45520,1.35178)%
(0.48236,1.37172)(0.54264,1.41362)(0.59731,1.44854)(0.63601,1.47126)(0.64726,1.47750)%
(0.69239,1.50081)(0.71682,1.51207)(0.73194,1.51846)(0.76617,1.53095)(0.79324,1.53818)%
(0.80337,1.53997)(0.81199,1.54091)(0.81894,1.54097)(0.81960,1.54091)%
\polyline(-1.58905,1.76317)(-1.51761,1.62566)(-1.48666,1.56408)(-1.43736,1.46439)%
(-1.40094,1.39011)(-1.34718,1.28069)(-1.32366,1.23328)(-1.25215,1.09278)(-1.24463,1.07842)%
(-1.18598,0.97026)(-1.12623,0.86833)(-1.06147,0.76948)(-1.00211,0.69104)(-0.98721,0.67332)%
(-0.94490,0.62749)(-0.88868,0.57712)(-0.83261,0.53896)(-0.77776,0.51316)(-0.72283,0.49836)%
(-0.66734,0.49403)(-0.65601,0.49440)(-0.61207,0.49956)(-0.55618,0.51423)(-0.49907,0.53765)%
(-0.44167,0.56867)(-0.38285,0.60710)(-0.34881,0.63196)(-0.32264,0.65221)(-0.26107,0.70313)%
(-0.21911,0.74001)(-0.19747,0.75959)(-0.13195,0.82063)(-0.10257,0.84865)(-0.06384,0.88595)%
(0.00284,0.95064)(0.08190,1.02705)(0.09966,1.04405)(0.16135,1.10230)(0.18867,1.12761)%
(0.24710,1.18056)(0.27076,1.20147)(0.34196,1.26234)(0.41727,1.32292)(0.45520,1.35178)%
(0.48236,1.37172)(0.54264,1.41362)(0.59731,1.44854)(0.63601,1.47126)(0.64726,1.47750)%
(0.69239,1.50081)(0.71682,1.51207)(0.73194,1.51846)(0.76617,1.53095)(0.79324,1.53818)%
(0.80337,1.53997)(0.81199,1.54091)(0.81894,1.54097)(0.81960,1.54091)%
\polyline(0.67308,1.69233)(0.67266,1.69294)(0.63947,1.74872)(0.62763,1.77003)(0.57900,1.86412)%
(0.52672,1.97472)(0.47088,2.09954)(0.44236,2.16451)(0.41253,2.23248)(0.35227,2.36680)%
(0.30651,2.46330)(0.28983,2.49674)(0.22664,2.61258)(0.18305,2.68043)(0.16218,2.70888)%
(0.09760,2.77855)(0.06322,2.80351)(0.03277,2.81823)(-0.03157,2.82604)(-0.05768,2.82021)%
(-0.09526,2.80294)(-0.15768,2.75232)(-0.18523,2.72197)(-0.21862,2.67928)(-0.27743,2.59069)%
(-0.33210,2.49658)(-0.38722,2.39488)(-0.43708,2.30119)(-0.48285,2.21747)(-0.52390,2.14728)%
(-0.55941,2.09281)(-0.58827,2.05485)(-0.59929,2.04242)(-0.60596,2.03579)%
\polyline(-0.60596,2.03579)(-0.60845,2.03331)(-0.61704,2.02629)%
\polyline(0.33522,1.25657)(0.36370,1.22607)(0.38155,1.20765)(0.42147,1.16781)(0.47746,1.11581)%
(0.53190,1.07055)(0.54793,1.05837)(0.58518,1.03231)(0.63729,1.00168)(0.68835,0.97905)%
(0.73858,0.96469)(0.76487,0.96054)(0.78827,0.95889)(0.83753,0.96200)(0.88644,0.97427)%
(0.93523,0.99588)(0.98411,1.02709)(1.02863,1.06375)(1.08305,1.11886)(1.13331,1.17923)%
(1.18418,1.24869)(1.23569,1.32647)(1.25669,1.36001)(1.28787,1.41151)(1.34039,1.50186)%
(1.39285,1.59525)(1.44489,1.68931)(1.45881,1.71450)(1.49634,1.78212)(1.54646,1.87095)%
(1.59468,1.95376)(1.64054,2.02917)(1.68359,2.09626)(1.72342,2.15458)(1.75960,2.20403)%
(1.77312,2.22117)%
\polyline(1.77312,2.22117)(1.79168,2.24473)%
\polyline(1.79168,2.24473)(1.81869,2.27640)%
\polyline(1.81869,2.27640)(1.81918,2.27697)(1.82381,2.28214)(1.83938,2.29883)%
\polyline(-2.12002,2.63309)(-2.07910,2.58389)(-2.03819,2.53109)(-1.99727,2.47467)%
(-1.95636,2.41469)(-1.91544,2.35127)(-1.87453,2.28460)(-1.83361,2.21498)(-1.79270,2.14277)%
(-1.75178,2.06844)(-1.71087,1.99251)(-1.66995,1.91558)(-1.62904,1.83832)(-1.58926,1.76357)%
\polyline(-1.58926,1.76357)(-1.58812,1.76144)(-1.54721,1.68569)(-1.54514,1.68197)%
\polyline(-1.54514,1.68197)(-1.51440,1.62647)%
\polyline(-1.51440,1.62647)(-1.50629,1.61186)(-1.47954,1.56533)%
\polyline(-1.47954,1.56533)(-1.46538,1.54070)(-1.42446,1.47301)(-1.38355,1.40952)%
(-1.34263,1.35091)(-1.30172,1.29785)(-1.26080,1.25088)(-1.21989,1.21050)(-1.17897,1.17708)%
(-1.13806,1.15090)(-1.09714,1.13215)(-1.05623,1.12088)(-1.01531,1.11704)(-0.97440,1.12044)%
(-0.93348,1.13080)(-0.89257,1.14775)(-0.85165,1.17080)(-0.81074,1.19938)(-0.76982,1.23286)%
(-0.72891,1.27053)(-0.68799,1.31166)(-0.64708,1.35548)(-0.60616,1.40120)(-0.56525,1.44806)%
(-0.52433,1.49530)(-0.48342,1.54221)(-0.44250,1.58812)(-0.40159,1.63243)(-0.36067,1.67462)%
(-0.31976,1.71423)(-0.27884,1.75087)(-0.23793,1.78428)(-0.19701,1.81424)(-0.18739,1.82045)%
\polyline(-0.18739,1.82045)(-0.15610,1.84063)(-0.11518,1.86342)(-0.07427,1.88260)%
\polyline(2.12002,2.58228)(2.07613,2.54069)(2.03225,2.49761)(1.98836,2.45333)(1.94448,2.40819)%
(1.90059,2.36256)(1.85671,2.31688)(1.84131,2.30099)%
\polyline(1.84131,2.30099)(1.82322,2.28231)%
\polyline(1.82322,2.28231)(1.81282,2.27157)(1.76893,2.22708)(1.76837,2.22653)%
\polyline(1.76837,2.22653)(1.72505,2.18383)(1.68116,2.14220)(1.63728,2.10255)(1.59339,2.06516)%
(1.54950,2.03022)(1.50575,1.99795)%
\polyline(1.50575,1.99795)(1.50562,1.99785)(1.46644,1.97127)%
\polyline(1.46644,1.97127)(1.46173,1.96809)(1.41785,1.94088)(1.37396,1.91608)(1.33008,1.89349)%
(1.28619,1.87287)(1.24230,1.85392)(1.19842,1.83637)(1.15453,1.81995)(1.15247,1.81921)%
\polyline(1.15247,1.81921)(1.11065,1.80440)(1.06676,1.78955)(1.02287,1.77530)(0.97899,1.76157)%
(0.93510,1.74843)(0.89122,1.73598)(0.84733,1.72442)(0.80345,1.71398)(0.75956,1.70494)%
(0.71567,1.69757)(0.67308,1.69233)%
\polyline(0.67308,1.69233)(0.67179,1.69218)(0.62790,1.68899)(0.58402,1.68822)(0.54013,1.68999)%
(0.49625,1.69437)(0.45236,1.70132)(0.40847,1.71073)(0.36459,1.72239)(0.32070,1.73603)%
(0.27682,1.75130)(0.23293,1.76780)(0.18904,1.78512)(0.14516,1.80281)(0.10127,1.82045)%
(0.05739,1.83761)(0.01350,1.85390)(-0.03038,1.86899)(-0.07427,1.88260)%
\polyline(0.07427,3.33277)(0.11518,3.28358)(0.15610,3.23077)(0.19701,3.17436)(0.23793,3.11438)%
(0.27884,3.05095)(0.31976,2.98429)(0.36067,2.91467)(0.40159,2.84246)(0.44250,2.76812)%
(0.48342,2.69219)(0.52433,2.61527)(0.56525,2.53800)(0.60616,2.46112)(0.64708,2.38538)%
(0.68799,2.31155)(0.72891,2.24039)(0.76982,2.17270)(0.81074,2.10920)(0.85165,2.05060)%
(0.89257,1.99753)(0.93348,1.95057)(0.97440,1.91018)(1.01531,1.87677)(1.05623,1.85059)%
(1.09714,1.83184)(1.13806,1.82057)(1.15247,1.81921)%
\polyline(1.46538,1.97021)(1.50629,2.01135)(1.54721,2.05517)(1.58812,2.10089)(1.62904,2.14775)%
(1.66995,2.19499)(1.71087,2.24190)(1.75178,2.28781)(1.79270,2.33212)(1.83361,2.37431)%
(1.87453,2.41392)(1.91544,2.45056)(1.95636,2.48397)(1.99727,2.51393)(2.03819,2.54032)%
(2.07910,2.56311)(2.12002,2.58228)%
\polyline(0.07427,3.33277)(0.03038,3.29118)(-0.01350,3.24810)(-0.05739,3.20382)(-0.10127,3.15868)%
(-0.14516,3.11306)(-0.18904,3.06738)(-0.23293,3.02206)(-0.27682,2.97758)(-0.32070,2.93432)%
(-0.36459,2.89269)(-0.40847,2.85304)(-0.45236,2.81565)(-0.49625,2.78071)(-0.54013,2.74835)%
(-0.58402,2.71858)(-0.62790,2.69137)(-0.67179,2.66657)(-0.71567,2.64398)(-0.75956,2.62336)%
(-0.80345,2.60441)(-0.84733,2.58687)(-0.89122,2.57044)(-0.93510,2.55489)(-0.97899,2.54005)%
(-1.02287,2.52579)(-1.06676,2.51206)(-1.11065,2.49892)(-1.15453,2.48648)(-1.19842,2.47492)%
(-1.24230,2.46447)(-1.28619,2.45543)(-1.33008,2.44806)(-1.37396,2.44267)(-1.41785,2.43948)%
(-1.46173,2.43871)(-1.50562,2.44048)(-1.54950,2.44486)(-1.59339,2.45181)(-1.63728,2.46122)%
(-1.68116,2.47288)(-1.72505,2.48652)(-1.76893,2.50179)(-1.81282,2.51829)(-1.85671,2.53561)%
(-1.90059,2.55331)(-1.94448,2.57094)(-1.98836,2.58810)(-2.03225,2.60439)(-2.07613,2.61948)%
(-2.12002,2.63309)%
\polyline(2.19429,0.69969)(-2.19429,-0.69969)%
\polyline(-2.04575,0.75049)(2.04575,-0.75049)%
\polyline(0.00000,-1.54000)(0.00000,0.94788)%
\polyline(0.00000,2.82221)(0.00000,3.26136)%
\polyline(0.00000,3.26136)(0.00000,3.75877)%
\end{picture}}%

%% file: fig/mant2.tex
{\unitlength=1cm%
\begin{picture}%
(6.2,5.86)(-2.96,-1.54)%
\linethickness{0.008in}
\settowidth{\Width}{$x$}\setlength{\Width}{-0.5\Width}%
\settoheight{\Height}{$x$}\settodepth{\Depth}{$x$}\setlength{\Height}{-0.5\Height}\setlength{\Depth}{0.5\Depth}\addtolength{\Height}{\Depth}%
\put(-2.3800000,-0.7600000){\hspace*{\Width}\raisebox{\Height}{$x$}}%
\settowidth{\Width}{$y$}\setlength{\Width}{-0.5\Width}%
\settoheight{\Height}{$y$}\settodepth{\Depth}{$y$}\setlength{\Height}{-0.5\Height}\setlength{\Depth}{0.5\Depth}\addtolength{\Height}{\Depth}%
\put(2.2200000,-0.8200000){\hspace*{\Width}\raisebox{\Height}{$y$}}%
\settowidth{\Width}{$z$}\setlength{\Width}{-0.5\Width}%
\settoheight{\Height}{$z$}\settodepth{\Depth}{$z$}\setlength{\Height}{-0.5\Height}\setlength{\Depth}{0.5\Depth}\addtolength{\Height}{\Depth}%
\put(0.0000000,3.9500000){\hspace*{\Width}\raisebox{\Height}{$z$}}%
\polyline(-1.58905,1.76317)(-1.54846,1.68504)%
\polyline(-1.54846,1.68504)(-1.51761,1.62566)(-1.50925,1.60901)%
\polyline(-1.50925,1.60901)(-1.48666,1.56408)(-1.47226,1.53494)%
\polyline(-1.47226,1.53494)(-1.43736,1.46439)(-1.40094,1.39011)(-1.34718,1.28069)%
(-1.32366,1.23328)(-1.25215,1.09278)(-1.24463,1.07842)(-1.18598,0.97026)(-1.12623,0.86833)%
(-1.06147,0.76948)(-1.00211,0.69104)(-0.98721,0.67332)(-0.94490,0.62749)(-0.88868,0.57712)%
(-0.83261,0.53896)(-0.77776,0.51316)(-0.72283,0.49836)(-0.66734,0.49403)(-0.65601,0.49440)%
(-0.61207,0.49956)(-0.55618,0.51423)(-0.49907,0.53765)(-0.44167,0.56867)(-0.38285,0.60710)%
(-0.34881,0.63196)(-0.32264,0.65221)(-0.26107,0.70313)(-0.21911,0.74001)(-0.19747,0.75959)%
(-0.13195,0.82063)(-0.10257,0.84865)(-0.06384,0.88595)(0.00284,0.95064)(0.08190,1.02705)%
(0.09966,1.04405)(0.16135,1.10230)(0.18867,1.12761)(0.24710,1.18056)(0.27076,1.20147)%
(0.33522,1.25657)%
\polyline(0.33522,1.25657)(0.34196,1.26234)(0.41727,1.32292)(0.45520,1.35178)(0.48236,1.37172)%
(0.54264,1.41362)(0.59731,1.44854)(0.63601,1.47126)(0.64726,1.47750)(0.69239,1.50081)%
(0.71682,1.51207)(0.73194,1.51846)(0.76617,1.53095)(0.79324,1.53818)(0.80179,1.53970)%
\polyline(-1.58905,1.76317)(-1.51761,1.62566)(-1.48666,1.56408)(-1.43736,1.46439)%
(-1.40094,1.39011)(-1.34718,1.28069)(-1.32366,1.23328)(-1.25215,1.09278)(-1.24463,1.07842)%
(-1.18598,0.97026)(-1.12623,0.86833)(-1.06147,0.76948)(-1.00211,0.69104)(-0.98721,0.67332)%
(-0.94490,0.62749)(-0.88868,0.57712)(-0.83261,0.53896)(-0.77776,0.51316)(-0.72283,0.49836)%
(-0.66734,0.49403)(-0.65601,0.49440)(-0.61207,0.49956)(-0.55618,0.51423)(-0.49907,0.53765)%
(-0.44167,0.56867)(-0.38285,0.60710)(-0.34881,0.63196)(-0.32264,0.65221)(-0.26107,0.70313)%
(-0.21911,0.74001)(-0.19747,0.75959)(-0.13195,0.82063)(-0.10257,0.84865)(-0.06384,0.88595)%
(0.00284,0.95064)(0.08190,1.02705)(0.09966,1.04405)(0.16135,1.10230)(0.18867,1.12761)%
(0.24710,1.18056)(0.27076,1.20147)(0.34196,1.26234)(0.41727,1.32292)(0.45520,1.35178)%
(0.48236,1.37172)(0.54264,1.41362)(0.59731,1.44854)(0.63601,1.47126)(0.64726,1.47750)%
(0.69239,1.50081)(0.71682,1.51207)(0.73194,1.51846)(0.76617,1.53095)(0.79324,1.53818)%
(0.80337,1.53997)(0.81199,1.54091)(0.81894,1.54097)(0.81960,1.54091)%
\polyline(-1.58905,1.76317)(-1.51761,1.62566)(-1.48666,1.56408)(-1.43736,1.46439)%
(-1.40094,1.39011)(-1.34718,1.28069)(-1.32366,1.23328)(-1.25215,1.09278)(-1.24463,1.07842)%
(-1.18598,0.97026)(-1.12623,0.86833)(-1.06147,0.76948)(-1.00211,0.69104)(-0.98721,0.67332)%
(-0.94490,0.62749)(-0.88868,0.57712)(-0.83261,0.53896)(-0.77776,0.51316)(-0.72283,0.49836)%
(-0.66734,0.49403)(-0.65601,0.49440)(-0.61207,0.49956)(-0.55618,0.51423)(-0.49907,0.53765)%
(-0.44167,0.56867)(-0.38285,0.60710)(-0.34881,0.63196)(-0.32264,0.65221)(-0.26107,0.70313)%
(-0.21911,0.74001)(-0.19747,0.75959)(-0.13195,0.82063)(-0.10257,0.84865)(-0.06384,0.88595)%
(0.00284,0.95064)(0.08190,1.02705)(0.09966,1.04405)(0.16135,1.10230)(0.18867,1.12761)%
(0.24710,1.18056)(0.27076,1.20147)(0.34196,1.26234)(0.41727,1.32292)(0.45520,1.35178)%
(0.48236,1.37172)(0.54264,1.41362)(0.59731,1.44854)(0.63601,1.47126)(0.64726,1.47750)%
(0.69239,1.50081)(0.71682,1.51207)(0.73194,1.51846)(0.76617,1.53095)(0.79324,1.53818)%
(0.80337,1.53997)(0.81199,1.54091)(0.81894,1.54097)(0.81960,1.54091)%
\polyline(0.67308,1.69233)(0.67266,1.69294)(0.63947,1.74872)(0.62763,1.77003)(0.57900,1.86412)%
(0.52672,1.97472)(0.47088,2.09954)(0.44236,2.16451)(0.41253,2.23248)(0.35227,2.36680)%
(0.30651,2.46330)(0.28983,2.49674)(0.22664,2.61258)(0.18305,2.68043)(0.16218,2.70888)%
(0.09760,2.77855)(0.06322,2.80351)(0.03277,2.81823)(-0.03157,2.82604)(-0.05768,2.82021)%
(-0.09526,2.80294)(-0.15768,2.75232)(-0.18523,2.72197)(-0.21862,2.67928)(-0.27743,2.59069)%
(-0.33210,2.49658)(-0.38722,2.39488)(-0.43708,2.30119)(-0.48285,2.21747)(-0.52390,2.14728)%
(-0.55941,2.09281)(-0.58827,2.05485)(-0.59929,2.04242)(-0.60596,2.03579)%
\polyline(-0.60596,2.03579)(-0.60845,2.03331)(-0.61704,2.02629)%
\polyline(0.33522,1.25657)(0.36370,1.22607)(0.38155,1.20765)(0.42147,1.16781)(0.47746,1.11581)%
(0.53190,1.07055)(0.54793,1.05837)(0.58518,1.03231)(0.63729,1.00168)(0.68835,0.97905)%
(0.73858,0.96469)(0.76487,0.96054)(0.78827,0.95889)(0.83753,0.96200)(0.88644,0.97427)%
(0.93523,0.99588)(0.98411,1.02709)(1.02863,1.06375)(1.08305,1.11886)(1.13331,1.17923)%
(1.18418,1.24869)(1.23569,1.32647)(1.25669,1.36001)(1.28787,1.41151)(1.34039,1.50186)%
(1.39285,1.59525)(1.44489,1.68931)(1.45881,1.71450)(1.49634,1.78212)(1.54646,1.87095)%
(1.59468,1.95376)(1.64054,2.02917)(1.68359,2.09626)(1.72342,2.15458)(1.75960,2.20403)%
(1.77312,2.22117)%
\polyline(1.77312,2.22117)(1.79168,2.24473)%
\polyline(1.79168,2.24473)(1.81869,2.27640)%
\polyline(1.81869,2.27640)(1.81918,2.27697)(1.82381,2.28214)(1.83938,2.29883)%
\polyline(-2.12002,2.63309)(-2.07910,2.58389)(-2.03819,2.53109)(-1.99727,2.47467)%
(-1.95636,2.41469)(-1.91544,2.35127)(-1.87453,2.28460)(-1.83361,2.21498)(-1.79270,2.14277)%
(-1.75178,2.06844)(-1.71087,1.99251)(-1.66995,1.91558)(-1.62904,1.83832)(-1.58926,1.76357)%
\polyline(-1.58926,1.76357)(-1.58812,1.76144)(-1.54721,1.68569)(-1.54514,1.68197)%
\polyline(-1.54514,1.68197)(-1.51440,1.62647)%
\polyline(-1.51440,1.62647)(-1.50629,1.61186)(-1.47954,1.56533)%
\polyline(-1.47954,1.56533)(-1.46538,1.54070)(-1.42446,1.47301)(-1.38355,1.40952)%
(-1.34263,1.35091)(-1.30172,1.29785)(-1.26080,1.25088)(-1.21989,1.21050)(-1.17897,1.17708)%
(-1.13806,1.15090)(-1.09714,1.13215)(-1.05623,1.12088)(-1.01531,1.11704)(-0.97440,1.12044)%
(-0.93348,1.13080)(-0.89257,1.14775)(-0.85165,1.17080)(-0.81074,1.19938)(-0.76982,1.23286)%
(-0.72891,1.27053)(-0.68799,1.31166)(-0.64708,1.35548)(-0.60616,1.40120)(-0.56525,1.44806)%
(-0.52433,1.49530)(-0.48342,1.54221)(-0.44250,1.58812)(-0.40159,1.63243)(-0.36067,1.67462)%
(-0.31976,1.71423)(-0.27884,1.75087)(-0.23793,1.78428)(-0.19701,1.81424)(-0.18739,1.82045)%
\polyline(-0.18739,1.82045)(-0.15610,1.84063)(-0.11518,1.86342)(-0.07427,1.88260)%
\polyline(2.12002,2.58228)(2.07613,2.54069)(2.03225,2.49761)(1.98836,2.45333)(1.94448,2.40819)%
(1.90059,2.36256)(1.85671,2.31688)(1.84131,2.30099)%
\polyline(1.84131,2.30099)(1.82322,2.28231)%
\polyline(1.82322,2.28231)(1.81282,2.27157)(1.76893,2.22708)(1.76837,2.22653)%
\polyline(1.76837,2.22653)(1.72505,2.18383)(1.68116,2.14220)(1.63728,2.10255)(1.59339,2.06516)%
(1.54950,2.03022)(1.50575,1.99795)%
\polyline(1.50575,1.99795)(1.50562,1.99785)(1.46644,1.97127)%
\polyline(1.46644,1.97127)(1.46173,1.96809)(1.41785,1.94088)(1.37396,1.91608)(1.33008,1.89349)%
(1.28619,1.87287)(1.24230,1.85392)(1.19842,1.83637)(1.15453,1.81995)(1.15247,1.81921)%
\polyline(1.15247,1.81921)(1.11065,1.80440)(1.06676,1.78955)(1.02287,1.77530)(0.97899,1.76157)%
(0.93510,1.74843)(0.89122,1.73598)(0.84733,1.72442)(0.80345,1.71398)(0.75956,1.70494)%
(0.71567,1.69757)(0.67308,1.69233)%
\polyline(0.67308,1.69233)(0.67179,1.69218)(0.62790,1.68899)(0.58402,1.68822)(0.54013,1.68999)%
(0.49625,1.69437)(0.45236,1.70132)(0.40847,1.71073)(0.36459,1.72239)(0.32070,1.73603)%
(0.27682,1.75130)(0.23293,1.76780)(0.18904,1.78512)(0.14516,1.80281)(0.10127,1.82045)%
(0.05739,1.83761)(0.01350,1.85390)(-0.03038,1.86899)(-0.07427,1.88260)%
\polyline(0.07427,3.33277)(0.11518,3.28358)(0.15610,3.23077)(0.19701,3.17436)(0.23793,3.11438)%
(0.27884,3.05095)(0.31976,2.98429)(0.36067,2.91467)(0.40159,2.84246)(0.44250,2.76812)%
(0.48342,2.69219)(0.52433,2.61527)(0.56525,2.53800)(0.60616,2.46112)(0.64708,2.38538)%
(0.68799,2.31155)(0.72891,2.24039)(0.76982,2.17270)(0.81074,2.10920)(0.85165,2.05060)%
(0.89257,1.99753)(0.93348,1.95057)(0.97440,1.91018)(1.01531,1.87677)(1.05623,1.85059)%
(1.09714,1.83184)(1.13806,1.82057)(1.15247,1.81921)%
\polyline(1.46538,1.97021)(1.50629,2.01135)(1.54721,2.05517)(1.58812,2.10089)(1.62904,2.14775)%
(1.66995,2.19499)(1.71087,2.24190)(1.75178,2.28781)(1.79270,2.33212)(1.83361,2.37431)%
(1.87453,2.41392)(1.91544,2.45056)(1.95636,2.48397)(1.99727,2.51393)(2.03819,2.54032)%
(2.07910,2.56311)(2.12002,2.58228)%
\polyline(0.07427,3.33277)(0.03038,3.29118)(-0.01350,3.24810)(-0.05739,3.20382)(-0.10127,3.15868)%
(-0.14516,3.11306)(-0.18904,3.06738)(-0.23293,3.02206)(-0.27682,2.97758)(-0.32070,2.93432)%
(-0.36459,2.89269)(-0.40847,2.85304)(-0.45236,2.81565)(-0.49625,2.78071)(-0.54013,2.74835)%
(-0.58402,2.71858)(-0.62790,2.69137)(-0.67179,2.66657)(-0.71567,2.64398)(-0.75956,2.62336)%
(-0.80345,2.60441)(-0.84733,2.58687)(-0.89122,2.57044)(-0.93510,2.55489)(-0.97899,2.54005)%
(-1.02287,2.52579)(-1.06676,2.51206)(-1.11065,2.49892)(-1.15453,2.48648)(-1.19842,2.47492)%
(-1.24230,2.46447)(-1.28619,2.45543)(-1.33008,2.44806)(-1.37396,2.44267)(-1.41785,2.43948)%
(-1.46173,2.43871)(-1.50562,2.44048)(-1.54950,2.44486)(-1.59339,2.45181)(-1.63728,2.46122)%
(-1.68116,2.47288)(-1.72505,2.48652)(-1.76893,2.50179)(-1.81282,2.51829)(-1.85671,2.53561)%
(-1.90059,2.55331)(-1.94448,2.57094)(-1.98836,2.58810)(-2.03225,2.60439)(-2.07613,2.61948)%
(-2.12002,2.63309)%
\polyline(2.19429,0.69969)(-2.19429,-0.69969)%
\polyline(-2.04575,0.75049)(2.04575,-0.75049)%
\polyline(0.00000,-1.54000)(0.00000,0.94788)%
\polyline(0.00000,2.82221)(0.00000,3.26136)%
\polyline(0.00000,3.26136)(0.00000,3.75877)%
\polyline(0.81960,1.54091)(0.80870,1.54449)(0.80092,1.54834)(0.77806,1.56386)(0.74826,1.59173)%
(0.72920,1.61343)(0.71299,1.63403)(0.67308,1.69233)%
\polyline(-0.61704,2.02629)(-0.60762,2.03018)(-0.60222,2.03171)(-0.58522,2.03453)%
(-0.55353,2.03440)(-0.51444,2.02762)(-0.46896,2.01296)(-0.41754,1.98949)(-0.36013,1.95618)%
(-0.32884,1.93529)(-0.29683,1.91214)(-0.22734,1.85627)(-0.21405,1.84479)(-0.15126,1.78746)%
(-0.12082,1.75798)(-0.06755,1.70412)(-0.03767,1.67278)(0.02493,1.60498)(0.03881,1.58961)%
(0.11015,1.50930)(0.12813,1.48882)(0.17772,1.43215)(0.24225,1.35875)(0.30411,1.28986)%
(0.33522,1.25657)%
\polyline(1.15247,1.81921)(1.17897,1.81673)(1.21989,1.82013)(1.26080,1.83049)(1.30172,1.84744)%
(1.34263,1.87048)(1.38355,1.89907)(1.42446,1.93254)(1.46538,1.97021)%
\polyline(0.00000,0.94788)(0.00000,2.82221)%
\linethickness{0.003in}
\put(-0.00412,3.25505){\circle*{0.120000}}\put(1.13984,1.80891){\circle*{0.120000}}
\put(1.47799,1.97439){\circle*{0.120000}}\put(0.67218,1.68660){\circle*{0.120000}}
\put(-0.19119,1.82330){\circle*{0.120000}}\put(0.33403,1.24772){\circle*{0.120000}}
\put(0.00307,2.81617){\circle*{0.120000}}\put(-0.02705,0.94323){\circle*{0.120000}}
\put(-0.61896,2.04416){\circle*{0.120000}}\put(0.81458,1.57179){\circle*{0.120000}}
\linethickness{0.008in}
\end{picture}}%

%% file: fig/p001.tex
{\unitlength=1cm%
\begin{picture}%
(9.94,5.4)(-5.19,-2.08)%
\linethickness{0.008in}
\polyline(-1.46071,-0.83234)(-1.47531,-0.84067)(-1.48992,-0.84899)(-1.50453,-0.85731)%
(-1.51913,-0.86564)(-1.53374,-0.87396)(-1.54835,-0.88229)(-1.56296,-0.89061)(-1.57756,-0.89893)%
(-1.59217,-0.90726)(-1.60678,-0.91558)(-1.62138,-0.92390)(-1.63599,-0.93223)(-1.65060,-0.94055)%
(-1.66521,-0.94887)(-1.67981,-0.95720)(-1.69442,-0.96552)(-1.70903,-0.97384)(-1.72363,-0.98217)%
(-1.73824,-0.99049)(-1.75285,-0.99881)(-1.76745,-1.00714)(-1.78206,-1.01546)(-1.79667,-1.02378)%
(-1.81128,-1.03211)(-1.82588,-1.04043)(-1.84049,-1.04875)(-1.85510,-1.05708)(-1.86970,-1.06540)%
(-1.88431,-1.07372)(-1.89892,-1.08205)(-1.91353,-1.09037)(-1.92813,-1.09869)(-1.94274,-1.10702)%
(-1.95735,-1.11534)(-1.97195,-1.12366)(-1.98656,-1.13199)(-2.00117,-1.14031)(-2.01577,-1.14864)%
(-2.03038,-1.15696)(-2.04499,-1.16528)(-2.05960,-1.17361)(-2.07420,-1.18193)(-2.08881,-1.19025)%
(-2.10342,-1.19858)(-2.11802,-1.20690)(-2.13263,-1.21522)(-2.14724,-1.22355)(-2.16185,-1.23187)%
(-2.17645,-1.24019)(-2.19106,-1.24852)%
\polyline(-1.46071,-0.83234)(-1.47531,-0.84067)(-1.48992,-0.84899)(-1.50453,-0.85731)%
(-1.51913,-0.86564)(-1.53374,-0.87396)(-1.54835,-0.88229)(-1.56296,-0.89061)(-1.57756,-0.89893)%
(-1.59217,-0.90726)(-1.60678,-0.91558)(-1.62138,-0.92390)(-1.63599,-0.93223)(-1.65060,-0.94055)%
(-1.66521,-0.94887)(-1.67981,-0.95720)(-1.69442,-0.96552)(-1.70903,-0.97384)(-1.72363,-0.98217)%
(-1.73824,-0.99049)(-1.75285,-0.99881)(-1.76745,-1.00714)(-1.78206,-1.01546)(-1.79667,-1.02378)%
(-1.81128,-1.03211)(-1.82588,-1.04043)(-1.84049,-1.04875)(-1.85510,-1.05708)(-1.86970,-1.06540)%
(-1.88431,-1.07372)(-1.89892,-1.08205)(-1.91353,-1.09037)(-1.92813,-1.09869)(-1.94274,-1.10702)%
(-1.95735,-1.11534)(-1.97195,-1.12366)(-1.98656,-1.13199)(-2.00117,-1.14031)(-2.01577,-1.14864)%
(-2.03038,-1.15696)(-2.04499,-1.16528)(-2.05960,-1.17361)(-2.07420,-1.18193)(-2.08881,-1.19025)%
(-2.10342,-1.19858)(-2.11802,-1.20690)(-2.13263,-1.21522)(-2.14724,-1.22355)(-2.16185,-1.23187)%
(-2.17645,-1.24019)(-2.19106,-1.24852)%
\end{picture}}%

%% file: fig/p004.tex
{\unitlength=1cm%
\begin{picture}%
(9.94,5.4)(-5.19,-2.08)%
\linethickness{0.008in}
\polyline(0.71054,-1.14160)(0.70165,-1.14273)(0.69211,-1.14392)(0.69188,-1.14395)%
\polyline(0.69188,-1.14395)(0.68265,-1.14507)(0.67466,-1.14601)%
\polyline(0.67466,-1.14601)(0.67325,-1.14618)(0.66392,-1.14725)(0.66373,-1.14727)%
\polyline(0.66373,-1.14727)(0.65466,-1.14828)(0.64548,-1.14927)(0.63643,-1.15021)%
(0.62747,-1.15112)(0.61857,-1.15199)(0.60973,-1.15281)(0.60097,-1.15360)(0.59227,-1.15435)%
(0.58365,-1.15507)(0.57508,-1.15574)(0.56659,-1.15638)(0.55816,-1.15698)(0.54980,-1.15754)%
(0.54150,-1.15807)(0.53328,-1.15856)(0.52512,-1.15901)(0.51702,-1.15943)(0.50903,-1.15980)%
(0.50113,-1.16013)(0.49328,-1.16044)(0.48551,-1.16071)(0.47855,-1.16092)%
\polyline(0.47855,-1.16092)(0.47779,-1.16094)(0.47015,-1.16113)(0.46369,-1.16127)%
\polyline(0.46369,-1.16127)(0.46257,-1.16129)(0.45504,-1.16141)(0.44832,-1.16149)%
\polyline(0.44832,-1.16149)(0.44759,-1.16150)(0.44020,-1.16155)(0.43323,-1.16156)%
\polyline(0.43323,-1.16156)(0.43287,-1.16157)(0.42560,-1.16154)(0.41841,-1.16149)%
(0.41269,-1.16141)%
\polyline(0.41269,-1.16141)(0.41126,-1.16140)(0.40420,-1.16127)(0.39768,-1.16112)%
\polyline(0.39768,-1.16112)(0.39028,-1.16092)(0.38478,-1.16073)%
\polyline(0.38478,-1.16073)(0.38344,-1.16068)(0.37665,-1.16042)(0.36993,-1.16013)%
(0.36427,-1.15984)%
\polyline(0.36427,-1.15984)(0.36328,-1.15979)(0.35668,-1.15942)(0.35177,-1.15912)%
\polyline(0.35177,-1.15912)(0.35014,-1.15902)(0.34367,-1.15859)(0.33885,-1.15823)%
\polyline(0.33885,-1.15823)(0.33726,-1.15812)(0.33091,-1.15761)(0.32619,-1.15721)%
\polyline(0.71126,-1.14150)(0.70165,-1.14273)(0.69211,-1.14392)(0.68265,-1.14507)%
(0.67325,-1.14618)(0.66392,-1.14725)(0.65466,-1.14828)(0.64548,-1.14927)(0.63643,-1.15021)%
(0.62747,-1.15112)(0.61857,-1.15199)(0.60973,-1.15281)(0.60097,-1.15360)(0.59227,-1.15435)%
(0.58365,-1.15507)(0.57508,-1.15574)(0.56659,-1.15638)(0.55816,-1.15698)(0.54980,-1.15754)%
(0.54150,-1.15807)(0.53328,-1.15856)(0.52512,-1.15901)(0.51702,-1.15943)(0.50903,-1.15980)%
(0.50113,-1.16013)(0.49328,-1.16044)(0.48551,-1.16071)(0.47779,-1.16094)(0.47015,-1.16113)%
(0.46257,-1.16129)(0.45504,-1.16141)(0.44759,-1.16150)(0.44020,-1.16155)(0.43287,-1.16157)%
(0.42560,-1.16154)(0.41841,-1.16149)(0.41126,-1.16140)(0.40420,-1.16127)(0.39768,-1.16112)%
(0.39028,-1.16092)(0.38344,-1.16068)(0.37665,-1.16042)(0.36993,-1.16013)(0.36328,-1.15979)%
(0.35668,-1.15942)(0.35014,-1.15902)(0.34367,-1.15859)(0.33726,-1.15812)(0.33091,-1.15761)%
(0.32462,-1.15708)(0.31840,-1.15650)%
\polyline(-2.19106,-1.24852)(-2.12330,-1.25243)(-2.05499,-1.25602)(-1.98615,-1.25927)%
(-1.91680,-1.26220)(-1.84695,-1.26480)(-1.77663,-1.26704)(-1.70586,-1.26897)(-1.63465,-1.27054)%
(-1.56303,-1.27179)(-1.49102,-1.27269)(-1.41864,-1.27323)(-1.34589,-1.27343)(-1.27283,-1.27329)%
(-1.19945,-1.27280)(-1.12578,-1.27196)(-1.05184,-1.27078)(-0.97765,-1.26923)(-0.90322,-1.26734)%
(-0.82861,-1.26510)(-0.75379,-1.26251)(-0.67882,-1.25955)(-0.60370,-1.25625)(-0.52846,-1.25259)%
(-0.45312,-1.24858)(-0.37770,-1.24422)(-0.30221,-1.23949)(-0.22670,-1.23442)(-0.15116,-1.22899)%
(-0.07562,-1.22321)(-0.00011,-1.21708)(0.07534,-1.21059)(0.15073,-1.20376)(0.22604,-1.19657)%
(0.30123,-1.18904)(0.37629,-1.18116)(0.45120,-1.17292)(0.47471,-1.17023)%
\polyline(0.47471,-1.17023)(0.51102,-1.16607)%
\polyline(0.51102,-1.16607)(0.52592,-1.16436)%
\polyline(0.52594,-1.16436)(0.60047,-1.15544)%
\polyline(0.60047,-1.15544)(0.67481,-1.14618)(0.71030,-1.14158)%
\polyline(0.71030,-1.14158)(0.74890,-1.13658)(0.82274,-1.12664)(0.89630,-1.11636)%
(0.96956,-1.10575)(1.04251,-1.09481)(1.11513,-1.08352)(1.18738,-1.07192)(1.25925,-1.05999)%
(1.33074,-1.04773)(1.40180,-1.03515)(1.47243,-1.02225)%
\polyline(1.01352,-1.20070)(1.02270,-1.19712)(1.03188,-1.19356)(1.04105,-1.18998)%
(1.05023,-1.18642)(1.05941,-1.18285)(1.06859,-1.17928)(1.07777,-1.17571)(1.08695,-1.17215)%
(1.09613,-1.16857)(1.10531,-1.16501)(1.11448,-1.16144)(1.12366,-1.15787)(1.13284,-1.15431)%
(1.14201,-1.15073)(1.15119,-1.14716)(1.16037,-1.14359)(1.16956,-1.14003)(1.17873,-1.13646)%
(1.18791,-1.13289)(1.19709,-1.12932)(1.20626,-1.12576)(1.21544,-1.12218)(1.22462,-1.11861)%
(1.23380,-1.11504)(1.24297,-1.11147)(1.25215,-1.10791)(1.26133,-1.10434)(1.27051,-1.10077)%
(1.27969,-1.09720)(1.28887,-1.09363)(1.29805,-1.09007)(1.30722,-1.08649)(1.31640,-1.08293)%
(1.32558,-1.07936)(1.33475,-1.07579)(1.34393,-1.07222)(1.35311,-1.06865)(1.36229,-1.06508)%
(1.37147,-1.06152)(1.38065,-1.05794)(1.38983,-1.05438)(1.39900,-1.05080)(1.40818,-1.04724)%
(1.41736,-1.04367)(1.42653,-1.04010)(1.43572,-1.03653)(1.44489,-1.03297)(1.45407,-1.02939)%
(1.46325,-1.02583)(1.47243,-1.02225)%
\polyline(-1.46071,-0.83234)(-1.41555,-0.84504)(-1.37006,-0.85753)(-1.32426,-0.86983)%
(-1.27813,-0.88194)(-1.23169,-0.89383)(-1.18498,-0.90553)(-1.13796,-0.91701)(-1.09067,-0.92830)%
(-1.04312,-0.93937)(-0.99530,-0.95023)(-0.94724,-0.96089)(-0.89894,-0.97134)(-0.85041,-0.98157)%
(-0.80167,-0.99158)(-0.75271,-1.00138)(-0.70356,-1.01096)(-0.65422,-1.02032)(-0.60469,-1.02947)%
(-0.55499,-1.03839)(-0.50515,-1.04709)(-0.45515,-1.05556)(-0.40500,-1.06382)(-0.35474,-1.07185)%
(-0.30434,-1.07964)(-0.25384,-1.08722)(-0.20325,-1.09457)(-0.15256,-1.10169)(-0.10179,-1.10857)%
(-0.05096,-1.11523)(-0.00008,-1.12166)(0.05086,-1.12785)(0.10183,-1.13380)(0.15283,-1.13954)%
(0.20384,-1.14503)(0.25486,-1.15028)(0.30587,-1.15531)(0.33122,-1.15769)%
\polyline(0.47471,-1.17023)(0.50965,-1.17303)%
\polyline(0.50967,-1.17304)(0.56048,-1.17687)(0.61124,-1.18048)(0.66192,-1.18383)%
(0.71250,-1.18696)(0.76299,-1.18985)(0.81336,-1.19249)(0.86361,-1.19490)(0.91373,-1.19707)%
(0.96371,-1.19900)(1.01352,-1.20070)%
\polyline(-1.46071,-0.83234)(-1.47531,-0.84067)(-1.48992,-0.84899)(-1.50453,-0.85731)%
(-1.51913,-0.86564)(-1.53374,-0.87396)(-1.54835,-0.88229)(-1.56296,-0.89061)(-1.57756,-0.89893)%
(-1.59217,-0.90726)(-1.60678,-0.91558)(-1.62138,-0.92390)(-1.63599,-0.93223)(-1.65060,-0.94055)%
(-1.66521,-0.94887)(-1.67981,-0.95720)(-1.69442,-0.96552)(-1.70903,-0.97384)(-1.72363,-0.98217)%
(-1.73824,-0.99049)(-1.75285,-0.99881)(-1.76745,-1.00714)(-1.78206,-1.01546)(-1.79667,-1.02378)%
(-1.81128,-1.03211)(-1.82588,-1.04043)(-1.84049,-1.04875)(-1.85510,-1.05708)(-1.86970,-1.06540)%
(-1.88431,-1.07372)(-1.89892,-1.08205)(-1.91353,-1.09037)(-1.92813,-1.09869)(-1.94274,-1.10702)%
(-1.95735,-1.11534)(-1.97195,-1.12366)(-1.98656,-1.13199)(-2.00117,-1.14031)(-2.01577,-1.14864)%
(-2.03038,-1.15696)(-2.04499,-1.16528)(-2.05960,-1.17361)(-2.07420,-1.18193)(-2.08881,-1.19025)%
(-2.10342,-1.19858)(-2.11802,-1.20690)(-2.13263,-1.21522)(-2.14724,-1.22355)(-2.16185,-1.23187)%
(-2.17645,-1.24019)(-2.19106,-1.24852)%
\polyline(-1.58243,-0.90171)(-1.53351,-0.91293)(-1.48422,-0.92395)(-1.43456,-0.93474)%
(-1.38458,-0.94531)(-1.33424,-0.95565)(-1.28359,-0.96578)(-1.23261,-0.97567)(-1.18134,-0.98533)%
(-1.12977,-0.99477)(-1.07792,-1.00398)(-1.02580,-1.01295)(-0.97343,-1.02168)(-0.92082,-1.03019)%
(-0.86797,-1.03845)(-0.81489,-1.04648)(-0.76161,-1.05426)(-0.70812,-1.06181)(-0.65445,-1.06911)%
(-0.60060,-1.07617)(-0.54659,-1.08299)(-0.49242,-1.08956)(-0.43812,-1.09589)(-0.38369,-1.10197)%
(-0.32914,-1.10781)(-0.27449,-1.11339)(-0.21973,-1.11872)(-0.16492,-1.12381)(-0.11002,-1.12864)%
(-0.05508,-1.13323)(-0.00009,-1.13756)(0.05494,-1.14164)(0.10998,-1.14547)(0.16502,-1.14904)%
(0.22007,-1.15237)(0.27509,-1.15543)(0.33009,-1.15824)(0.36872,-1.16004)(0.38351,-1.16073)%
\polyline(0.51359,-1.16578)(0.54953,-1.16696)(0.56319,-1.16734)(0.60419,-1.16850)%
(0.65874,-1.16979)(0.71319,-1.17082)(0.76750,-1.17161)(0.82168,-1.17212)(0.87570,-1.17240)%
(0.92955,-1.17242)(0.98323,-1.17218)(1.03672,-1.17170)(1.09001,-1.17095)%
\polyline(-1.70416,-0.97107)(-1.65147,-0.98084)(-1.59838,-0.99036)(-1.54488,-0.99965)%
(-1.49101,-1.00869)(-1.43679,-1.01748)(-1.38219,-1.02604)(-1.32726,-1.03433)(-1.27200,-1.04238)%
(-1.21642,-1.05017)(-1.16054,-1.05771)(-1.10437,-1.06500)(-1.04793,-1.07203)(-0.99122,-1.07881)%
(-0.93426,-1.08532)(-0.87707,-1.09158)(-0.81965,-1.09757)(-0.76202,-1.10330)(-0.70420,-1.10876)%
(-0.64619,-1.11396)(-0.58803,-1.11890)(-0.52970,-1.12356)(-0.47123,-1.12796)(-0.41264,-1.13210)%
(-0.35394,-1.13596)(-0.29513,-1.13955)(-0.23624,-1.14288)(-0.17727,-1.14593)(-0.11825,-1.14872)%
(-0.05919,-1.15122)(-0.00009,-1.15346)(0.05902,-1.15544)(0.11813,-1.15713)(0.17723,-1.15855)%
(0.23630,-1.15970)(0.29534,-1.16058)(0.35431,-1.16118)(0.41322,-1.16151)(0.42923,-1.16153)%
\polyline(0.53078,-1.16137)(0.58940,-1.16088)(0.64790,-1.16013)(0.70626,-1.15911)%
(0.76447,-1.15780)(0.82250,-1.15624)(0.88036,-1.15441)(0.93803,-1.15230)(0.99550,-1.14993)%
(1.05273,-1.14729)(1.10973,-1.14439)(1.16649,-1.14122)%
\polyline(-1.82588,-1.04043)(-1.76942,-1.04873)(-1.71252,-1.05678)(-1.65520,-1.06455)%
(-1.59746,-1.07206)(-1.53932,-1.07931)(-1.48080,-1.08629)(-1.42191,-1.09299)(-1.36267,-1.09942)%
(-1.30307,-1.10558)(-1.24317,-1.11146)(-1.18293,-1.11706)(-1.12242,-1.12239)(-1.06162,-1.12743)%
(-1.00056,-1.13219)(-0.93925,-1.13667)(-0.87770,-1.14087)(-0.81593,-1.14478)(-0.75396,-1.14841)%
(-0.69180,-1.15175)(-0.62947,-1.15480)(-0.56698,-1.15756)(-0.50435,-1.16003)(-0.44159,-1.16222)%
(-0.37873,-1.16411)(-0.31577,-1.16572)(-0.25273,-1.16703)(-0.18963,-1.16805)(-0.12648,-1.16878)%
(-0.06329,-1.16922)(-0.00009,-1.16937)(0.06310,-1.16922)(0.12628,-1.16878)(0.18943,-1.16805)%
(0.25253,-1.16703)(0.31557,-1.16572)(0.37853,-1.16412)(0.44140,-1.16222)(0.49016,-1.16053)%
(0.50152,-1.16013)%
\polyline(0.56679,-1.15757)(0.62928,-1.15481)(0.67364,-1.15263)(0.69161,-1.15175)%
(0.75377,-1.14842)(0.81574,-1.14479)(0.87751,-1.14088)(0.93906,-1.13669)(1.00037,-1.13221)%
(1.06143,-1.12744)(1.12223,-1.12240)(1.18276,-1.11708)(1.24297,-1.11148)%
\polyline(-1.94761,-1.10979)(-1.88739,-1.11663)(-1.82668,-1.12319)(-1.76552,-1.12946)%
(-1.70390,-1.13544)(-1.64186,-1.14114)(-1.57941,-1.14654)(-1.51655,-1.15165)(-1.45333,-1.15646)%
(-1.38973,-1.16098)(-1.32578,-1.16520)(-1.26151,-1.16912)(-1.19691,-1.17274)(-1.13202,-1.17605)%
(-1.06686,-1.17907)(-1.00142,-1.18177)(-0.93574,-1.18417)(-0.86983,-1.18626)(-0.80372,-1.18805)%
(-0.73741,-1.18953)(-0.67091,-1.19070)(-0.60427,-1.19155)(-0.53747,-1.19211)(-0.47055,-1.19234)%
(-0.40353,-1.19227)(-0.33641,-1.19189)(-0.26923,-1.19118)(-0.20198,-1.19017)(-0.13470,-1.18885)%
(-0.06741,-1.18722)(-0.00011,-1.18527)(0.06718,-1.18301)(0.13443,-1.18044)(0.20164,-1.17756)%
(0.26877,-1.17437)(0.33581,-1.17086)(0.40276,-1.16706)(0.42173,-1.16588)(0.46958,-1.16294)%
(0.53626,-1.15850)(0.56089,-1.15676)(0.57551,-1.15572)%
\polyline(0.63113,-1.15162)(0.66915,-1.14873)(0.73532,-1.14338)(0.80127,-1.13773)%
(0.86701,-1.13178)(0.93251,-1.12552)(0.99774,-1.11896)(1.06271,-1.11211)(1.12738,-1.10496)%
(1.19174,-1.09751)(1.25577,-1.08977)(1.31946,-1.08173)%
\polyline(-2.06933,-1.17915)(-2.00534,-1.18453)(-1.94083,-1.18960)(-1.87583,-1.19437)%
(-1.81035,-1.19882)(-1.74440,-1.20297)(-1.67802,-1.20679)(-1.61120,-1.21031)(-1.54399,-1.21351)%
(-1.47638,-1.21638)(-1.40840,-1.21894)(-1.34007,-1.22117)(-1.27140,-1.22309)(-1.20243,-1.22467)%
(-1.13315,-1.22593)(-1.06360,-1.22687)(-0.99380,-1.22748)(-0.92374,-1.22775)(-0.85347,-1.22770)%
(-0.78301,-1.22732)(-0.71235,-1.22660)(-0.64154,-1.22556)(-0.57058,-1.22417)(-0.49951,-1.22247)%
(-0.42832,-1.22042)(-0.35706,-1.21805)(-0.28572,-1.21534)(-0.21434,-1.21229)(-0.14293,-1.20892)%
(-0.07151,-1.20521)(-0.00011,-1.20118)(0.07127,-1.19680)(0.14258,-1.19210)(0.21383,-1.18707)%
(0.28500,-1.18170)(0.35605,-1.17602)(0.42698,-1.17000)(0.44843,-1.16807)(0.47312,-1.16585)%
(0.49775,-1.16365)(0.56837,-1.15698)(0.63688,-1.15017)(0.65512,-1.14828)%
\polyline(0.67508,-1.14620)(0.68886,-1.14476)(0.70903,-1.14266)(0.77903,-1.13501)%
(0.84879,-1.12705)(0.91828,-1.11876)(0.98751,-1.11016)(1.05643,-1.10125)(1.12504,-1.09201)%
(1.19331,-1.08247)(1.26123,-1.07262)(1.32878,-1.06246)(1.39594,-1.05200)%
\polyline(0.15992,-1.14032)(0.16146,-1.14142)(0.16299,-1.14252)(0.16452,-1.14363)%
(0.16605,-1.14473)(0.16758,-1.14583)(0.16912,-1.14694)(0.17064,-1.14805)(0.17218,-1.14915)%
(0.17371,-1.15026)(0.17524,-1.15136)(0.17677,-1.15246)(0.17831,-1.15357)(0.17983,-1.15468)%
(0.18137,-1.15578)(0.18290,-1.15688)(0.18444,-1.15799)(0.18596,-1.15909)(0.18750,-1.16020)%
(0.18903,-1.16130)(0.19056,-1.16241)(0.19209,-1.16351)(0.19362,-1.16462)(0.19515,-1.16572)%
(0.19669,-1.16682)(0.19822,-1.16793)(0.19975,-1.16904)(0.20128,-1.17014)(0.20281,-1.17124)%
(0.20434,-1.17235)(0.20587,-1.17345)(0.20741,-1.17456)(0.20894,-1.17566)(0.21047,-1.17677)%
(0.21200,-1.17787)(0.21353,-1.17898)(0.21506,-1.18008)(0.21660,-1.18118)(0.21813,-1.18229)%
(0.21966,-1.18340)(0.22119,-1.18450)(0.22272,-1.18560)(0.22425,-1.18671)(0.22579,-1.18781)%
(0.22732,-1.18892)(0.22885,-1.19003)(0.23038,-1.19112)(0.23192,-1.19223)(0.23344,-1.19334)%
(0.23498,-1.19445)(0.23651,-1.19554)%
\end{picture}}%

%% file: fig/p007.tex
{\unitlength=1cm%
\begin{picture}%
(9.94,5.4)(-5.19,-2.08)%
\linethickness{0.008in}
\polyline(0.71067,-1.14158)(0.70185,-1.14271)(0.69237,-1.14389)(0.69213,-1.14391)%
\polyline(0.69213,-1.14391)(0.68296,-1.14504)(0.67625,-1.14582)%
\polyline(0.67625,-1.14582)(0.67362,-1.14614)(0.66435,-1.14720)(0.65515,-1.14822)%
(0.64602,-1.14921)(0.63696,-1.15016)(0.62797,-1.15106)(0.61905,-1.15194)(0.61018,-1.15277)%
(0.60140,-1.15356)(0.59267,-1.15432)(0.58402,-1.15503)(0.57544,-1.15572)(0.56692,-1.15635)%
(0.55846,-1.15697)(0.55008,-1.15753)(0.54177,-1.15805)(0.53352,-1.15855)(0.52533,-1.15900)%
(0.51722,-1.15941)(0.50918,-1.15979)(0.50119,-1.16013)(0.49328,-1.16044)(0.48620,-1.16068)%
\polyline(0.48620,-1.16068)(0.48544,-1.16071)(0.47776,-1.16094)(0.47089,-1.16111)%
\polyline(0.47089,-1.16111)(0.47015,-1.16113)(0.46260,-1.16129)(0.45585,-1.16140)%
\polyline(0.45585,-1.16140)(0.45511,-1.16141)(0.44769,-1.16150)(0.44105,-1.16154)%
\polyline(0.44105,-1.16154)(0.44032,-1.16155)(0.43303,-1.16156)(0.42579,-1.16154)%
\polyline(0.42579,-1.16154)(0.41862,-1.16149)(0.41151,-1.16140)(0.40551,-1.16130)%
\polyline(0.40551,-1.16130)(0.40446,-1.16128)(0.39747,-1.16111)(0.39055,-1.16092)%
\polyline(0.39055,-1.16092)(0.38369,-1.16069)(0.38028,-1.16056)%
\polyline(0.38028,-1.16056)(0.37690,-1.16044)(0.37016,-1.16013)%
\polyline(0.37016,-1.16013)(0.36349,-1.15980)(0.35687,-1.15944)(0.35033,-1.15903)%
(0.34384,-1.15860)(0.33742,-1.15813)(0.33455,-1.15790)%
\polyline(0.33455,-1.15790)(0.33106,-1.15763)(0.32476,-1.15708)(0.32086,-1.15673)%
\polyline(-2.19106,-1.24852)(-2.05586,-1.25597)(-1.91855,-1.26213)(-1.77931,-1.26697)%
(-1.63826,-1.27047)(-1.49559,-1.27263)(-1.35143,-1.27343)(-1.20597,-1.27287)(-1.05934,-1.27091)%
(-0.91172,-1.26758)(-0.76327,-1.26286)(-0.61417,-1.25673)(-0.46456,-1.24921)(-0.31464,-1.24029)%
(-0.16454,-1.22998)(-0.01445,-1.21827)(0.13547,-1.20518)(0.28506,-1.19070)(0.43415,-1.17484)%
(0.47445,-1.17016)%
\polyline(0.47445,-1.17016)(0.50101,-1.16708)%
\polyline(0.50101,-1.16708)(0.58256,-1.15761)%
\polyline(0.58257,-1.15761)(0.67625,-1.14582)%
\polyline(0.67625,-1.14582)(0.69236,-1.14380)%
\polyline(0.69236,-1.14380)(0.70994,-1.14159)%
\polyline(0.70994,-1.14159)(0.73016,-1.13904)(0.87678,-1.11912)(1.02223,-1.09788)%
(1.16638,-1.07533)(1.30907,-1.05149)(1.45012,-1.02637)(1.58941,-1.00000)(1.72679,-0.97241)%
(1.86208,-0.94360)(1.99516,-0.91361)(2.12590,-0.88248)(2.25413,-0.85022)(2.37973,-0.81686)%
(2.50258,-0.78244)(2.62253,-0.74697)(2.73947,-0.71052)(2.85327,-0.67309)(2.96383,-0.63475)%
(3.07102,-0.59550)(3.17476,-0.55540)(3.27492,-0.51448)(3.37140,-0.47278)(3.46413,-0.43036)%
(3.55300,-0.38723)(3.63793,-0.34345)(3.71885,-0.29906)(3.79568,-0.25410)(3.86835,-0.20862)%
(3.93681,-0.16267)(4.00098,-0.11626)(4.01671,-0.10396)%
\polyline(4.01671,-0.10396)(4.03820,-0.08716)%
\polyline(4.03820,-0.08716)(4.06080,-0.06949)%
\polyline(3.05420,-0.99973)(3.07432,-0.98112)(3.09446,-0.96252)(3.11459,-0.94391)%
(3.13472,-0.92532)(3.15486,-0.90671)(3.17499,-0.88810)(3.19512,-0.86950)(3.21525,-0.85089)%
(3.23538,-0.83229)(3.25551,-0.81369)(3.27565,-0.79508)(3.29578,-0.77647)(3.31591,-0.75787)%
(3.33604,-0.73926)(3.35617,-0.72066)(3.37631,-0.70206)(3.39644,-0.68345)(3.41657,-0.66484)%
(3.43670,-0.64624)(3.45684,-0.62763)(3.47697,-0.60902)(3.49710,-0.59043)(3.51724,-0.57182)%
(3.53736,-0.55322)(3.55750,-0.53461)(3.57764,-0.51600)(3.59776,-0.49740)(3.61790,-0.47879)%
(3.63802,-0.46019)(3.65816,-0.44159)(3.67830,-0.42298)(3.69842,-0.40437)(3.71856,-0.38577)%
(3.73869,-0.36716)(3.75882,-0.34857)(3.77896,-0.32996)(3.79909,-0.31135)(3.81922,-0.29275)%
(3.83935,-0.27414)(3.85948,-0.25553)(3.87962,-0.23693)(3.89975,-0.21833)(3.91988,-0.19972)%
(3.94001,-0.18112)(3.96014,-0.16251)(3.98027,-0.14390)(4.00041,-0.12531)(4.01964,-0.10753)%
\polyline(4.01964,-0.10753)(4.02054,-0.10670)(4.03931,-0.08935)%
\polyline(4.03931,-0.08935)(4.04067,-0.08809)(4.06080,-0.06949)%
\polyline(-1.46071,-0.83234)(-1.37064,-0.85738)(-1.27930,-0.88163)(-1.18675,-0.90509)%
(-1.09308,-0.92773)(-0.99834,-0.94956)(-0.90262,-0.97055)(-0.80600,-0.99070)(-0.70854,-1.01000)%
(-0.61034,-1.02843)(-0.51146,-1.04600)(-0.41200,-1.06268)(-0.31200,-1.07848)(-0.21158,-1.09338)%
(-0.11080,-1.10737)(-0.00975,-1.12045)(0.09150,-1.13262)(0.19286,-1.14386)(0.29425,-1.15418)%
(0.32508,-1.15703)%
\polyline(0.47445,-1.17016)(0.49677,-1.17202)%
\polyline(0.49679,-1.17202)(0.59775,-1.17954)(0.69843,-1.18612)(0.79871,-1.19174)%
(0.89852,-1.19644)(0.99777,-1.20019)(1.09639,-1.20299)(1.19428,-1.20485)(1.29136,-1.20577)%
(1.38754,-1.20576)(1.48275,-1.20480)(1.57690,-1.20292)(1.66990,-1.20010)(1.76165,-1.19636)%
(1.85211,-1.19171)(1.94116,-1.18614)(2.02872,-1.17967)(2.11473,-1.17230)(2.19909,-1.16405)%
(2.28172,-1.15491)(2.36255,-1.14490)(2.44149,-1.13404)(2.51845,-1.12232)(2.59338,-1.10978)%
(2.66618,-1.09641)(2.73677,-1.08223)(2.80510,-1.06726)(2.87108,-1.05151)(2.93463,-1.03499)%
(2.99570,-1.01773)(3.05420,-0.99973)%
\polyline(-1.46071,-0.83234)(-1.47531,-0.84067)(-1.48992,-0.84899)(-1.50453,-0.85731)%
(-1.51913,-0.86564)(-1.53374,-0.87396)(-1.54835,-0.88229)(-1.56296,-0.89061)(-1.57756,-0.89893)%
(-1.59217,-0.90726)(-1.60678,-0.91558)(-1.62138,-0.92390)(-1.63599,-0.93223)(-1.65060,-0.94055)%
(-1.66521,-0.94887)(-1.67981,-0.95720)(-1.69442,-0.96552)(-1.70903,-0.97384)(-1.72363,-0.98217)%
(-1.73824,-0.99049)(-1.75285,-0.99881)(-1.76745,-1.00714)(-1.78206,-1.01546)(-1.79667,-1.02378)%
(-1.81128,-1.03211)(-1.82588,-1.04043)(-1.84049,-1.04875)(-1.85510,-1.05708)(-1.86970,-1.06540)%
(-1.88431,-1.07372)(-1.89892,-1.08205)(-1.91353,-1.09037)(-1.92813,-1.09869)(-1.94274,-1.10702)%
(-1.95735,-1.11534)(-1.97195,-1.12366)(-1.98656,-1.13199)(-2.00117,-1.14031)(-2.01577,-1.14864)%
(-2.03038,-1.15696)(-2.04499,-1.16528)(-2.05960,-1.17361)(-2.07420,-1.18193)(-2.08881,-1.19025)%
(-2.10342,-1.19858)(-2.11802,-1.20690)(-2.13263,-1.21522)(-2.14724,-1.22355)(-2.16185,-1.23187)%
(-2.17645,-1.24019)(-2.19106,-1.24852)%
\polyline(-1.58243,-0.90171)(-1.48484,-0.92381)(-1.38584,-0.94504)(-1.28552,-0.96540)%
(-1.18394,-0.98485)(-1.08121,-1.00340)(-0.97742,-1.02103)(-0.87266,-1.03773)(-0.76701,-1.05348)%
(-0.66057,-1.06829)(-0.55343,-1.08214)(-0.44568,-1.09503)(-0.33743,-1.10694)(-0.22876,-1.11786)%
(-0.11975,-1.12781)(-0.01053,-1.13676)(0.09883,-1.14471)(0.20823,-1.15167)(0.31739,-1.15762)%
(0.31756,-1.15762)(0.35499,-1.15932)(0.37017,-1.16001)%
\polyline(0.51233,-1.16568)(0.53567,-1.16653)(0.64426,-1.16947)(0.75239,-1.17141)%
(0.85999,-1.17235)(0.96694,-1.17228)(1.07316,-1.17122)(1.17856,-1.16916)(1.28303,-1.16611)%
(1.38647,-1.16207)(1.48882,-1.15706)(1.58994,-1.15109)(1.68976,-1.14414)(1.78819,-1.13623)%
(1.88514,-1.12737)(1.98051,-1.11758)(2.07420,-1.10687)(2.16615,-1.09525)(2.25624,-1.08271)%
(2.34441,-1.06928)(2.43056,-1.05498)(2.51460,-1.03984)(2.59647,-1.02384)(2.67606,-1.00700)%
(2.75330,-0.98936)(2.82814,-0.97092)(2.90045,-0.95171)(2.97019,-0.93173)(3.03729,-0.91102)%
(3.10166,-0.88960)(3.16324,-0.86748)(3.22195,-0.84469)%
\polyline(-1.70416,-0.97107)(-1.59905,-0.99025)(-1.49239,-1.00846)(-1.38427,-1.02572)%
(-1.27481,-1.04198)(-1.16409,-1.05725)(-1.05222,-1.07151)(-0.93931,-1.08476)(-0.82548,-1.09698)%
(-0.71080,-1.10815)(-0.59540,-1.11828)(-0.47939,-1.12736)(-0.36286,-1.13539)(-0.24593,-1.14235)%
(-0.12872,-1.14824)(-0.01132,-1.15306)(0.10615,-1.15680)(0.22360,-1.15947)(0.34087,-1.16107)%
(0.42579,-1.16144)%
\polyline(0.55393,-1.16113)(0.57458,-1.16103)(0.69076,-1.15941)(0.80635,-1.15671)%
(0.92127,-1.15294)(1.03537,-1.14812)(1.14856,-1.14225)(1.26073,-1.13533)(1.37178,-1.12737)%
(1.48160,-1.11839)(1.59008,-1.10837)(1.69713,-1.09736)(1.80263,-1.08535)(1.90650,-1.07236)%
(2.00863,-1.05839)(2.10891,-1.04346)(2.20726,-1.02760)(2.30357,-1.01081)(2.39776,-0.99311)%
(2.48973,-0.97453)(2.57940,-0.95507)(2.66666,-0.93476)(2.75145,-0.91362)(2.83368,-0.89167)%
(2.91325,-0.86893)(2.99010,-0.84542)(3.06414,-0.82117)(3.13529,-0.79621)(3.20350,-0.77055)%
(3.26869,-0.74421)(3.33078,-0.71724)(3.38973,-0.68965)%
\polyline(-1.82588,-1.04043)(-1.71325,-1.05668)(-1.59893,-1.07188)(-1.48303,-1.08603)%
(-1.36567,-1.09910)(-1.24696,-1.11109)(-1.12703,-1.12199)(-1.00598,-1.13178)(-0.88394,-1.14046)%
(-0.76103,-1.14801)(-0.63737,-1.15443)(-0.51308,-1.15971)(-0.38828,-1.16385)(-0.26311,-1.16683)%
(-0.13767,-1.16867)(-0.01210,-1.16936)(0.11348,-1.16890)(0.23896,-1.16728)(0.36420,-1.16451)%
(0.39800,-1.16345)(0.48907,-1.16059)%
\polyline(0.59224,-1.15640)(0.61060,-1.15565)(0.61347,-1.15553)(0.73726,-1.14933)%
(0.86033,-1.14200)(0.98255,-1.13354)(1.10379,-1.12396)(1.22395,-1.11328)(1.34290,-1.10150)%
(1.46053,-1.08863)(1.57672,-1.07469)(1.69136,-1.05969)(1.80432,-1.04364)(1.91551,-1.02657)%
(2.02481,-1.00848)(2.13212,-0.98940)(2.23731,-0.96934)(2.34031,-0.94833)(2.44100,-0.92638)%
(2.53928,-0.90352)(2.63505,-0.87977)(2.72824,-0.85515)(2.81873,-0.82969)(2.90645,-0.80341)%
(2.99129,-0.77634)(3.07319,-0.74851)(3.15205,-0.71993)(3.22781,-0.69065)(3.30039,-0.66068)%
(3.36971,-0.63006)(3.43572,-0.59882)(3.49834,-0.56700)(3.55750,-0.53461)%
\polyline(-1.94761,-1.10979)(-1.82745,-1.12310)(-1.70547,-1.13530)(-1.58179,-1.14634)%
(-1.45653,-1.15622)(-1.32984,-1.16494)(-1.20183,-1.17247)(-1.07264,-1.17881)(-0.94240,-1.18394)%
(-0.81126,-1.18787)(-0.67934,-1.19057)(-0.54677,-1.19205)(-0.41371,-1.19230)(-0.28028,-1.19132)%
(-0.14662,-1.18911)(-0.01288,-1.18567)(0.12081,-1.18099)(0.25432,-1.17508)(0.38752,-1.16795)%
(0.42354,-1.16568)(0.47385,-1.16252)(0.52024,-1.15960)(0.54091,-1.15810)(0.56310,-1.15649)%
(0.57585,-1.15558)(0.58883,-1.15464)%
\polyline(0.62979,-1.15167)(0.65237,-1.15004)(0.78376,-1.13926)(0.91429,-1.12729)%
(1.04383,-1.11414)(1.17221,-1.09980)(1.29934,-1.08431)(1.42507,-1.06766)(1.54928,-1.04989)%
(1.67184,-1.03099)(1.79262,-1.01100)(1.91151,-0.98992)(2.02838,-0.96779)(2.14312,-0.94460)%
(2.25560,-0.92041)(2.36571,-0.89522)(2.47336,-0.86906)(2.57843,-0.84196)(2.68080,-0.81393)%
(2.78038,-0.78502)(2.87707,-0.75523)(2.97079,-0.72462)(3.06143,-0.69321)(3.14890,-0.66101)%
(3.23313,-0.62808)(3.31401,-0.59443)(3.39149,-0.56012)(3.46549,-0.52515)(3.53593,-0.48958)%
(3.60275,-0.45344)(3.66588,-0.41676)(3.69369,-0.39934)(3.72526,-0.37957)%
\polyline(-2.06933,-1.17915)(-1.94165,-1.18954)(-1.81202,-1.19872)(-1.68055,-1.20666)%
(-1.54740,-1.21335)(-1.41271,-1.21879)(-1.27663,-1.22295)(-1.13930,-1.22583)(-1.00087,-1.22743)%
(-0.86149,-1.22772)(-0.72131,-1.22671)(-0.58047,-1.22439)(-0.43914,-1.22075)(-0.29746,-1.21581)%
(-0.15559,-1.20954)(-0.01366,-1.20197)(0.12814,-1.19308)(0.26969,-1.18289)(0.41083,-1.17139)%
(0.44921,-1.16791)(0.47203,-1.16583)(0.49602,-1.16365)(0.55140,-1.15861)(0.60279,-1.15344)%
(0.61904,-1.15181)(0.63334,-1.15037)(0.64646,-1.14905)(0.65987,-1.14769)(0.69127,-1.14454)%
(0.83027,-1.12920)(0.96826,-1.11259)(1.10511,-1.09473)(1.24064,-1.07565)(1.37473,-1.05534)%
(1.50724,-1.03383)(1.63803,-1.01115)(1.76697,-0.98730)(1.89390,-0.96231)(2.01871,-0.93619)%
(2.14125,-0.90900)(2.26143,-0.88073)(2.37909,-0.85142)(2.49413,-0.82110)(2.60642,-0.78979)%
(2.71585,-0.75752)(2.82232,-0.72434)(2.92570,-0.69026)(3.02592,-0.65532)(3.12286,-0.61955)%
(3.21641,-0.58299)(3.30652,-0.54569)(3.39306,-0.50766)(3.47597,-0.46894)(3.55517,-0.42959)%
(3.63059,-0.38963)(3.70214,-0.34910)(3.76978,-0.30805)(3.83342,-0.26652)(3.85492,-0.25137)%
(3.87815,-0.23501)(3.89303,-0.22453)%
\polyline(0.15992,-1.14032)(0.16146,-1.14142)(0.16299,-1.14252)(0.16452,-1.14363)%
(0.16605,-1.14473)(0.16758,-1.14583)(0.16912,-1.14694)(0.17064,-1.14805)(0.17218,-1.14915)%
(0.17371,-1.15026)(0.17524,-1.15136)(0.17677,-1.15246)(0.17831,-1.15357)(0.17983,-1.15468)%
(0.18137,-1.15578)(0.18290,-1.15688)(0.18444,-1.15799)(0.18596,-1.15909)(0.18750,-1.16020)%
(0.18903,-1.16130)(0.19056,-1.16241)(0.19209,-1.16351)(0.19362,-1.16462)(0.19515,-1.16572)%
(0.19669,-1.16682)(0.19822,-1.16793)(0.19975,-1.16904)(0.20128,-1.17014)(0.20281,-1.17124)%
(0.20434,-1.17235)(0.20587,-1.17345)(0.20741,-1.17456)(0.20894,-1.17566)(0.21047,-1.17677)%
(0.21200,-1.17787)(0.21353,-1.17898)(0.21506,-1.18008)(0.21660,-1.18118)(0.21813,-1.18229)%
(0.21966,-1.18340)(0.22119,-1.18450)(0.22272,-1.18560)(0.22425,-1.18671)(0.22579,-1.18781)%
(0.22732,-1.18892)(0.22885,-1.19003)(0.23038,-1.19112)(0.23192,-1.19223)(0.23344,-1.19334)%
(0.23498,-1.19445)(0.23651,-1.19554)%
\polyline(0.16299,-1.14252)(0.16452,-1.14363)(0.16605,-1.14473)(0.16758,-1.14583)%
(0.16912,-1.14694)(0.17064,-1.14805)(0.17186,-1.14892)%
\polyline(1.79371,-1.19471)(1.80852,-1.18634)(1.82354,-1.17785)(1.83857,-1.16936)%
(1.85360,-1.16086)(1.86863,-1.15237)(1.88366,-1.14387)(1.89869,-1.13538)(1.91371,-1.12688)%
(1.92874,-1.11839)(1.94377,-1.10990)(1.95880,-1.10140)(1.97383,-1.09291)(1.98886,-1.08441)%
(2.00389,-1.07592)(2.01892,-1.06742)(2.03394,-1.05893)(2.04897,-1.05044)(2.06400,-1.04194)%
(2.07903,-1.03345)(2.09406,-1.02495)(2.10909,-1.01646)(2.12412,-1.00797)(2.13914,-0.99947)%
(2.15417,-0.99098)(2.16920,-0.98248)(2.18423,-0.97399)(2.19926,-0.96549)(2.21429,-0.95700)%
(2.22932,-0.94851)(2.24435,-0.94001)(2.25937,-0.93152)(2.27440,-0.92302)(2.28943,-0.91453)%
(2.30446,-0.90604)(2.31949,-0.89754)(2.33452,-0.88905)(2.34954,-0.88055)(2.36458,-0.87206)%
(2.37960,-0.86357)(2.39463,-0.85507)(2.40966,-0.84658)(2.42469,-0.83808)(2.43972,-0.82959)%
(2.45475,-0.82109)(2.46977,-0.81260)(2.48481,-0.80410)(2.49983,-0.79561)(2.51486,-0.78712)%
(2.52989,-0.77862)(2.54491,-0.77014)%
\polyline(3.05564,-0.99927)(3.07375,-0.98252)(3.07577,-0.98065)(3.09389,-0.96389)%
(3.09591,-0.96203)(3.11401,-0.94528)(3.11604,-0.94341)(3.13415,-0.92666)(3.13617,-0.92479)%
(3.15427,-0.90805)(3.15630,-0.90617)(3.17441,-0.88943)(3.17644,-0.88755)(3.19453,-0.87081)%
(3.19657,-0.86893)(3.21467,-0.85219)(3.21670,-0.85031)(3.23480,-0.83358)(3.23683,-0.83169)%
(3.25493,-0.81496)(3.25697,-0.81307)(3.27506,-0.79634)(3.27710,-0.79445)(3.29519,-0.77772)%
(3.29723,-0.77584)(3.31532,-0.75911)(3.31736,-0.75722)(3.33545,-0.74049)(3.33750,-0.73859)%
(3.35558,-0.72187)(3.35763,-0.71997)(3.37570,-0.70326)(3.37776,-0.70136)(3.39584,-0.68464)%
(3.39789,-0.68274)(3.41596,-0.66602)(3.41803,-0.66412)(3.43610,-0.64740)(3.43816,-0.64549)%
(3.45622,-0.62879)(3.45829,-0.62688)(3.47636,-0.61017)(3.47842,-0.60826)(3.49649,-0.59156)%
(3.49856,-0.58964)(3.51662,-0.57293)(3.51869,-0.57102)(3.53675,-0.55432)(3.53882,-0.55240)%
(3.55688,-0.53570)(3.55895,-0.53378)(3.57701,-0.51709)(3.57909,-0.51516)(3.59714,-0.49846)%
(3.59921,-0.49655)(3.61727,-0.47985)(3.61935,-0.47793)(3.63739,-0.46124)(3.63948,-0.45930)%
(3.65753,-0.44262)(3.65962,-0.44068)(3.67765,-0.42400)(3.67974,-0.42207)(3.69779,-0.40538)%
(3.69988,-0.40345)(3.71791,-0.38677)(3.72001,-0.38483)(3.73805,-0.36815)(3.74015,-0.36620)%
(3.75818,-0.34953)(3.76027,-0.34759)(3.77831,-0.33091)(3.78041,-0.32897)(3.79844,-0.31230)%
(3.80054,-0.31035)(3.81857,-0.29368)(3.82068,-0.29173)(3.83870,-0.27506)(3.84081,-0.27311)%
(3.85883,-0.25644)(3.86094,-0.25449)(3.87896,-0.23783)(3.88107,-0.23587)(3.89909,-0.21921)%
(3.90121,-0.21725)(3.91922,-0.20059)(3.92134,-0.19864)(3.93934,-0.18198)(3.94147,-0.18001)%
(3.95948,-0.16336)(3.96160,-0.16139)(3.97960,-0.14475)(3.98174,-0.14277)(3.99974,-0.12612)%
(4.00187,-0.12416)(4.01693,-0.11023)(4.02200,-0.10554)(4.03693,-0.09173)(4.04214,-0.08691)%
(4.06013,-0.07028)%
\end{picture}}%

%% file: fig/p010.tex
{\unitlength=1cm%
\begin{picture}%
(9.94,5.4)(-5.19,-2.08)%
\linethickness{0.008in}
\polyline(0.71001,-1.14166)(0.70249,-1.14263)(0.69286,-1.14383)(0.68330,-1.14499)%
(0.68099,-1.14527)%
\polyline(0.68099,-1.14527)(0.67382,-1.14612)(0.66440,-1.14719)(0.65507,-1.14824)%
(0.64580,-1.14923)(0.63660,-1.15019)(0.62748,-1.15112)(0.62168,-1.15168)(0.60977,-1.15281)%
(0.60108,-1.15360)(0.59245,-1.15434)(0.58391,-1.15505)(0.57541,-1.15572)(0.56699,-1.15635)%
(0.55862,-1.15695)(0.55032,-1.15751)(0.54209,-1.15804)(0.53392,-1.15852)(0.52582,-1.15898)%
(0.51778,-1.15939)(0.50980,-1.15977)(0.50189,-1.16011)(0.49404,-1.16041)(0.48780,-1.16063)%
\polyline(0.48780,-1.16063)(0.48626,-1.16069)(0.47855,-1.16092)(0.47279,-1.16106)%
\polyline(0.47279,-1.16106)(0.47089,-1.16111)(0.46330,-1.16128)(0.45764,-1.16137)%
\polyline(0.45764,-1.16137)(0.45577,-1.16140)(0.44831,-1.16149)(0.44238,-1.16154)%
\polyline(0.44238,-1.16154)(0.44091,-1.16155)(0.43358,-1.16157)(0.42738,-1.16155)%
\polyline(0.42738,-1.16155)(0.42630,-1.16155)(0.41910,-1.16149)(0.41230,-1.16142)%
\polyline(0.41230,-1.16142)(0.41195,-1.16141)(0.40487,-1.16129)(0.39786,-1.16113)%
(0.39090,-1.16093)(0.38402,-1.16071)(0.37719,-1.16045)(0.37043,-1.16014)(0.36373,-1.15982)%
(0.35710,-1.15944)(0.35052,-1.15905)(0.34401,-1.15861)(0.33940,-1.15827)%
\polyline(0.33940,-1.15827)(0.33756,-1.15814)(0.33118,-1.15764)(0.32501,-1.15711)%
\polyline(0.71220,-1.14139)(0.70249,-1.14263)(0.69286,-1.14383)(0.68330,-1.14499)%
(0.67382,-1.14612)(0.66440,-1.14719)(0.65507,-1.14824)(0.64580,-1.14923)(0.63660,-1.15019)%
(0.62748,-1.15112)(0.62168,-1.15168)(0.60977,-1.15281)(0.60108,-1.15360)(0.59245,-1.15434)%
(0.58391,-1.15505)(0.57541,-1.15572)(0.56699,-1.15635)(0.55862,-1.15695)(0.55032,-1.15751)%
(0.54209,-1.15804)(0.53392,-1.15852)(0.52582,-1.15898)(0.51778,-1.15939)(0.50980,-1.15977)%
(0.50189,-1.16011)(0.49404,-1.16041)(0.48626,-1.16069)(0.47855,-1.16092)(0.47089,-1.16111)%
(0.46330,-1.16128)(0.45577,-1.16140)(0.44831,-1.16149)(0.44091,-1.16155)(0.43358,-1.16157)%
(0.42630,-1.16155)(0.41910,-1.16149)(0.41195,-1.16141)(0.40487,-1.16129)(0.39786,-1.16113)%
(0.39090,-1.16093)(0.38402,-1.16071)(0.37719,-1.16045)(0.37043,-1.16014)(0.36373,-1.15982)%
(0.35710,-1.15944)(0.35052,-1.15905)(0.34401,-1.15861)(0.33756,-1.15814)(0.33118,-1.15764)%
(0.32486,-1.15710)(0.31861,-1.15652)%
\polyline(4.17260,0.03025)(4.16236,0.02012)(4.16206,0.01982)%
\polyline(4.16206,0.01982)(4.15186,0.00967)(4.15172,0.00952)%
\polyline(4.15172,0.00952)(4.14708,0.00488)(4.14097,-0.00125)(4.14047,-0.00175)%
\polyline(4.14047,-0.00175)(4.12973,-0.01258)(4.12946,-0.01286)%
\polyline(4.12946,-0.01286)(4.11851,-0.02397)(4.10729,-0.03544)(4.10701,-0.03572)%
\polyline(4.10701,-0.03572)(4.09606,-0.04697)(4.09578,-0.04726)%
\polyline(4.09578,-0.04726)(4.08484,-0.05858)(4.08456,-0.05888)%
\polyline(4.08456,-0.05888)(4.07362,-0.07027)(4.07337,-0.07052)%
\polyline(4.07337,-0.07052)(4.06336,-0.08103)(4.05050,-0.09458)(4.03865,-0.10717)%
(4.02681,-0.11983)(4.01497,-0.13259)(4.00312,-0.14542)(3.99128,-0.15834)(3.97946,-0.17135)%
(3.96714,-0.18497)(3.95479,-0.19875)(3.94242,-0.21263)(3.93007,-0.22661)(3.91771,-0.24068)%
(3.90536,-0.25484)(3.89302,-0.26910)(3.88028,-0.28395)(3.86750,-0.29894)(3.85472,-0.31404)%
(3.84196,-0.32925)(3.82921,-0.34457)(3.81646,-0.35999)(3.80480,-0.37420)(3.79058,-0.39167)%
(3.77748,-0.40789)(3.76440,-0.42422)(3.75134,-0.44067)(3.73828,-0.45723)(3.72524,-0.47392)%
(3.71771,-0.48362)(3.71209,-0.49088)(3.69877,-0.50821)(3.68547,-0.52565)(3.67219,-0.54322)%
(3.65894,-0.56090)(3.64571,-0.57871)(3.63249,-0.59663)(3.61908,-0.61499)(3.60566,-0.63353)%
(3.59292,-0.65127)%
\polyline(3.59292,-0.65127)(3.59226,-0.65220)(3.57956,-0.67004)%
\polyline(3.57956,-0.67004)(3.57891,-0.67098)(3.56623,-0.68894)%
\polyline(3.56623,-0.68894)(3.56557,-0.68988)(3.55260,-0.70843)%
\polyline(-2.19106,-1.24852)(-1.98702,-1.25923)(-1.77842,-1.26699)(-1.56575,-1.27175)%
(-1.34959,-1.27343)(-1.13045,-1.27203)(-0.90889,-1.26750)(-0.68547,-1.25983)(-0.46076,-1.24900)%
(-0.23529,-1.23501)(-0.00967,-1.21788)(0.21556,-1.19759)(0.43983,-1.17421)%
\polyline(0.43983,-1.17421)(0.66259,-1.14773)%
\polyline(0.66259,-1.14772)(0.68099,-1.14527)%
\polyline(0.68099,-1.14527)(0.70828,-1.14161)%
\polyline(0.70828,-1.14161)(0.88329,-1.11821)(1.10137,-1.08569)(1.31629,-1.05024)%
(1.52754,-1.01191)(1.73460,-0.97079)(1.93696,-0.92694)(2.13414,-0.88046)(2.32566,-0.83144)%
(2.51107,-0.77998)(2.68993,-0.72620)(2.86184,-0.67020)(3.02638,-0.61210)(3.18319,-0.55204)%
(3.33190,-0.49014)(3.47221,-0.42654)(3.53211,-0.39688)%
\polyline(3.53211,-0.39688)(3.60379,-0.36137)(3.72638,-0.29479)(3.83970,-0.22693)%
(3.94354,-0.15795)(4.03770,-0.08801)(4.07985,-0.05262)%
\polyline(4.07985,-0.05262)(4.10436,-0.03205)%
\polyline(4.10436,-0.03205)(4.12199,-0.01724)(4.12839,-0.01109)%
\polyline(4.12839,-0.01109)(4.13986,-0.00007)%
\polyline(4.13986,-0.00007)(4.15096,0.01060)%
\polyline(4.15096,0.01060)(4.17085,0.02972)%
\polyline(4.17085,0.02972)(4.19629,0.05417)(4.26045,0.12609)(4.31439,0.19835)(4.35806,0.27079)%
(4.39139,0.34325)(4.41440,0.41557)(4.42709,0.48759)(4.42951,0.55917)(4.42172,0.63014)%
(4.40383,0.70034)(4.37595,0.76965)(4.33823,0.83789)(4.29084,0.90493)(4.23397,0.97063)%
(4.16785,1.03486)(4.09270,1.09747)%
\polyline(3.51323,-0.37249)(3.52483,-0.34310)(3.53641,-0.31370)(3.54800,-0.28430)%
(3.55960,-0.25490)(3.57118,-0.22550)(3.58277,-0.19610)(3.59436,-0.16670)(3.60594,-0.13731)%
(3.61754,-0.10790)(3.62913,-0.07850)(3.64071,-0.04910)(3.65230,-0.01971)(3.66390,0.00969)%
(3.67548,0.03910)(3.68707,0.06850)(3.69866,0.09789)(3.71025,0.12729)(3.72184,0.15669)%
(3.73343,0.18610)(3.74502,0.21549)(3.75661,0.24489)(3.76819,0.27429)(3.77979,0.30369)%
(3.79138,0.33309)(3.80296,0.36249)(3.81455,0.39189)(3.82614,0.42128)(3.83773,0.45069)%
(3.84932,0.48009)(3.86091,0.50949)(3.87250,0.53888)(3.88409,0.56828)(3.89567,0.59769)%
(3.90726,0.62709)(3.91886,0.65648)(3.93044,0.68588)(3.94203,0.71528)(3.95363,0.74469)%
(3.96521,0.77408)(3.97680,0.80348)(3.98839,0.83288)(3.99998,0.86229)(4.01157,0.89168)%
(4.02316,0.92108)(4.03474,0.95048)(4.04633,0.97988)(4.05793,1.00927)(4.06951,1.03868)%
(4.08110,1.06808)(4.09270,1.09747)%
\polyline(-1.46071,-0.83234)(-1.32484,-0.86968)(-1.18616,-0.90523)(-1.04493,-0.93895)%
(-0.90139,-0.97081)(-0.75582,-1.00076)(-0.60846,-1.02878)(-0.45958,-1.05482)(-0.30945,-1.07886)%
(-0.15834,-1.10089)(-0.00652,-1.12085)(0.14572,-1.13876)(0.29812,-1.15455)(0.32518,-1.15699)%
\polyline(0.32518,-1.15699)(0.33940,-1.15827)%
\polyline(0.45039,-1.16827)(0.60224,-1.17986)(0.75341,-1.18932)(0.90359,-1.19665)%
(1.05251,-1.20186)(1.19986,-1.20493)(1.34539,-1.20588)(1.48878,-1.20471)(1.62976,-1.20144)%
(1.76803,-1.19607)(1.90332,-1.18863)(2.03535,-1.17914)(2.16383,-1.16762)(2.28848,-1.15411)%
(2.40904,-1.13864)(2.52521,-1.12124)(2.63676,-1.10196)(2.74340,-1.08083)(2.84489,-1.05792)%
(2.94097,-1.03327)(3.03141,-1.00692)(3.11595,-0.97896)(3.19438,-0.94943)(3.26647,-0.91841)%
(3.33202,-0.88596)(3.39082,-0.85217)(3.44268,-0.81710)(3.48743,-0.78083)(3.52489,-0.74345)%
(3.55277,-0.70781)%
\polyline(3.55277,-0.70781)(3.55492,-0.70507)(3.56465,-0.68801)%
\polyline(3.53211,-0.39687)(3.51323,-0.37249)%
\polyline(-1.46071,-0.83234)(-1.47531,-0.84067)(-1.48992,-0.84899)(-1.50453,-0.85731)%
(-1.51913,-0.86564)(-1.53374,-0.87396)(-1.54835,-0.88229)(-1.56296,-0.89061)(-1.57756,-0.89893)%
(-1.59217,-0.90726)(-1.60678,-0.91558)(-1.62138,-0.92390)(-1.63599,-0.93223)(-1.65060,-0.94055)%
(-1.66521,-0.94887)(-1.67981,-0.95720)(-1.69442,-0.96552)(-1.70903,-0.97384)(-1.72363,-0.98217)%
(-1.73824,-0.99049)(-1.75285,-0.99881)(-1.76745,-1.00714)(-1.78206,-1.01546)(-1.79667,-1.02378)%
(-1.81128,-1.03211)(-1.82588,-1.04043)(-1.84049,-1.04875)(-1.85510,-1.05708)(-1.86970,-1.06540)%
(-1.88431,-1.07372)(-1.89892,-1.08205)(-1.91353,-1.09037)(-1.92813,-1.09869)(-1.94274,-1.10702)%
(-1.95735,-1.11534)(-1.97195,-1.12366)(-1.98656,-1.13199)(-2.00117,-1.14031)(-2.01577,-1.14864)%
(-2.03038,-1.15696)(-2.04499,-1.16528)(-2.05960,-1.17361)(-2.07420,-1.18193)(-2.08881,-1.19025)%
(-2.10342,-1.19858)(-2.11802,-1.20690)(-2.13263,-1.21522)(-2.14724,-1.22355)(-2.16185,-1.23187)%
(-2.17645,-1.24019)(-2.19106,-1.24852)%
\polyline(-1.58243,-0.90171)(-1.43520,-0.93461)(-1.28487,-0.96552)(-1.13173,-0.99442)%
(-0.97610,-1.02125)(-0.81825,-1.04598)(-0.65853,-1.06857)(-0.49722,-1.08900)(-0.33467,-1.10723)%
(-0.17116,-1.12324)(-0.00704,-1.13702)(0.15737,-1.14856)(0.32128,-1.15781)(0.32174,-1.15783)%
(0.35134,-1.15909)(0.37043,-1.15991)(0.38437,-1.16051)(0.40016,-1.16119)%
\polyline(0.48576,-1.16484)(0.64909,-1.16957)(0.81140,-1.17204)(0.97238,-1.17226)%
(1.13167,-1.17020)(1.28899,-1.16591)(1.44398,-1.15939)(1.59634,-1.15067)(1.74574,-1.13977)%
(1.89187,-1.12672)(2.03443,-1.11156)(2.17310,-1.09432)(2.30759,-1.07504)(2.43760,-1.05377)%
(2.56285,-1.03056)(2.68305,-1.00545)(2.79793,-0.97853)(2.90723,-0.94983)(3.01069,-0.91942)%
(3.10807,-0.88738)(3.19912,-0.85377)(3.28362,-0.81867)(3.36135,-0.78216)(3.43214,-0.74432)%
(3.49575,-0.70524)(3.55202,-0.66500)(3.60080,-0.62370)(3.61337,-0.61077)(3.64192,-0.58143)%
(3.64469,-0.57785)(3.66223,-0.55515)%
\polyline(3.71700,-0.29989)(3.70593,-0.26669)(3.68220,-0.22029)(3.65016,-0.17386)%
(3.60981,-0.12750)%
\polyline(-1.70416,-0.97107)(-1.54557,-0.99953)(-1.38358,-1.02582)(-1.21854,-1.04988)%
(-1.05079,-1.07169)(-0.88070,-1.09118)(-0.70860,-1.10835)(-0.53488,-1.12316)(-0.35988,-1.13558)%
(-0.18399,-1.14559)(-0.00757,-1.15319)(0.16900,-1.15836)(0.34534,-1.16110)(0.37345,-1.16116)%
(0.40117,-1.16121)(0.41230,-1.16123)(0.42673,-1.16125)(0.44098,-1.16127)(0.46118,-1.16131)%
\polyline(0.52113,-1.16142)(0.69593,-1.15930)(0.86940,-1.15477)(1.04116,-1.14785)%
(1.21085,-1.13854)(1.37811,-1.12688)(1.54258,-1.11290)(1.70390,-1.09662)(1.86172,-1.07811)%
(2.01571,-1.05738)(2.16553,-1.03448)(2.31085,-1.00949)(2.45135,-0.98245)(2.58672,-0.95342)%
(2.71666,-0.92247)(2.84088,-0.88968)(2.95910,-0.85510)(3.07106,-0.81882)(3.17649,-0.78093)%
(3.27516,-0.74149)(3.36684,-0.70061)(3.45130,-0.65839)(3.52835,-0.61489)(3.59780,-0.57024)%
(3.65948,-0.52453)(3.70561,-0.48446)%
\polyline(3.71270,-0.47831)(3.71323,-0.47785)(3.72108,-0.46968)(3.75009,-0.43950)%
(3.75892,-0.43031)(3.76394,-0.42385)(3.77740,-0.40652)%
\polyline(3.86039,-0.21319)(3.86267,-0.18360)(3.85798,-0.13323)(3.84472,-0.08278)%
(3.82291,-0.03237)(3.79256,0.01790)(3.75370,0.06788)(3.70639,0.11749)%
\polyline(-1.82588,-1.04043)(-1.65593,-1.06446)(-1.48229,-1.08611)(-1.30534,-1.10535)%
(-1.12549,-1.12212)(-0.94313,-1.13640)(-0.75867,-1.14814)(-0.57252,-1.15733)(-0.38510,-1.16393)%
(-0.19682,-1.16795)(-0.00810,-1.16936)(0.18064,-1.16817)(0.36895,-1.16438)(0.45348,-1.16150)%
(0.47278,-1.16085)%
\polyline(0.49515,-1.16008)(0.51020,-1.15957)(0.52220,-1.15916)%
\polyline(0.55649,-1.15799)(0.74276,-1.14903)(0.92738,-1.13750)(1.10994,-1.12345)%
(1.29002,-1.10688)(1.46724,-1.08786)(1.64118,-1.06641)(1.81146,-1.04258)(1.97771,-1.01644)%
(2.13955,-0.98802)(2.29663,-0.95741)(2.44859,-0.92467)(2.59510,-0.88986)(2.73584,-0.85308)%
(2.87047,-0.81439)(2.99871,-0.77389)(3.12027,-0.73166)(3.23488,-0.68781)(3.34229,-0.64243)%
(3.44226,-0.59561)(3.53455,-0.54746)(3.61897,-0.49810)(3.69533,-0.44763)(3.76346,-0.39616)%
(3.79769,-0.36617)(3.82321,-0.34381)(3.82586,-0.34106)(3.83893,-0.32750)(3.85326,-0.31265)%
(3.87444,-0.29069)(3.88102,-0.28238)%
\polyline(3.97188,-0.13690)(3.97599,-0.12793)(3.99221,-0.07294)(3.99953,-0.01780)%
(3.99795,0.03739)(3.98746,0.09249)(3.96810,0.14739)(3.93990,0.20196)(3.90291,0.25608)%
(3.85723,0.30963)(3.80296,0.36249)%
\polyline(-1.94761,-1.10979)(-1.76629,-1.12939)(-1.58099,-1.14641)(-1.39215,-1.16082)%
(-1.20019,-1.17256)(-1.00557,-1.18161)(-0.80874,-1.18793)(-0.61017,-1.19150)(-0.41032,-1.19229)%
(-0.20964,-1.19031)(-0.00862,-1.18554)(0.19228,-1.17798)(0.39259,-1.16766)(0.42205,-1.16572)%
(0.52104,-1.15922)(0.55031,-1.15730)(0.56740,-1.15617)(0.58054,-1.15532)%
\polyline(0.59186,-1.15457)(0.78961,-1.13876)(0.98538,-1.12023)(1.17872,-1.09904)%
(1.36920,-1.07523)(1.55636,-1.04884)(1.73977,-1.01991)(1.91902,-0.98854)(2.09370,-0.95477)%
(2.26339,-0.91867)(2.42774,-0.88034)(2.58634,-0.83985)(2.73886,-0.79728)(2.88495,-0.75273)%
(3.02428,-0.70631)(3.15654,-0.65810)(3.28144,-0.60823)(3.39872,-0.55680)(3.50810,-0.50393)%
(3.60935,-0.44973)(3.70227,-0.39432)(3.78665,-0.33782)(3.86232,-0.28036)(3.91294,-0.23620)%
(3.92913,-0.22208)(3.94040,-0.21056)(3.95273,-0.19799)(3.96630,-0.18414)(3.98415,-0.16594)%
(3.98694,-0.16309)(3.99156,-0.15743)%
\polyline(4.05867,-0.06857)(4.07516,-0.04354)(4.10541,0.01677)(4.12636,0.07725)(4.13797,0.13776)%
(4.14027,0.19818)(4.13324,0.25838)(4.11696,0.31821)(4.09148,0.37756)(4.05688,0.43629)%
(4.01326,0.49426)(3.96077,0.55138)(3.89954,0.60748)%
\polyline(-2.06933,-1.17915)(-1.87666,-1.19431)(-1.67970,-1.20671)(-1.47895,-1.21628)%
(-1.27489,-1.22300)(-1.06801,-1.22682)(-0.85882,-1.22771)(-0.64783,-1.22566)(-0.43554,-1.22064)%
(-0.22247,-1.21266)(-0.00915,-1.20171)(0.20392,-1.18779)(0.41622,-1.17093)(0.44842,-1.16791)%
(0.60161,-1.15355)(0.62723,-1.15115)(0.64625,-1.14909)(0.67532,-1.14594)(0.83644,-1.12849)%
(1.04337,-1.10296)(1.24751,-1.07464)(1.44836,-1.04357)(1.64548,-1.00981)(1.83836,-0.97343)%
(2.02658,-0.93450)(2.20968,-0.89311)(2.38723,-0.84933)(2.55884,-0.80326)(2.72409,-0.75502)%
(2.88262,-0.70469)(3.03407,-0.65239)(3.17809,-0.59822)(3.31438,-0.54233)(3.44262,-0.48480)%
(3.56254,-0.42579)(3.67390,-0.36543)(3.77645,-0.30383)(3.86998,-0.24115)(3.95432,-0.17753)%
(3.98486,-0.15128)(4.01109,-0.12875)(4.02930,-0.11309)(4.03701,-0.10543)(4.04923,-0.09329)%
(4.07238,-0.07027)(4.08205,-0.06066)(4.09479,-0.04799)%
\polyline(4.11841,-0.02026)(4.12599,-0.01135)(4.13710,0.00170)%
\polyline(4.14755,0.01397)(4.15067,0.01763)(4.15625,0.02561)(4.16309,0.03539)(4.19685,0.08363)%
(4.23327,0.14986)(4.25991,0.21617)(4.27672,0.28241)(4.28374,0.34846)(4.28099,0.41416)%
(4.26854,0.47937)(4.24645,0.54393)(4.21486,0.60772)(4.17386,0.67061)(4.12362,0.73245)%
(4.06431,0.79311)(3.99611,0.85248)%
\polyline(0.15992,-1.14032)(0.16146,-1.14142)(0.16299,-1.14252)(0.16452,-1.14363)%
(0.16605,-1.14473)(0.16758,-1.14583)(0.16912,-1.14694)(0.17064,-1.14805)(0.17218,-1.14915)%
(0.17371,-1.15026)(0.17524,-1.15136)(0.17677,-1.15246)(0.17831,-1.15357)(0.17983,-1.15468)%
(0.18137,-1.15578)(0.18290,-1.15688)(0.18444,-1.15799)(0.18596,-1.15909)(0.18750,-1.16020)%
(0.18903,-1.16130)(0.19056,-1.16241)(0.19209,-1.16351)(0.19362,-1.16462)(0.19515,-1.16572)%
(0.19669,-1.16682)(0.19822,-1.16793)(0.19975,-1.16904)(0.20128,-1.17014)(0.20281,-1.17124)%
(0.20434,-1.17235)(0.20587,-1.17345)(0.20741,-1.17456)(0.20894,-1.17566)(0.21047,-1.17677)%
(0.21200,-1.17787)(0.21353,-1.17898)(0.21506,-1.18008)(0.21660,-1.18118)(0.21813,-1.18229)%
(0.21966,-1.18340)(0.22119,-1.18450)(0.22272,-1.18560)(0.22425,-1.18671)(0.22579,-1.18781)%
(0.22732,-1.18892)(0.22885,-1.19003)(0.23038,-1.19112)(0.23192,-1.19223)(0.23344,-1.19334)%
(0.23498,-1.19445)(0.23651,-1.19554)%
\polyline(1.79382,-1.19465)(1.80852,-1.18634)(1.80992,-1.18555)(1.82354,-1.17785)%
(1.83857,-1.16936)(1.85360,-1.16086)(1.86863,-1.15237)(1.88366,-1.14387)(1.89869,-1.13538)%
(1.91371,-1.12688)(1.92874,-1.11839)(1.94377,-1.10990)(1.95880,-1.10140)(1.97383,-1.09291)%
(1.98886,-1.08441)(2.00389,-1.07592)(2.01892,-1.06742)(2.03394,-1.05893)(2.04897,-1.05044)%
(2.06400,-1.04194)(2.07903,-1.03345)(2.09406,-1.02495)(2.10909,-1.01646)(2.12412,-1.00797)%
(2.13914,-0.99947)(2.15417,-0.99098)(2.16920,-0.98248)(2.18423,-0.97399)(2.19926,-0.96549)%
(2.21429,-0.95700)(2.22932,-0.94851)(2.24435,-0.94001)(2.25937,-0.93152)(2.27440,-0.92302)%
(2.28943,-0.91453)(2.30446,-0.90604)(2.31949,-0.89754)(2.33452,-0.88905)(2.34954,-0.88055)%
(2.36458,-0.87206)(2.37960,-0.86357)(2.39463,-0.85507)(2.40966,-0.84658)(2.42469,-0.83808)%
(2.43972,-0.82959)(2.45475,-0.82109)(2.46977,-0.81260)(2.48481,-0.80410)(2.49983,-0.79561)%
(2.51486,-0.78712)(2.52989,-0.77862)(2.54492,-0.77013)%
\polyline(3.05625,-0.99870)(3.07577,-0.98065)(3.09591,-0.96203)(3.11604,-0.94341)%
(3.13617,-0.92479)(3.15630,-0.90617)(3.17644,-0.88755)(3.19657,-0.86893)(3.21670,-0.85031)%
(3.23683,-0.83169)(3.25697,-0.81307)(3.27710,-0.79445)(3.29723,-0.77584)(3.31736,-0.75722)%
(3.33750,-0.73859)(3.35763,-0.71997)(3.37776,-0.70136)(3.39789,-0.68274)(3.41803,-0.66412)%
(3.43816,-0.64549)(3.45829,-0.62688)(3.47842,-0.60826)(3.49856,-0.58964)(3.51869,-0.57102)%
(3.53882,-0.55240)(3.55895,-0.53378)(3.57909,-0.51516)(3.59921,-0.49655)(3.61935,-0.47793)%
(3.63948,-0.45930)(3.65962,-0.44068)(3.67974,-0.42207)(3.69988,-0.40345)(3.72001,-0.38483)%
(3.74015,-0.36620)(3.76027,-0.34759)(3.78041,-0.32897)(3.80054,-0.31035)(3.82068,-0.29173)%
(3.84081,-0.27311)(3.86094,-0.25449)(3.88107,-0.23587)(3.90121,-0.21725)(3.92134,-0.19864)%
(3.94147,-0.18001)(3.96160,-0.16139)(3.98174,-0.14277)(4.00081,-0.12514)(4.00187,-0.12416)%
(4.02107,-0.10640)(4.02200,-0.10554)(4.04152,-0.08748)(4.04214,-0.08691)(4.06167,-0.06884)%
\polyline(3.81650,-0.24083)(3.81930,-0.23613)(3.83528,-0.20939)(3.85126,-0.18266)%
(3.86724,-0.15592)(3.88323,-0.12919)(3.89920,-0.10246)(3.91518,-0.07572)(3.93117,-0.04899)%
(3.94714,-0.02226)(3.96313,0.00447)(3.97911,0.03121)(3.99508,0.05795)(4.01107,0.08468)%
(4.02704,0.11142)(4.04303,0.13815)(4.05901,0.16488)(4.07498,0.19161)(4.09097,0.21835)%
(4.10695,0.24508)(4.12293,0.27182)(4.13891,0.29856)(4.15489,0.32529)(4.17087,0.35202)%
(4.18685,0.37875)(4.20283,0.40549)(4.21881,0.43222)(4.23479,0.45895)(4.25077,0.48569)%
(4.26675,0.51243)(4.28273,0.53916)(4.29871,0.56589)(4.31469,0.59263)(4.33068,0.61936)%
(4.34665,0.64609)(4.36263,0.67282)(4.37862,0.69956)(4.39388,0.72509)%
\end{picture}}%

%% file: fig/p013.tex
{\unitlength=1cm%
\begin{picture}%
(9.94,5.4)(-5.19,-2.08)%
\linethickness{0.008in}
\polyline(0.71004,-1.14165)(0.70223,-1.14266)(0.69283,-1.14383)(0.69260,-1.14386)%
\polyline(0.69260,-1.14386)(0.68350,-1.14497)(0.67423,-1.14607)(0.66503,-1.14713)%
(0.65590,-1.14815)(0.64683,-1.14913)%
\polyline(0.64683,-1.14913)(0.63783,-1.15007)(0.62890,-1.15098)(0.62003,-1.15184)%
(0.61122,-1.15268)(0.60249,-1.15347)(0.59382,-1.15423)(0.58521,-1.15494)(0.57667,-1.15562)%
(0.56819,-1.15626)(0.55979,-1.15687)(0.55145,-1.15744)(0.54317,-1.15797)(0.53496,-1.15846)%
(0.52681,-1.15892)(0.51873,-1.15934)(0.51072,-1.15972)(0.50277,-1.16007)(0.49488,-1.16038)%
(0.48707,-1.16065)(0.48630,-1.16068)%
\polyline(0.48630,-1.16068)(0.47930,-1.16089)(0.47623,-1.16097)%
\polyline(0.47623,-1.16097)(0.47162,-1.16110)(0.46589,-1.16122)%
\polyline(0.46589,-1.16122)(0.46399,-1.16127)(0.45644,-1.16139)(0.44893,-1.16148)%
(0.44596,-1.16150)%
\polyline(0.71170,-1.14145)(0.70223,-1.14266)(0.69283,-1.14383)(0.68350,-1.14497)%
(0.67423,-1.14607)(0.66503,-1.14713)(0.65590,-1.14815)(0.64683,-1.14913)(0.63783,-1.15007)%
(0.62890,-1.15098)(0.62003,-1.15184)(0.61122,-1.15268)(0.60249,-1.15347)(0.59382,-1.15423)%
(0.58521,-1.15494)(0.57667,-1.15562)(0.56819,-1.15626)(0.55979,-1.15687)(0.55145,-1.15744)%
(0.54317,-1.15797)(0.53496,-1.15846)(0.52681,-1.15892)(0.51873,-1.15934)(0.51072,-1.15972)%
(0.50277,-1.16007)(0.49488,-1.16038)(0.48707,-1.16065)(0.47930,-1.16089)(0.47162,-1.16110)%
(0.46399,-1.16127)(0.45644,-1.16139)(0.44893,-1.16148)(0.44151,-1.16154)(0.43414,-1.16157)%
(0.42684,-1.16155)(0.41960,-1.16150)(0.41244,-1.16141)(0.40532,-1.16129)(0.39828,-1.16114)%
(0.39131,-1.16095)(0.38439,-1.16072)(0.37754,-1.16046)(0.37075,-1.16017)(0.36403,-1.15983)%
(0.35737,-1.15947)(0.35078,-1.15906)(0.34425,-1.15863)(0.33779,-1.15816)(0.33138,-1.15765)%
(0.32505,-1.15712)(0.31877,-1.15654)%
\polyline(0.43597,-1.16155)(0.43414,-1.16157)(0.42684,-1.16155)(0.41960,-1.16150)%
(0.41565,-1.16146)%
\polyline(0.71170,-1.14145)(0.70223,-1.14266)(0.69283,-1.14383)(0.68350,-1.14497)%
(0.67423,-1.14607)(0.66503,-1.14713)(0.65590,-1.14815)(0.64683,-1.14913)(0.63783,-1.15007)%
(0.62890,-1.15098)(0.62003,-1.15184)(0.61122,-1.15268)(0.60249,-1.15347)(0.59382,-1.15423)%
(0.58521,-1.15494)(0.57667,-1.15562)(0.56819,-1.15626)(0.55979,-1.15687)(0.55145,-1.15744)%
(0.54317,-1.15797)(0.53496,-1.15846)(0.52681,-1.15892)(0.51873,-1.15934)(0.51072,-1.15972)%
(0.50277,-1.16007)(0.49488,-1.16038)(0.48707,-1.16065)(0.47930,-1.16089)(0.47162,-1.16110)%
(0.46399,-1.16127)(0.45644,-1.16139)(0.44893,-1.16148)(0.44151,-1.16154)(0.43414,-1.16157)%
(0.42684,-1.16155)(0.41960,-1.16150)(0.41244,-1.16141)(0.40532,-1.16129)(0.39828,-1.16114)%
(0.39131,-1.16095)(0.38439,-1.16072)(0.37754,-1.16046)(0.37075,-1.16017)(0.36403,-1.15983)%
(0.35737,-1.15947)(0.35078,-1.15906)(0.34425,-1.15863)(0.33779,-1.15816)(0.33138,-1.15765)%
(0.32505,-1.15712)(0.31877,-1.15654)%
\polyline(0.40567,-1.16130)(0.40532,-1.16129)(0.39828,-1.16114)(0.39131,-1.16095)%
(0.38542,-1.16075)%
\polyline(0.38542,-1.16075)(0.38439,-1.16072)(0.37754,-1.16046)(0.37075,-1.16017)%
(0.36403,-1.15983)(0.35737,-1.15947)(0.35078,-1.15906)(0.34610,-1.15875)%
\polyline(0.34610,-1.15875)(0.34425,-1.15863)(0.33779,-1.15816)(0.33154,-1.15767)%
\polyline(0.33154,-1.15767)(0.33138,-1.15765)(0.32505,-1.15712)(0.31892,-1.15655)%
\polyline(4.17294,0.03060)(4.16216,0.01993)(4.16189,0.01965)%
\polyline(4.16189,0.01965)(4.15110,0.00890)(4.15082,0.00864)%
\polyline(4.15082,0.00864)(4.14004,-0.00218)(4.13977,-0.00245)%
\polyline(4.13977,-0.00245)(4.12898,-0.01334)(4.12871,-0.01361)%
\polyline(4.12871,-0.01361)(4.11794,-0.02456)(4.11765,-0.02484)%
\polyline(4.11765,-0.02484)(4.10689,-0.03584)(4.10660,-0.03613)%
\polyline(4.10660,-0.03613)(4.09583,-0.04722)(4.09556,-0.04749)%
\polyline(4.09556,-0.04749)(4.08477,-0.05865)(4.08446,-0.05898)%
\polyline(4.08446,-0.05898)(4.07373,-0.07016)(4.07362,-0.07026)%
\polyline(4.07362,-0.07026)(4.06980,-0.07426)(4.06214,-0.08230)(4.05024,-0.09486)%
(4.03835,-0.10749)(4.02646,-0.12020)(4.01456,-0.13301)(4.00268,-0.14590)(3.99080,-0.15888)%
(3.97891,-0.17194)(3.96704,-0.18510)(3.95616,-0.19722)(3.94254,-0.21250)(3.92999,-0.22670)%
(3.91744,-0.24099)(3.90489,-0.25538)(3.89236,-0.26988)(3.87983,-0.28448)(3.86731,-0.29918)%
(3.85478,-0.31397)(3.84228,-0.32887)(3.82928,-0.34448)(3.81625,-0.36025)(3.80322,-0.37613)%
(3.79021,-0.39213)(3.77720,-0.40824)(3.76421,-0.42446)(3.75124,-0.44079)(3.73828,-0.45724)%
(3.72534,-0.47380)(3.71197,-0.49103)(3.69863,-0.50839)(3.68531,-0.52586)(3.67201,-0.54346)%
(3.65873,-0.56118)(3.64547,-0.57902)(3.63225,-0.59696)(3.61905,-0.61504)(3.60936,-0.62840)%
(3.60579,-0.63336)(3.59310,-0.65102)%
\polyline(3.59310,-0.65102)(3.59229,-0.65216)(3.58018,-0.66919)%
\polyline(3.58018,-0.66919)(3.57883,-0.67107)(3.56607,-0.68917)%
\polyline(3.56607,-0.68917)(3.56540,-0.69012)(3.55402,-0.70641)%
\polyline(-2.19106,-1.24852)(-1.91855,-1.26213)(-1.63826,-1.27047)(-1.35143,-1.27343)%
(-1.05934,-1.27091)(-0.76327,-1.26286)(-0.46456,-1.24921)(-0.16454,-1.22998)(0.13547,-1.20518)%
(0.43415,-1.17484)%
\polyline(0.43415,-1.17484)(0.64609,-1.14921)%
\polyline(0.64609,-1.14921)(0.67421,-1.14581)%
\polyline(0.67421,-1.14581)(0.69279,-1.14356)%
\polyline(0.69279,-1.14356)(0.70789,-1.14173)%
\polyline(0.70789,-1.14173)(0.73016,-1.13904)(1.02223,-1.09788)(1.30907,-1.05149)%
(1.58941,-1.00000)(1.86208,-0.94360)(2.12590,-0.88248)(2.37973,-0.81686)(2.62253,-0.74697)%
(2.85327,-0.67309)(3.07102,-0.59550)(3.27492,-0.51448)(3.46413,-0.43036)(3.53184,-0.39650)%
\polyline(3.53184,-0.39650)(3.63793,-0.34345)(3.79568,-0.25410)(3.93681,-0.16267)%
(4.06080,-0.06949)(4.07973,-0.05268)%
\polyline(4.07973,-0.05268)(4.10331,-0.03175)%
\polyline(4.10331,-0.03175)(4.12680,-0.01089)%
\polyline(4.12680,-0.01089)(4.15021,0.00989)%
\polyline(4.15021,0.00989)(4.16728,0.02505)(4.17187,0.03000)%
\polyline(4.17187,0.03000)(4.25592,0.12059)(4.32648,0.21675)(4.37881,0.31317)(4.41285,0.40945)%
(4.42863,0.50523)(4.42625,0.60015)(4.40591,0.69384)(4.36788,0.78595)(4.31250,0.87611)%
(4.24018,0.96402)(4.15144,1.04934)(4.04681,1.13175)(3.92695,1.21097)(3.79253,1.28672)%
(3.64427,1.35873)(3.48300,1.42675)(3.30954,1.49057)(3.12477,1.54998)(2.92961,1.60479)%
(2.72499,1.65484)(2.51191,1.69997)(2.29135,1.74007)(2.06432,1.77505)(1.83184,1.80480)%
\polyline(1.83126,0.27436)(1.83127,0.30497)(1.83128,0.33558)(1.83129,0.36619)(1.83130,0.39680)%
(1.83132,0.42741)(1.83133,0.45801)(1.83134,0.48862)(1.83135,0.51923)(1.83136,0.54984)%
(1.83137,0.58045)(1.83139,0.61106)(1.83140,0.64167)(1.83141,0.67228)(1.83142,0.70288)%
(1.83143,0.73349)(1.83145,0.76410)(1.83145,0.79471)(1.83147,0.82532)(1.83148,0.85593)%
(1.83149,0.88654)(1.83150,0.91714)(1.83151,0.94775)(1.83152,0.97836)(1.83154,1.00897)%
(1.83155,1.03958)(1.83156,1.07019)(1.83157,1.10080)(1.83158,1.13140)(1.83160,1.16201)%
(1.83161,1.19262)(1.83162,1.22323)(1.83163,1.25384)(1.83164,1.28445)(1.83165,1.31506)%
(1.83166,1.34567)(1.83168,1.37628)(1.83169,1.40688)(1.83170,1.43749)(1.83171,1.46810)%
(1.83172,1.49871)(1.83173,1.52932)(1.83175,1.55993)(1.83176,1.59053)(1.83177,1.62114)%
(1.83178,1.65175)(1.83179,1.68236)(1.83181,1.71297)(1.83182,1.74358)(1.83183,1.77419)%
(1.83184,1.80480)%
\polyline(-1.46071,-0.83234)(-1.27930,-0.88163)(-1.09308,-0.92773)(-0.90262,-0.97055)%
(-0.70854,-1.01000)(-0.51146,-1.04600)(-0.31200,-1.07848)(-0.11080,-1.10737)(0.09150,-1.13262)%
(0.29425,-1.15418)(0.31972,-1.15642)%
\polyline(0.31972,-1.15642)(0.33171,-1.15748)%
\polyline(0.33171,-1.15748)(0.34609,-1.15875)%
\polyline(0.47475,-1.17008)(0.49677,-1.17202)(0.69843,-1.18612)(0.89852,-1.19644)%
(1.09639,-1.20299)(1.29136,-1.20577)(1.48275,-1.20480)(1.66990,-1.20010)(1.85211,-1.19171)%
(2.02872,-1.17967)(2.19909,-1.16405)(2.36255,-1.14490)(2.51845,-1.12232)(2.66618,-1.09641)%
(2.80510,-1.06726)(2.93463,-1.03499)(3.05420,-0.99973)(3.16323,-0.96163)(3.26120,-0.92083)%
(3.34761,-0.87749)(3.42200,-0.83181)(3.48392,-0.78395)(3.53296,-0.73411)(3.55132,-0.70767)%
\polyline(3.53185,-0.39650)(3.49100,-0.34833)(3.42783,-0.29080)(3.35021,-0.23342)%
(3.25826,-0.17646)(3.15212,-0.12017)(3.03202,-0.06484)(2.89822,-0.01071)(2.75107,0.04195)%
(2.59096,0.09288)(2.41838,0.14183)(2.23383,0.18856)(2.03791,0.23280)(1.83126,0.27436)%
\polyline(-1.46071,-0.83234)(-1.47531,-0.84067)(-1.48992,-0.84899)(-1.50453,-0.85731)%
(-1.51913,-0.86564)(-1.53374,-0.87396)(-1.54835,-0.88229)(-1.56296,-0.89061)(-1.57756,-0.89893)%
(-1.59217,-0.90726)(-1.60678,-0.91558)(-1.62138,-0.92390)(-1.63599,-0.93223)(-1.65060,-0.94055)%
(-1.66521,-0.94887)(-1.67981,-0.95720)(-1.69442,-0.96552)(-1.70903,-0.97384)(-1.72363,-0.98217)%
(-1.73824,-0.99049)(-1.75285,-0.99881)(-1.76745,-1.00714)(-1.78206,-1.01546)(-1.79667,-1.02378)%
(-1.81128,-1.03211)(-1.82588,-1.04043)(-1.84049,-1.04875)(-1.85510,-1.05708)(-1.86970,-1.06540)%
(-1.88431,-1.07372)(-1.89892,-1.08205)(-1.91353,-1.09037)(-1.92813,-1.09869)(-1.94274,-1.10702)%
(-1.95735,-1.11534)(-1.97195,-1.12366)(-1.98656,-1.13199)(-2.00117,-1.14031)(-2.01577,-1.14864)%
(-2.03038,-1.15696)(-2.04499,-1.16528)(-2.05960,-1.17361)(-2.07420,-1.18193)(-2.08881,-1.19025)%
(-2.10342,-1.19858)(-2.11802,-1.20690)(-2.13263,-1.21522)(-2.14724,-1.22355)(-2.16185,-1.23187)%
(-2.17358,-1.23855)%
\polyline(-2.17358,-1.23855)(-2.17645,-1.24019)(-2.19106,-1.24852)%
\polyline(-1.58243,-0.90171)(-1.38584,-0.94504)(-1.18394,-0.98485)(-0.97742,-1.02103)%
(-0.76701,-1.05348)(-0.55343,-1.08214)(-0.33743,-1.10694)(-0.11975,-1.12781)(0.09883,-1.14471)%
(0.31725,-1.15761)(0.31756,-1.15762)(0.34760,-1.15886)%
\polyline(0.39166,-1.16065)(0.40890,-1.16135)%
\polyline(0.51116,-1.16552)(0.53567,-1.16653)(0.75239,-1.17141)(0.96694,-1.17228)%
(1.17856,-1.16916)(1.38647,-1.16207)(1.58994,-1.15109)(1.78819,-1.13623)(1.98051,-1.11758)%
(2.16615,-1.09525)(2.34441,-1.06928)(2.51460,-1.03984)(2.67606,-1.00700)(2.82814,-0.97092)%
(2.97019,-0.93173)(3.10166,-0.88960)(3.22196,-0.84469)(3.33056,-0.79718)(3.42698,-0.74726)%
(3.51076,-0.69512)(3.58147,-0.64098)(3.61366,-0.60955)(3.63874,-0.58505)(3.64414,-0.57792)%
(3.66321,-0.55271)%
\polyline(3.71640,-0.29901)(3.71361,-0.28670)(3.68500,-0.22501)(3.64165,-0.16324)%
(3.58364,-0.10165)(3.51101,-0.04050)(3.42393,0.01993)(3.32260,0.07940)(3.20727,0.13764)%
(3.07827,0.19440)(2.93597,0.24940)(2.78083,0.30242)(2.61330,0.35321)(2.43396,0.40152)%
(2.24342,0.44714)(2.04231,0.48985)(1.83135,0.52944)%
\polyline(-1.70416,-0.97107)(-1.49239,-1.00846)(-1.27481,-1.04198)(-1.05222,-1.07151)%
(-0.82548,-1.09698)(-0.59540,-1.11828)(-0.36286,-1.13539)(-0.12872,-1.14824)(0.10615,-1.15680)%
(0.34085,-1.16107)(0.37599,-1.16106)(0.39531,-1.16105)(0.40567,-1.16106)%
\polyline(0.44611,-1.16105)(0.45681,-1.16105)(0.47347,-1.16104)%
\polyline(0.54831,-1.16103)(0.57458,-1.16103)(0.80635,-1.15671)(1.03537,-1.14812)%
(1.26073,-1.13533)(1.48160,-1.11839)(1.69713,-1.09736)(1.90650,-1.07236)(2.10891,-1.04346)%
(2.30357,-1.01081)(2.48973,-0.97453)(2.66666,-0.93476)(2.83368,-0.89167)(2.99010,-0.84542)%
(3.13529,-0.79621)(3.26869,-0.74421)(3.38973,-0.68965)(3.49791,-0.63274)(3.59278,-0.57368)%
(3.67390,-0.51274)(3.70486,-0.48384)(3.72072,-0.46903)(3.74094,-0.45015)(3.74912,-0.44020)%
(3.76261,-0.42379)%
\polyline(3.85973,-0.21260)(3.86266,-0.18828)(3.85560,-0.12126)(3.83339,-0.05414)%
(3.79603,0.01280)(3.74361,0.07927)(3.67627,0.14503)(3.59420,0.20979)(3.49765,0.27329)%
(3.38693,0.33527)(3.26241,0.39547)(3.12453,0.45364)(2.97373,0.50952)(2.81058,0.56290)%
(2.63564,0.61353)(2.44955,0.66122)(2.25300,0.70573)(2.04671,0.74688)(1.83145,0.78450)%
\polyline(-1.82588,-1.04043)(-1.59893,-1.07188)(-1.36567,-1.09910)(-1.12703,-1.12199)%
(-0.88394,-1.14046)(-0.63737,-1.15443)(-0.38828,-1.16385)(-0.13767,-1.16867)(0.11348,-1.16890)%
(0.36415,-1.16451)(0.44799,-1.16149)(0.46588,-1.16085)%
\polyline(0.49193,-1.15991)(0.50314,-1.15950)(0.51911,-1.15893)(0.53918,-1.15821)%
\polyline(0.58547,-1.15654)(0.61347,-1.15553)(0.86033,-1.14200)(1.10379,-1.12396)%
(1.34290,-1.10150)(1.57672,-1.07469)(1.80432,-1.04364)(2.02481,-1.00848)(2.23731,-0.96934)%
(2.44100,-0.92638)(2.63505,-0.87977)(2.81873,-0.82969)(2.99129,-0.77634)(3.15205,-0.71993)%
(3.30039,-0.66068)(3.43572,-0.59882)(3.55750,-0.53461)(3.66525,-0.46829)(3.75856,-0.40012)%
(3.81052,-0.35394)(3.83705,-0.33037)(3.84021,-0.32682)(3.86595,-0.29796)%
\polyline(3.97065,-0.13723)(3.98080,-0.11444)(3.99752,-0.04118)(3.99848,0.03225)(3.98367,0.10554)%
(3.95317,0.17842)(3.90706,0.25060)(3.84557,0.32179)(3.76891,0.39171)(3.67738,0.46008)%
(3.57137,0.52665)(3.45127,0.59113)(3.31756,0.65329)(3.17078,0.71287)(3.01150,0.76964)%
(2.84034,0.82337)(2.65798,0.87386)(2.46514,0.92090)(2.26259,0.96432)(2.05111,1.00393)%
(1.83155,1.03958)%
\polyline(-1.94761,-1.10979)(-1.70547,-1.13530)(-1.45653,-1.15622)(-1.20183,-1.17247)%
(-0.94240,-1.18394)(-0.67934,-1.19057)(-0.41371,-1.19230)(-0.14662,-1.18911)(0.12081,-1.18099)%
(0.38752,-1.16795)(0.47423,-1.16209)(0.50623,-1.15992)(0.52679,-1.15853)(0.54313,-1.15743)%
(0.55641,-1.15653)(0.56858,-1.15570)(0.58562,-1.15455)(0.60499,-1.15324)%
\polyline(0.62265,-1.15204)(0.65237,-1.15004)(0.91429,-1.12729)(1.17221,-1.09980)%
(1.42507,-1.06766)(1.67184,-1.03099)(1.91151,-0.98992)(2.14312,-0.94460)(2.36571,-0.89522)%
(2.57843,-0.84196)(2.78038,-0.78502)(2.97079,-0.72462)(3.14890,-0.66101)(3.31401,-0.59443)%
(3.46549,-0.52515)(3.60275,-0.45344)(3.72527,-0.37957)(3.83260,-0.30384)(3.89862,-0.24822)%
(3.92434,-0.22655)(3.92713,-0.22366)(3.94002,-0.21031)(3.96566,-0.18376)(3.98214,-0.16669)%
(3.99083,-0.15769)(4.00018,-0.14800)(4.00264,-0.14473)%
\polyline(4.05987,-0.06850)(4.10321,0.01165)(4.13008,0.09212)(4.14043,0.17260)(4.13429,0.25278)%
(4.11174,0.33234)(4.07294,0.41098)(4.01811,0.48840)(3.94752,0.56430)(3.86154,0.63839)%
(3.76058,0.71037)(3.64509,0.78000)(3.51560,0.84700)(3.37271,0.91111)(3.21703,0.97210)%
(3.04925,1.02975)(2.87010,1.08384)(2.68031,1.13419)(2.48073,1.18059)(2.27217,1.22291)%
(2.05552,1.26097)(1.83165,1.29466)%
\polyline(-2.06933,-1.17915)(-1.81202,-1.19872)(-1.54740,-1.21335)(-1.27663,-1.22295)%
(-1.00087,-1.22743)(-0.72131,-1.22671)(-0.43914,-1.22075)(-0.15559,-1.20954)(0.12814,-1.19308)%
(0.41083,-1.17139)(0.47431,-1.16531)(0.57136,-1.15602)(0.60245,-1.15304)(0.61998,-1.15136)%
(0.63197,-1.15022)%
\polyline(0.64502,-1.14897)(0.65989,-1.14754)(0.69127,-1.14454)(0.96826,-1.11259)%
(0.97519,-1.11165)(1.24064,-1.07565)(1.50724,-1.03383)(1.76697,-0.98730)(2.01871,-0.93619)%
(2.26143,-0.88073)(2.49413,-0.82110)(2.71585,-0.75752)(2.92570,-0.69026)(3.12286,-0.61955)%
(3.30652,-0.54569)(3.47597,-0.46894)(3.63059,-0.38963)(3.76978,-0.30805)(3.89303,-0.22453)%
(3.98470,-0.15153)(3.99994,-0.13939)(4.01088,-0.12891)(4.02319,-0.11712)(4.03592,-0.10492)%
(4.06088,-0.08101)(4.07276,-0.06962)(4.08148,-0.06126)%
\polyline(4.12865,-0.00701)(4.13617,0.00197)(4.14576,0.01342)(4.15598,0.02562)(4.16333,0.03438)%
(4.16592,0.03845)(4.21934,0.12234)(4.25803,0.21055)(4.27936,0.29868)(4.28335,0.38637)%
(4.27011,0.47331)(4.23981,0.55914)(4.19272,0.64355)(4.12914,0.72621)(4.04948,0.80681)%
(3.95417,0.88507)(3.84376,0.96067)(3.71880,1.03336)(3.57994,1.10286)(3.42786,1.16893)%
(3.26329,1.23134)(3.08701,1.28987)(2.89984,1.34432)(2.70265,1.39452)(2.49632,1.44029)%
(2.28176,1.48149)(2.05992,1.51801)(1.83174,1.54972)%
\polyline(0.15992,-1.14032)(0.16146,-1.14142)(0.16299,-1.14252)(0.16452,-1.14363)%
(0.16605,-1.14473)(0.16758,-1.14583)(0.16912,-1.14694)(0.17064,-1.14805)(0.17218,-1.14915)%
(0.17371,-1.15026)(0.17524,-1.15136)(0.17677,-1.15246)(0.17831,-1.15357)(0.17983,-1.15468)%
(0.18137,-1.15578)(0.18290,-1.15688)(0.18444,-1.15799)(0.18596,-1.15909)(0.18750,-1.16020)%
(0.18903,-1.16130)(0.19056,-1.16241)(0.19209,-1.16351)(0.19362,-1.16462)(0.19515,-1.16572)%
(0.19669,-1.16682)(0.19822,-1.16793)(0.19975,-1.16904)(0.20128,-1.17014)(0.20281,-1.17124)%
(0.20434,-1.17235)(0.20587,-1.17345)(0.20741,-1.17456)(0.20894,-1.17566)(0.21047,-1.17677)%
(0.21200,-1.17787)(0.21353,-1.17898)(0.21506,-1.18008)(0.21660,-1.18118)(0.21813,-1.18229)%
(0.21966,-1.18340)(0.22119,-1.18450)(0.22272,-1.18560)(0.22425,-1.18671)(0.22579,-1.18781)%
(0.22732,-1.18892)(0.22885,-1.19003)(0.23038,-1.19112)(0.23192,-1.19223)(0.23344,-1.19334)%
(0.23498,-1.19445)(0.23651,-1.19554)%
\polyline(0.16299,-1.14252)(0.16452,-1.14363)(0.16605,-1.14473)(0.16758,-1.14583)%
(0.16912,-1.14694)(0.17064,-1.14805)(0.17218,-1.14915)(0.17281,-1.14961)(0.17371,-1.15026)%
(0.17524,-1.15136)(0.17677,-1.15246)(0.17831,-1.15357)(0.17983,-1.15468)(0.18137,-1.15578)%
(0.18290,-1.15688)(0.18444,-1.15799)(0.18596,-1.15909)(0.18750,-1.16020)(0.18903,-1.16130)%
(0.19056,-1.16241)(0.19209,-1.16351)(0.19362,-1.16462)(0.19515,-1.16572)(0.19669,-1.16682)%
(0.19822,-1.16793)(0.19975,-1.16904)(0.20128,-1.17014)(0.20281,-1.17124)(0.20434,-1.17235)%
(0.20587,-1.17345)(0.20741,-1.17456)(0.20894,-1.17566)(0.21047,-1.17677)(0.21200,-1.17787)%
(0.21353,-1.17898)(0.21506,-1.18008)(0.21660,-1.18118)(0.21813,-1.18229)(0.21966,-1.18340)%
(0.22119,-1.18450)(0.22272,-1.18560)(0.22425,-1.18671)(0.22579,-1.18781)(0.22732,-1.18892)%
(0.22885,-1.19003)(0.23038,-1.19112)(0.23192,-1.19223)(0.23344,-1.19334)(0.23498,-1.19445)%
(0.23574,-1.19499)%
\polyline(1.79431,-1.19437)(1.80852,-1.18634)(1.80992,-1.18555)(1.82354,-1.17785)%
(1.83857,-1.16936)(1.85360,-1.16086)(1.86863,-1.15237)(1.88366,-1.14387)(1.89869,-1.13538)%
(1.91371,-1.12688)(1.92874,-1.11839)(1.94377,-1.10990)(1.95880,-1.10140)(1.97383,-1.09291)%
(1.98886,-1.08441)(2.00389,-1.07592)(2.01892,-1.06742)(2.03394,-1.05893)(2.04897,-1.05044)%
(2.06400,-1.04194)(2.07903,-1.03345)(2.09406,-1.02495)(2.10909,-1.01646)(2.12412,-1.00797)%
(2.13914,-0.99947)(2.15417,-0.99098)(2.16920,-0.98248)(2.18423,-0.97399)(2.19926,-0.96549)%
(2.21429,-0.95700)(2.22932,-0.94851)(2.24435,-0.94001)(2.25937,-0.93152)(2.27440,-0.92302)%
(2.28943,-0.91453)(2.30446,-0.90604)(2.31949,-0.89754)(2.33452,-0.88905)(2.34954,-0.88055)%
(2.36458,-0.87206)(2.37960,-0.86357)(2.39463,-0.85507)(2.40966,-0.84658)(2.42469,-0.83808)%
(2.43972,-0.82959)(2.45475,-0.82109)(2.46977,-0.81260)(2.48481,-0.80410)(2.49983,-0.79561)%
(2.51486,-0.78712)(2.52925,-0.77898)(2.52989,-0.77862)(2.54492,-0.77013)%
\polyline(3.05564,-0.99927)(3.07577,-0.98065)(3.09591,-0.96203)(3.11604,-0.94341)%
(3.13617,-0.92479)(3.15630,-0.90617)(3.17644,-0.88755)(3.19657,-0.86893)(3.21670,-0.85031)%
(3.23683,-0.83169)(3.25697,-0.81307)(3.27710,-0.79445)(3.29723,-0.77584)(3.31736,-0.75722)%
(3.33750,-0.73859)(3.35763,-0.71997)(3.37776,-0.70136)(3.39789,-0.68274)(3.41803,-0.66412)%
(3.43816,-0.64549)(3.45829,-0.62688)(3.47842,-0.60826)(3.49856,-0.58964)(3.51869,-0.57102)%
(3.53882,-0.55240)(3.55895,-0.53378)(3.57909,-0.51516)(3.59921,-0.49655)(3.61935,-0.47793)%
(3.63948,-0.45930)(3.65962,-0.44068)(3.67974,-0.42207)(3.69988,-0.40345)(3.72001,-0.38483)%
(3.74015,-0.36620)(3.76027,-0.34759)(3.78041,-0.32897)(3.80054,-0.31035)(3.82068,-0.29173)%
(3.84081,-0.27311)(3.86094,-0.25449)(3.88107,-0.23587)(3.90121,-0.21725)(3.92134,-0.19864)%
(3.94147,-0.18001)(3.96160,-0.16139)(3.98174,-0.14277)(4.00100,-0.12496)(4.00187,-0.12416)%
(4.01956,-0.10779)(4.02200,-0.10554)(4.04015,-0.08875)(4.04214,-0.08691)(4.06196,-0.06858)%
\polyline(3.81670,-0.24048)(3.81930,-0.23613)(3.83528,-0.20939)(3.85126,-0.18266)%
(3.86724,-0.15592)(3.88323,-0.12919)(3.89920,-0.10246)(3.91518,-0.07572)(3.93117,-0.04899)%
(3.94714,-0.02226)(3.96313,0.00447)(3.97911,0.03121)(3.99508,0.05795)(4.01107,0.08468)%
(4.02704,0.11142)(4.04303,0.13815)(4.05901,0.16488)(4.07498,0.19161)(4.09097,0.21835)%
(4.10695,0.24508)(4.12293,0.27182)(4.13891,0.29856)(4.15489,0.32529)(4.17087,0.35202)%
(4.18685,0.37875)(4.20283,0.40549)(4.21881,0.43222)(4.23479,0.45895)(4.25077,0.48569)%
(4.26675,0.51243)(4.28273,0.53916)(4.29871,0.56589)(4.31469,0.59263)(4.33068,0.61936)%
(4.34665,0.64609)(4.36263,0.67282)(4.37862,0.69956)(4.39336,0.72423)%
\polyline(3.18677,-0.13766)(3.19373,-0.10679)(3.20069,-0.07591)(3.20765,-0.04504)%
(3.21461,-0.01417)(3.22157,0.01670)(3.22852,0.04758)(3.23549,0.07845)(3.24244,0.10932)%
(3.24940,0.14020)(3.25635,0.17107)(3.26332,0.20194)(3.27028,0.23280)(3.27724,0.26368)%
(3.28419,0.29455)(3.29115,0.32542)(3.29811,0.35629)(3.30507,0.38717)(3.31202,0.41804)%
(3.31899,0.44891)(3.32594,0.47979)(3.33291,0.51066)(3.33986,0.54153)(3.34682,0.57240)%
(3.35378,0.60328)(3.36074,0.63415)(3.36769,0.66502)(3.37465,0.69589)(3.38161,0.72677)%
(3.38857,0.75764)(3.39553,0.78851)(3.40249,0.81938)(3.40945,0.85026)(3.41641,0.88113)%
(3.42337,0.91200)(3.43032,0.94288)(3.43728,0.97375)(3.44424,1.00462)(3.45120,1.03550)%
(3.45816,1.06636)(3.46512,1.09723)(3.47208,1.12810)(3.47904,1.15898)(3.48599,1.18985)%
(3.49295,1.22072)(3.49991,1.25159)(3.50687,1.28247)(3.51382,1.31334)(3.52078,1.34421)%
(3.52774,1.37508)(3.53449,1.40503)%
\polyline(1.82588,0.27538)(1.82588,0.30530)(1.82588,0.30598)(1.82588,0.33591)(1.82588,0.33658)%
(1.82588,0.36652)(1.82588,0.36719)(1.82588,0.39713)(1.82588,0.39779)(1.82588,0.42773)%
(1.82588,0.42839)(1.82588,0.45834)(1.82588,0.45899)(1.82588,0.48894)(1.82588,0.48959)%
(1.82588,0.51955)(1.82588,0.52020)(1.82588,0.55016)(1.82588,0.55080)(1.82588,0.58077)%
(1.82588,0.58140)(1.82588,0.61137)(1.82588,0.61200)(1.82588,0.64198)(1.82588,0.64260)%
(1.82588,0.67258)(1.82588,0.67321)(1.82588,0.70319)(1.82588,0.70381)(1.82588,0.73380)%
(1.82588,0.73441)(1.82588,0.76441)(1.82588,0.76501)(1.82588,0.79501)(1.82588,0.79561)%
(1.82588,0.82561)(1.82588,0.82622)(1.82588,0.85623)(1.82588,0.85682)(1.82588,0.88683)%
(1.82588,0.88742)(1.82588,0.91744)(1.82588,0.91802)(1.82588,0.94804)(1.82588,0.94862)%
(1.82588,0.97865)(1.82588,0.97923)(1.82588,1.00925)(1.82588,1.00983)(1.82588,1.03987)%
(1.82588,1.04043)(1.82588,1.07047)(1.82588,1.07103)(1.82588,1.10108)(1.82588,1.10163)%
(1.82588,1.13168)(1.82588,1.13224)(1.82588,1.16229)(1.82588,1.16284)(1.82588,1.19290)%
(1.82588,1.19344)(1.82588,1.22350)(1.82588,1.22404)(1.82588,1.25411)(1.82588,1.25464)%
(1.82588,1.28471)(1.82588,1.28525)(1.82588,1.31532)(1.82588,1.31585)(1.82588,1.34592)%
(1.82588,1.34645)(1.82588,1.37654)(1.82588,1.37705)(1.82588,1.40714)(1.82588,1.40766)%
(1.82588,1.43775)(1.82588,1.43826)(1.82588,1.46835)(1.82588,1.46886)(1.82588,1.49896)%
(1.82588,1.49946)(1.82588,1.52957)(1.82588,1.53006)(1.82588,1.56018)(1.82588,1.56067)%
(1.82588,1.59078)(1.82588,1.59127)(1.82588,1.62138)(1.82588,1.62187)(1.82588,1.65199)%
(1.82588,1.65247)(1.82588,1.68260)(1.82588,1.68307)(1.82588,1.71321)(1.82588,1.71368)%
(1.82588,1.74381)(1.82588,1.74428)(1.82588,1.77442)(1.82588,1.77488)(1.82588,1.80548)%
\end{picture}}%

%% file: fig/p016.tex
{\unitlength=1cm%
\begin{picture}%
(9.94,5.4)(-5.19,-2.08)%
\linethickness{0.008in}
\polyline(0.70986,-1.14168)(0.70255,-1.14262)(0.69318,-1.14380)(0.69294,-1.14382)%
\polyline(0.69294,-1.14382)(0.68387,-1.14492)(0.67463,-1.14602)(0.67440,-1.14604)%
\polyline(0.67440,-1.14604)(0.66546,-1.14708)(0.65636,-1.14810)(0.64731,-1.14908)%
(0.63833,-1.15002)(0.63113,-1.15075)%
\polyline(0.63113,-1.15075)(0.62942,-1.15092)(0.62056,-1.15179)(0.61178,-1.15262)%
(0.60307,-1.15342)(0.59442,-1.15417)(0.58583,-1.15489)(0.57731,-1.15557)(0.56885,-1.15622)%
(0.56047,-1.15683)(0.55213,-1.15739)(0.54387,-1.15793)(0.53567,-1.15842)(0.52754,-1.15888)%
(0.51947,-1.15931)(0.51147,-1.15969)(0.50352,-1.16004)(0.49565,-1.16035)(0.48785,-1.16063)%
(0.48631,-1.16068)%
\polyline(0.48631,-1.16068)(0.48010,-1.16087)(0.47355,-1.16105)%
\polyline(0.47355,-1.16105)(0.47241,-1.16108)(0.46479,-1.16125)(0.46101,-1.16131)%
\polyline(0.46101,-1.16131)(0.45724,-1.16138)(0.44975,-1.16148)(0.44864,-1.16148)%
\polyline(0.44864,-1.16148)(0.44233,-1.16154)(0.43569,-1.16156)%
\polyline(0.43569,-1.16156)(0.43497,-1.16156)(0.42767,-1.16155)(0.42333,-1.16152)%
\polyline(0.42333,-1.16152)(0.42044,-1.16151)(0.41327,-1.16143)(0.40616,-1.16131)%
(0.39913,-1.16116)(0.39214,-1.16097)(0.38523,-1.16075)(0.37838,-1.16049)(0.37159,-1.16020)%
(0.36487,-1.15987)(0.35851,-1.15953)%
\polyline(0.35851,-1.15953)(0.35821,-1.15951)(0.35162,-1.15912)(0.34703,-1.15882)%
\polyline(0.34703,-1.15882)(0.34508,-1.15869)(0.33861,-1.15822)(0.33237,-1.15773)%
\polyline(0.33237,-1.15773)(0.33221,-1.15772)(0.32586,-1.15719)(0.31958,-1.15662)%
\polyline(4.17199,0.02966)(4.16206,0.01982)(4.16181,0.01958)%
\polyline(4.16181,0.01958)(4.15187,0.00968)(4.15161,0.00941)%
\polyline(4.15161,0.00941)(4.14074,-0.00148)(4.14045,-0.00177)%
\polyline(4.14045,-0.00177)(4.12932,-0.01299)(4.12904,-0.01329)%
\polyline(4.12904,-0.01329)(4.11790,-0.02458)(4.11762,-0.02487)%
\polyline(4.11762,-0.02487)(4.10649,-0.03625)(4.10621,-0.03654)%
\polyline(4.10621,-0.03654)(4.09508,-0.04798)(4.09480,-0.04828)%
\polyline(4.09480,-0.04828)(4.08366,-0.05980)(4.07225,-0.07170)(4.06085,-0.08366)%
(4.04943,-0.09570)(4.03803,-0.10783)(4.02662,-0.12003)(4.01521,-0.13231)(4.00927,-0.13875)%
(4.00335,-0.14517)(3.99101,-0.15864)(3.97867,-0.17222)(3.96632,-0.18589)(3.95398,-0.19965)%
(3.94165,-0.21351)(3.92931,-0.22746)(3.91699,-0.24151)(3.90466,-0.25565)(3.89235,-0.26990)%
(3.88003,-0.28424)(3.86773,-0.29868)(3.85483,-0.31394)(3.84181,-0.32943)(3.82881,-0.34505)%
(3.81582,-0.36076)(3.80284,-0.37660)(3.78987,-0.39254)(3.77692,-0.40859)(3.76398,-0.42476)%
(3.75105,-0.44104)(3.73813,-0.45742)(3.72524,-0.47392)(3.72065,-0.47981)(3.71201,-0.49100)%
(3.69859,-0.50843)(3.68521,-0.52599)(3.67185,-0.54367)(3.65851,-0.56147)(3.64519,-0.57939)%
(3.63191,-0.59743)(3.61865,-0.61558)(3.60542,-0.63386)(3.59353,-0.65043)%
\polyline(3.58985,-0.65559)(3.57929,-0.67041)%
\polyline(3.57772,-0.67264)(3.57592,-0.67520)(3.56688,-0.68802)%
\polyline(3.56688,-0.68802)(3.56559,-0.68985)(3.55342,-0.70726)%
\polyline(-2.19106,-1.24852)(-1.84872,-1.26473)(-1.49467,-1.27264)(-1.13138,-1.27204)%
(-0.76138,-1.26278)(-0.38725,-1.24478)(-0.01158,-1.21803)(0.36300,-1.18258)(0.47143,-1.16971)%
\polyline(0.47143,-1.16971)(0.63112,-1.15075)%
\polyline(0.63112,-1.15075)(0.67459,-1.14560)%
\polyline(0.67459,-1.14560)(0.69313,-1.14339)%
\polyline(0.69313,-1.14339)(0.70723,-1.14172)%
\polyline(0.70723,-1.14172)(0.73391,-1.13855)(1.09861,-1.08612)(1.45459,-1.02555)%
(1.79944,-0.95714)(2.13084,-0.88127)(2.44659,-0.79835)(2.74460,-0.70887)(3.02297,-0.61335)%
(3.27992,-0.51238)(3.51386,-0.40654)(3.53112,-0.39747)%
\polyline(3.53112,-0.39747)(3.72337,-0.29650)(3.90725,-0.18291)(4.06447,-0.06649)%
(4.07784,-0.05428)%
\polyline(4.07784,-0.05428)(4.10251,-0.03174)%
\polyline(4.10251,-0.03174)(4.12646,-0.00986)%
\polyline(4.12646,-0.00986)(4.15005,0.01168)%
\polyline(4.15005,0.01168)(4.16943,0.02939)%
\polyline(4.16943,0.02939)(4.19423,0.05205)(4.29593,0.17199)(4.36917,0.29260)(4.41379,0.41312)%
(4.42981,0.53285)(4.41746,0.65105)(4.37719,0.76702)(4.30962,0.88008)(4.21558,0.98958)%
(4.09606,1.09486)(3.95219,1.19534)(3.78530,1.29048)(3.59682,1.37974)(3.38833,1.46266)%
(3.16150,1.53880)(2.91812,1.60780)(2.66002,1.66931)(2.38914,1.72308)(2.10743,1.76887)%
(1.81687,1.80651)(1.51950,1.83588)(1.21732,1.85692)(0.91235,1.86961)(0.60653,1.87398)%
(0.30181,1.87011)(0.00007,1.85812)(-0.29688,1.83819)(-0.58728,1.81052)(-0.86947,1.77536)%
(-1.14191,1.73300)%
\polyline(-1.33193,0.49110)(-1.32813,0.51594)(-1.32433,0.54078)(-1.32053,0.56561)%
(-1.31673,0.59045)(-1.31293,0.61529)(-1.30912,0.64013)(-1.30533,0.66497)(-1.30153,0.68980)%
(-1.29773,0.71464)(-1.29393,0.73948)(-1.29013,0.76432)(-1.28632,0.78916)(-1.28253,0.81400)%
(-1.27873,0.83884)(-1.27492,0.86367)(-1.27113,0.88851)(-1.26732,0.91334)(-1.26352,0.93818)%
(-1.25973,0.96303)(-1.25592,0.98786)(-1.25212,1.01270)(-1.24832,1.03754)(-1.24452,1.06238)%
(-1.24072,1.08721)(-1.23692,1.11205)(-1.23312,1.13690)(-1.22932,1.16173)(-1.22552,1.18657)%
(-1.22172,1.21140)(-1.21791,1.23624)(-1.21412,1.26108)(-1.21032,1.28592)(-1.20651,1.31076)%
(-1.20272,1.33559)(-1.19891,1.36044)(-1.19511,1.38527)(-1.19132,1.41011)(-1.18751,1.43495)%
(-1.18371,1.45979)(-1.17991,1.48463)(-1.17611,1.50946)(-1.17231,1.53430)(-1.16852,1.55914)%
(-1.16471,1.58398)(-1.16090,1.60881)(-1.15710,1.63365)(-1.15331,1.65850)(-1.14951,1.68333)%
(-1.14570,1.70817)(-1.14191,1.73300)%
\polyline(-1.46071,-0.83234)(-1.23288,-0.89353)(-0.99773,-0.94970)(-0.75644,-1.00064)%
(-0.51020,-1.04622)(-0.26024,-1.08628)(-0.00781,-1.12070)(0.24581,-1.14937)(0.32624,-1.15662)%
\polyline(0.32624,-1.15662)(0.33898,-1.15777)%
\polyline(0.33898,-1.15777)(0.35850,-1.15953)%
\polyline(0.47143,-1.16971)(0.49935,-1.17223)(0.75149,-1.18921)(1.00092,-1.20029)%
(1.24633,-1.20546)(1.48637,-1.20475)(1.71973,-1.19819)(1.94508,-1.18588)(2.16114,-1.16789)%
(2.36661,-1.14437)(2.56022,-1.11548)(2.74075,-1.08140)(2.90701,-1.04235)(3.05783,-0.99857)%
(3.19214,-0.95033)(3.30891,-0.89795)(3.40717,-0.84175)(3.48603,-0.78207)(3.54473,-0.71933)%
(3.55154,-0.70753)%
\polyline(3.53113,-0.39747)(3.51481,-0.37433)(3.44166,-0.30234)(3.34587,-0.23050)%
(3.22764,-0.15935)(3.08728,-0.08938)(2.92527,-0.02112)(2.74223,0.04491)(2.53896,0.10822)%
(2.31637,0.16831)(2.07556,0.22469)(1.81775,0.27691)(1.54431,0.32454)(1.25677,0.36716)%
(0.95676,0.40439)(0.64604,0.43591)(0.32649,0.46140)(0.00007,0.48061)(-0.33114,0.49332)%
(-0.66501,0.49938)(-0.99936,0.49867)(-1.33193,0.49110)%
\polyline(-1.46071,-0.83234)(-1.47531,-0.84067)(-1.48992,-0.84899)(-1.50453,-0.85731)%
(-1.51913,-0.86564)(-1.53374,-0.87396)(-1.54835,-0.88229)(-1.56296,-0.89061)(-1.57756,-0.89893)%
(-1.59217,-0.90726)(-1.60678,-0.91558)(-1.62138,-0.92390)(-1.63599,-0.93223)(-1.65060,-0.94055)%
(-1.66521,-0.94887)(-1.67981,-0.95720)(-1.69442,-0.96552)(-1.70903,-0.97384)(-1.72363,-0.98217)%
(-1.73824,-0.99049)(-1.75285,-0.99881)(-1.76745,-1.00714)(-1.78206,-1.01546)(-1.79667,-1.02378)%
(-1.81128,-1.03211)(-1.82588,-1.04043)(-1.84049,-1.04875)(-1.85510,-1.05708)(-1.86970,-1.06540)%
(-1.88431,-1.07372)(-1.89892,-1.08205)(-1.91353,-1.09037)(-1.92813,-1.09869)(-1.94274,-1.10702)%
(-1.95735,-1.11534)(-1.97195,-1.12366)(-1.98656,-1.13199)(-2.00117,-1.14031)(-2.01577,-1.14864)%
(-2.03038,-1.15696)(-2.04499,-1.16528)(-2.05960,-1.17361)(-2.07420,-1.18193)(-2.08881,-1.19025)%
(-2.10342,-1.19858)(-2.11802,-1.20690)(-2.13263,-1.21522)(-2.14724,-1.22355)(-2.16185,-1.23187)%
(-2.17645,-1.24019)(-2.19106,-1.24852)%
\polyline(-1.58243,-0.90171)(-1.33552,-0.95539)(-1.08055,-1.00352)(-0.81893,-1.04587)%
(-0.55206,-1.08232)(-0.28141,-1.11270)(-0.00844,-1.13692)(0.26289,-1.15474)(0.26535,-1.15491)%
(0.33380,-1.15784)(0.34512,-1.15833)(0.35853,-1.15890)(0.37331,-1.15954)(0.38561,-1.16006)%
(0.39950,-1.16066)(0.41907,-1.16150)%
\polyline(0.50837,-1.16532)(0.53845,-1.16661)(0.80934,-1.17203)(1.07654,-1.17117)%
(1.33851,-1.16408)(1.59379,-1.15083)(1.84087,-1.13155)(2.07835,-1.10637)(2.30478,-1.07546)%
(2.51883,-1.03903)(2.71916,-0.99732)(2.90452,-0.95058)(3.07371,-0.89910)(3.22560,-0.84322)%
(3.35916,-0.78327)(3.47341,-0.71963)(3.56750,-0.65268)(3.61264,-0.60962)(3.64066,-0.58287)%
(3.64331,-0.57916)(3.66204,-0.55293)%
\polyline(3.71605,-0.30034)(3.71265,-0.28396)(3.67369,-0.20672)(3.61168,-0.12947)%
(3.52675,-0.05272)(3.41912,0.02299)(3.28917,0.09717)(3.13746,0.16929)(2.96464,0.23886)%
(2.77155,0.30539)(2.55914,0.36840)(2.32850,0.42744)(2.08087,0.48205)(1.81760,0.53185)%
(1.54018,0.57642)(1.25020,0.61545)(0.94935,0.64859)(0.63945,0.67558)(0.32237,0.69619)%
(0.00007,0.71019)(-0.32543,0.71747)(-0.65206,0.71791)(-0.97771,0.71144)(-1.30026,0.69808)%
\polyline(-1.70416,-0.97107)(-1.43816,-1.01726)(-1.16338,-1.05734)(-0.88142,-1.09111)%
(-0.59392,-1.11841)(-0.30257,-1.13911)(-0.00907,-1.15314)(0.28284,-1.16039)(0.28488,-1.16044)%
(0.38192,-1.16062)(0.40617,-1.16067)(0.42232,-1.16071)(0.43568,-1.16073)(0.45743,-1.16077)%
(0.47143,-1.16080)(0.48190,-1.16082)%
\polyline(0.54531,-1.16094)(0.57754,-1.16101)(0.86720,-1.15485)(1.15215,-1.14204)%
(1.43070,-1.12269)(1.70120,-1.09692)(1.96201,-1.06491)(2.21160,-1.02687)(2.44842,-0.98304)%
(2.67105,-0.93371)(2.87810,-0.87916)(3.06829,-0.81977)(3.24042,-0.75587)(3.39338,-0.68787)%
(3.52618,-0.61620)(3.63791,-0.54130)(3.70502,-0.48334)(3.72784,-0.46363)(3.73481,-0.45536)%
\polyline(3.86084,-0.21892)(3.85831,-0.13515)(3.83204,-0.05115)(3.78207,0.03254)(3.70856,0.11540)%
(3.61183,0.19689)(3.49235,0.27649)(3.35070,0.35368)(3.18763,0.42797)(3.00401,0.49885)%
(2.80087,0.56587)(2.57931,0.62859)(2.34062,0.68656)(2.08618,0.73942)(1.81746,0.78678)%
(1.53605,0.82832)(1.24363,0.86374)(0.94195,0.89280)(0.63287,0.91526)(0.31826,0.93097)%
(0.00007,0.93978)(-0.31972,0.94162)(-0.63910,0.93642)(-0.95606,0.92423)(-1.26858,0.90507)%
\polyline(-1.82588,-1.04043)(-1.54080,-1.07913)(-1.24621,-1.11117)(-0.94391,-1.13634)%
(-0.63579,-1.15450)(-0.32374,-1.16553)(-0.00970,-1.16936)(0.30441,-1.16598)(0.43464,-1.16156)%
(0.46095,-1.16067)%
\polyline(0.47111,-1.16032)(0.48628,-1.15981)(0.51184,-1.15894)(0.52793,-1.15840)%
(0.54622,-1.15778)%
\polyline(0.58225,-1.15655)(0.61663,-1.15539)(0.92505,-1.13767)(1.22776,-1.11292)%
(1.52288,-1.08130)(1.80861,-1.04301)(2.08315,-0.99827)(2.34484,-0.94737)(2.59206,-0.89062)%
(2.82327,-0.82837)(3.03704,-0.76101)(3.23206,-0.68895)(3.40713,-0.61263)(3.56115,-0.53253)%
(3.69319,-0.44914)(3.80242,-0.36298)(3.81033,-0.35482)(3.83794,-0.32636)(3.86620,-0.29722)%
(3.88008,-0.28292)(3.88816,-0.27458)(3.89222,-0.26865)%
\polyline(3.96939,-0.13690)(3.98726,-0.09324)(4.00000,-0.00143)(3.98803,0.09039)(3.95144,0.18166)%
(3.89045,0.27180)(3.80544,0.36026)(3.69692,0.44650)(3.56558,0.52999)(3.41223,0.61019)%
(3.23781,0.68664)(3.04339,0.75884)(2.83018,0.82636)(2.59949,0.88877)(2.35275,0.94569)%
(2.09149,0.99678)(1.81731,1.04171)(1.53190,1.08021)(1.23705,1.11204)(0.93455,1.13700)%
(0.62628,1.15495)(0.31414,1.16576)(0.00007,1.16937)(-0.31400,1.16576)(-0.62615,1.15495)%
(-0.93442,1.13701)(-1.23692,1.11205)%
\polyline(-1.94761,-1.10979)(-1.64344,-1.14100)(-1.32903,-1.16499)(-1.00640,-1.18157)%
(-0.67766,-1.19059)(-0.34491,-1.19195)(-0.01032,-1.18558)(0.32394,-1.17151)(0.46550,-1.16223)%
(0.49265,-1.16045)(0.52749,-1.15818)(0.54381,-1.15711)(0.55748,-1.15621)(0.56922,-1.15544)%
(0.58622,-1.15433)(0.61061,-1.15273)(0.65573,-1.14977)(0.98291,-1.12048)(1.30337,-1.08380)%
(1.61507,-1.03992)(1.91601,-0.98909)(2.20431,-0.93163)(2.47811,-0.86787)(2.73570,-0.79820)%
(2.97549,-0.72304)(3.19598,-0.64285)(3.39583,-0.55813)(3.57383,-0.46939)(3.72892,-0.37719)%
(3.86020,-0.28208)(3.89775,-0.24781)(3.92506,-0.22287)(3.95176,-0.19849)(3.96692,-0.18466)%
(3.97777,-0.17148)%
\polyline(4.05648,-0.07125)(4.10453,0.01472)(4.13478,0.11546)(4.13915,0.21607)(4.11776,0.31594)%
(4.07084,0.41446)(3.99883,0.51106)(3.90231,0.60513)(3.78201,0.69612)(3.63882,0.78348)%
(3.47376,0.86671)(3.28798,0.94531)(3.08276,1.01883)(2.85949,1.08684)(2.61967,1.14895)%
(2.36488,1.20483)(2.09680,1.25414)(1.81716,1.29665)(1.52777,1.33211)(1.23048,1.36034)%
(0.92715,1.38121)(0.61970,1.39463)(0.31003,1.40054)(0.00007,1.39895)(-0.30830,1.38990)%
(-0.61319,1.37348)(-0.91277,1.34980)(-1.20525,1.31903)%
\polyline(-2.06933,-1.17915)(-1.74608,-1.20287)(-1.41185,-1.21881)(-1.06889,-1.22681)%
(-0.71952,-1.22668)(-0.36608,-1.21836)(-0.01096,-1.20181)(0.34347,-1.17704)(0.46682,-1.16550)%
(0.55703,-1.15705)(0.58578,-1.15437)(0.60301,-1.15276)(0.62050,-1.15112)(0.63103,-1.15013)%
(0.65625,-1.14777)(0.67293,-1.14621)(0.69482,-1.14417)(1.04076,-1.10331)(1.37898,-1.05467)%
(1.70725,-0.99852)(2.02343,-0.93518)(2.32544,-0.86499)(2.61135,-0.78837)(2.87933,-0.70578)%
(3.12770,-0.61771)(3.35492,-0.52469)(3.55960,-0.42732)(3.74054,-0.32615)(3.89669,-0.22184)%
(3.99799,-0.13893)(4.01139,-0.12796)(4.02722,-0.11501)(4.03607,-0.10579)(4.04754,-0.09382)%
\polyline(4.13741,0.00215)(4.14633,0.01482)(4.15447,0.02636)(4.16094,0.03557)(4.20884,0.10353)%
(4.25916,0.21392)(4.28229,0.32415)(4.27831,0.43356)(4.24747,0.54148)(4.19023,0.64728)%
(4.10721,0.75032)(3.99918,0.84999)(3.86710,0.94573)(3.71206,1.03698)(3.53529,1.12322)%
(3.33815,1.20399)(3.12213,1.27882)(2.88881,1.34732)(2.63985,1.40913)(2.37701,1.46395)%
(2.10211,1.51151)(1.81702,1.55158)(1.52364,1.58400)(1.22390,1.60863)(0.91975,1.62542)%
(0.61311,1.63431)(0.30592,1.63532)(0.00007,1.62854)(-0.30258,1.61404)(-0.60023,1.59199)%
(-0.89112,1.56258)(-1.17358,1.52602)%
\polyline(0.16146,-1.14142)(0.16299,-1.14252)(0.16452,-1.14363)(0.16605,-1.14473)%
(0.16758,-1.14583)(0.16912,-1.14694)(0.17020,-1.14772)(0.17064,-1.14805)(0.17218,-1.14915)%
(0.17371,-1.15026)(0.17524,-1.15136)(0.17677,-1.15246)(0.17831,-1.15357)(0.17983,-1.15468)%
(0.18137,-1.15578)(0.18290,-1.15688)(0.18444,-1.15799)(0.18596,-1.15909)(0.18750,-1.16020)%
(0.18903,-1.16130)(0.19056,-1.16241)(0.19209,-1.16351)(0.19362,-1.16462)(0.19515,-1.16572)%
(0.19669,-1.16682)(0.19822,-1.16793)(0.19975,-1.16904)(0.20128,-1.17014)(0.20281,-1.17124)%
(0.20434,-1.17235)(0.20587,-1.17345)(0.20741,-1.17456)(0.20894,-1.17566)(0.21047,-1.17677)%
(0.21200,-1.17787)(0.21353,-1.17898)(0.21506,-1.18008)(0.21660,-1.18118)(0.21813,-1.18229)%
(0.21966,-1.18340)(0.22119,-1.18450)(0.22272,-1.18560)(0.22425,-1.18671)(0.22579,-1.18781)%
(0.22732,-1.18892)(0.22885,-1.19003)(0.23038,-1.19112)(0.23192,-1.19223)(0.23344,-1.19334)%
(0.23498,-1.19445)(0.23530,-1.19467)%
\polyline(1.79481,-1.19409)(1.80852,-1.18634)(1.82354,-1.17785)(1.83857,-1.16936)%
(1.85360,-1.16086)(1.86863,-1.15237)(1.88366,-1.14387)(1.89869,-1.13538)(1.91371,-1.12688)%
(1.92874,-1.11839)(1.94377,-1.10990)(1.95880,-1.10140)(1.97383,-1.09291)(1.98886,-1.08441)%
(2.00389,-1.07592)(2.01892,-1.06742)(2.03394,-1.05893)(2.04897,-1.05044)(2.06400,-1.04194)%
(2.07903,-1.03345)(2.09406,-1.02495)(2.10909,-1.01646)(2.12412,-1.00797)(2.13914,-0.99947)%
(2.15417,-0.99098)(2.16920,-0.98248)(2.18423,-0.97399)(2.19926,-0.96549)(2.21429,-0.95700)%
(2.22932,-0.94851)(2.24435,-0.94001)(2.25937,-0.93152)(2.27440,-0.92302)(2.28943,-0.91453)%
(2.30446,-0.90604)(2.31949,-0.89754)(2.33452,-0.88905)(2.34954,-0.88055)(2.36458,-0.87206)%
(2.37960,-0.86357)(2.39463,-0.85507)(2.40966,-0.84658)(2.42469,-0.83808)(2.43972,-0.82959)%
(2.45475,-0.82109)(2.46977,-0.81260)(2.48481,-0.80410)(2.49983,-0.79561)(2.51486,-0.78712)%
(2.52989,-0.77862)(2.54492,-0.77013)%
\polyline(3.05564,-0.99927)(3.07577,-0.98065)(3.09591,-0.96203)(3.11604,-0.94341)%
(3.13617,-0.92479)(3.15630,-0.90617)(3.17644,-0.88755)(3.19657,-0.86893)(3.21670,-0.85031)%
(3.23683,-0.83169)(3.25697,-0.81307)(3.27710,-0.79445)(3.29723,-0.77584)(3.31736,-0.75722)%
(3.33750,-0.73859)(3.35763,-0.71997)(3.37776,-0.70136)(3.39789,-0.68274)(3.41803,-0.66412)%
(3.43816,-0.64549)(3.45829,-0.62688)(3.47842,-0.60826)(3.49856,-0.58964)(3.51869,-0.57102)%
(3.53882,-0.55240)(3.55895,-0.53378)(3.57909,-0.51516)(3.59921,-0.49655)(3.61935,-0.47793)%
(3.63948,-0.45930)(3.65962,-0.44068)(3.67974,-0.42207)(3.69988,-0.40345)(3.72001,-0.38483)%
(3.74015,-0.36620)(3.76027,-0.34759)(3.78041,-0.32897)(3.80054,-0.31035)(3.82068,-0.29173)%
(3.84081,-0.27311)(3.86094,-0.25449)(3.88107,-0.23587)(3.90121,-0.21725)(3.92134,-0.19864)%
(3.94147,-0.18001)(3.96160,-0.16139)(3.98174,-0.14277)(3.99932,-0.12651)(4.00187,-0.12416)%
(4.01958,-0.10777)(4.02200,-0.10554)(4.04214,-0.08691)(4.06227,-0.06829)%
\polyline(3.81825,-0.23789)(3.81930,-0.23613)(3.83528,-0.20939)(3.85126,-0.18266)%
(3.86724,-0.15592)(3.88323,-0.12919)(3.89920,-0.10246)(3.91518,-0.07572)(3.93117,-0.04899)%
(3.94714,-0.02226)(3.96313,0.00447)(3.97911,0.03121)(3.99508,0.05795)(4.01107,0.08468)%
(4.02704,0.11142)(4.04303,0.13815)(4.05901,0.16488)(4.07498,0.19161)(4.09097,0.21835)%
(4.10695,0.24508)(4.12293,0.27182)(4.13891,0.29856)(4.15489,0.32529)(4.17087,0.35202)%
(4.18685,0.37875)(4.20283,0.40549)(4.21881,0.43222)(4.23479,0.45895)(4.25077,0.48569)%
(4.26675,0.51243)(4.28273,0.53916)(4.29871,0.56589)(4.31469,0.59263)(4.33068,0.61936)%
(4.34665,0.64609)(4.36263,0.67282)(4.37862,0.69956)(4.39253,0.72284)%
\polyline(3.18677,-0.13766)(3.19373,-0.10679)(3.20069,-0.07591)(3.20765,-0.04504)%
(3.21461,-0.01417)(3.22157,0.01670)(3.22852,0.04758)(3.23549,0.07845)(3.24244,0.10932)%
(3.24940,0.14020)(3.25635,0.17107)(3.26332,0.20194)(3.27028,0.23280)(3.27724,0.26368)%
(3.28419,0.29455)(3.29115,0.32542)(3.29811,0.35629)(3.30507,0.38717)(3.31202,0.41804)%
(3.31899,0.44891)(3.32594,0.47979)(3.33291,0.51066)(3.33986,0.54153)(3.34682,0.57240)%
(3.35378,0.60328)(3.36074,0.63415)(3.36769,0.66502)(3.37465,0.69589)(3.38161,0.72677)%
(3.38857,0.75764)(3.39553,0.78851)(3.40249,0.81938)(3.40945,0.85026)(3.41641,0.88113)%
(3.42337,0.91200)(3.43032,0.94288)(3.43728,0.97375)(3.44424,1.00462)(3.45120,1.03550)%
(3.45816,1.06636)(3.46512,1.09723)(3.47208,1.12810)(3.47904,1.15898)(3.48599,1.18985)%
(3.49295,1.22072)(3.49991,1.25159)(3.50687,1.28247)(3.51382,1.31334)(3.52078,1.34421)%
(3.52774,1.37508)(3.53439,1.40457)%
\polyline(1.82588,0.27538)(1.82588,0.30598)(1.82588,0.33658)(1.82588,0.36719)(1.82588,0.39779)%
(1.82588,0.42839)(1.82588,0.45899)(1.82588,0.48959)(1.82588,0.52020)(1.82588,0.55080)%
(1.82588,0.58140)(1.82588,0.61200)(1.82588,0.64260)(1.82588,0.67321)(1.82588,0.70381)%
(1.82588,0.73441)(1.82588,0.76501)(1.82588,0.79561)(1.82588,0.82622)(1.82588,0.85682)%
(1.82588,0.88742)(1.82588,0.91802)(1.82588,0.94862)(1.82588,0.97923)(1.82588,1.00983)%
(1.82588,1.04043)(1.82588,1.07103)(1.82588,1.10163)(1.82588,1.13224)(1.82588,1.16284)%
(1.82588,1.19344)(1.82588,1.22404)(1.82588,1.25464)(1.82588,1.28525)(1.82588,1.31585)%
(1.82588,1.34645)(1.82588,1.37705)(1.82588,1.40766)(1.82588,1.43826)(1.82588,1.46886)%
(1.82588,1.49946)(1.82588,1.53006)(1.82588,1.56067)(1.82588,1.59127)(1.82588,1.62187)%
(1.82588,1.65247)(1.82588,1.68307)(1.82588,1.71368)(1.82588,1.74428)(1.82588,1.77488)%
(1.82588,1.80535)%
\polyline(-0.20848,0.48941)(-0.20806,0.51654)(-0.20765,0.54369)(-0.20724,0.57083)%
(-0.20683,0.59797)(-0.20642,0.62511)(-0.20601,0.65225)(-0.20560,0.67939)(-0.20519,0.70654)%
(-0.20479,0.73368)(-0.20437,0.76082)(-0.20396,0.78796)(-0.20355,0.81509)(-0.20314,0.84223)%
(-0.20273,0.86937)(-0.20232,0.89652)(-0.20191,0.92366)(-0.20150,0.95080)(-0.20109,0.97794)%
(-0.20068,1.00508)(-0.20026,1.03222)(-0.19985,1.05936)(-0.19944,1.08651)(-0.19903,1.11365)%
(-0.19862,1.14079)(-0.19822,1.16793)(-0.19781,1.19507)(-0.19740,1.22221)(-0.19699,1.24935)%
(-0.19657,1.27649)(-0.19616,1.30364)(-0.19575,1.33078)(-0.19534,1.35792)(-0.19493,1.38506)%
(-0.19452,1.41220)(-0.19411,1.43935)(-0.19370,1.46648)(-0.19329,1.49363)(-0.19287,1.52077)%
(-0.19246,1.54790)(-0.19206,1.57504)(-0.19165,1.60218)(-0.19124,1.62932)(-0.19083,1.65647)%
(-0.19042,1.68361)(-0.19001,1.71075)(-0.18960,1.73789)(-0.18918,1.76503)(-0.18877,1.79217)%
(-0.18836,1.81931)(-0.18797,1.84550)%
\end{picture}}%

%% file: fig/p019.tex
{\unitlength=1cm%
\begin{picture}%
(9.94,5.4)(-5.19,-2.08)%
\linethickness{0.008in}
\polyline(0.70922,-1.14176)(0.70706,-1.14205)(0.69750,-1.14326)(0.69726,-1.14328)%
\polyline(0.69726,-1.14328)(0.68802,-1.14443)(0.67861,-1.14555)(0.67859,-1.14555)%
\polyline(0.67859,-1.14555)(0.66926,-1.14664)(0.65998,-1.14769)(0.65975,-1.14772)%
\polyline(0.49122,-1.16051)(0.48886,-1.16060)(0.48102,-1.16084)(0.47595,-1.16098)%
\polyline(0.47595,-1.16098)(0.47324,-1.16106)(0.46554,-1.16123)(0.46094,-1.16132)%
\polyline(0.46094,-1.16132)(0.45789,-1.16137)(0.45032,-1.16147)(0.44580,-1.16151)%
\polyline(0.44580,-1.16151)(0.44282,-1.16154)(0.43538,-1.16157)(0.42800,-1.16155)%
(0.42070,-1.16151)(0.41345,-1.16143)(0.40628,-1.16132)(0.39918,-1.16116)(0.39214,-1.16097)%
(0.38516,-1.16074)(0.37825,-1.16049)(0.37141,-1.16020)(0.36463,-1.15986)(0.35792,-1.15949)%
\polyline(0.35792,-1.15949)(0.35128,-1.15910)(0.34470,-1.15866)(0.33835,-1.15820)%
\polyline(0.33835,-1.15820)(0.33818,-1.15819)(0.33174,-1.15768)(0.32535,-1.15714)%
(0.31919,-1.15658)%
\polyline(4.17054,0.02824)(4.16000,0.01780)(4.15974,0.01753)%
\polyline(4.15974,0.01753)(4.14921,0.00703)(4.14894,0.00676)%
\polyline(4.14894,0.00676)(4.13841,-0.00382)(4.13814,-0.00409)%
\polyline(4.13814,-0.00409)(4.12761,-0.01472)(4.12734,-0.01499)%
\polyline(4.12734,-0.01499)(4.11681,-0.02569)(4.10601,-0.03673)(4.10575,-0.03700)%
\polyline(4.10575,-0.03700)(4.09523,-0.04783)(4.09496,-0.04810)%
\polyline(4.09496,-0.04810)(4.08444,-0.05901)(4.08412,-0.05934)%
\polyline(4.08412,-0.05934)(4.07364,-0.07025)(4.07354,-0.07035)%
\polyline(4.07354,-0.07035)(4.06980,-0.07426)(4.06202,-0.08243)(4.04996,-0.09515)%
(4.03788,-0.10798)(4.02582,-0.12089)(4.01376,-0.13388)(4.00171,-0.14697)(3.98965,-0.16014)%
(3.97758,-0.17340)(3.96554,-0.18675)(3.95349,-0.20019)(3.94145,-0.21373)(3.92939,-0.22736)%
(3.91735,-0.24109)(3.90533,-0.25490)(3.89867,-0.26258)(3.89287,-0.26929)(3.87993,-0.28437)%
(3.86698,-0.29955)(3.85405,-0.31485)(3.84112,-0.33025)(3.82821,-0.34576)(3.81531,-0.36138)%
(3.80242,-0.37711)(3.78954,-0.39295)(3.77667,-0.40890)(3.76381,-0.42496)(3.75097,-0.44113)%
(3.73814,-0.45741)(3.72533,-0.47381)(3.71189,-0.49115)(3.69845,-0.50862)(3.68504,-0.52623)%
(3.67164,-0.54394)(3.65828,-0.56178)(3.64494,-0.57974)(3.63162,-0.59782)(3.61834,-0.61602)%
(3.60507,-0.63435)(3.59316,-0.65095)%
\polyline(3.59316,-0.65095)(3.59184,-0.65278)(3.58062,-0.66856)%
\polyline(3.58062,-0.66856)(3.57864,-0.67134)(3.56547,-0.69002)(3.55256,-0.70851)%
\polyline(-2.19106,-1.24852)(-1.77931,-1.26697)(-1.35143,-1.27343)(-0.91172,-1.26758)%
(-0.46456,-1.24921)(-0.01445,-1.21827)(0.43415,-1.17484)%
\polyline(0.43415,-1.17484)(0.60744,-1.15302)%
\polyline(0.60744,-1.15302)(0.65063,-1.14759)%
\polyline(0.65063,-1.14759)(0.66908,-1.14527)%
\polyline(0.66908,-1.14527)(0.68783,-1.14291)%
\polyline(0.68783,-1.14291)(0.70494,-1.14075)%
\polyline(0.70494,-1.14075)(0.71649,-1.13930)%
\polyline(0.71649,-1.13930)(0.87678,-1.11912)(1.30907,-1.05149)(1.72679,-0.97241)%
(2.12590,-0.88248)(2.50258,-0.78244)(2.85327,-0.67309)(3.17476,-0.55540)(3.46413,-0.43036)%
(3.52997,-0.39642)%
\polyline(3.52997,-0.39642)(3.71885,-0.29906)(3.93681,-0.16267)(4.07854,-0.05185)%
\polyline(4.07854,-0.05185)(4.10277,-0.03291)%
\polyline(4.10277,-0.03291)(4.11626,-0.02237)(4.12570,-0.01269)%
\polyline(4.12570,-0.01269)(4.14705,0.00916)%
\polyline(4.14705,0.00916)(4.16833,0.03093)%
\polyline(4.16833,0.03093)(4.25592,0.12059)(4.35492,0.26495)(4.41285,0.40945)(4.42971,0.55283)%
(4.40591,0.69384)(4.34234,0.83129)(4.24018,0.96402)(4.10107,1.09092)(3.92695,1.21097)%
(3.72008,1.32320)(3.48300,1.42675)(3.21851,1.52085)(2.92961,1.60479)(2.61945,1.67803)%
(2.29135,1.74007)(1.94869,1.79057)(1.59495,1.82928)(1.23356,1.85602)(0.86800,1.87076)%
(0.50163,1.87357)(0.13778,1.86461)(-0.22037,1.84411)(-0.56976,1.81242)(-0.90752,1.76996)%
(-1.23096,1.71721)(-1.53758,1.65474)(-1.82512,1.58316)(-2.09155,1.50312)(-2.33505,1.41535)%
(-2.55410,1.32059)(-2.74735,1.21959)(-2.91375,1.11317)(-3.05249,1.00214)%
\polyline(-4.05668,0.07040)(-4.03660,0.08903)(-4.01651,0.10767)(-3.99643,0.12630)%
(-3.97634,0.14494)(-3.95625,0.16357)(-3.93617,0.18221)(-3.91608,0.20084)(-3.89600,0.21947)%
(-3.87592,0.23811)(-3.85583,0.25674)(-3.83575,0.27538)(-3.81567,0.29401)(-3.79559,0.31265)%
(-3.77550,0.33128)(-3.75542,0.34992)(-3.73533,0.36855)(-3.71525,0.38718)(-3.69517,0.40582)%
(-3.67508,0.42445)(-3.65500,0.44309)(-3.63492,0.46172)(-3.61483,0.48036)(-3.59475,0.49899)%
(-3.57467,0.51764)(-3.55458,0.53627)(-3.53450,0.55490)(-3.51442,0.57354)(-3.49433,0.59217)%
(-3.47425,0.61081)(-3.45416,0.62944)(-3.43408,0.64808)(-3.41400,0.66671)(-3.39392,0.68535)%
(-3.37383,0.70398)(-3.35375,0.72261)(-3.33367,0.74125)(-3.31357,0.75988)(-3.29349,0.77852)%
(-3.27341,0.79715)(-3.25332,0.81579)(-3.23324,0.83442)(-3.21315,0.85306)(-3.19307,0.87169)%
(-3.17299,0.89032)(-3.15290,0.90896)(-3.13282,0.92759)(-3.11274,0.94623)(-3.09266,0.96486)%
(-3.07257,0.98350)(-3.05249,1.00214)%
\polyline(-1.46071,-0.83234)(-1.18675,-0.90509)(-0.90262,-0.97055)(-0.61034,-1.02843)%
(-0.31200,-1.07848)(-0.00975,-1.12045)(0.29425,-1.15418)(0.32571,-1.15681)%
\polyline(0.32571,-1.15681)(0.33853,-1.15788)%
\polyline(0.33853,-1.15788)(0.35769,-1.15948)%
\polyline(0.47500,-1.16929)(0.59775,-1.17954)(0.89852,-1.19644)(1.19428,-1.20485)%
(1.48275,-1.20480)(1.76165,-1.19636)(2.02872,-1.17967)(2.28172,-1.15491)(2.51845,-1.12232)%
(2.73677,-1.08223)(2.93463,-1.03499)(3.11005,-0.98102)(3.26120,-0.92083)(3.38634,-0.85493)%
(3.48392,-0.78395)(3.55254,-0.70852)%
\polyline(3.55254,-0.70852)(3.56254,-0.68798)%
\polyline(3.52997,-0.39642)(3.51716,-0.37707)(3.42783,-0.29080)(3.30602,-0.20487)%
(3.15212,-0.12017)(2.96681,-0.03761)(2.75107,0.04195)(2.50620,0.11762)(2.23383,0.18856)%
(1.93588,0.25394)(1.61459,0.31299)(1.27247,0.36501)(0.91232,0.40932)(0.53717,0.44536)%
(0.15029,0.47261)(-0.24483,0.49067)(-0.64459,0.49921)(-1.04521,0.49803)(-1.44284,0.48703)%
(-1.83358,0.46621)(-2.21350,0.43569)(-2.57873,0.39572)(-2.92547,0.34664)(-3.25003,0.28894)%
(-3.54891,0.22320)(-3.81881,0.15009)(-4.05668,0.07040)%
\polyline(-1.46071,-0.83234)(-1.47531,-0.84067)(-1.48992,-0.84899)(-1.50453,-0.85731)%
(-1.51913,-0.86564)(-1.53374,-0.87396)(-1.54835,-0.88229)(-1.56296,-0.89061)(-1.57756,-0.89893)%
(-1.59217,-0.90726)(-1.60678,-0.91558)(-1.62138,-0.92390)(-1.63599,-0.93223)(-1.65060,-0.94055)%
(-1.66521,-0.94887)(-1.67981,-0.95720)(-1.69442,-0.96552)(-1.70903,-0.97384)(-1.72363,-0.98217)%
(-1.73824,-0.99049)(-1.75285,-0.99881)(-1.76745,-1.00714)(-1.78206,-1.01546)(-1.79667,-1.02378)%
(-1.81128,-1.03211)(-1.82588,-1.04043)(-1.84049,-1.04875)(-1.85510,-1.05708)(-1.86970,-1.06540)%
(-1.88431,-1.07372)(-1.89892,-1.08205)(-1.91353,-1.09037)(-1.92813,-1.09869)(-1.94274,-1.10702)%
(-1.95735,-1.11534)(-1.97195,-1.12366)(-1.98656,-1.13199)(-2.00117,-1.14031)(-2.01577,-1.14864)%
(-2.03038,-1.15696)(-2.04499,-1.16528)(-2.05960,-1.17361)(-2.07420,-1.18193)(-2.08881,-1.19025)%
(-2.10342,-1.19858)(-2.11802,-1.20690)(-2.13263,-1.21522)(-2.14724,-1.22355)(-2.16185,-1.23187)%
(-2.17645,-1.24019)(-2.18376,-1.24435)%
\polyline(-1.58243,-0.90171)(-1.28552,-0.96540)(-0.97742,-1.02103)(-0.66057,-1.06829)%
(-0.33743,-1.10694)(-0.01053,-1.13676)(0.31664,-1.15756)(0.31756,-1.15762)(0.34330,-1.15856)%
\polyline(0.36465,-1.15933)(0.37519,-1.15972)(0.38553,-1.16008)(0.39954,-1.16060)%
(0.42552,-1.16154)%
\polyline(0.51423,-1.16476)(0.64426,-1.16947)(0.96694,-1.17228)(1.28303,-1.16611)%
(1.58994,-1.15109)(1.88514,-1.12737)(2.16615,-1.09525)(2.43056,-1.05498)(2.67606,-1.00700)%
(2.90045,-0.95171)(3.10166,-0.88960)(3.27775,-0.82125)(3.42698,-0.74726)(3.54777,-0.66829)%
(3.61204,-0.60947)(3.63874,-0.58505)(3.64330,-0.57846)%
\polyline(3.71563,-0.30072)(3.68500,-0.22501)(3.61447,-0.13240)(3.51101,-0.04050)%
(3.37503,0.04980)(3.20727,0.13764)(3.00876,0.22213)(2.78083,0.30242)(2.52508,0.37769)%
(2.24342,0.44714)(1.93802,0.51004)(1.61132,0.56571)(1.26598,0.61352)(0.90493,0.65289)%
(0.53124,0.68340)(0.14821,0.70461)(-0.24076,0.71624)(-0.63212,0.71807)(-1.02227,0.71002)%
(-1.40753,0.69206)(-1.78425,0.66429)(-2.14877,0.62693)(-2.49754,0.58028)(-2.82706,0.52476)%
(-3.13405,0.46088)(-3.41532,0.38926)(-3.66796,0.31060)(-3.88930,0.22569)%
\polyline(-1.70416,-0.97107)(-1.38427,-1.02572)(-1.05222,-1.07151)(-0.71080,-1.10815)%
(-0.36286,-1.13539)(-0.01132,-1.15306)(0.34084,-1.16107)(0.38098,-1.16087)(0.40035,-1.16079)%
(0.41345,-1.16072)(0.43057,-1.16064)(0.44296,-1.16058)(0.45803,-1.16051)(0.47337,-1.16044)%
\polyline(0.55149,-1.16006)(0.69076,-1.15941)(1.03537,-1.14812)(1.37178,-1.12737)%
(1.69713,-1.09736)(2.00863,-1.05839)(2.30357,-1.01081)(2.57940,-0.95507)(2.83368,-0.89167)%
(3.06414,-0.82117)(3.26869,-0.74421)(3.44546,-0.66147)(3.59278,-0.57368)(3.70920,-0.48164)%
(3.72042,-0.46893)(3.73401,-0.45357)(3.76126,-0.42271)(3.77554,-0.40654)(3.78905,-0.39125)%
(3.79356,-0.38614)(3.79976,-0.37433)%
\polyline(3.85851,-0.21166)(3.86266,-0.18828)(3.84639,-0.08769)(3.79603,0.01280)(3.71179,0.11226)%
(3.59420,0.20979)(3.44404,0.30448)(3.26241,0.39547)(3.05071,0.48187)(2.81058,0.56290)%
(2.54395,0.63776)(2.25300,0.70573)(1.94015,0.76615)(1.60803,0.81842)(1.25950,0.86201)%
(0.89754,0.89647)(0.52532,0.92143)(0.14612,0.93661)(-0.23668,0.94181)(-0.61965,0.93695)%
(-0.99932,0.92201)(-1.37221,0.89709)(-1.73491,0.86239)(-2.08404,0.81817)(-2.41634,0.76485)%
(-2.72867,0.70288)(-3.01805,0.63282)(-3.28172,0.55533)(-3.51712,0.47112)(-3.72194,0.38098)%
\polyline(-1.82588,-1.04043)(-1.48303,-1.08603)(-1.12703,-1.12199)(-0.76103,-1.14801)%
(-0.38828,-1.16385)(-0.01210,-1.16936)(0.36414,-1.16451)(0.40459,-1.16287)(0.41491,-1.16245)%
(0.43666,-1.16156)(0.46090,-1.16058)(0.47590,-1.15996)(0.49093,-1.15935)(0.58875,-1.15537)%
(0.73726,-1.14933)(1.10379,-1.12396)(1.46053,-1.08863)(1.80432,-1.04364)(2.13212,-0.98940)%
(2.44100,-0.92638)(2.72824,-0.85515)(2.99129,-0.77634)(3.22781,-0.69065)(3.43572,-0.59882)%
(3.61316,-0.50169)(3.75856,-0.40012)(3.82226,-0.34037)(3.85185,-0.31261)(3.87063,-0.29499)%
(3.87888,-0.28356)%
\polyline(3.96730,-0.13882)(3.99113,-0.07785)(3.99848,0.03225)(3.97038,0.14205)(3.90706,0.25060)%
(3.80911,0.35693)(3.67738,0.46008)(3.51305,0.55917)(3.31756,0.65329)(3.09266,0.74162)%
(2.84034,0.82337)(2.56282,0.89782)(2.26259,0.96432)(1.94229,1.02226)(1.60476,1.07113)%
(1.25301,1.11051)(0.89016,1.14004)(0.51939,1.15947)(0.14404,1.16861)(-0.23260,1.16739)%
(-0.60717,1.15582)(-0.97636,1.13400)(-1.33690,1.10212)(-1.68558,1.06047)(-2.01932,1.00942)%
(-2.33514,0.94942)(-2.63026,0.88100)(-2.90206,0.80477)(-3.14813,0.72140)(-3.36628,0.63163)%
(-3.53091,0.54826)(-3.55458,0.53627)%
\polyline(-1.94761,-1.10979)(-1.58179,-1.14634)(-1.20183,-1.17247)(-0.81126,-1.18787)%
(-0.41371,-1.19230)(-0.01288,-1.18567)(0.38752,-1.16795)(0.43472,-1.16453)(0.48922,-1.16059)%
(0.59958,-1.15260)(0.62606,-1.15068)(0.78376,-1.13926)(1.17221,-1.09980)(1.54928,-1.04989)%
(1.91151,-0.98992)(2.25560,-0.92041)(2.57843,-0.84196)(2.87707,-0.75523)(3.14890,-0.66101)%
(3.39149,-0.56012)(3.60275,-0.45344)(3.78085,-0.34192)(3.89866,-0.24720)(3.92434,-0.22655)%
(3.92645,-0.22424)(3.93859,-0.21092)%
\polyline(4.05434,-0.07077)(4.05975,-0.06164)(4.10321,0.01165)(4.13732,0.13238)(4.13429,0.25278)%
(4.09436,0.37180)(4.01811,0.48840)(3.90643,0.60159)(3.76058,0.71037)(3.58206,0.81385)%
(3.37271,0.91111)(3.13461,1.00136)(2.87010,1.08384)(2.58169,1.15789)(2.27217,1.22291)%
(1.94442,1.27836)(1.60149,1.32385)(1.24652,1.35901)(0.88277,1.38361)(0.51348,1.39751)%
(0.14195,1.40061)(-0.22852,1.39296)(-0.59471,1.37468)(-0.95342,1.34598)(-1.30159,1.30715)%
(-1.63625,1.25856)(-1.95458,1.20067)(-2.25394,1.13399)(-2.53186,1.05911)(-2.78608,0.97671)%
(-3.01454,0.88746)(-3.21544,0.79214)(-3.34857,0.71418)(-3.38721,0.69156)%
\polyline(-2.06933,-1.17915)(-1.68055,-1.20666)(-1.27663,-1.22295)(-0.86149,-1.22772)%
(-0.43914,-1.22075)(-0.01366,-1.20197)(0.41083,-1.17139)(0.45649,-1.16680)(0.60305,-1.15205)%
(0.66336,-1.14599)(0.68415,-1.14389)%
\polyline(0.69733,-1.14257)(0.83027,-1.12920)(1.24064,-1.07565)(1.63803,-1.01115)%
(2.01871,-0.93619)(2.37909,-0.85142)(2.71585,-0.75752)(3.02592,-0.65532)(3.30652,-0.54569)%
(3.55517,-0.42959)(3.76978,-0.30805)(3.94856,-0.18214)(3.98220,-0.15146)(4.00821,-0.12772)%
(4.03422,-0.10399)(4.06011,-0.08036)(4.07242,-0.06914)(4.08351,-0.05901)%
\polyline(4.09506,-0.04672)(4.10447,-0.03477)%
\polyline(4.11411,-0.02252)(4.12290,-0.01136)%
\polyline(4.13301,0.00149)(4.14237,0.01338)%
\polyline(4.14505,0.01678)(4.15166,0.02517)(4.19349,0.07831)(4.25803,0.21055)(4.28351,0.34260)%
(4.27011,0.47331)(4.21835,0.60155)(4.12914,0.72621)(4.00375,0.84625)(3.84376,0.96067)%
(3.65107,1.06853)(3.42786,1.16893)(3.17656,1.26110)(2.89984,1.34432)(2.60057,1.41796)%
(2.28176,1.48149)(1.94656,1.53447)(1.59822,1.57655)(1.24005,1.60751)(0.87538,1.62719)%
(0.50756,1.63554)(0.13986,1.63260)(-0.22445,1.61854)(-0.58223,1.59356)(-0.93047,1.55797)%
(-1.26627,1.51218)(-1.58691,1.45665)(-1.88985,1.39191)(-2.17275,1.31855)(-2.43346,1.23723)%
(-2.67008,1.14865)(-2.88094,1.05353)(-3.06460,0.95266)(-3.18890,0.86794)(-3.21985,0.84684)%
\polyline(0.15992,-1.14032)(0.16146,-1.14142)(0.16299,-1.14252)(0.16452,-1.14363)%
(0.16605,-1.14473)(0.16758,-1.14583)(0.16912,-1.14694)(0.17064,-1.14805)(0.17218,-1.14915)%
(0.17371,-1.15026)(0.17524,-1.15136)(0.17677,-1.15246)(0.17831,-1.15357)(0.17983,-1.15468)%
(0.18137,-1.15578)(0.18290,-1.15688)(0.18444,-1.15799)(0.18596,-1.15909)(0.18750,-1.16020)%
(0.18903,-1.16130)(0.19056,-1.16241)(0.19209,-1.16351)(0.19362,-1.16462)(0.19515,-1.16572)%
(0.19669,-1.16682)(0.19822,-1.16793)(0.19975,-1.16904)(0.20128,-1.17014)(0.20281,-1.17124)%
(0.20434,-1.17235)(0.20587,-1.17345)(0.20741,-1.17456)(0.20894,-1.17566)(0.21047,-1.17677)%
(0.21200,-1.17787)(0.21353,-1.17898)(0.21506,-1.18008)(0.21660,-1.18118)(0.21813,-1.18229)%
(0.21966,-1.18340)(0.22119,-1.18450)(0.22272,-1.18560)(0.22425,-1.18671)(0.22579,-1.18781)%
(0.22732,-1.18892)(0.22885,-1.19003)(0.23038,-1.19112)(0.23192,-1.19223)(0.23344,-1.19334)%
(0.23498,-1.19445)(0.23651,-1.19554)%
\polyline(0.16299,-1.14252)(0.16452,-1.14363)(0.16605,-1.14473)(0.16758,-1.14583)%
(0.16912,-1.14694)(0.17064,-1.14805)(0.17157,-1.14872)%
\polyline(0.17929,-1.15427)(0.17983,-1.15468)(0.18137,-1.15578)(0.18290,-1.15688)%
(0.18444,-1.15799)(0.18596,-1.15909)(0.18750,-1.16020)(0.18903,-1.16130)(0.19056,-1.16241)%
(0.19209,-1.16351)(0.19362,-1.16462)(0.19515,-1.16572)(0.19669,-1.16682)(0.19822,-1.16793)%
(0.19975,-1.16904)(0.20128,-1.17014)(0.20281,-1.17124)(0.20434,-1.17235)(0.20587,-1.17345)%
(0.20741,-1.17456)(0.20894,-1.17566)(0.21047,-1.17677)(0.21200,-1.17787)(0.21353,-1.17898)%
(0.21506,-1.18008)(0.21660,-1.18118)(0.21813,-1.18229)(0.21966,-1.18340)(0.22119,-1.18450)%
(0.22272,-1.18560)(0.22425,-1.18671)(0.22579,-1.18781)(0.22732,-1.18892)(0.22885,-1.19003)%
(0.23038,-1.19112)(0.23192,-1.19223)(0.23344,-1.19334)(0.23459,-1.19416)%
\polyline(1.79442,-1.19431)(1.80852,-1.18634)(1.80992,-1.18555)(1.82354,-1.17785)%
(1.83857,-1.16936)(1.85360,-1.16086)(1.86863,-1.15237)(1.88366,-1.14387)(1.89869,-1.13538)%
(1.91371,-1.12688)(1.92874,-1.11839)(1.94377,-1.10990)(1.95880,-1.10140)(1.97383,-1.09291)%
(1.98886,-1.08441)(2.00389,-1.07592)(2.01892,-1.06742)(2.03394,-1.05893)(2.04897,-1.05044)%
(2.06400,-1.04194)(2.07903,-1.03345)(2.09406,-1.02495)(2.10909,-1.01646)(2.12412,-1.00797)%
(2.13914,-0.99947)(2.15417,-0.99098)(2.16920,-0.98248)(2.18423,-0.97399)(2.19926,-0.96549)%
(2.21429,-0.95700)(2.22932,-0.94851)(2.24435,-0.94001)(2.25937,-0.93152)(2.27440,-0.92302)%
(2.28943,-0.91453)(2.30446,-0.90604)(2.31949,-0.89754)(2.33452,-0.88905)(2.34954,-0.88055)%
(2.36458,-0.87206)(2.37960,-0.86357)(2.39463,-0.85507)(2.40966,-0.84658)(2.42469,-0.83808)%
(2.43972,-0.82959)(2.45475,-0.82109)(2.46977,-0.81260)(2.48481,-0.80410)(2.49983,-0.79561)%
(2.51486,-0.78712)(2.52989,-0.77862)(2.54492,-0.77013)%
\polyline(3.05808,-0.99701)(3.07577,-0.98065)(3.09591,-0.96203)(3.11604,-0.94341)%
(3.13617,-0.92479)(3.15630,-0.90617)(3.17644,-0.88755)(3.19657,-0.86893)(3.21670,-0.85031)%
(3.23683,-0.83169)(3.25697,-0.81307)(3.27710,-0.79445)(3.29723,-0.77584)(3.31736,-0.75722)%
(3.33750,-0.73859)(3.35763,-0.71997)(3.37776,-0.70136)(3.39789,-0.68274)(3.41803,-0.66412)%
(3.43816,-0.64549)(3.45829,-0.62688)(3.47842,-0.60826)(3.49856,-0.58964)(3.51869,-0.57102)%
(3.53882,-0.55240)(3.55895,-0.53378)(3.57909,-0.51516)(3.59921,-0.49655)(3.61935,-0.47793)%
(3.63948,-0.45930)(3.65962,-0.44068)(3.67974,-0.42207)(3.69988,-0.40345)(3.72001,-0.38483)%
(3.74015,-0.36620)(3.76027,-0.34759)(3.78041,-0.32897)(3.80054,-0.31035)(3.82068,-0.29173)%
(3.84081,-0.27311)(3.86094,-0.25449)(3.88107,-0.23587)(3.90121,-0.21725)(3.92134,-0.19864)%
(3.94147,-0.18001)(3.96160,-0.16139)(3.98174,-0.14277)(4.00086,-0.12509)(4.00187,-0.12416)%
(4.02062,-0.10681)(4.02200,-0.10554)(4.04068,-0.08826)(4.04214,-0.08691)(4.06074,-0.06971)%
\polyline(3.81923,-0.23624)(3.81930,-0.23613)(3.83528,-0.20939)(3.85126,-0.18266)%
(3.86724,-0.15592)(3.88323,-0.12919)(3.89920,-0.10246)(3.91518,-0.07572)(3.93117,-0.04899)%
(3.94714,-0.02226)(3.96313,0.00447)(3.97911,0.03121)(3.99508,0.05795)(4.01107,0.08468)%
(4.02704,0.11142)(4.04303,0.13815)(4.05901,0.16488)(4.07498,0.19161)(4.09097,0.21835)%
(4.10695,0.24508)(4.12293,0.27182)(4.13891,0.29856)(4.15489,0.32529)(4.17087,0.35202)%
(4.18685,0.37875)(4.20283,0.40549)(4.21881,0.43222)(4.23479,0.45895)(4.25077,0.48569)%
(4.26675,0.51243)(4.28273,0.53916)(4.29871,0.56589)(4.31469,0.59263)(4.33068,0.61936)%
(4.34665,0.64609)(4.36263,0.67282)(4.37862,0.69956)(4.39252,0.72281)%
\polyline(3.18677,-0.13766)(3.19373,-0.10679)(3.20069,-0.07591)(3.20765,-0.04504)%
(3.21461,-0.01417)(3.22157,0.01670)(3.22852,0.04758)(3.23549,0.07845)(3.24244,0.10932)%
(3.24940,0.14020)(3.25635,0.17107)(3.26332,0.20194)(3.27028,0.23280)(3.27724,0.26368)%
(3.28419,0.29455)(3.29115,0.32542)(3.29811,0.35629)(3.30507,0.38717)(3.31202,0.41804)%
(3.31899,0.44891)(3.32594,0.47979)(3.33291,0.51066)(3.33986,0.54153)(3.34682,0.57240)%
(3.35378,0.60328)(3.36074,0.63415)(3.36769,0.66502)(3.37465,0.69589)(3.38161,0.72677)%
(3.38857,0.75764)(3.39553,0.78851)(3.40249,0.81938)(3.40945,0.85026)(3.41641,0.88113)%
(3.42337,0.91200)(3.43032,0.94288)(3.43728,0.97375)(3.44424,1.00462)(3.45120,1.03550)%
(3.45816,1.06636)(3.46512,1.09723)(3.47208,1.12810)(3.47904,1.15898)(3.48599,1.18985)%
(3.49295,1.22072)(3.49991,1.25159)(3.50687,1.28247)(3.51382,1.31334)(3.52078,1.34421)%
(3.52774,1.37508)(3.53434,1.40433)%
\polyline(1.82588,0.27538)(1.82588,0.30598)(1.82588,0.33658)(1.82588,0.36719)(1.82588,0.39779)%
(1.82588,0.42839)(1.82588,0.45899)(1.82588,0.48959)(1.82588,0.52020)(1.82588,0.55080)%
(1.82588,0.58140)(1.82588,0.61200)(1.82588,0.64260)(1.82588,0.67321)(1.82588,0.70381)%
(1.82588,0.73441)(1.82588,0.76501)(1.82588,0.79561)(1.82588,0.82622)(1.82588,0.85682)%
(1.82588,0.88742)(1.82588,0.91802)(1.82588,0.94862)(1.82588,0.97923)(1.82588,1.00983)%
(1.82588,1.04043)(1.82588,1.07103)(1.82588,1.10163)(1.82588,1.13224)(1.82588,1.16284)%
(1.82588,1.19344)(1.82588,1.22404)(1.82588,1.25464)(1.82588,1.28525)(1.82588,1.31585)%
(1.82588,1.34645)(1.82588,1.37705)(1.82588,1.40766)(1.82588,1.43826)(1.82588,1.46886)%
(1.82588,1.49946)(1.82588,1.53006)(1.82588,1.56067)(1.82588,1.59127)(1.82588,1.62187)%
(1.82588,1.65247)(1.82588,1.68307)(1.82588,1.71368)(1.82588,1.74428)(1.82588,1.77488)%
(1.82588,1.80401)%
\polyline(-0.20848,0.48941)(-0.20806,0.51654)(-0.20765,0.54369)(-0.20724,0.57083)%
(-0.20683,0.59797)(-0.20642,0.62511)(-0.20601,0.65225)(-0.20560,0.67939)(-0.20519,0.70654)%
(-0.20479,0.73368)(-0.20437,0.76082)(-0.20396,0.78796)(-0.20355,0.81509)(-0.20314,0.84223)%
(-0.20273,0.86937)(-0.20232,0.89652)(-0.20191,0.92366)(-0.20150,0.95080)(-0.20109,0.97794)%
(-0.20068,1.00508)(-0.20026,1.03222)(-0.19985,1.05936)(-0.19944,1.08651)(-0.19903,1.11365)%
(-0.19862,1.14079)(-0.19822,1.16793)(-0.19781,1.19507)(-0.19740,1.22221)(-0.19699,1.24935)%
(-0.19657,1.27649)(-0.19616,1.30364)(-0.19575,1.33078)(-0.19534,1.35792)(-0.19493,1.38506)%
(-0.19452,1.41220)(-0.19411,1.43935)(-0.19370,1.46648)(-0.19329,1.49363)(-0.19287,1.52077)%
(-0.19246,1.54790)(-0.19206,1.57504)(-0.19165,1.60218)(-0.19124,1.62932)(-0.19083,1.65647)%
(-0.19042,1.68361)(-0.19001,1.71075)(-0.18960,1.73789)(-0.18918,1.76503)(-0.18877,1.79217)%
(-0.18836,1.81931)(-0.18796,1.84597)%
\polyline(-2.38612,0.41818)(-2.37744,0.44075)(-2.36877,0.46333)(-2.36009,0.48590)%
(-2.35142,0.50847)(-2.34274,0.53104)(-2.33406,0.55361)(-2.32539,0.57618)(-2.31671,0.59875)%
(-2.30803,0.62132)(-2.29936,0.64390)(-2.29068,0.66647)(-2.28200,0.68904)(-2.27333,0.71162)%
(-2.26465,0.73419)(-2.25598,0.75676)(-2.24730,0.77933)(-2.23862,0.80191)(-2.22994,0.82448)%
(-2.22126,0.84705)(-2.21258,0.86963)(-2.20391,0.89220)(-2.19523,0.91477)(-2.18655,0.93733)%
(-2.17788,0.95991)(-2.16920,0.98248)(-2.16052,1.00505)(-2.15185,1.02763)(-2.14317,1.05020)%
(-2.13450,1.07277)(-2.12582,1.09534)(-2.11714,1.11792)(-2.10847,1.14049)(-2.09979,1.16306)%
(-2.09111,1.18564)(-2.08244,1.20821)(-2.07376,1.23078)(-2.06508,1.25335)(-2.05641,1.27593)%
(-2.04773,1.29849)(-2.03905,1.32107)(-2.03037,1.34364)(-2.02169,1.36621)(-2.01302,1.38878)%
(-2.00434,1.41135)(-1.99566,1.43393)(-1.98699,1.45650)(-1.97831,1.47907)(-1.96963,1.50164)%
(-1.96096,1.52422)(-1.95308,1.54472)%
\polyline(-4.06227,0.06829)(-4.04279,0.08631)(-4.04214,0.08691)(-4.02200,0.10554)%
(-4.00188,0.12414)(-4.00187,0.12416)(-3.98176,0.14275)(-3.98174,0.14277)(-3.96164,0.16136)%
(-3.96160,0.16139)(-3.94151,0.17997)(-3.94147,0.18001)(-3.92140,0.19858)(-3.92134,0.19864)%
(-3.90128,0.21718)(-3.90121,0.21725)(-3.88116,0.23579)(-3.88107,0.23587)(-3.86104,0.25440)%
(-3.86094,0.25449)(-3.84092,0.27301)(-3.84081,0.27311)(-3.82080,0.29161)(-3.82068,0.29173)%
(-3.80068,0.31022)(-3.80054,0.31035)(-3.78056,0.32883)(-3.78041,0.32897)(-3.76044,0.34744)%
(-3.76027,0.34759)(-3.74032,0.36605)(-3.74015,0.36620)(-3.72020,0.38465)(-3.72001,0.38483)%
(-3.70008,0.40326)(-3.69988,0.40345)(-3.67996,0.42187)(-3.67974,0.42207)(-3.65984,0.44048)%
(-3.65962,0.44068)(-3.63972,0.45908)(-3.63948,0.45930)(-3.61960,0.47769)(-3.61935,0.47793)%
(-3.59948,0.49630)(-3.59921,0.49655)(-3.57936,0.51491)(-3.57909,0.51516)(-3.55924,0.53352)%
(-3.55895,0.53378)(-3.53912,0.55212)(-3.53882,0.55240)(-3.51900,0.57074)(-3.51869,0.57102)%
(-3.49888,0.58934)(-3.49856,0.58964)(-3.47876,0.60795)(-3.47842,0.60826)(-3.45864,0.62656)%
(-3.45829,0.62688)(-3.43852,0.64517)(-3.43816,0.64549)(-3.41840,0.66377)(-3.41803,0.66412)%
(-3.39828,0.68238)(-3.39789,0.68274)(-3.37816,0.70099)(-3.37776,0.70136)(-3.35804,0.71959)%
(-3.35763,0.71997)(-3.33792,0.73821)(-3.33750,0.73859)(-3.31780,0.75681)(-3.31736,0.75722)%
(-3.29767,0.77543)(-3.29723,0.77584)(-3.27756,0.79403)(-3.27710,0.79445)(-3.25744,0.81264)%
(-3.25697,0.81307)(-3.23731,0.83125)(-3.23683,0.83169)(-3.21720,0.84985)(-3.21670,0.85031)%
(-3.19708,0.86846)(-3.19657,0.86893)(-3.17696,0.88706)(-3.17644,0.88755)(-3.15684,0.90568)%
(-3.15630,0.90617)(-3.13672,0.92428)(-3.13617,0.92479)(-3.11659,0.94290)(-3.11604,0.94341)%
(-3.09648,0.96150)(-3.09591,0.96203)(-3.07636,0.98011)(-3.07577,0.98065)(-3.05623,0.99872)%
\end{picture}}%

%% file: fig/p022.tex
{\unitlength=1cm%
\begin{picture}%
(9.94,5.4)(-5.19,-2.08)%
\linethickness{0.008in}
\polyline(0.70919,-1.14177)(0.70847,-1.14187)(0.69894,-1.14308)(0.69870,-1.14311)%
\polyline(0.69870,-1.14311)(0.68947,-1.14425)(0.68007,-1.14539)(0.67983,-1.14541)%
\polyline(0.67983,-1.14541)(0.67074,-1.14648)(0.66148,-1.14754)(0.66125,-1.14756)%
\polyline(0.47519,-1.16101)(0.47482,-1.16102)(0.46711,-1.16120)(0.45947,-1.16134)%
(0.45188,-1.16145)(0.44437,-1.16153)(0.43693,-1.16156)(0.42954,-1.16156)(0.42223,-1.16153)%
(0.41498,-1.16146)(0.40780,-1.16135)(0.40068,-1.16120)(0.39363,-1.16102)(0.38664,-1.16080)%
(0.37973,-1.16054)(0.37327,-1.16028)%
\polyline(0.37327,-1.16028)(0.37287,-1.16026)(0.36608,-1.15994)(0.35936,-1.15958)%
(0.35919,-1.15957)%
\polyline(0.35919,-1.15957)(0.35270,-1.15919)(0.34611,-1.15876)(0.34595,-1.15874)%
\polyline(0.34595,-1.15874)(0.33958,-1.15829)(0.33312,-1.15780)(0.32689,-1.15727)%
\polyline(0.71808,-1.14062)(0.70847,-1.14187)(0.69894,-1.14308)(0.68947,-1.14425)%
(0.68007,-1.14539)(0.67074,-1.14648)(0.66148,-1.14754)(0.65228,-1.14855)(0.64315,-1.14952)%
(0.63409,-1.15046)(0.62510,-1.15136)(0.61618,-1.15222)(0.60733,-1.15304)(0.59854,-1.15382)%
(0.58982,-1.15457)(0.58116,-1.15527)(0.57259,-1.15594)(0.56407,-1.15657)(0.55562,-1.15716)%
(0.54724,-1.15772)(0.53892,-1.15824)(0.53068,-1.15871)(0.52250,-1.15916)(0.51438,-1.15956)%
(0.50634,-1.15992)(0.49835,-1.16026)(0.49045,-1.16054)(0.48260,-1.16081)(0.47482,-1.16102)%
(0.46711,-1.16120)(0.45947,-1.16134)(0.45188,-1.16145)(0.44437,-1.16153)(0.43693,-1.16156)%
(0.42954,-1.16156)(0.42223,-1.16153)(0.41498,-1.16146)(0.40780,-1.16135)(0.40068,-1.16120)%
(0.39363,-1.16102)(0.38664,-1.16080)(0.37973,-1.16054)(0.37287,-1.16026)(0.36608,-1.15994)%
(0.35936,-1.15958)(0.35270,-1.15919)(0.34611,-1.15876)(0.33958,-1.15829)(0.33312,-1.15780)%
(0.32672,-1.15727)(0.32039,-1.15669)%
\polyline(4.17137,0.02905)(4.16176,0.01953)(4.16152,0.01929)%
\polyline(4.16152,0.01929)(4.15189,0.00970)(4.15154,0.00934)%
\polyline(4.15154,0.00934)(4.14052,-0.00169)(4.14024,-0.00198)%
\polyline(4.14024,-0.00198)(4.12891,-0.01341)(4.12864,-0.01369)%
\polyline(4.12864,-0.01369)(4.11732,-0.02518)(4.11695,-0.02555)%
\polyline(4.11695,-0.02555)(4.10571,-0.03704)(4.10543,-0.03733)%
\polyline(4.10543,-0.03733)(4.09412,-0.04898)(4.09382,-0.04927)%
\polyline(4.09382,-0.04927)(4.08251,-0.06099)(4.08223,-0.06128)%
\polyline(4.08223,-0.06128)(4.07092,-0.07309)(4.05932,-0.08525)(4.04772,-0.09750)%
(4.03613,-0.10984)(4.02454,-0.12225)(4.01295,-0.13475)(4.00136,-0.14733)(3.98978,-0.15999)%
(3.97819,-0.17273)(3.96661,-0.18556)(3.95502,-0.19849)(3.94242,-0.21263)(3.92965,-0.22707)%
(3.91688,-0.24163)(3.90412,-0.25629)(3.89136,-0.27104)(3.87861,-0.28590)(3.86586,-0.30087)%
(3.85312,-0.31594)(3.84039,-0.33111)(3.82768,-0.34640)(3.81497,-0.36179)(3.80227,-0.37729)%
(3.78959,-0.39289)(3.77691,-0.40860)(3.76424,-0.42443)(3.75159,-0.44035)(3.73818,-0.45737)%
(3.72474,-0.47456)(3.71133,-0.49187)(3.69793,-0.50930)(3.68456,-0.52684)(3.67121,-0.54452)%
(3.65788,-0.56230)(3.64458,-0.58022)(3.63130,-0.59825)(3.61806,-0.61640)(3.60484,-0.63467)%
(3.59165,-0.65305)(3.58046,-0.66878)%
\polyline(3.58046,-0.66878)(3.57849,-0.67156)(3.56536,-0.69019)%
\polyline(3.56536,-0.69019)(3.55357,-0.70706)%
\polyline(-3.06418,0.99143)(-3.08553,0.97184)%
\polyline(-3.08553,0.97184)(-3.08686,0.97062)(-3.11362,0.94622)(-3.11948,0.94092)%
\polyline(-3.11948,0.94092)(-3.13761,0.92449)(-3.13907,0.92317)%
\polyline(-3.13907,0.92317)(-3.16783,0.89729)(-3.19522,0.87281)(-3.22266,0.84845)%
(-3.25016,0.82419)(-3.27772,0.80003)(-3.30534,0.77597)(-3.33301,0.75200)(-3.36074,0.72813)%
(-3.38853,0.70436)(-3.41636,0.68069)(-3.44425,0.65711)(-3.47219,0.63363)(-3.50018,0.61023)%
(-3.52822,0.58692)(-3.55631,0.56371)(-3.58445,0.54058)(-3.61263,0.51754)(-3.64086,0.49457)%
(-3.66913,0.47171)(-3.69745,0.44893)(-3.72580,0.42622)(-3.75421,0.40360)(-3.78264,0.38106)%
(-3.79399,0.37210)(-3.81135,0.35842)(-3.84022,0.33576)(-3.86912,0.31318)(-3.89804,0.29071)%
(-3.92698,0.26832)(-3.95595,0.24601)(-3.98494,0.22381)(-4.01395,0.20167)(-4.04298,0.17962)%
(-4.07204,0.15765)(-4.10111,0.13577)(-4.13019,0.11396)(-4.15930,0.09224)(-4.18842,0.07060)%
(-4.21756,0.04902)(-4.24672,0.02753)(-4.27590,0.00611)(-4.30509,-0.01524)(-4.33429,-0.03651)%
(-4.36264,-0.05708)%
\polyline(-4.36264,-0.05708)(-4.36351,-0.05771)(-4.39201,-0.07831)%
\polyline(-4.39201,-0.07831)(-4.39274,-0.07884)(-4.42186,-0.09980)%
\polyline(-4.42186,-0.09980)(-4.42198,-0.09989)(-4.45050,-0.12037)%
\polyline(-4.45050,-0.12037)(-4.45124,-0.12089)(-4.47978,-0.14129)%
\polyline(-2.19106,-1.24852)(-1.70856,-1.26890)(-1.20503,-1.27286)(-0.68736,-1.25991)%
(-0.16263,-1.22984)(0.36205,-1.18268)(0.47372,-1.16888)%
\polyline(0.47372,-1.16888)(0.58907,-1.15463)%
\polyline(0.58907,-1.15463)(0.63393,-1.14909)%
\polyline(0.63393,-1.14909)(0.65207,-1.14685)%
\polyline(0.65207,-1.14685)(0.67050,-1.14457)%
\polyline(0.67050,-1.14457)(0.68922,-1.14225)%
\polyline(0.68922,-1.14225)(0.70440,-1.14038)%
\polyline(0.70440,-1.14038)(0.71785,-1.13872)%
\polyline(0.71785,-1.13872)(0.87957,-1.11874)(1.38295,-1.03852)(1.86550,-0.94285)%
(2.32087,-0.83271)(2.74314,-0.70934)(3.12697,-0.57415)(3.46759,-0.42872)(3.52907,-0.39646)%
\polyline(3.52907,-0.39646)(3.76096,-0.27479)(4.00373,-0.11419)(4.07494,-0.05210)%
\polyline(4.07494,-0.05210)(4.09989,-0.03035)%
\polyline(4.09989,-0.03035)(4.12475,-0.00867)%
\polyline(4.12475,-0.00867)(4.14926,0.01270)%
\polyline(4.14926,0.01270)(4.17018,0.03094)%
\polyline(4.17018,0.03094)(4.19334,0.05114)(4.32803,0.21921)(4.40684,0.38801)(4.42963,0.55555)%
(4.39704,0.71983)(4.31045,0.87895)(4.17202,1.03110)(3.98453,1.17455)(3.75141,1.30775)%
(3.47659,1.42927)(3.16450,1.53787)(2.81994,1.63249)(2.44801,1.71225)(2.05406,1.77648)%
(1.64354,1.82472)(1.22196,1.85667)(0.79483,1.87228)(0.36752,1.87164)(-0.05475,1.85505)%
(-0.46699,1.82299)(-0.86453,1.77605)(-1.24299,1.71500)(-1.59842,1.64072)(-1.92726,1.55420)%
(-2.22634,1.45652)(-2.49298,1.34883)(-2.72498,1.23234)(-2.92053,1.10832)(-3.05915,0.99387)%
\polyline(-3.05915,0.99387)(-3.07832,0.97804)(-3.08577,0.96959)%
\polyline(-3.08577,0.96959)(-3.13083,0.91843)%
\polyline(-3.13083,0.91843)(-3.19745,0.84279)(-3.27748,0.70388)(-3.31832,0.56256)%
(-3.32030,0.42011)(-3.28466,0.27993)%
\polyline(-3.28466,0.27993)(-3.28411,0.27774)(-3.21076,0.13665)(-3.10156,-0.00204)%
\polyline(-4.50832,-0.71939)(-4.48018,-0.70505)(-4.45204,-0.69069)(-4.42391,-0.67635)%
(-4.39577,-0.66201)(-4.36764,-0.64765)(-4.33950,-0.63331)(-4.31137,-0.61896)(-4.28323,-0.60461)%
(-4.25509,-0.59027)(-4.22696,-0.57592)(-4.19883,-0.56157)(-4.17069,-0.54723)(-4.14256,-0.53288)%
(-4.11442,-0.51853)(-4.08628,-0.50419)(-4.05816,-0.48984)(-4.03002,-0.47550)(-4.00188,-0.46114)%
(-3.97374,-0.44680)(-3.94561,-0.43246)(-3.91747,-0.41810)(-3.88934,-0.40376)(-3.86121,-0.38941)%
(-3.83307,-0.37506)(-3.80493,-0.36072)(-3.77680,-0.34637)(-3.74866,-0.33202)(-3.72053,-0.31768)%
(-3.69240,-0.30333)(-3.66426,-0.28898)(-3.63612,-0.27463)(-3.60798,-0.26029)(-3.57985,-0.24594)%
(-3.55172,-0.23159)(-3.52358,-0.21725)(-3.49544,-0.20291)(-3.46731,-0.18855)(-3.43917,-0.17421)%
(-3.41104,-0.15986)(-3.38291,-0.14551)(-3.35477,-0.13117)(-3.32663,-0.11682)(-3.29849,-0.10247)%
(-3.27036,-0.08813)(-3.24223,-0.07378)(-3.21409,-0.05943)(-3.18596,-0.04508)(-3.15782,-0.03074)%
(-3.12968,-0.01639)(-3.10156,-0.00204)%
\polyline(-1.46071,-0.83234)(-1.13975,-0.91658)(-0.80538,-0.99082)(-0.46084,-1.05461)%
(-0.10951,-1.10754)(0.24517,-1.14931)(0.32079,-1.15579)%
\polyline(0.32079,-1.15579)(0.33352,-1.15687)%
\polyline(0.33352,-1.15687)(0.34650,-1.15798)%
\polyline(0.34650,-1.15798)(0.37328,-1.16028)%
\polyline(0.47372,-1.16888)(0.59968,-1.17967)(0.95043,-1.19852)(1.29382,-1.20578)%
(1.62620,-1.20155)(1.94397,-1.18595)(2.24352,-1.15927)(2.52135,-1.12186)(2.77405,-1.07422)%
(2.99836,-1.01694)(3.19118,-0.95072)(3.34966,-0.87636)(3.47122,-0.79480)(3.55357,-0.70703)%
\polyline(3.55357,-0.70703)(3.56215,-0.68772)%
\polyline(3.52908,-0.39646)(3.45900,-0.31754)(3.32531,-0.21704)(3.14779,-0.11804)%
(2.92746,-0.02197)(2.66597,0.06976)(2.36552,0.15576)(2.02892,0.23472)(1.65952,0.30536)%
(1.26126,0.36654)(0.83854,0.41719)(0.39625,0.45641)(-0.06031,0.48344)(-0.52551,0.49767)%
(-0.99341,0.49873)(-1.45793,0.48641)(-1.91283,0.46071)(-2.35191,0.42184)(-2.76899,0.37025)%
(-3.15812,0.30658)(-3.28467,0.27993)%
\polyline(-4.47643,-0.14145)(-4.49999,-0.15610)(-4.61784,-0.26770)(-4.67862,-0.38173)%
(-4.68088,-0.49630)(-4.62403,-0.60950)(-4.50832,-0.71939)%
\polyline(-1.46071,-0.83234)(-1.47531,-0.84067)(-1.48992,-0.84899)(-1.50453,-0.85731)%
(-1.51913,-0.86564)(-1.53374,-0.87396)(-1.54835,-0.88229)(-1.56296,-0.89061)(-1.57756,-0.89893)%
(-1.59217,-0.90726)(-1.60678,-0.91558)(-1.62138,-0.92390)(-1.63599,-0.93223)(-1.65060,-0.94055)%
(-1.66521,-0.94887)(-1.67981,-0.95720)(-1.69442,-0.96552)(-1.70903,-0.97384)(-1.72363,-0.98217)%
(-1.73824,-0.99049)(-1.75285,-0.99881)(-1.76745,-1.00714)(-1.78206,-1.01546)(-1.79667,-1.02378)%
(-1.81128,-1.03211)(-1.82588,-1.04043)(-1.84049,-1.04875)(-1.85510,-1.05708)(-1.86970,-1.06540)%
(-1.88431,-1.07372)(-1.89892,-1.08205)(-1.91353,-1.09037)(-1.92813,-1.09869)(-1.94274,-1.10702)%
(-1.95735,-1.11534)(-1.97195,-1.12366)(-1.98656,-1.13199)(-2.00117,-1.14031)(-2.01577,-1.14864)%
(-2.03038,-1.15696)(-2.04499,-1.16528)(-2.05960,-1.17361)(-2.07420,-1.18193)(-2.08881,-1.19025)%
(-2.10342,-1.19858)(-2.11802,-1.20690)(-2.13263,-1.21522)(-2.14724,-1.22355)(-2.16185,-1.23187)%
(-2.17645,-1.24019)(-2.19106,-1.24852)%
\polyline(-1.58243,-0.90171)(-1.23456,-0.97530)(-0.87199,-1.03783)(-0.49860,-1.08883)%
(-0.11836,-1.12793)(0.26409,-1.15482)(0.26465,-1.15487)(0.32639,-1.15723)(0.34598,-1.15799)%
\polyline(0.38011,-1.15930)(0.39402,-1.15984)(0.40819,-1.16037)(0.43891,-1.16156)%
\polyline(0.51072,-1.16432)(0.64633,-1.16951)(1.02253,-1.17185)(1.38910,-1.16197)%
(1.74198,-1.14007)(2.07716,-1.10652)(2.39076,-1.06175)(2.67906,-1.00634)(2.93854,-0.94098)%
(3.16592,-0.86648)(3.35821,-0.78374)(3.51273,-0.69376)(3.62715,-0.59766)(3.62854,-0.59572)%
(3.64139,-0.57779)%
\polyline(3.71499,-0.29892)(3.71293,-0.28474)(3.65229,-0.17661)(3.54659,-0.06887)%
(3.39633,0.03709)(3.20259,0.13984)(2.96696,0.23800)(2.69163,0.33020)(2.37927,0.41518)%
(2.03310,0.49168)(1.65686,0.55859)(1.25471,0.61490)(0.83125,0.65971)(0.39146,0.69228)%
(-0.05939,0.71204)(-0.51576,0.71856)(-0.97194,0.71162)(-1.42210,0.69118)(-1.86044,0.65737)%
(-2.28113,0.61057)(-2.67855,0.55130)(-3.04727,0.48029)(-3.31940,0.41381)%
\polyline(-4.21365,0.04638)(-4.23714,0.02956)(-4.27477,0.00260)(-4.29624,-0.01278)%
(-4.32541,-0.04515)(-4.40126,-0.12933)(-4.45223,-0.24809)(-4.44809,-0.36729)(-4.38849,-0.48514)%
(-4.27386,-0.59983)%
\polyline(-1.70416,-0.97107)(-1.32935,-1.03402)(-0.93859,-1.08483)(-0.53636,-1.12304)%
(-0.12722,-1.14830)(0.28322,-1.16040)(0.28413,-1.16043)(0.37153,-1.16021)(0.40068,-1.16012)%
\polyline(0.41498,-1.16009)(0.42660,-1.16006)(0.44474,-1.16001)(0.45984,-1.15997)%
(0.47484,-1.15993)(0.50811,-1.15984)%
\polyline(0.54771,-1.15974)(0.69297,-1.15936)(1.09460,-1.14518)(1.48438,-1.11814)%
(1.85775,-1.07860)(2.21036,-1.02708)(2.53800,-0.96423)(2.83677,-0.89082)(3.10302,-0.80774)%
(3.33348,-0.71603)(3.52524,-0.61676)(3.67578,-0.51117)(3.71792,-0.46773)(3.74776,-0.43696)%
(3.76200,-0.42227)(3.77675,-0.40706)(3.78309,-0.40053)(3.78832,-0.39096)%
\polyline(3.85619,-0.21179)(3.86222,-0.16953)(3.83243,-0.05200)(3.75623,0.06493)(3.63418,0.17982)%
(3.46735,0.29123)(3.25739,0.39773)(3.00647,0.49798)(2.71729,0.59066)(2.39302,0.67459)%
(2.03730,0.74864)(1.65419,0.81182)(1.24816,0.86325)(0.82397,0.90222)(0.38668,0.92815)%
(-0.05845,0.94065)(-0.50600,0.93945)(-0.95046,0.92451)(-1.38628,0.89594)(-1.80803,0.85404)%
(-2.21035,0.79930)(-2.58810,0.73234)(-2.93641,0.65400)(-3.25073,0.56525)(-3.31862,0.54116)%
\polyline(-3.95366,0.24597)(-3.98166,0.22270)(-4.01109,0.19823)(-4.06682,0.15190)%
(-4.09249,0.13056)(-4.18467,0.00906)(-4.22585,-0.11445)(-4.21529,-0.23828)(-4.15294,-0.36078)%
(-4.03939,-0.48028)%
\polyline(-1.82588,-1.04043)(-1.42415,-1.09274)(-1.00520,-1.13184)(-0.57411,-1.15726)%
(-0.13607,-1.16869)(0.30361,-1.16599)(0.41932,-1.16154)(0.44433,-1.16057)(0.45942,-1.15999)%
(0.47475,-1.15940)(0.48763,-1.15891)(0.49868,-1.15848)(0.51472,-1.15786)(0.53103,-1.15723)%
\polyline(0.54762,-1.15660)(0.58031,-1.15534)(0.73962,-1.14920)(1.16669,-1.11852)%
(1.57967,-1.07432)(1.97353,-1.01713)(2.34356,-0.94764)(2.68524,-0.86671)(2.99447,-0.77529)%
(3.26751,-0.67451)(3.50104,-0.56557)(3.69226,-0.44979)(3.82320,-0.34152)(3.83885,-0.32857)%
(3.85065,-0.31382)%
\polyline(3.96530,-0.13962)(3.99160,-0.07574)(3.99592,0.05281)(3.95194,0.18073)(3.86018,0.30647)%
(3.72176,0.42850)(3.53837,0.54535)(3.31219,0.65561)(3.04598,0.75795)(2.74295,0.85112)%
(2.40677,0.93401)(2.04149,1.00560)(1.65153,1.06504)(1.24162,1.11161)(0.81669,1.14473)%
(0.38189,1.16403)(-0.05753,1.16925)(-0.49625,1.16033)(-0.92897,1.13739)(-1.35046,1.10071)%
(-1.75563,1.05071)(-2.13958,0.98802)(-2.49766,0.91338)(-2.82556,0.82771)(-3.11929,0.73202)%
(-3.28942,0.66257)%
\polyline(-3.75228,0.40399)(-3.76239,0.39703)(-3.79110,0.36904)(-3.83291,0.32830)%
(-3.88874,0.27388)(-3.96808,0.14743)(-3.99946,0.01919)(-3.98250,-0.10928)(-3.91739,-0.23643)%
(-3.80493,-0.36072)%
\polyline(-1.94761,-1.10979)(-1.51895,-1.15146)(-1.07182,-1.17885)(-0.61186,-1.19148)%
(-0.14493,-1.18908)(0.32308,-1.17156)(0.43758,-1.16352)(0.47232,-1.16109)(0.50625,-1.15870)%
(0.52239,-1.15756)(0.53879,-1.15642)(0.55633,-1.15519)(0.56818,-1.15435)(0.58147,-1.15342)%
(0.59886,-1.15220)(0.61654,-1.15096)(0.65476,-1.14827)(0.78626,-1.13905)(1.23878,-1.09185)%
(1.67494,-1.03050)(2.08931,-0.95566)(2.47675,-0.86821)(2.83249,-0.76919)(3.15218,-0.65976)%
(3.43199,-0.54127)(3.66860,-0.41510)(3.85929,-0.28281)(3.89619,-0.24741)(3.92381,-0.22092)%
(3.95124,-0.19460)(3.97619,-0.17066)(3.98874,-0.15862)(4.00153,-0.14634)(4.00191,-0.14598)%
(4.01103,-0.13230)%
\polyline(4.04988,-0.07395)(4.06037,-0.05820)(4.09497,-0.00626)(4.13761,0.13468)(4.12962,0.27516)%
(4.07144,0.41347)(3.96413,0.54801)(3.80936,0.67719)(3.60938,0.79948)(3.36699,0.91350)%
(3.08548,1.01792)(2.76861,1.11158)(2.42052,1.19342)(2.04568,1.26256)(1.64887,1.31827)%
(1.23507,1.35996)(0.80940,1.38725)(0.37709,1.39990)(-0.05660,1.39785)(-0.48650,1.38122)%
(-0.90749,1.35028)(-1.31464,1.30547)(-1.70324,1.24739)(-2.06880,1.17674)(-2.40721,1.09443)%
(-2.71470,1.00141)(-2.98786,0.89880)(-3.22372,0.78777)(-3.23203,0.78276)%
\polyline(-3.51519,0.59294)(-3.55439,0.56143)(-3.57408,0.54561)(-3.58120,0.53737)%
(-3.68498,0.41721)(-3.75150,0.28580)(-3.77307,0.15283)(-3.74970,0.01973)(-3.68184,-0.11207)%
(-3.57048,-0.24115)%
\polyline(-2.06933,-1.17915)(-1.61375,-1.21018)(-1.13842,-1.22585)(-0.64961,-1.22569)%
(-0.15378,-1.20945)(0.34257,-1.17711)(0.45825,-1.16574)(0.52853,-1.15883)(0.57245,-1.15451)%
(0.58964,-1.15282)(0.60714,-1.15110)(0.62535,-1.14931)(0.64342,-1.14753)(0.65872,-1.14602)%
(0.67104,-1.14481)(0.68981,-1.14297)(0.83292,-1.12889)(1.31087,-1.06519)(1.77022,-0.98667)%
(2.20509,-0.89419)(2.60994,-0.78877)(2.97972,-0.67167)(3.30989,-0.54424)(3.59648,-0.40803)%
(3.83617,-0.26465)(3.98178,-0.15069)(4.00863,-0.12969)(4.02632,-0.11584)(4.03363,-0.10779)%
(4.05650,-0.08265)%
\polyline(4.12423,-0.00818)(4.13520,0.00389)%
\polyline(4.16399,0.03553)(4.16498,0.03661)(4.25091,0.19088)(4.28362,0.34511)(4.26333,0.49749)%
(4.19095,0.64621)(4.06808,0.78955)(3.89695,0.92587)(3.68040,1.05362)(3.42179,1.17139)%
(3.12499,1.27790)(2.79427,1.37204)(2.43426,1.45284)(2.04987,1.51952)(1.64620,1.57149)%
(1.22851,1.60831)(0.80211,1.62976)(0.37231,1.63577)(-0.05567,1.62645)(-0.47674,1.60211)%
(-0.88601,1.56317)(-1.27882,1.51024)(-1.65083,1.44406)(-1.99803,1.36548)(-2.31677,1.27547)%
(-2.60384,1.17512)(-2.85642,1.06557)(-3.07213,0.94805)(-3.14563,0.89644)(-3.15760,0.88804)%
\polyline(-3.21819,0.84550)(-3.24907,0.82382)(-3.27471,0.79951)(-3.28533,0.78944)%
(-3.32911,0.74792)(-3.38161,0.69815)(-3.38577,0.69421)(-3.48123,0.56055)(-3.53490,0.42419)%
(-3.54669,0.28647)(-3.51691,0.14874)(-3.44630,0.01228)(-3.33601,-0.12160)%
\polyline(0.16146,-1.14142)(0.16299,-1.14252)(0.16452,-1.14363)(0.16605,-1.14473)%
(0.16758,-1.14583)(0.16912,-1.14694)(0.17064,-1.14805)(0.17186,-1.14893)%
\polyline(0.15992,-1.14032)(0.16146,-1.14142)(0.16299,-1.14252)(0.16452,-1.14363)%
(0.16605,-1.14473)(0.16758,-1.14583)(0.16912,-1.14694)(0.17064,-1.14805)(0.17218,-1.14915)%
(0.17371,-1.15026)(0.17524,-1.15136)(0.17677,-1.15246)(0.17831,-1.15357)(0.17983,-1.15468)%
(0.18137,-1.15578)(0.18290,-1.15688)(0.18444,-1.15799)(0.18596,-1.15909)(0.18750,-1.16020)%
(0.18903,-1.16130)(0.19056,-1.16241)(0.19209,-1.16351)(0.19362,-1.16462)(0.19515,-1.16572)%
(0.19669,-1.16682)(0.19822,-1.16793)(0.19975,-1.16904)(0.20128,-1.17014)(0.20281,-1.17124)%
(0.20434,-1.17235)(0.20587,-1.17345)(0.20741,-1.17456)(0.20894,-1.17566)(0.21047,-1.17677)%
(0.21200,-1.17787)(0.21353,-1.17898)(0.21506,-1.18008)(0.21660,-1.18118)(0.21813,-1.18229)%
(0.21966,-1.18340)(0.22119,-1.18450)(0.22272,-1.18560)(0.22425,-1.18671)(0.22579,-1.18781)%
(0.22732,-1.18892)(0.22885,-1.19003)(0.23038,-1.19112)(0.23192,-1.19223)(0.23344,-1.19334)%
(0.23498,-1.19445)(0.23651,-1.19554)%
\polyline(0.17957,-1.15449)(0.17983,-1.15468)(0.18137,-1.15578)(0.18290,-1.15688)%
(0.18444,-1.15799)(0.18596,-1.15909)(0.18750,-1.16020)(0.18903,-1.16130)(0.19056,-1.16241)%
(0.19209,-1.16351)(0.19362,-1.16462)(0.19515,-1.16572)(0.19669,-1.16682)(0.19822,-1.16793)%
(0.19975,-1.16904)(0.20128,-1.17014)(0.20281,-1.17124)(0.20434,-1.17235)(0.20587,-1.17345)%
(0.20741,-1.17456)(0.20894,-1.17566)(0.21047,-1.17677)(0.21200,-1.17787)(0.21353,-1.17898)%
(0.21506,-1.18008)(0.21660,-1.18118)(0.21813,-1.18229)(0.21966,-1.18340)(0.22119,-1.18450)%
(0.22272,-1.18560)(0.22425,-1.18671)(0.22579,-1.18781)(0.22732,-1.18892)(0.22885,-1.19003)%
(0.23038,-1.19112)(0.23192,-1.19223)(0.23344,-1.19334)(0.23456,-1.19414)%
\polyline(1.79640,-1.19319)(1.80852,-1.18634)(1.81110,-1.18489)(1.82354,-1.17785)%
(1.83857,-1.16936)(1.85360,-1.16086)(1.86863,-1.15237)(1.88366,-1.14387)(1.89869,-1.13538)%
(1.91371,-1.12688)(1.92874,-1.11839)(1.94377,-1.10990)(1.95880,-1.10140)(1.97383,-1.09291)%
(1.98886,-1.08441)(2.00389,-1.07592)(2.01892,-1.06742)(2.03394,-1.05893)(2.04897,-1.05044)%
(2.06400,-1.04194)(2.07903,-1.03345)(2.09406,-1.02495)(2.10909,-1.01646)(2.12412,-1.00797)%
(2.13914,-0.99947)(2.15417,-0.99098)(2.16920,-0.98248)(2.18423,-0.97399)(2.19926,-0.96549)%
(2.21429,-0.95700)(2.22932,-0.94851)(2.24435,-0.94001)(2.25937,-0.93152)(2.27440,-0.92302)%
(2.28943,-0.91453)(2.30446,-0.90604)(2.31949,-0.89754)(2.33452,-0.88905)(2.34954,-0.88055)%
(2.36458,-0.87206)(2.37960,-0.86357)(2.39463,-0.85507)(2.40966,-0.84658)(2.42469,-0.83808)%
(2.43972,-0.82959)(2.45475,-0.82109)(2.46977,-0.81260)(2.48481,-0.80410)(2.49983,-0.79561)%
(2.51486,-0.78712)(2.52989,-0.77862)(2.54116,-0.77226)%
\polyline(3.05908,-0.99609)(3.07577,-0.98065)(3.09591,-0.96203)(3.11604,-0.94341)%
(3.13617,-0.92479)(3.15630,-0.90617)(3.17644,-0.88755)(3.19657,-0.86893)(3.21670,-0.85031)%
(3.23683,-0.83169)(3.25697,-0.81307)(3.27710,-0.79445)(3.29723,-0.77584)(3.31736,-0.75722)%
(3.33750,-0.73859)(3.35763,-0.71997)(3.37776,-0.70136)(3.39789,-0.68274)(3.41803,-0.66412)%
(3.43816,-0.64549)(3.45829,-0.62688)(3.47842,-0.60826)(3.49856,-0.58964)(3.51869,-0.57102)%
(3.53882,-0.55240)(3.55895,-0.53378)(3.57909,-0.51516)(3.59921,-0.49655)(3.61935,-0.47793)%
(3.63948,-0.45930)(3.65962,-0.44068)(3.67974,-0.42207)(3.69988,-0.40345)(3.72001,-0.38483)%
(3.74015,-0.36620)(3.76027,-0.34759)(3.78041,-0.32897)(3.80054,-0.31035)(3.82068,-0.29173)%
(3.84081,-0.27311)(3.86094,-0.25449)(3.88107,-0.23587)(3.90121,-0.21725)(3.92134,-0.19864)%
(3.94147,-0.18001)(3.96160,-0.16139)(3.98174,-0.14277)(4.00187,-0.12416)(4.02177,-0.10575)%
(4.02200,-0.10554)(4.04022,-0.08868)(4.04214,-0.08691)(4.05467,-0.07532)(4.06227,-0.06829)%
\polyline(3.81923,-0.23625)(3.81930,-0.23613)(3.83528,-0.20939)(3.85126,-0.18266)%
(3.86724,-0.15592)(3.88323,-0.12919)(3.89920,-0.10246)(3.91518,-0.07572)(3.93117,-0.04899)%
(3.94714,-0.02226)(3.96313,0.00447)(3.97911,0.03121)(3.99508,0.05795)(4.01107,0.08468)%
(4.02704,0.11142)(4.04303,0.13815)(4.05901,0.16488)(4.07498,0.19161)(4.09097,0.21835)%
(4.10695,0.24508)(4.12293,0.27182)(4.13891,0.29856)(4.15489,0.32529)(4.17087,0.35202)%
(4.18685,0.37875)(4.20283,0.40549)(4.21881,0.43222)(4.23479,0.45895)(4.25077,0.48569)%
(4.26675,0.51243)(4.28273,0.53916)(4.29871,0.56589)(4.31469,0.59263)(4.33068,0.61936)%
(4.34665,0.64609)(4.36263,0.67282)(4.37862,0.69956)(4.39403,0.72535)%
\polyline(3.18677,-0.13766)(3.19373,-0.10679)(3.20069,-0.07591)(3.20765,-0.04504)%
(3.21461,-0.01417)(3.22157,0.01670)(3.22852,0.04758)(3.23549,0.07845)(3.24244,0.10932)%
(3.24940,0.14020)(3.25635,0.17107)(3.26332,0.20194)(3.27028,0.23280)(3.27724,0.26368)%
(3.28419,0.29455)(3.29115,0.32542)(3.29811,0.35629)(3.30507,0.38717)(3.31202,0.41804)%
(3.31899,0.44891)(3.32594,0.47979)(3.33291,0.51066)(3.33986,0.54153)(3.34682,0.57240)%
(3.35378,0.60328)(3.36074,0.63415)(3.36769,0.66502)(3.37465,0.69589)(3.38161,0.72677)%
(3.38857,0.75764)(3.39553,0.78851)(3.40249,0.81938)(3.40945,0.85026)(3.41641,0.88113)%
(3.42337,0.91200)(3.43032,0.94288)(3.43728,0.97375)(3.44424,1.00462)(3.45120,1.03550)%
(3.45816,1.06636)(3.46512,1.09723)(3.47208,1.12810)(3.47904,1.15898)(3.48599,1.18985)%
(3.49295,1.22072)(3.49991,1.25159)(3.50687,1.28247)(3.51382,1.31334)(3.52078,1.34421)%
(3.52774,1.37508)(3.53422,1.40379)%
\polyline(1.82588,0.27538)(1.82588,0.30598)(1.82588,0.33658)(1.82588,0.36719)(1.82588,0.39779)%
(1.82588,0.42839)(1.82588,0.45899)(1.82588,0.48959)(1.82588,0.52020)(1.82588,0.55080)%
(1.82588,0.58140)(1.82588,0.61200)(1.82588,0.64260)(1.82588,0.67321)(1.82588,0.70381)%
(1.82588,0.73441)(1.82588,0.76501)(1.82588,0.79561)(1.82588,0.82622)(1.82588,0.85682)%
(1.82588,0.88742)(1.82588,0.91802)(1.82588,0.94862)(1.82588,0.97923)(1.82588,1.00983)%
(1.82588,1.04043)(1.82588,1.07103)(1.82588,1.10163)(1.82588,1.13224)(1.82588,1.16284)%
(1.82588,1.19344)(1.82588,1.22404)(1.82588,1.25464)(1.82588,1.28525)(1.82588,1.31585)%
(1.82588,1.34645)(1.82588,1.37705)(1.82588,1.40766)(1.82588,1.43826)(1.82588,1.46886)%
(1.82588,1.49946)(1.82588,1.53006)(1.82588,1.56067)(1.82588,1.59127)(1.82588,1.62187)%
(1.82588,1.65247)(1.82588,1.68307)(1.82588,1.71368)(1.82588,1.74428)(1.82588,1.77488)%
(1.82588,1.80329)%
\polyline(-0.20848,0.48941)(-0.20806,0.51654)(-0.20765,0.54369)(-0.20724,0.57083)%
(-0.20683,0.59797)(-0.20642,0.62511)(-0.20601,0.65225)(-0.20560,0.67939)(-0.20519,0.70654)%
(-0.20479,0.73368)(-0.20437,0.76082)(-0.20396,0.78796)(-0.20355,0.81509)(-0.20314,0.84223)%
(-0.20273,0.86937)(-0.20232,0.89652)(-0.20191,0.92366)(-0.20150,0.95080)(-0.20109,0.97794)%
(-0.20068,1.00508)(-0.20026,1.03222)(-0.19985,1.05936)(-0.19944,1.08651)(-0.19903,1.11365)%
(-0.19862,1.14079)(-0.19822,1.16793)(-0.19781,1.19507)(-0.19740,1.22221)(-0.19699,1.24935)%
(-0.19657,1.27649)(-0.19616,1.30364)(-0.19575,1.33078)(-0.19534,1.35792)(-0.19493,1.38506)%
(-0.19452,1.41220)(-0.19411,1.43935)(-0.19370,1.46648)(-0.19329,1.49363)(-0.19287,1.52077)%
(-0.19246,1.54790)(-0.19206,1.57504)(-0.19165,1.60218)(-0.19124,1.62932)(-0.19083,1.65647)%
(-0.19042,1.68361)(-0.19001,1.71075)(-0.18960,1.73789)(-0.18918,1.76503)(-0.18877,1.79217)%
(-0.18836,1.81931)(-0.18798,1.84469)%
\polyline(-2.38612,0.41818)(-2.37744,0.44075)(-2.36877,0.46333)(-2.36009,0.48590)%
(-2.35142,0.50847)(-2.34274,0.53104)(-2.33406,0.55361)(-2.32539,0.57618)(-2.31671,0.59875)%
(-2.30803,0.62132)(-2.29936,0.64390)(-2.29068,0.66647)(-2.28200,0.68904)(-2.27333,0.71162)%
(-2.26465,0.73419)(-2.25598,0.75676)(-2.24730,0.77933)(-2.23862,0.80191)(-2.22994,0.82448)%
(-2.22126,0.84705)(-2.21258,0.86963)(-2.20391,0.89220)(-2.19523,0.91477)(-2.18655,0.93733)%
(-2.17788,0.95991)(-2.16920,0.98248)(-2.16052,1.00505)(-2.15185,1.02763)(-2.14317,1.05020)%
(-2.13450,1.07277)(-2.12582,1.09534)(-2.11714,1.11792)(-2.10847,1.14049)(-2.09979,1.16306)%
(-2.09111,1.18564)(-2.08244,1.20821)(-2.07376,1.23078)(-2.06508,1.25335)(-2.05641,1.27593)%
(-2.04773,1.29849)(-2.03905,1.32107)(-2.03037,1.34364)(-2.02169,1.36621)(-2.01302,1.38878)%
(-2.00434,1.41135)(-1.99566,1.43393)(-1.98699,1.45650)(-1.97831,1.47907)(-1.96963,1.50164)%
(-1.96096,1.52422)(-1.95262,1.54592)%
\polyline(-3.13114,0.92944)(-3.11914,0.94054)(-3.09591,0.96203)(-3.09092,0.96664)%
\polyline(-3.07886,0.97779)(-3.07577,0.98065)(-3.06975,0.98622)(-3.05806,0.99703)%
\polyline(-4.67989,-0.44645)(-4.65937,-0.43481)(-4.63170,-0.41911)(-4.60402,-0.40340)%
(-4.57634,-0.38770)(-4.54866,-0.37200)(-4.52098,-0.35629)(-4.49330,-0.34059)(-4.46562,-0.32489)%
(-4.43794,-0.30919)(-4.41027,-0.29349)(-4.38259,-0.27778)(-4.35491,-0.26208)(-4.32723,-0.24638)%
(-4.29955,-0.23067)(-4.27187,-0.21497)(-4.24419,-0.19927)(-4.21651,-0.18357)(-4.18884,-0.16787)%
(-4.16116,-0.15216)(-4.13348,-0.13646)(-4.10580,-0.12076)(-4.07812,-0.10505)(-4.05044,-0.08935)%
(-4.02276,-0.07365)(-3.99508,-0.05795)(-3.96740,-0.04224)(-3.93973,-0.02654)(-3.91205,-0.01084)%
(-3.88437,0.00486)(-3.85669,0.02056)(-3.82901,0.03627)(-3.80134,0.05197)(-3.77365,0.06767)%
(-3.74597,0.08338)(-3.71830,0.09908)(-3.69062,0.11478)(-3.66294,0.13048)(-3.63526,0.14618)%
(-3.60758,0.16189)(-3.57990,0.17759)(-3.55223,0.19329)(-3.52454,0.20900)(-3.49687,0.22470)%
(-3.46919,0.24040)(-3.44151,0.25610)(-3.41383,0.27180)(-3.38615,0.28751)(-3.35847,0.30321)%
(-3.33080,0.31891)(-3.30322,0.33456)%
\end{picture}}%

%% file: fig/p025.tex
{\unitlength=1cm%
\begin{picture}%
(9.94,5.4)(-5.19,-2.08)%
\linethickness{0.008in}
\polyline(0.71106,-1.14153)(0.70289,-1.14257)(0.69363,-1.14373)(0.69340,-1.14377)%
\polyline(0.69340,-1.14377)(0.68441,-1.14486)(0.67527,-1.14595)(0.67504,-1.14597)%
\polyline(0.43569,-1.16157)(0.42987,-1.16156)(0.42265,-1.16153)(0.41549,-1.16146)%
\polyline(0.41549,-1.16146)(0.40840,-1.16136)(0.40137,-1.16122)(0.39440,-1.16104)%
(0.38749,-1.16083)(0.38065,-1.16059)(0.37806,-1.16048)%
\polyline(0.37806,-1.16048)(0.37387,-1.16031)(0.36715,-1.15999)(0.36698,-1.15998)%
\polyline(0.36698,-1.15998)(0.36050,-1.15964)(0.35406,-1.15927)%
\polyline(0.35406,-1.15927)(0.35390,-1.15927)(0.34736,-1.15885)(0.34089,-1.15840)%
(0.34072,-1.15838)%
\polyline(0.34072,-1.15838)(0.33447,-1.15791)(0.32828,-1.15740)%
\polyline(0.71222,-1.14137)(0.70289,-1.14257)(0.69363,-1.14373)(0.68441,-1.14486)%
(0.67527,-1.14595)(0.66618,-1.14699)(0.65717,-1.14800)(0.64820,-1.14899)(0.63930,-1.14992)%
(0.63047,-1.15082)(0.62170,-1.15168)(0.61299,-1.15251)(0.60434,-1.15331)(0.59576,-1.15406)%
(0.58724,-1.15478)(0.57878,-1.15546)(0.57038,-1.15610)(0.56204,-1.15671)(0.55378,-1.15728)%
(0.54556,-1.15782)(0.53741,-1.15832)(0.52933,-1.15878)(0.52131,-1.15921)(0.51334,-1.15960)%
(0.50544,-1.15997)(0.49761,-1.16028)(0.48982,-1.16057)(0.48212,-1.16081)(0.47447,-1.16103)%
(0.46688,-1.16120)(0.45935,-1.16135)(0.45189,-1.16146)(0.44449,-1.16153)(0.43714,-1.16157)%
(0.42987,-1.16156)(0.42265,-1.16153)(0.41549,-1.16146)(0.40840,-1.16136)(0.40137,-1.16122)%
(0.39440,-1.16104)(0.38749,-1.16083)(0.38065,-1.16059)(0.37387,-1.16031)(0.36715,-1.15999)%
(0.36050,-1.15964)(0.35390,-1.15927)(0.34736,-1.15885)(0.34089,-1.15840)(0.33447,-1.15791)%
(0.32813,-1.15738)(0.32184,-1.15683)%
\polyline(4.17266,0.03032)(4.16147,0.01925)(4.16119,0.01896)%
\polyline(4.16119,0.01896)(4.15000,0.00782)(4.14972,0.00753)%
\polyline(4.14972,0.00753)(4.13853,-0.00369)(4.13818,-0.00404)%
\polyline(4.13818,-0.00404)(4.12705,-0.01527)(4.12678,-0.01556)%
\polyline(4.12678,-0.01556)(4.11559,-0.02692)(4.11531,-0.02721)%
\polyline(4.11531,-0.02721)(4.10412,-0.03866)(4.10384,-0.03894)%
\polyline(4.10384,-0.03894)(4.09266,-0.05046)(4.09238,-0.05076)%
\polyline(4.09238,-0.05076)(4.08120,-0.06234)(4.08092,-0.06263)%
\polyline(4.08092,-0.06263)(4.06974,-0.07430)(4.06946,-0.07459)%
\polyline(4.06946,-0.07459)(4.05828,-0.08634)(4.04681,-0.09846)(4.03537,-0.11065)%
(4.02390,-0.12292)(4.01245,-0.13528)(4.00100,-0.14772)(3.98955,-0.16024)(3.97810,-0.17284)%
(3.96664,-0.18553)(3.95616,-0.19722)(3.94225,-0.21283)(3.92942,-0.22734)(3.91659,-0.24195)%
(3.90378,-0.25666)(3.89098,-0.27148)(3.87818,-0.28640)(3.86538,-0.30143)(3.85259,-0.31656)%
(3.83982,-0.33180)(3.82705,-0.34714)(3.81430,-0.36260)(3.80155,-0.37815)(3.78882,-0.39383)%
(3.77609,-0.40961)(3.76338,-0.42550)(3.75069,-0.44149)(3.73801,-0.45760)(3.72533,-0.47381)%
(3.71180,-0.49125)(3.69828,-0.50884)(3.68478,-0.52656)(3.67131,-0.54439)(3.65785,-0.56235)%
(3.64442,-0.58042)(3.63103,-0.59862)(3.61765,-0.61694)(3.60431,-0.63538)(3.59300,-0.65117)%
\polyline(3.59300,-0.65117)(3.59100,-0.65394)(3.57971,-0.66982)%
\polyline(3.57971,-0.66982)(3.57772,-0.67263)(3.57309,-0.67920)%
\polyline(3.57309,-0.67920)(3.56646,-0.68861)%
\polyline(3.56646,-0.68861)(3.56448,-0.69144)(3.55275,-0.70824)%
\polyline(-3.06784,0.98807)(-3.08793,0.96965)(-3.09195,0.96597)%
\polyline(-3.09195,0.96597)(-3.11478,0.94516)(-3.11613,0.94394)%
\polyline(-3.11613,0.94394)(-3.13903,0.92321)%
\polyline(-3.13903,0.92321)(-3.14172,0.92077)(-3.14712,0.91592)%
\polyline(-3.14712,0.91592)(-3.16873,0.89649)(-3.19583,0.87228)(-3.22298,0.84817)%
(-3.25021,0.82415)(-3.27749,0.80022)(-3.30486,0.77638)(-3.33229,0.75263)(-3.35978,0.72895)%
(-3.37847,0.71295)(-3.38749,0.70525)(-3.41554,0.68138)(-3.44366,0.65762)(-3.47180,0.63395)%
(-3.49999,0.61039)(-3.52822,0.58692)(-3.55650,0.56355)(-3.58482,0.54027)(-3.61317,0.51709)%
(-3.64156,0.49401)(-3.66999,0.47101)(-3.69846,0.44810)(-3.72697,0.42529)(-3.75551,0.40256)%
(-3.78408,0.37991)(-3.81269,0.35736)(-3.84133,0.33489)(-3.87000,0.31250)(-3.89871,0.29019)%
(-3.92744,0.26796)(-3.95621,0.24582)(-3.98500,0.22376)(-4.01382,0.20177)(-4.04267,0.17986)%
(-4.07154,0.15802)(-4.10044,0.13627)(-4.12938,0.11458)(-4.15832,0.09297)(-4.18730,0.07143)%
(-4.19911,0.06268)(-4.21653,0.04979)(-4.24593,0.02810)(-4.27535,0.00651)(-4.30477,-0.01500)%
(-4.33420,-0.03645)(-4.36363,-0.05780)(-4.39234,-0.07855)%
\polyline(-4.39234,-0.07855)(-4.39307,-0.07908)(-4.42252,-0.10028)(-4.45124,-0.12089)%
\polyline(-4.45124,-0.12089)(-4.45197,-0.12141)(-4.48120,-0.14231)%
\polyline(-2.19106,-1.24852)(-1.63826,-1.27047)(-1.05934,-1.27091)(-0.46456,-1.24921)%
(0.13547,-1.20518)(0.46305,-1.16875)%
\polyline(0.46305,-1.16875)(0.59308,-1.15429)%
\polyline(0.59308,-1.15429)(0.65708,-1.14717)%
\polyline(0.65708,-1.14717)(0.67519,-1.14516)%
\polyline(0.67519,-1.14516)(0.69356,-1.14312)%
\polyline(0.69356,-1.14312)(0.70658,-1.14166)%
\polyline(0.70658,-1.14166)(0.73016,-1.13904)(1.30907,-1.05149)(1.86208,-0.94360)%
(2.37973,-0.81686)(2.85327,-0.67309)(3.27492,-0.51448)(3.52575,-0.39630)%
\polyline(3.52575,-0.39630)(3.63793,-0.34345)(3.93681,-0.16267)(4.05996,-0.06236)%
\polyline(4.05996,-0.06236)(4.08535,-0.04167)%
\polyline(4.08535,-0.04167)(4.11066,-0.02106)%
\polyline(4.11066,-0.02106)(4.13530,-0.00099)%
\polyline(4.13530,-0.00099)(4.14842,0.00969)%
\polyline(4.14842,0.00969)(4.16728,0.02505)(4.17233,0.03113)%
\polyline(4.17233,0.03113)(4.32648,0.21675)(4.41285,0.40945)(4.42625,0.60015)(4.36788,0.78595)%
(4.24018,0.96402)(4.04681,1.13175)(3.79253,1.28672)(3.48300,1.42675)(3.12477,1.54998)%
(2.72499,1.65484)(2.29135,1.74007)(1.83184,1.80480)(1.35465,1.84843)(0.86800,1.87076)%
(0.37991,1.87188)(-0.10181,1.85221)(-0.56976,1.81242)(-1.01705,1.75349)(-1.43738,1.67661)%
(-1.82512,1.58316)(-2.17534,1.47469)(-2.48388,1.35291)(-2.74735,1.21959)(-2.96311,1.07664)%
(-3.05792,0.99067)%
\polyline(-3.05792,0.99067)(-3.11164,0.94195)%
\polyline(-3.11164,0.94195)(-3.12930,0.92593)(-3.13646,0.91624)%
\polyline(-3.13646,0.91624)(-3.24482,0.76941)(-3.30925,0.60899)(-3.32285,0.44655)%
(-3.28648,0.28394)%
\polyline(-3.28648,0.28394)(-3.20159,0.12290)(-3.07019,-0.03485)(-2.89469,-0.18772)%
(-2.67801,-0.33419)(-2.42340,-0.47287)(-2.13449,-0.60246)(-1.81521,-0.72178)(-1.46976,-0.82976)%
\polyline(-2.20465,-1.24768)(-2.18995,-1.23932)(-2.18031,-1.23384)%
\polyline(-2.18031,-1.23384)(-2.17526,-1.23097)(-2.16564,-1.22550)%
\polyline(-2.16564,-1.22550)(-2.16056,-1.22261)(-2.15097,-1.21716)%
\polyline(-2.15097,-1.21716)(-2.14586,-1.21426)(-2.13630,-1.20882)%
\polyline(-2.13630,-1.20882)(-2.13117,-1.20590)(-2.12212,-1.20075)%
\polyline(-2.12212,-1.20075)(-2.11647,-1.19754)(-2.10745,-1.19241)%
\polyline(-2.10745,-1.19241)(-2.10177,-1.18918)(-2.09278,-1.18407)%
\polyline(-2.09278,-1.18407)(-2.08707,-1.18082)(-2.07811,-1.17573)%
\polyline(-2.07811,-1.17573)(-2.07237,-1.17247)(-2.06344,-1.16738)%
\polyline(-2.06344,-1.16738)(-2.05768,-1.16410)(-2.04877,-1.15904)%
\polyline(-2.04877,-1.15904)(-2.04298,-1.15574)(-2.03410,-1.15069)%
\polyline(-2.03410,-1.15069)(-2.02828,-1.14738)(-2.01943,-1.14235)%
\polyline(-2.01943,-1.14235)(-2.01358,-1.13903)(-2.00476,-1.13401)%
\polyline(-2.00476,-1.13401)(-1.99888,-1.13067)(-1.99058,-1.12595)%
\polyline(-1.99058,-1.12595)(-1.98419,-1.12231)(-1.97591,-1.11760)%
\polyline(-1.97591,-1.11760)(-1.96949,-1.11395)(-1.96124,-1.10926)%
\polyline(-1.96124,-1.10926)(-1.95479,-1.10560)(-1.94657,-1.10092)%
\polyline(-1.94657,-1.10092)(-1.94009,-1.09724)(-1.93190,-1.09258)%
\polyline(-1.93190,-1.09258)(-1.92539,-1.08888)(-1.91723,-1.08423)%
\polyline(-1.91723,-1.08423)(-1.91070,-1.08051)(-1.90256,-1.07589)%
\polyline(-1.90256,-1.07589)(-1.89600,-1.07216)(-1.88789,-1.06755)%
\polyline(-1.88789,-1.06755)(-1.88131,-1.06380)(-1.87371,-1.05949)%
\polyline(-1.87371,-1.05949)(-1.86661,-1.05544)(-1.85904,-1.05114)%
\polyline(-1.85904,-1.05114)(-1.85190,-1.04708)(-1.84437,-1.04279)%
\polyline(-1.84437,-1.04279)(-1.83721,-1.03873)(-1.82970,-1.03445)%
\polyline(-1.82970,-1.03445)(-1.82251,-1.03037)(-1.81503,-1.02611)%
\polyline(-1.81503,-1.02611)(-1.80782,-1.02201)(-1.80037,-1.01777)%
\polyline(-1.80037,-1.01777)(-1.79312,-1.01365)(-1.78569,-1.00943)%
\polyline(-1.78569,-1.00943)(-1.77842,-1.00529)(-1.77102,-1.00109)%
\polyline(-1.77102,-1.00109)(-1.76372,-0.99694)(-1.75635,-0.99274)%
\polyline(-1.75635,-0.99274)(-1.74902,-0.98857)(-1.74217,-0.98468)%
\polyline(-1.74217,-0.98468)(-1.73433,-0.98021)(-1.72750,-0.97634)%
\polyline(-1.72750,-0.97634)(-1.71963,-0.97185)(-1.71284,-0.96799)%
\polyline(-1.71284,-0.96799)(-1.70493,-0.96350)(-1.69817,-0.95965)%
\polyline(-1.69817,-0.95965)(-1.69023,-0.95514)(-1.68350,-0.95131)%
\polyline(-1.68350,-0.95131)(-1.67553,-0.94678)(-1.66882,-0.94297)%
\polyline(-1.66882,-0.94297)(-1.66084,-0.93842)(-1.65415,-0.93463)%
\polyline(-1.65415,-0.93463)(-1.64614,-0.93007)(-1.63948,-0.92627)%
\polyline(-1.63948,-0.92627)(-1.63144,-0.92171)(-1.62482,-0.91793)%
\polyline(-1.62482,-0.91793)(-1.61674,-0.91335)(-1.61063,-0.90987)%
\polyline(-1.61063,-0.90987)(-1.60205,-0.90498)(-1.59597,-0.90153)%
\polyline(-1.59597,-0.90153)(-1.58735,-0.89663)(-1.58130,-0.89318)%
\polyline(-1.58130,-0.89318)(-1.57265,-0.88827)(-1.56662,-0.88484)%
\polyline(-1.56662,-0.88484)(-1.55795,-0.87991)(-1.55195,-0.87650)%
\polyline(-1.55195,-0.87650)(-1.54325,-0.87155)(-1.53728,-0.86816)%
\polyline(-1.53728,-0.86816)(-1.52856,-0.86320)(-1.52261,-0.85982)%
\polyline(-1.52261,-0.85982)(-1.51386,-0.85484)(-1.50795,-0.85147)%
\polyline(-1.50795,-0.85147)(-1.49916,-0.84648)(-1.49328,-0.84313)%
\polyline(-1.49328,-0.84313)(-1.48446,-0.83812)(-1.47910,-0.83507)%
\polyline(-1.47910,-0.83507)(-1.46980,-0.82979)%
\polyline(-1.46071,-0.83234)(-1.09308,-0.92773)(-0.70854,-1.01000)(-0.31200,-1.07848)%
(0.09150,-1.13262)(0.32861,-1.15568)%
\polyline(0.32861,-1.15568)(0.34136,-1.15691)%
\polyline(0.34136,-1.15691)(0.35434,-1.15818)%
\polyline(0.35434,-1.15818)(0.37806,-1.16048)%
\polyline(0.46306,-1.16875)(0.49677,-1.17202)(0.89852,-1.19644)(1.29136,-1.20577)%
(1.66990,-1.20010)(2.02872,-1.17967)(2.36255,-1.14490)(2.66618,-1.09641)(2.93463,-1.03499)%
(3.16323,-0.96163)(3.34761,-0.87749)(3.48392,-0.78395)(3.54859,-0.70663)%
\polyline(3.54859,-0.70663)(3.56326,-0.68910)%
\polyline(3.56326,-0.68910)(3.56878,-0.68250)(3.57252,-0.66938)%
\polyline(3.52576,-0.39630)(3.49100,-0.34833)(3.35021,-0.23342)(3.15212,-0.12017)%
(2.89822,-0.01071)(2.59096,0.09288)(2.23383,0.18856)(1.83126,0.27436)(1.38865,0.34850)%
(0.91232,0.40932)(0.40934,0.45544)(-0.11243,0.48569)(-0.64459,0.49921)(-1.17828,0.49546)%
(-1.70434,0.47423)(-2.21350,0.43569)(-2.69656,0.38034)(-3.14452,0.30910)(-3.28515,0.27923)%
\polyline(-4.41944,-0.10188)(-4.42565,-0.10510)(-4.44845,-0.12305)%
\polyline(-4.44845,-0.12305)(-4.47793,-0.14623)%
\polyline(-4.47793,-0.14623)(-4.58547,-0.23082)(-4.67173,-0.36047)(-4.68201,-0.49131)%
(-4.61526,-0.62046)(-4.47189,-0.74506)(-4.25378,-0.86225)(-3.96423,-0.96930)(-3.60796,-1.06360)%
(-3.19098,-1.14278)(-2.72047,-1.20475)(-2.20465,-1.24768)%
\polyline(-1.46675,-0.83579)(-1.47531,-0.84067)(-1.48139,-0.84413)%
\polyline(-1.48139,-0.84413)(-1.48992,-0.84899)(-1.49651,-0.85275)%
\polyline(-1.49651,-0.85275)(-1.50453,-0.85731)(-1.51115,-0.86109)%
\polyline(-1.51115,-0.86109)(-1.51913,-0.86564)(-1.52578,-0.86943)%
\polyline(-1.52578,-0.86943)(-1.53374,-0.87396)(-1.54042,-0.87777)%
\polyline(-1.54042,-0.87777)(-1.54835,-0.88229)(-1.55505,-0.88611)%
\polyline(-1.55505,-0.88611)(-1.56296,-0.89061)(-1.56969,-0.89445)%
\polyline(-1.56969,-0.89445)(-1.57756,-0.89893)(-1.58432,-0.90278)%
\polyline(-1.58432,-0.90278)(-1.59217,-0.90726)(-1.59896,-0.91112)%
\polyline(-1.59896,-0.91112)(-1.60678,-0.91558)(-1.61359,-0.91946)%
\polyline(-1.61359,-0.91946)(-1.62138,-0.92390)(-1.62871,-0.92808)%
\polyline(-1.62871,-0.92808)(-1.63599,-0.93223)(-1.64335,-0.93642)%
\polyline(-1.64335,-0.93642)(-1.65060,-0.94055)(-1.65798,-0.94476)%
\polyline(-1.65798,-0.94476)(-1.66521,-0.94887)(-1.67262,-0.95310)%
\polyline(-1.67262,-0.95310)(-1.67981,-0.95720)(-1.68725,-0.96144)%
\polyline(-1.68725,-0.96144)(-1.69442,-0.96552)(-1.70189,-0.96977)%
\polyline(-1.70189,-0.96977)(-1.70903,-0.97384)(-1.71652,-0.97811)%
\polyline(-1.71652,-0.97811)(-1.72363,-0.98217)(-1.73116,-0.98645)%
\polyline(-1.73116,-0.98645)(-1.73824,-0.99049)(-1.74628,-0.99507)%
\polyline(-1.74628,-0.99507)(-1.75285,-0.99881)(-1.76091,-1.00341)%
\polyline(-1.76091,-1.00341)(-1.76745,-1.00714)(-1.77555,-1.01175)%
\polyline(-1.77555,-1.01175)(-1.78206,-1.01546)(-1.79019,-1.02009)%
\polyline(-1.79019,-1.02009)(-1.79667,-1.02378)(-1.80482,-1.02843)%
\polyline(-1.80482,-1.02843)(-1.81128,-1.03211)(-1.81946,-1.03677)%
\polyline(-1.81946,-1.03677)(-1.82588,-1.04043)(-1.83409,-1.04511)%
\polyline(-1.83409,-1.04511)(-1.84049,-1.04875)(-1.84872,-1.05345)%
\polyline(-1.84872,-1.05345)(-1.85510,-1.05708)(-1.86336,-1.06179)%
\polyline(-1.86336,-1.06179)(-1.86970,-1.06540)(-1.87848,-1.07040)%
\polyline(-1.88436,-1.07375)(-1.89312,-1.07874)%
\polyline(-1.89902,-1.08211)(-1.90775,-1.08708)%
\polyline(-1.91368,-1.09046)(-1.92239,-1.09542)%
\polyline(-2.00117,-1.14031)(-2.01068,-1.14573)%
\polyline(-2.01583,-1.14867)(-2.02531,-1.15407)%
\polyline(-2.03050,-1.15702)(-2.03995,-1.16241)%
\polyline(-2.04516,-1.16538)(-2.05459,-1.17075)%
\polyline(-2.05982,-1.17373)(-2.06922,-1.17909)%
\polyline(-2.07449,-1.18209)(-2.08386,-1.18743)%
\polyline(-2.08915,-1.19044)(-2.09849,-1.19577)%
\polyline(-2.10381,-1.19880)(-2.11313,-1.20411)%
\polyline(-2.11847,-1.20715)(-2.12776,-1.21245)%
\polyline(-2.13264,-1.21523)(-2.14288,-1.22106)%
\polyline(-2.14731,-1.22359)(-2.15752,-1.22940)%
\polyline(-2.16197,-1.23194)(-2.17215,-1.23774)%
\polyline(-2.17663,-1.24030)(-2.18679,-1.24608)%
\polyline(-1.58243,-0.90171)(-1.18394,-0.98485)(-0.76701,-1.05348)(-0.33743,-1.10694)%
(0.09883,-1.14471)(0.32819,-1.15616)(0.33952,-1.15673)(0.35399,-1.15745)(0.36708,-1.15811)%
(0.37731,-1.15862)(0.39483,-1.15949)(0.40882,-1.16019)(0.43635,-1.16157)%
\polyline(0.49934,-1.16471)(0.53567,-1.16653)(0.96694,-1.17228)(1.38647,-1.16207)%
(1.78819,-1.13623)(2.16615,-1.09525)(2.51460,-1.03984)(2.82814,-0.97092)(3.10166,-0.88960)%
(3.33056,-0.79718)(3.51076,-0.69512)(3.62747,-0.59475)(3.63874,-0.58505)(3.64253,-0.57900)%
\polyline(3.71072,-0.29942)(3.68500,-0.22501)(3.58364,-0.10165)(3.42393,0.01993)(3.20727,0.13764)%
(2.93597,0.24940)(2.61330,0.35321)(2.24342,0.44714)(1.83135,0.52944)(1.38299,0.59849)%
(0.90493,0.65289)(0.40444,0.69151)(-0.11066,0.71344)(-0.63212,0.71807)(-1.15141,0.70512)%
(-1.65985,0.67462)(-2.14877,0.62693)(-2.60968,0.56274)(-3.03442,0.48306)(-3.31552,0.41384)%
\polyline(-4.21458,0.04914)(-4.22885,0.04065)(-4.23713,0.03307)(-4.24476,0.02611)%
(-4.29795,-0.02250)(-4.37276,-0.09085)(-4.44692,-0.22597)(-4.44942,-0.36210)(-4.37965,-0.49657)%
(-4.23828,-0.62669)(-4.02727,-0.74983)(-3.74986,-0.86344)(-3.41054,-0.96514)(-3.01490,-1.05272)%
(-2.56959,-1.12425)(-2.08217,-1.17803)%
\polyline(-1.70416,-0.97107)(-1.27481,-1.04198)(-0.82548,-1.09698)(-0.36286,-1.13539)%
(0.10615,-1.15680)(0.34975,-1.15901)(0.36417,-1.15913)(0.37807,-1.15926)%
\polyline(0.42375,-1.15967)(0.43542,-1.15978)(0.45974,-1.15999)(0.47487,-1.16013)%
(0.49648,-1.16032)%
\polyline(0.53562,-1.16068)(0.57458,-1.16103)(1.03537,-1.14812)(1.48160,-1.11839)%
(1.90650,-1.07236)(2.30357,-1.01081)(2.66666,-0.93476)(2.99010,-0.84542)(3.26869,-0.74421)%
(3.49791,-0.63274)(3.67390,-0.51274)(3.70293,-0.48204)(3.73170,-0.45160)(3.75967,-0.42200)%
(3.77430,-0.40653)(3.78813,-0.39190)(3.79356,-0.38614)(3.79859,-0.37534)%
\polyline(3.85493,-0.21220)(3.85560,-0.12126)(3.79603,0.01280)(3.67627,0.14503)(3.49765,0.27329)%
(3.26241,0.39547)(2.97373,0.50952)(2.63564,0.61353)(2.25300,0.70573)(1.83145,0.78450)%
(1.37732,0.84848)(0.89754,0.89647)(0.39953,0.92759)(-0.10889,0.94119)(-0.61965,0.93695)%
(-1.12454,0.91480)(-1.61536,0.87503)(-2.08404,0.81817)(-2.52282,0.74513)(-2.92431,0.65704)%
(-3.28172,0.55533)(-3.31477,0.54311)%
\polyline(-3.99688,0.21056)(-4.03205,0.18641)(-4.06480,0.15127)(-4.16006,0.04912)%
(-4.22210,-0.09147)(-4.21684,-0.23290)(-4.14404,-0.37267)(-4.00465,-0.50832)(-3.80075,-0.63741)%
(-3.53549,-0.75760)(-3.21311,-0.86669)(-2.83882,-0.96268)(-2.41871,-1.04375)(-1.95969,-1.10837)%
\polyline(-1.82588,-1.04043)(-1.36567,-1.09910)(-0.88394,-1.14046)(-0.38828,-1.16385)%
(0.11348,-1.16890)(0.40108,-1.16121)(0.43566,-1.16029)(0.45931,-1.15965)(0.47443,-1.15925)%
(0.49017,-1.15883)(0.50578,-1.15841)(0.52168,-1.15799)(0.53448,-1.15764)(0.55746,-1.15703)%
\polyline(0.61347,-1.15553)(1.10379,-1.12396)(1.57672,-1.07469)(2.02481,-1.00848)%
(2.44100,-0.92638)(2.81873,-0.82969)(3.15205,-0.71993)(3.43572,-0.59882)(3.66525,-0.46829)%
(3.80709,-0.35442)(3.83705,-0.33037)(3.84946,-0.31441)%
\polyline(3.96415,-0.14040)(3.99752,-0.04118)(3.98367,0.10554)(3.90706,0.25060)(3.76891,0.39171)%
(3.57137,0.52665)(3.31756,0.65329)(3.01150,0.76964)(2.65798,0.87386)(2.26259,0.96432)%
(1.83155,1.03958)(1.37165,1.09846)(0.89016,1.14004)(0.39463,1.16366)(-0.10712,1.16895)%
(-0.60717,1.15582)(-1.09766,1.12448)(-1.57087,1.07542)(-2.01932,1.00942)(-2.43595,0.92752)%
(-2.81420,0.83100)(-3.14813,0.72140)(-3.28801,0.66188)%
\polyline(-3.69388,0.44508)(-3.72232,0.42236)(-3.75081,0.39960)(-3.78070,0.37573)%
(-3.83524,0.33216)(-3.94736,0.18909)(-3.99729,0.04304)(-3.98424,-0.10369)(-3.90843,-0.24878)%
(-3.77104,-0.38995)(-3.57423,-0.52498)(-3.32112,-0.65174)(-3.01569,-0.76823)(-2.66273,-0.87262)%
(-2.26784,-0.96326)(-1.83721,-1.03873)%
\polyline(-1.94761,-1.10979)(-1.45653,-1.15622)(-0.94240,-1.18394)(-0.41371,-1.19230)%
(0.12081,-1.18099)(0.45771,-1.16137)(0.50536,-1.15859)(0.52121,-1.15767)(0.53732,-1.15674)%
(0.55369,-1.15578)(0.56529,-1.15511)(0.57914,-1.15430)(0.59615,-1.15331)(0.60818,-1.15261)%
(0.65237,-1.15004)(1.17221,-1.09980)(1.67184,-1.03099)(2.14312,-0.94460)(2.57843,-0.84196)%
(2.97079,-0.72462)(3.31401,-0.59443)(3.60275,-0.45344)(3.83260,-0.30384)(3.89452,-0.24626)%
(3.92265,-0.22009)(3.95200,-0.19280)(3.97595,-0.17054)(3.98899,-0.15841)(4.00018,-0.14800)%
(4.00104,-0.14667)(4.01012,-0.13261)%
\polyline(4.05045,-0.07011)(4.06183,-0.05247)(4.10321,0.01165)(4.14043,0.17260)(4.11174,0.33234)%
(4.01811,0.48840)(3.86154,0.63839)(3.64509,0.78000)(3.37271,0.91111)(3.04925,1.02975)%
(2.68031,1.13419)(2.27217,1.22291)(1.83165,1.29466)(1.36599,1.34845)(0.88277,1.38361)%
(0.38971,1.39974)(-0.10535,1.39670)(-0.59471,1.37468)(-1.07079,1.33415)(-1.52637,1.27581)%
(-1.95458,1.20067)(-2.34908,1.10991)(-2.70410,1.00497)(-3.01454,0.88746)(-3.23898,0.77733)%
\polyline(-3.46864,0.63265)(-3.48490,0.62197)(-3.49739,0.61026)(-3.51773,0.59117)%
(-3.52578,0.58361)(-3.57876,0.53390)(-3.63843,0.47791)(-3.73466,0.32906)(-3.77247,0.17754)%
(-3.75165,0.02552)(-3.67282,-0.12489)(-3.53741,-0.27159)(-3.34772,-0.41256)(-3.10675,-0.54589)%
(-2.81826,-0.66978)(-2.48665,-0.78256)(-2.11696,-0.88277)(-1.72597,-0.96666)(-1.71481,-0.96905)%
\polyline(-2.06933,-1.17915)(-1.54740,-1.21335)(-1.00087,-1.22743)(-0.43914,-1.22075)%
(0.12814,-1.19308)(0.45392,-1.16500)(0.50557,-1.16055)(0.52061,-1.15925)(0.56196,-1.15569)%
(0.57253,-1.15477)(0.59565,-1.15278)(0.61289,-1.15130)(0.62906,-1.14990)(0.63966,-1.14899)%
(0.65757,-1.14745)(0.67624,-1.14583)(0.69127,-1.14454)(1.24064,-1.07565)(1.76697,-0.98730)%
(2.26143,-0.88073)(2.71585,-0.75752)(3.12286,-0.61955)(3.47597,-0.46894)(3.76978,-0.30805)%
(3.98273,-0.15201)(3.99994,-0.13939)(4.00820,-0.13061)(4.01959,-0.11849)%
\polyline(4.15492,0.02544)(4.16333,0.03438)(4.25803,0.21055)(4.28335,0.38637)(4.23981,0.55914)%
(4.12914,0.72621)(3.95417,0.88507)(3.71880,1.03336)(3.42786,1.16893)(3.08701,1.28987)%
(2.70265,1.39452)(2.28176,1.48149)(1.83174,1.54972)(1.36032,1.59844)(0.87538,1.62719)%
(0.38481,1.63581)(-0.10358,1.62445)(-0.58223,1.59356)(-1.04392,1.54383)(-1.48188,1.47622)%
(-1.88985,1.39191)(-2.26222,1.29230)(-2.59399,1.17894)(-2.88094,1.05353)(-3.11955,0.91790)%
(-3.13097,0.90914)(-3.15572,0.89014)%
\polyline(-3.27434,0.79909)(-3.28969,0.78732)(-3.30710,0.77395)(-3.32898,0.74951)%
(-3.36925,0.70453)(-3.37705,0.69582)(-3.44163,0.62366)(-3.52196,0.46903)(-3.54766,0.31205)%
(-3.51907,0.15473)(-3.43720,-0.00100)(-3.30380,-0.15322)(-3.12120,-0.30013)(-2.89238,-0.44004)%
(-2.62083,-0.57132)(-2.31057,-0.69252)(-1.96608,-0.80227)(-1.59226,-0.89941)%
\polyline(0.16146,-1.14142)(0.16299,-1.14252)(0.16452,-1.14363)(0.16605,-1.14473)%
(0.16758,-1.14583)(0.16912,-1.14694)(0.16964,-1.14731)%
\polyline(0.15992,-1.14032)(0.16146,-1.14142)(0.16299,-1.14252)(0.16452,-1.14363)%
(0.16605,-1.14473)(0.16758,-1.14583)(0.16912,-1.14694)(0.17064,-1.14805)(0.17218,-1.14915)%
(0.17371,-1.15026)(0.17524,-1.15136)(0.17677,-1.15246)(0.17831,-1.15357)(0.17983,-1.15468)%
(0.18137,-1.15578)(0.18290,-1.15688)(0.18444,-1.15799)(0.18596,-1.15909)(0.18750,-1.16020)%
(0.18903,-1.16130)(0.19056,-1.16241)(0.19209,-1.16351)(0.19362,-1.16462)(0.19515,-1.16572)%
(0.19669,-1.16682)(0.19822,-1.16793)(0.19975,-1.16904)(0.20128,-1.17014)(0.20281,-1.17124)%
(0.20434,-1.17235)(0.20587,-1.17345)(0.20741,-1.17456)(0.20894,-1.17566)(0.21047,-1.17677)%
(0.21200,-1.17787)(0.21353,-1.17898)(0.21506,-1.18008)(0.21660,-1.18118)(0.21813,-1.18229)%
(0.21966,-1.18340)(0.22119,-1.18450)(0.22272,-1.18560)(0.22425,-1.18671)(0.22579,-1.18781)%
(0.22732,-1.18892)(0.22885,-1.19003)(0.23038,-1.19112)(0.23192,-1.19223)(0.23344,-1.19334)%
(0.23498,-1.19445)(0.23651,-1.19554)%
\polyline(0.17736,-1.15289)(0.17831,-1.15357)(0.17983,-1.15468)(0.18137,-1.15578)%
(0.18290,-1.15688)(0.18444,-1.15799)(0.18596,-1.15909)(0.18750,-1.16020)(0.18903,-1.16130)%
(0.19056,-1.16241)(0.19209,-1.16351)(0.19362,-1.16462)(0.19515,-1.16572)(0.19669,-1.16682)%
(0.19822,-1.16793)(0.19975,-1.16904)(0.20128,-1.17014)(0.20281,-1.17124)(0.20434,-1.17235)%
(0.20587,-1.17345)(0.20741,-1.17456)(0.20894,-1.17566)(0.21047,-1.17677)(0.21200,-1.17787)%
(0.21353,-1.17898)(0.21506,-1.18008)(0.21660,-1.18118)(0.21813,-1.18229)(0.21966,-1.18340)%
(0.22119,-1.18450)(0.22272,-1.18560)(0.22425,-1.18671)(0.22579,-1.18781)(0.22732,-1.18892)%
(0.22885,-1.19003)(0.23038,-1.19112)(0.23192,-1.19223)(0.23344,-1.19334)(0.23458,-1.19416)%
\polyline(1.79698,-1.19286)(1.80852,-1.18634)(1.80992,-1.18555)(1.82354,-1.17785)%
(1.83857,-1.16936)(1.85360,-1.16086)(1.86863,-1.15237)(1.88366,-1.14387)(1.89869,-1.13538)%
(1.91371,-1.12688)(1.92874,-1.11839)(1.94377,-1.10990)(1.95880,-1.10140)(1.97383,-1.09291)%
(1.98886,-1.08441)(2.00389,-1.07592)(2.01892,-1.06742)(2.03394,-1.05893)(2.04897,-1.05044)%
(2.06400,-1.04194)(2.07903,-1.03345)(2.09406,-1.02495)(2.10909,-1.01646)(2.12412,-1.00797)%
(2.13914,-0.99947)(2.15417,-0.99098)(2.16920,-0.98248)(2.18423,-0.97399)(2.19926,-0.96549)%
(2.21429,-0.95700)(2.22932,-0.94851)(2.24435,-0.94001)(2.25937,-0.93152)(2.27440,-0.92302)%
(2.28943,-0.91453)(2.30446,-0.90604)(2.31949,-0.89754)(2.33452,-0.88905)(2.34954,-0.88055)%
(2.36458,-0.87206)(2.37960,-0.86357)(2.39463,-0.85507)(2.40966,-0.84658)(2.42469,-0.83808)%
(2.43972,-0.82959)(2.45475,-0.82109)(2.46977,-0.81260)(2.48481,-0.80410)(2.49983,-0.79561)%
(2.51486,-0.78712)(2.52926,-0.77898)(2.52989,-0.77862)(2.54492,-0.77013)%
\polyline(3.06080,-0.99450)(3.07577,-0.98065)(3.09591,-0.96203)(3.11604,-0.94341)%
(3.13617,-0.92479)(3.15630,-0.90617)(3.17644,-0.88755)(3.19657,-0.86893)(3.21670,-0.85031)%
(3.23683,-0.83169)(3.25697,-0.81307)(3.27710,-0.79445)(3.29723,-0.77584)(3.31736,-0.75722)%
(3.33750,-0.73859)(3.35763,-0.71997)(3.37776,-0.70136)(3.39789,-0.68274)(3.41803,-0.66412)%
(3.43816,-0.64549)(3.45829,-0.62688)(3.47842,-0.60826)(3.49856,-0.58964)(3.51869,-0.57102)%
(3.53882,-0.55240)(3.55895,-0.53378)(3.57909,-0.51516)(3.59921,-0.49655)(3.61935,-0.47793)%
(3.63948,-0.45930)(3.65962,-0.44068)(3.67974,-0.42207)(3.69988,-0.40345)(3.72001,-0.38483)%
(3.74015,-0.36620)(3.76027,-0.34759)(3.78041,-0.32897)(3.80054,-0.31035)(3.82068,-0.29173)%
(3.84081,-0.27311)(3.86094,-0.25449)(3.88107,-0.23587)(3.90121,-0.21725)(3.92134,-0.19864)%
(3.94147,-0.18001)(3.96160,-0.16139)(3.98174,-0.14277)(4.00051,-0.12541)(4.00187,-0.12416)%
(4.01949,-0.10786)(4.02200,-0.10554)(4.04121,-0.08777)(4.04214,-0.08691)(4.04962,-0.07999)%
(4.06227,-0.06829)%
\polyline(3.81810,-0.23814)(3.81930,-0.23613)(3.83528,-0.20939)(3.85126,-0.18266)%
(3.86724,-0.15592)(3.88323,-0.12919)(3.89920,-0.10246)(3.91518,-0.07572)(3.93117,-0.04899)%
(3.94714,-0.02226)(3.96313,0.00447)(3.97911,0.03121)(3.99508,0.05795)(4.01107,0.08468)%
(4.02704,0.11142)(4.04303,0.13815)(4.05901,0.16488)(4.07498,0.19161)(4.09097,0.21835)%
(4.10695,0.24508)(4.12293,0.27182)(4.13891,0.29856)(4.15489,0.32529)(4.17087,0.35202)%
(4.18685,0.37875)(4.20283,0.40549)(4.21881,0.43222)(4.23479,0.45895)(4.25077,0.48569)%
(4.26675,0.51243)(4.28273,0.53916)(4.29871,0.56589)(4.31469,0.59263)(4.33068,0.61936)%
(4.34665,0.64609)(4.36263,0.67282)(4.37862,0.69956)(4.38937,0.71755)(4.39459,0.72630)%
\polyline(3.18677,-0.13766)(3.19373,-0.10679)(3.20069,-0.07591)(3.20765,-0.04504)%
(3.21461,-0.01417)(3.22157,0.01670)(3.22852,0.04758)(3.23549,0.07845)(3.24244,0.10932)%
(3.24940,0.14020)(3.25635,0.17107)(3.26332,0.20194)(3.27028,0.23280)(3.27724,0.26368)%
(3.28419,0.29455)(3.29115,0.32542)(3.29811,0.35629)(3.30507,0.38717)(3.31202,0.41804)%
(3.31899,0.44891)(3.32594,0.47979)(3.33291,0.51066)(3.33986,0.54153)(3.34682,0.57240)%
(3.35378,0.60328)(3.36074,0.63415)(3.36769,0.66502)(3.37465,0.69589)(3.38161,0.72677)%
(3.38857,0.75764)(3.39553,0.78851)(3.40249,0.81938)(3.40945,0.85026)(3.41641,0.88113)%
(3.42337,0.91200)(3.43032,0.94288)(3.43728,0.97375)(3.44424,1.00462)(3.45120,1.03550)%
(3.45816,1.06636)(3.46512,1.09723)(3.47208,1.12810)(3.47904,1.15898)(3.48599,1.18985)%
(3.49295,1.22072)(3.49991,1.25159)(3.50687,1.28247)(3.51382,1.31334)(3.52078,1.34421)%
(3.52774,1.37508)(3.53417,1.40360)%
\polyline(1.82588,0.27538)(1.82588,0.30598)(1.82588,0.33658)(1.82588,0.36719)(1.82588,0.39779)%
(1.82588,0.42839)(1.82588,0.45899)(1.82588,0.48959)(1.82588,0.52020)(1.82588,0.55080)%
(1.82588,0.58140)(1.82588,0.61200)(1.82588,0.64260)(1.82588,0.67321)(1.82588,0.70381)%
(1.82588,0.73441)(1.82588,0.76501)(1.82588,0.79561)(1.82588,0.82622)(1.82588,0.85682)%
(1.82588,0.88742)(1.82588,0.91802)(1.82588,0.94862)(1.82588,0.97923)(1.82588,1.00983)%
(1.82588,1.04043)(1.82588,1.07103)(1.82588,1.10163)(1.82588,1.13224)(1.82588,1.16284)%
(1.82588,1.19344)(1.82588,1.22404)(1.82588,1.25464)(1.82588,1.28525)(1.82588,1.31585)%
(1.82588,1.34645)(1.82588,1.37705)(1.82588,1.40766)(1.82588,1.43826)(1.82588,1.46886)%
(1.82588,1.49946)(1.82588,1.53006)(1.82588,1.56067)(1.82588,1.59127)(1.82588,1.62187)%
(1.82588,1.65247)(1.82588,1.68307)(1.82588,1.71368)(1.82588,1.74428)(1.82588,1.77488)%
(1.82588,1.80534)%
\polyline(-0.20848,0.48941)(-0.20806,0.51654)(-0.20765,0.54369)(-0.20724,0.57083)%
(-0.20683,0.59797)(-0.20642,0.62511)(-0.20601,0.65225)(-0.20560,0.67939)(-0.20519,0.70654)%
(-0.20479,0.73368)(-0.20437,0.76082)(-0.20396,0.78796)(-0.20355,0.81509)(-0.20314,0.84223)%
(-0.20273,0.86937)(-0.20232,0.89652)(-0.20191,0.92366)(-0.20150,0.95080)(-0.20109,0.97794)%
(-0.20068,1.00508)(-0.20026,1.03222)(-0.19985,1.05936)(-0.19944,1.08651)(-0.19903,1.11365)%
(-0.19862,1.14079)(-0.19822,1.16793)(-0.19781,1.19507)(-0.19740,1.22221)(-0.19699,1.24935)%
(-0.19657,1.27649)(-0.19616,1.30364)(-0.19575,1.33078)(-0.19534,1.35792)(-0.19493,1.38506)%
(-0.19452,1.41220)(-0.19411,1.43935)(-0.19370,1.46648)(-0.19329,1.49363)(-0.19287,1.52077)%
(-0.19246,1.54790)(-0.19206,1.57504)(-0.19165,1.60218)(-0.19124,1.62932)(-0.19083,1.65647)%
(-0.19042,1.68361)(-0.19001,1.71075)(-0.18960,1.73789)(-0.18918,1.76503)(-0.18877,1.79217)%
(-0.18836,1.81931)(-0.18798,1.84488)%
\polyline(-2.38612,0.41818)(-2.37744,0.44075)(-2.36877,0.46333)(-2.36009,0.48590)%
(-2.35142,0.50847)(-2.34274,0.53104)(-2.33406,0.55361)(-2.32539,0.57618)(-2.31671,0.59875)%
(-2.30803,0.62132)(-2.29936,0.64390)(-2.29068,0.66647)(-2.28200,0.68904)(-2.27333,0.71162)%
(-2.26465,0.73419)(-2.25598,0.75676)(-2.24730,0.77933)(-2.23862,0.80191)(-2.22994,0.82448)%
(-2.22126,0.84705)(-2.21258,0.86963)(-2.20391,0.89220)(-2.19523,0.91477)(-2.18655,0.93733)%
(-2.17788,0.95991)(-2.16920,0.98248)(-2.16052,1.00505)(-2.15185,1.02763)(-2.14317,1.05020)%
(-2.13450,1.07277)(-2.12582,1.09534)(-2.11714,1.11792)(-2.10847,1.14049)(-2.09979,1.16306)%
(-2.09111,1.18564)(-2.08244,1.20821)(-2.07376,1.23078)(-2.06508,1.25335)(-2.05641,1.27593)%
(-2.04773,1.29849)(-2.03905,1.32107)(-2.03037,1.34364)(-2.02169,1.36621)(-2.01302,1.38878)%
(-2.00434,1.41135)(-1.99566,1.43393)(-1.98699,1.45650)(-1.97831,1.47907)(-1.96963,1.50164)%
(-1.96096,1.52422)(-1.95360,1.54337)%
\polyline(-3.19504,0.87034)(-3.17644,0.88755)(-3.16820,0.89517)(-3.15872,0.90393)%
(-3.14700,0.91478)%
\polyline(-3.12407,0.93598)(-3.11604,0.94341)(-3.11589,0.94355)(-3.09591,0.96203)%
(-3.09237,0.96530)(-3.07577,0.98065)(-3.06841,0.98745)(-3.05564,0.99927)%
\polyline(-4.67842,-0.44561)(-4.65937,-0.43481)(-4.63170,-0.41911)(-4.60402,-0.40340)%
(-4.57634,-0.38770)(-4.54866,-0.37200)(-4.52098,-0.35629)(-4.49330,-0.34059)(-4.46562,-0.32489)%
(-4.43794,-0.30919)(-4.41027,-0.29349)(-4.38259,-0.27778)(-4.35491,-0.26208)(-4.32723,-0.24638)%
(-4.29955,-0.23067)(-4.27187,-0.21497)(-4.24419,-0.19927)(-4.21651,-0.18357)(-4.18884,-0.16787)%
(-4.16116,-0.15216)(-4.13348,-0.13646)(-4.10580,-0.12076)(-4.07812,-0.10505)(-4.05044,-0.08935)%
(-4.02276,-0.07365)(-3.99508,-0.05795)(-3.96740,-0.04224)(-3.93973,-0.02654)(-3.91205,-0.01084)%
(-3.88437,0.00486)(-3.85669,0.02056)(-3.82901,0.03627)(-3.80134,0.05197)(-3.77365,0.06767)%
(-3.74597,0.08338)(-3.71830,0.09908)(-3.69062,0.11478)(-3.66294,0.13048)(-3.63526,0.14618)%
(-3.60758,0.16189)(-3.57990,0.17759)(-3.55223,0.19329)(-3.52454,0.20900)(-3.49687,0.22470)%
(-3.46919,0.24040)(-3.44151,0.25610)(-3.41383,0.27180)(-3.38615,0.28751)(-3.35847,0.30321)%
(-3.33080,0.31891)(-3.30334,0.33449)%
\polyline(-4.00733,-0.95336)(-3.98401,-0.94184)(-3.95804,-0.92902)(-3.93207,-0.91621)%
(-3.90610,-0.90339)(-3.88013,-0.89056)(-3.85416,-0.87774)(-3.82819,-0.86492)(-3.80222,-0.85210)%
(-3.77626,-0.83928)(-3.75028,-0.82646)(-3.72432,-0.81364)(-3.69834,-0.80082)(-3.67238,-0.78800)%
(-3.64640,-0.77517)(-3.62044,-0.76235)(-3.59446,-0.74954)(-3.56850,-0.73672)(-3.54252,-0.72390)%
(-3.51656,-0.71107)(-3.49059,-0.69825)(-3.46462,-0.68543)(-3.43864,-0.67262)(-3.41268,-0.65979)%
(-3.38670,-0.64697)(-3.36074,-0.63415)(-3.33476,-0.62133)(-3.30880,-0.60851)(-3.28283,-0.59568)%
(-3.25686,-0.58287)(-3.23088,-0.57005)(-3.20492,-0.55723)(-3.17894,-0.54440)(-3.15298,-0.53158)%
(-3.12701,-0.51876)(-3.10104,-0.50594)(-3.07507,-0.49312)(-3.04910,-0.48030)(-3.02313,-0.46748)%
(-2.99716,-0.45466)(-2.97119,-0.44184)(-2.94522,-0.42901)(-2.91925,-0.41619)(-2.89328,-0.40338)%
(-2.86731,-0.39056)(-2.84134,-0.37774)(-2.81537,-0.36491)(-2.78940,-0.35209)(-2.76344,-0.33927)%
(-2.73746,-0.32646)(-2.71149,-0.31363)%
\polyline(-2.18679,-1.24608)(-2.17711,-1.24057)(-2.17645,-1.24019)(-2.16245,-1.23221)%
(-2.16185,-1.23187)(-2.14779,-1.22386)(-2.14724,-1.22355)(-2.13312,-1.21550)(-2.13263,-1.21522)%
(-2.11881,-1.20735)(-2.11802,-1.20690)(-2.10418,-1.19901)(-2.10342,-1.19858)(-2.08956,-1.19068)%
(-2.08881,-1.19025)(-2.07493,-1.18235)(-2.07420,-1.18193)(-2.06400,-1.17612)(-2.05960,-1.17361)%
(-2.04560,-1.16563)(-2.04499,-1.16528)(-2.03094,-1.15728)(-2.03038,-1.15696)(-2.01627,-1.14892)%
(-2.01577,-1.14864)(-2.00146,-1.14048)%
\polyline(-1.99541,-1.13703)(-1.98663,-1.13203)(-1.98656,-1.13199)(-1.97270,-1.12409)%
(-1.97195,-1.12366)(-1.95808,-1.11576)(-1.94342,-1.10741)(-1.94274,-1.10702)(-1.93004,-1.09978)%
(-1.92813,-1.09869)(-1.91541,-1.09144)(-1.91353,-1.09037)(-1.89942,-1.08234)(-1.89892,-1.08205)%
(-1.88476,-1.07398)(-1.88431,-1.07372)(-1.87048,-1.06584)(-1.86970,-1.06540)(-1.85585,-1.05751)%
(-1.85510,-1.05708)(-1.84122,-1.04917)(-1.82741,-1.04130)(-1.81277,-1.03296)(-1.79813,-1.02462)%
(-1.79667,-1.02378)(-1.78258,-1.01575)(-1.78206,-1.01546)(-1.76885,-1.00793)(-1.76745,-1.00714)%
(-1.75421,-0.99959)(-1.75285,-0.99881)(-1.73900,-0.99092)(-1.72437,-0.98259)(-1.70972,-0.97424)%
(-1.69506,-0.96589)(-1.68040,-0.95753)(-1.66573,-0.94917)(-1.66521,-0.94887)(-1.65107,-0.94082)%
(-1.65060,-0.94055)(-1.63640,-0.93246)(-1.63599,-0.93223)(-1.62571,-0.92637)(-1.62138,-0.92390)%
(-1.61117,-0.91808)(-1.60678,-0.91558)(-1.59287,-0.90766)(-1.57821,-0.89930)(-1.56354,-0.89094)%
(-1.54888,-0.88259)(-1.53422,-0.87423)(-1.53374,-0.87396)(-1.51955,-0.86588)(-1.51913,-0.86564)%
(-1.50489,-0.85752)(-1.50453,-0.85731)(-1.49023,-0.84917)%
\polyline(-1.47531,-0.84067)(-1.46071,-0.83234)%
\end{picture}}%